\documentclass[a4paper,reqno,11pt]{amsart}

\usepackage{setspace, amssymb, amsmath, amsthm, graphicx, color, comment, booktabs, multirow, hyperref, mathtools, thmtools, thm-restate}
\usepackage[T1]{fontenc}
\usepackage[utf8]{inputenc}
\usepackage[inline]{enumitem}
\usepackage[dvipsnames,svgnames,table]{xcolor}
\usepackage[foot]{amsaddr}
\usepackage[capitalise, noabbrev]{cleveref}

\usepackage{tabularx}

\usepackage[longnamesfirst,numbers,sort&compress]{natbib}

\hypersetup{
    pdftitle={Centered colorings and weak coloring numbers in minor-closed graph classes},
    pdfauthor={Jędrzej Hodor, Hoang La, Piotr Micek, Clément Rambaud},
    colorlinks,
    linkcolor={RoyalBlue},
    citecolor={RubineRed},
    urlcolor={blue!80!black},
    nesting=false,
}

\input{style.toc}
\newcommand{\q}[1]{``#1''}
\newcommand{\defin}[1]{\emph{\textcolor{ForestGreen}{#1}}}
\newcommand{\defmath}[1]{\defin{\text{$#1$}}}

\DeclarePairedDelimiter\set{\{}{\}}

\newcommand{\td}{\operatorname{td}}
\newcommand{\ftd}{\operatorname{ftd}}
\newcommand{\rtd}{\operatorname{rtd}}
\newcommand{\srtd}{\operatorname{srtd}}
\newcommand{\pw}{\operatorname{pw}}
\newcommand{\tw}{\operatorname{tw}}

\newcommand{\ltw}{\operatorname{ltw}}

\newcommand{\Trace}{\operatorname{Trace}}

\newcommand{\cgF}{\mathcal{F}}

\newcommand{\calM}{\mathcal{M}}

\newcommand{\Oh}{\mathcal{O}}

\newcommand{\calA}{\mathcal{A}} 
\newcommand{\calB}{\mathcal{B}} 
\newcommand{\calC}{\mathcal{C}}
\newcommand{\calD}{\mathcal{D}}
\newcommand{\calE}{\mathcal{E}} 
\newcommand{\calF}{\mathcal{F}}
 
\newcommand{\calH}{\mathcal{H}}
\newcommand{\calI}{\mathcal{I}} 
 
\newcommand{\calK}{\mathcal{K}} 
\newcommand{\calL}{\mathcal{L}}

\newcommand{\calP}{\mathcal{P}} 
\newcommand{\calQ}{\mathcal{Q}} 
\newcommand{\calR}{\mathcal{R}}
\newcommand{\calS}{\mathcal{S}}
\newcommand{\calT}{\mathcal{T}}
\newcommand{\calU}{\mathcal{U}} 
\newcommand{\calV}{\mathcal{V}}
\newcommand{\calW}{\mathcal{W}}
\newcommand{\calX}{\mathcal{X}}  
\newcommand{\calY}{\mathcal{Y}}
\newcommand{\calZ}{\mathcal{Z}}

\newcommand{\Apex}{\mathbf{A}}
\newcommand{\Tree}{\mathbf{T}}
\newcommand{\Rt}{\mathcal{R}}
\newcommand{\SRt}{\mathcal{S}}
\newcommand{\edgeless}{\mathcal{E}}

\newcommand{\posint}{\NN_{1}}  
\newcommand{\nonnegint}{\NN}  
\newcommand{\nonnegrat}{\mathbb{Q}_{\geq 0}}

\newcommand{\rvar}[1]{\mathsf{#1}}
\newcommand{\param}{\mathrm{par}}
\newcommand{\paramktmf}{\mathrm{cen\text{-}k}}
\newcommand{\parambdtw}{\mathrm{cen\text{-}t}}
\newcommand{\cLRS}{c_{\mathrm{LRS}}}
\DeclareMathOperator\parent{p}
\let\root\relax
\DeclareMathOperator\root{r}
\DeclareMathOperator\frate{ftdfr}
\DeclareMathOperator\cen{cen}

\newcommand{\bigO}{\mathcal{O}}

\newcommand{\lca}{\mathrm{lca}}
\newcommand{\LCA}{\mathrm{LCA}}

\newcommand{\torso}{\mathrm{torso}}
\newcommand{\dist}{\mathrm{dist}}

\newcommand{\NN}{\mathbb{N}}

\let\leq\leqslant
\let\geq\geqslant

\let\subset\subseteq

\let\epsilon\varepsilon
\let\phi\varphi

\DeclareMathOperator\WReach{WReach}
\DeclareMathOperator\wcol{wcol}

\DeclareMathOperator\vc{vc}

\renewcommand{\setminus}{-}

\renewcommand{\root}{\mathrm{root}}

\newcommand{\stpath}[2]{$#1$--$#2$ path}

\newcommand{\stgeodesic}[2]{$#1$--$#2$ geodesic}

\newcommand{\ext}[3]{\mathrm{ext}(#1, #2, #3)}

\newcommand{\ttt}{\mathtt{t}}

\newcommand{\subtree}[2]{\mathrm{Sub}(#1,#2)}

\makeatletter
\newcommand{\pushright}[1]{\ifmeasuring@#1\else\omit\hfill$\displaystyle{#1}$\fi\ignorespaces}
\newcommand{\pushleft}[1]{\ifmeasuring@#1\else\omit$\displaystyle{#1}$\hfill\fi\ignorespaces}
\makeatother
\usepackage{tikz}
\usetikzlibrary{arrows}
\usetikzlibrary{calc}
\tikzset{myspacing/.style = {outer sep = 5pt, inner sep = 5pt}}
\tikzset{implies/.style = {-implies, double distance=3pt, thick}}
\tikzset{box/.style = {rectangle, draw, myspacing, align=center}}

\providecommand{\noopsort}[1]{}

\makeatletter
\def\thm@space@setup{
\thm@preskip=4mm
\thm@postskip=0mm
}
\makeatother

\crefformat{equation}{#2(#1)#3}
\let\eqref\cref
\crefformat{subsection}{Subsection #2#1#3}
\crefformat{subsubsection}{Subsubsection #2#1#3}

\theoremstyle{plain}
\newtheorem{thm}{Theorem}[section]
\newtheorem*{thm*}{Theorem}
\newtheorem{theorem}[thm]{Theorem}

\newtheorem{lemma}[thm]{Lemma}
\newtheorem*{lemma*}{Lemma}
\newtheorem{cor}[thm]{Corollary}
\newtheorem*{cor*}{Corollary}
\newtheorem{corollary}[thm]{Corollary}
\newtheorem*{corollary*}{Corollary}
\newtheorem{obs}[thm]{Observation}

\theoremstyle{remark}

\newtheorem*{problem*}{Open problem}

\newtheorem{claim}[thm]{Claim}
\newtheorem*{claim*}{Claim}

\crefname{obs}{Observation}{Observations}
\theoremstyle{definition}

\newtheorem*{conj*}{Conjecture}
\crefname{lem}{Lemma}{Lemmas}
\crefname{thm}{Theorem}{Theorems}
\crefname{cor}{Corollary}{Corollaries}

\newenvironment{proofclaim}[1][]
    {\let\oldqed\qedsymbol\renewcommand{\qedsymbol}{\ensuremath{\lozenge}}\begin{proof}[Proof of the claim] }{\end{proof}\renewcommand{\qedsymbol}{\oldqed}}

\setenumerate{label=\textup{(\roman*)}, noitemsep, topsep=3pt-\parskip, 
labelindent=.2em, leftmargin=*, widest=iii,}
\setitemize{noitemsep, topsep=-\parskip, labelindent=.2em, leftmargin=*, widest=iii,}

\newenvironment{enumerateOurAlph}{\begin{enumerate}
[label={\normalfont(\makebox[\mywidth]{\alph*})}]}{\end{enumerate}}

\newenvironment{enumerateOurAlphPrim}{\begin{enumerate}
[label={\normalfont(\makebox[\mywidthprim]{\alph*'})}]}{\end{enumerate}}

\newdimen\mywidth
\sbox0{m}%
\newdimen\mywidthprim
\sbox0{m}%
\newdimen\mywidthprimprim
\sbox0{m}%
\newdimen\mywidthA
\sbox0{m}%
\newdimen\mywidthAprim
\sbox0{m}%
\title[Centered colorings and weak coloring numbers]{Centered colorings and weak coloring numbers in minor-closed graph classes}

\thanks{The results of this paper for weak coloring numbers were announced at SODA 2025~\cite{wcol_paper}, see also the full version~\cite{wcol_paper_arxiv}; 
while the result for the centered chromatic numbers of $K_t$-minor-free graphs was announced at SODA 2026~\cite{centered_paper_soda}. A substantial portion of the material presented in this article is included in the PhD thesis of the fourth author~\cite{rambaud_thesis}.}

\begin{document}

\author[Hodor]{Jędrzej Hodor}
\address[J.~Hodor]{Theoretical Computer Science Department, 
Faculty of Mathematics and Computer Science and Doctoral School of Exact and Natural Sciences, Jagiellonian University, Krak\'ow, Poland}
\email{jedrzej.hodor@gmail.com}

\author[La]{Hoang La}
\address[H.~La]{LISN, Universit\'e Paris-Saclay, CNRS, Gif-sur-Yvette, France}
\email{hoang.la.research@gmail.com}

\author[Micek]{Piotr Micek}
\address[P.~Micek]{Theoretical Computer Science Department, 
Faculty of Mathematics and Computer Science, Jagiellonian University, Krak\'ow, Poland}
\email{piotr.micek@uj.edu.pl}

\author[Rambaud]{Cl\'ement Rambaud}
\address[C.~Rambaud]{Universit\'e C\^ote d'Azur, CNRS, Inria, I3S, Sophia Antipolis, France}
\email{clement.rambaud@normalesup.org}

\thanks{This research was funded by the National Science Center of Poland under grant UMO-2023/05/Y/ST6/00079 within the WEAVE-UNISONO program. 
Additionally, J.\ Hodor was supported by a Polish Ministry of Education and Science grant (Perły Nauki; PN/01/0265/2022). H.\ La  was supported by ANR PIA funding: ANR-20-IDEES-0002. C.\ Rambaud was supported by ANR DIGRAPHS funding: ANR-19-CE48-0013.}

\begin{abstract}
    Let $\mathcal{C}$ be a proper minor-closed class of graphs.
    Given the minors excluded in $\mathcal{C}$,
    we determine the maximum $q$-centered chromatic number and the maximum $q$th weak coloring number of graphs in $\mathcal{C}$ within an $\bigO(q)$-factor.
    Moreover, when $\mathcal{C}$ excludes a planar graph, we determine it within a constant factor.
    Our results imply that the $q$-centered chromatic number of $K_t$-minor-free graphs is in $\Oh(q^{t-1})$,
    improving on the previously known $\bigO(q^{h(t)})$ bound with a large and non-explicit function $h$.
    We include similar bounds for another family of parameters, the fractional treedepth fragility rates.
    All our bounds are proved via the same general framework.
\end{abstract}

\maketitle
\thispagestyle{empty}

\clearpage
%\newgeometry{top=1.5cm}
%\vspace*{-2cm}
{
\setstretch{0.9}
\tableofcontents
}
%\restoregeometry
\thispagestyle{empty}
\clearpage

\section{Introduction}\label{sec:intro}

One of the driving forces in graph theory continues to be the development of efficient algorithms for computationally hard problems for sparse graph classes.
Ne{\v{s}}et{\v{r}}il and Ossona de Mendez~\cite{sparsity} introduced the concepts of bounded expansion and nowhere denseness of classes of graphs.
These notions cover many well-studied classes of graphs such as planar graphs, graphs of bounded treewidth, graphs excluding a fixed minor, graphs of bounded book-thickness, and graphs that admit drawings with a bounded number of crossings per edge; see the paper of Ne{\v{s}}et{\v{r}}il, Ossona de Mendez, and Wood~\cite{NOdMW12} and the lecture notes of Pilipczuk, Pilipczuk, and Siebertz~\cite{notes}. 
Weak coloring numbers and centered chromatic numbers are key families of parameters capturing these concepts. 
In this work, we study bounds on these parameters in proper minor-closed classes of graphs. 

We start with an introduction of the main characters and some theorems that follow as corollaries from a bigger framework that we develop within this work.

Let $G$ be a graph, let $q$ be a positive integer, and let $C$ be a set of colors.
A coloring $\varphi\colon V(G)\to C$ of $G$ is \defin{$q$-centered} if for every connected\footnote{In this paper, connected graphs are nonnull, that is, they have at least one vertex.} subgraph $H$ of $G$, 
either $\varphi$ uses more than $q$ colors on $H$,
or there is a color that appears exactly once on $H$.
The \defin{$q$-centered chromatic number} of $G$, denoted by \defin{$\cen_q(G)$}, is the least nonnegative integer $k$ such that $G$ admits a $q$-centered coloring using $k$ colors.
The following statement is one of the simplest and most important outcomes of our work.

\begin{thm}\label{thm:Kt_minor_free-centered}
    Let $t$ be an integer with $t \geq 2$.
    There exists an integer $c$ such that
    for every $K_t$-minor-free graph $G$ and
    for every positive integer $q$,
    \[
        \cen_q(G) \leq c \cdot q^{t-1}.
    \]
\end{thm}
This improves on the work of~Pilipczuk and Siebertz~\cite{PS19} who proved that for every integer~$t$ with $t \geq 2$,
$K_t$-minor-free graphs have $q$-centered chromatic number upper bounded by a polynomial in $q$.
Contrary to \Cref{thm:Kt_minor_free-centered}, the degree of their polynomial is not explicitly given and arises from an application of the graph minor structure theorem
by Robertson and Seymour~\cite{GM16}.
On the other hand, D\k{e}bski, Micek, Schr\"{o}der, and Felsner~\cite{Dbski2021} showed that there exist $K_t$-minor-free graphs with $q$-centered chromatic number in $\Omega(q^{t-2})$.
Hence, \Cref{thm:Kt_minor_free-centered} is tight up to an $\bigO(q)$-factor. 

Let $G$ be a graph, let \defin{$\Pi(G)$} be the set of all vertex orderings of $G$, let $\sigma \in \Pi(G)$, and let $q$ be a nonnegative integer.
For all vertices $u$ and $v$ of $G$, we say that
$v$ is \defin{weakly $q$-reachable from $u$ in $(G,\sigma)$} if there exists
a path between $u$ and $v$ in $G$
containing at most $q$ edges 
such that for every vertex $w$ on the path, $v\leq_{\sigma} w$.
Let \defin{$\WReach_q[G, \sigma, u]$} be the set of vertices that are weakly $q$-reachable from $u$ in $(G,\sigma)$. 
The \defin{$q$th weak coloring number} of $G$ is defined as
\[
    \defin{\text{$\wcol_q(G)$}} = \min_{\sigma \in \Pi(G)}\max_{u \in V(G)}\ |\WReach_q[G, \sigma, u]|.
\]

In order to show the robustness of our framework, we include another family of parameters, introduced by Dvořák and Sereni~\cite{DS_2020}, also connected to concepts of sparsity, 
see the work of Dvořák~\cite{Dvorak_2016}. 
The \defin{treedepth} of a graph $G$, denoted by \defin{$\td(G)$}, is defined recursively as $0$ if $G$ is the null graph, $\min\{\td(G-v) + 1\mid v \in V(G)\}$ if $G$ is connected, and $\max\{ \td(C)\mid C \textrm{ component of }G\}$ otherwise.
For a nonempty set $\calY$, we say that $\lambda \colon \calY \rightarrow [0,1]$ is a \defin{probability distribution} on $\calY$ if $\sum_{Y \in \calY}\lambda(Y) = 1$.
Let $G$ be a graph and let $q$ be a positive integer.
The \defin{$q$th fractional treedepth fragility rate} of $G$ is the minimum positive integer $k$ such that there exists a family $\calY$ of subsets of $V(G)$ such that
\begin{enumerate}
    \item $\td(G - Y) \leq k$ for every $Y \in \calY$ and
    \item there exists a probability distribution $\lambda$ on $\calY$ such that for every $u \in V(G)$, we have $\sum_{u \in Y \in \calY}\lambda(Y) \leq \frac{1}{q}$.
\end{enumerate}
We denote the $q$th fractional treedepth fragility rate of $G$ by \defin{$\frate_q(G)$}.

The state of the art for weak coloring numbers and fractional treedepth fragility rates already includes results analogous to~\Cref{thm:Kt_minor_free-centered}. 
Namely, Van den Heuvel, Ossona de Mendez, Quiroz, Rabinovich, and Siebertz~\cite{vdHetal17} showed that for every integer $t$ with $t \geq 2$, there exists an integer $c$ such that for every $K_t$-minor-free graph $G$, and for every positive integer $q$, we have 
\[
\wcol_q(G)\leq c\cdot q^{t-1}.
\]

Dvořák and Sereni~\cite[Theorem~28]{DS_2020} gave bounds on the fractional treedepth fragility rates for graphs of bounded treewidth.
This, together with an idea of Esperet and Norin~\cite[Theorem~8]{Esperet-Norin} of applying a result of DeVos, Ding, Oporowski, Sanders, Reed, Seymour, and Vertigan~\cite{DDOSRSV04}; and a product structure theorem by Illingworth, Scott, and Wood~\cite[Theorem~4]{ISW22} implies
that for every integer $t$ with $t \geq 2$, there exists an integer $c$ such that for every $K_t$-minor-free graph $G$ and for every positive integer $q$, we have 
\[
    \frate_q(G) \leq c \cdot q^{t-1}.
\]
Since this statement was never proved explicitly in the literature, we chose to include its proof.
It also serves as an illustration of our methods.
See~\Cref{thm:Kt_minor_free-fragility} in~\Cref{sec:kt-minor-free-ftdfr-proof}.

Centered colorings and weak coloring numbers are crucial tools in designing parameterized algorithms in 
classes of graphs of bounded expansion.
For example,
Pilipczuk and Siebertz~\cite{PS19} showed that if $\calC$ is a class of graphs excluding a fixed minor,
then it can be decided whether a given $q$-vertex graph $H$ is a subgraph of an $n$-vertex
graph $G$ in $\calC$ in time $2^{\Oh(q\log q)}\cdot n^{\Oh(1)}$ and space $n^{\Oh(1)}$.
This algorithm relies on the fact that the union of any $q$ color classes in a $q$-centered coloring induces a subgraph of treedepth at most $q$.
Therefore, finding a $q$-centered coloring that uses  $q^{\Oh(1)}$ colors allows us to reduce the problem to graphs of bounded treedepth, 
on which the subgraph isomorphism problem can be solved efficiently. 
The running times of algorithms based on $q$-centered colorings 
heavily depend on the number of colors used.
Even more algorithmic problems were solved using the weak coloring numbers characterization of sparse graphs.
Dvo\v{r}\'ak~\cite{Dvorak13} showed a constant-factor approximation for distance versions of domination number and independence number,
with further applications in fixed-parameter algorithms and kernelization by Eickmeyer, Giannopoulou, Kreutzer, Kwon, Pilipczuk, Rabinovich, and Siebertz~\cite{EGKKPRS17}.
Grohe, Kreutzer, and Siebertz~\cite{GroheKreutzerSiebertz17} proved that deciding first-order properties is fixed-parameter tractable in nowhere dense graph classes. 
Reidl and Sullivan~\cite{ReidlSullivan23} presented an algorithm counting the number of occurrences of a fixed induced subgraph in sparse graphs. 
The time complexities of all these algorithms depend heavily on the asymptotics of $\wcol_q$ in respective graph classes.

The growth rates of centered chromatic numbers and weak coloring numbers have been extensively studied.
Grohe, Kreutzer, Rabinovich, Siebertz, and Stavropoulos~\cite{Grohe15} proved that
if $\tw(G)\leq t$, then $\wcol_q(G)\leq \binom{q+t}{t}$.\footnote{For a graph $G$, let \defin{$\tw(G)$}, \defin{$\pw(G)$}, and \defin{$\vc(G)$} stand for the treewidth, pathwidth, and vertex cover number of $G$, respectively.}
This is tight as for all nonnegative 
integers $q,t$ they constructed a graph $G_{q,t}$ 
with $\tw(G_{q,t}) = t$ and $\wcol_q(G_{q,t}) = \binom{q+t}{t}$.
Similarly, Pilipczuk and Siebertz~\cite{PS19} proved that if $\tw(G)\leq t$, then $\cen_q(G)\leq \binom{q+t}{t}$ which is again tight as proved by Dębski, Micek, Schr\"{o}der, and Felsner~\cite{Dbski2021}.
Also, Dvo\v{r}\'ak and Sereni~\cite{DS_2020} proved that if $\tw(G)\leq t$, then $\frate_q(G) = \Oh(q^t)$, which is also tight.
In the class of planar graphs, Dębski, Felsner, Micek, and Schr\"{o}der~\cite{Dbski2021} proved that $\cen_q(G) = \Oh(q^3 \log q)$; 
van den Heuvel, Ossona de Mendez, Quiroz, Rabinovich, and Siebertz~\cite{vdHetal17} proved that $\wcol_q(G) = \Oh(q^3)$; and
Dvo\v{r}\'ak and Sereni~\cite{DS_2020} proved that $\frate_q(G) = \Oh(q^3 \log q)$, while the best known lower bounds are in $\Omega(q^2\log q)$~\cite{Dbski2021,JM22,DS_2020}.
Interestingly, all these lower bound constructions in the planar case have bounded treewidth.

More generally, fix a nonnull graph $X$ and let $\param\in\set{\cen,\wcol,\frate}$. 
What is the growth rate, with respect to $q$, of the maximum of $\param_q(G)$ over all $X$-minor-free graphs $G$? 
Most previous work has focused on the case $\param = \wcol$.
Van den Heuvel et al.~\cite{vdHetal17} showed that $\wcol_q(G)=\Oh\left(q^{|V(X)|-1}\right)$. 
Subsequently, van den Heuvel and Wood~\cite{vandenHeuvel2018} proved that $\wcol_q(G)=\Oh\left(q^{\vc(X)+1}\right)$.
Dujmović, Hickingbotham, Hodor, Joret, La, Micek, Morin, Rambaud, and Wood~\cite{DHHJLMMRW24} proved that there exists an exponential function $g$ such that $\wcol_q(G)=\Oh\left(q^{g(\td(X))}\right)$. 
For $\param = \cen$, we only have 
$\cen_q(G)=\Oh(q^{h(|V(X)|)})$ for some large and non-explicit function $h$ as proved by Pilipczuk and Siebertz~\cite{PS19},
while for $\param = \frate$, we have $\frate_q(G)=\Oh(q^{|V(X)|-1})$ as follows from the aforementioned results~\cite{DS_2020,Esperet-Norin,ISW22}.

All these works 
can be viewed as attempts to understand the following graph parameters. 
For a given nonnull graph $X$ and $\param\in\set{\cen,\wcol,\frate}$, let
\begin{align*}
    \defin{\text{$f_{\param}(X)$}} = \inf\{\alpha \in \mathbb{R} \mid &\text{ there exists } c > 0 \text{ such that for every $X$-minor-free graph $G$}\\
    &\text{ and for every positive integer $q$, we have } \param_q(G) \leq c \cdot q^\alpha\}.
\end{align*}
The question is whether $f_{\param}$ is tied to\footnote{Two graph parameters $p,r$ are said to be \defin{tied} if there are two functions $\alpha,\beta$ such that $p(G) \leq \alpha(r(G))$ and $r(G) \leq \beta(p(G))$ for every graph $G$.} 
another well-established graph parameter. 
In other words, which property of $X$ governs the growth rate of $\param_q(G)$ for $X$-minor-free graphs $G$? 
Recall that for every graph $X$,
\[\tw(X) \leq \pw(X) \leq \td(X) - 1 \leq \vc(X) \leq |V(X)| - 1.\]
The aforementioned results imply that $\tw(X)-1 \leq f_{\wcol}(X)\leq g(\td(X))$.
However, 
$f_{\wcol}$ is not tied to any of these parameters. 
Indeed, neither pathwidth nor treedepth can lower bound $f_{\wcol}$.
For every positive integer $k$, let \defin{$T_k$} be a complete ternary tree of vertex-height $k$.
By Robertson-Seymour Excluded Tree Minor Theorem~\cite{GM1}, there is a constant depending on $k$ bounding pathwidth of $T_k$-minor-free graphs.
Also, we have $\wcol_q(G) = \Oh(\pw(G)^2 \cdot \log q)$ (we provide a proof in~\Cref{ssec:pathwidth} for completeness, see~\Cref{lemma:wcol_path_log}),
thus, $f_{\wcol}(T_k) =0$ while $\pw(T_k) = k$ and $\td(T_k) = k+1$.
Next, we argue that neither treewidth nor pathwidth can upper-bound $f_{\wcol}$.
For every positive integer $k$, let \defin{$L_k$} be a ladder with $k$ rungs. 
There is a graph \defin{$G_{q,t}$} (constructed by Grohe et al.~\cite{Grohe15}) such that $\wcol_q(G_{q,t})=\Omega(q^t)$, and if $k = \Omega(2^t)$, then $G_{q,t}$ excludes $L_k$ as a minor.
Therefore, $f_{\wcol}(L_k) = 2^{\Omega(k)}$, and $\tw(L_k) \leq \pw(L_k) \leq 2$.

Surprisingly, a key parameter in this context turns out to be \defin{$2$-treedepth}
as defined by Huynh, Joret, Micek, Seweryn, and Wollan~\cite{HJMSW22}.
They used this parameter to characterize classes of graphs excluding a fixed ladder as a minor.
In this paper, we prove that $f_\param$ is tied to $2$-treedepth for each $\param \in \{\cen,\wcol,\frate\}$.
More precisely, we determine $f_\param$ up to $\pm 1$ by introducing
two `rooted' versions of $2$-treedepth called $\rtd_2(\cdot)$ and $\srtd_2(\cdot)$.
They both admit treedepth-like recursive definitions that differ only in the base cases.
We now present these definitions.

Let $G$ be a graph.
A \defin{forest decomposition}
of $G$ is a pair $\mathcal{W} = \big(F,(W_x \mid x \in V(F))\big)$,
where $F$ is a forest and $W_x \subseteq V(G)$ for every $x \in V(F)$, satisfying the following conditions:
\begin{enumerate}
    \item for every $u \in V(G)$, the graph $F[\{x \in V(F) \mid u \in W_x\}]$ is a connected subgraph of $F$ and
    \item for every edge $uv \in E(G)$, there exists $x \in V(F)$ such that $u,v \in W_x$.
\end{enumerate}
We call $\calW$ a \defin{tree decomposition} when $F$ is a tree.
The sets $W_x$ are called the \defin{bags} of $\mathcal{W}$.
The sets $W_x \cap W_y$ for $xy \in E(F)$ are called the \defin{adhesions} of $\mathcal{W}$ and the \defin{adhesion} of $\calW$ is the maximum size of an adhesion of $\calW$. 
The \defin{width} of $\mathcal{W}$ is $\max_{x \in V(F)} |W_x|-1$,
and the \defin{treewidth} of $G$ is the minimum width of a tree decomposition of $G$.

A \defin{linear forest} is a disjoint union of paths.
A \defin{rooted forest} is a forest where every component is a rooted tree.
When $F$ is a rooted forest, for every $x \in V(F)$ that is not a root, let \defin{$\parent(F,x)$} be the parent of $x$ in $F$.
A forest decomposition $\big(F,(W_x\mid x\in V(F))\big)$ is \defin{rooted} if $F$ is a rooted forest.

\begin{figure}[tp]
    \centering
    \includegraphics{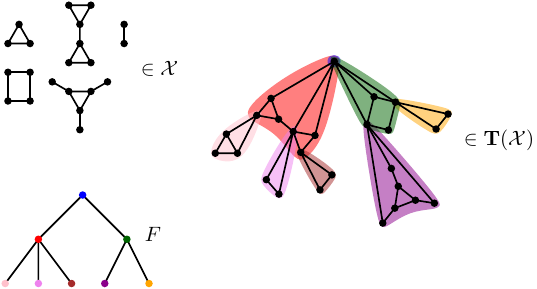} 
    \caption{An example of a graph in $\Tree(\mathcal{X})$ and a rooted tree decomposition indexed by $F$ witnessing this fact.
        The colored sets represent the bags of the tree decomposition.}
    \label{fig:T-operator}
\end{figure}

Let $\mathcal{X}$ be a class of graphs.
We define \defin{$\Tree(\mathcal{X})$} as the class of all the graphs $G$
such that there is a rooted forest decomposition $\big(F ,(W_x \mid x \in V(F))\big)$ of $G$
of adhesion at most $1$
such that for every $x \in V(F)$,
the graph $G[W_x \setminus W_{\parent(F,x)}] \in \mathcal{X}$ if $x$ is not a root,
and $|W_x| \leq 1$ if $x$ is a root. 
See Figure~\ref{fig:T-operator}.

Observe that $\calX \subset \Tree(\calX)$. 
The operator $\Tree$ is monotone in the sense that for all classes of graphs $\calX$ and $\calY$ with $\calX \subset \calY$, we have $\Tree(\calX) \subset \Tree(\calY)$.
We define recursively, for every nonnegative integer $t$, the classes
\begin{align*}
    \defin{\text{$\Rt_t$}} &= 
        \begin{cases}
            \text{only the null graph} & \textrm{if $t=0$,} \\
            \text{all edgeless graphs} & \textrm{if $t=1$,} \\
            \text{all forests} & \textrm{if $t=2$,} \\
            \Tree(\Rt_{t-1}) & \textrm{if $t\geq 3$,}
        \end{cases}
\ \ \text{ and } \ \
    \defin{\text{$\SRt_t$}} =
        \begin{cases}
            \text{only the null graph} & \textrm{if $t=0$,} \\
            \text{all edgeless graphs} & \textrm{if $t=1$,} \\
            \text{all linear forests} & \textrm{if $t=2$,} \\
            \Tree(\SRt_{t-1}) & \textrm{if $t \geq 3$.}
        \end{cases}
\end{align*}

\begin{figure}[tp]
    \centering
    \usetikzlibrary{arrows.meta}
    \vspace{1cm}
    \begin{tikzpicture}
        \node at (-2,0) {$\SRt_1=\Rt_1$};
        \node (E) at (2,0) {\includegraphics{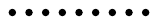}};

        \node at (-2,-2) {$\SRt_2$};
        \node (S2) at (2,-2) {\includegraphics{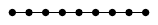}};

        \node at (-2,-4) {$\Rt_2$};
        \node (R2) at (2,-4) {\includegraphics{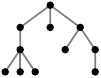}};
        
        \node at (-2,-6) {$\SRt_3$};
        \node (S3) at (2,-6) {\includegraphics{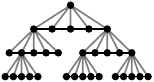}};

        \node at (-2,-8) {$\Rt_3$};
        \node (R3) at (2,-8) {\includegraphics{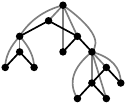}};

        \draw[bend left=67.5, |-{Computer Modern Rightarrow[length=1.25mm]}] (4,-0.25) to node[midway, right] {$\Tree$}(4,-3.75); 
        \draw[bend left=67.5, |-{Computer Modern Rightarrow[length=1.25mm]}] (4,-2.25) to node[midway, right] {$\Tree$}(4,-5.75);
        \draw[bend left=67.5, |-{Computer Modern Rightarrow[length=1.25mm]}] (4,-4.25) to node[midway, right] {$\Tree$}(4,-7.75);
    \end{tikzpicture}
    \caption{Examples of graphs in $\SRt_1=\Rt_1\subseteq \SRt_2\subseteq\Rt_2 \subseteq \SRt_3\subseteq\Rt_3$.}
    \label{fig:intro:small_values_of_rtd2_srtd2}
\end{figure}

See~\Cref{fig:intro:small_values_of_rtd2_srtd2}.
In the definition of $\Rt_t$, the base cases for $t \in \{1,2\}$ are redundant. 
Applying $\Tree$ to the class consisting only of the null graph yields the class of all edgeless graphs,
and applying $\Tree$ to the class of all edgeless graphs yields the class of all forests. 
However, we state these base cases explicitly to facilitate comparison with $\SRt_t$.

We are ready to define two key parameters, namely, \defin{rooted $2$-treedepth} and \defin{simple rooted $2$-treedepth}.
For every graph $G$, let\footnote{We denote by \defin{$\nonnegint$} the set of nonnegative integers and by \defin{$\posint$} the set of positive integers.} 
\begin{align*}
    \defin{\text{$\rtd_2(G)$}} &= \min \{ t \in \nonnegint \mid G \in \Rt_t\}, \\
    \defin{\text{$\srtd_2(G)$}} &= \min \{ t \in \nonnegint \mid G \in \SRt_t\}.
\end{align*}
It is easy to check that $\SRt_3$ contains all forests, hence, $\Rt_2 \subset \SRt_3$.
Thus, since $\Tree$ is a monotone operator, we obtain that for every nonnegative integer $t$,
\[\Rt_{t+1}\supseteq \SRt_{t+1} \supseteq \Rt_{t}.\]
Therefore, for every graph $G$, we have
\[\rtd_2(G) \leq \srtd_2(G) \leq \rtd_2(G)+1.\]
Both of these parameters are linearly tied to $2$-treedepth defined by Huynh et al.~\cite{HJMSW22}; see \Cref{lemma:rtd_2-and-td_2} in~\Cref{sec:td2_and_rtd2}.

We are now ready to present the main contributions of the paper. 
We give one statement for each of the three families of parameters. However, all proofs follow from a single abstract framework. 

\begin{theorem}\label{thm:centered}
For every integer $t$ with $t \geq 3$, for every graph $X$, there exists an integer $c$ such that for every $X$-minor-free graph $G$, for all integers $q$ with $q\geq 2$, if $\srtd_2(X)\leq t$, then
\begin{align*}
    \cen_q(G) &\leq c  \cdot q^{t-1},\\
    \cen_q(G) &\leq c  \cdot (\tw(G)+1) \cdot q^{t-2};\\
\intertext{and if $\rtd_2(X)\leq t$, then}
    \cen_q(G) &\leq c  \cdot q^{t-1} \log q,\\
    \cen_q(G) &\leq c  \cdot (\tw(G)+1) \cdot q^{t-2} \log q.\\
\end{align*}
\end{theorem}

\begin{theorem}\label{thm:wcol}
For every integer $t$ with $t \geq 2$, for every graph $X$, there exists an integer $c$ such that for every $X$-minor-free graph $G$, for all integers $q$ with $q\geq 2$, if $\srtd_2(X)\leq t$, then
\begin{align*}
    \wcol_q(G) &\leq c  \cdot q^{t-1},\\
    \wcol_q(G) &\leq c  \cdot (\tw(G)+1) \cdot q^{t-2};\\
    \intertext{and if $\rtd_2(X)\leq t$, then}
    \wcol_q(G) &\leq c  \cdot q^{t-1} \log q,\\
    \wcol_q(G) &\leq c  \cdot (\tw(G)+1) \cdot q^{t-2} \log q.\\
\end{align*} 
\end{theorem}
\begin{theorem}\label{thm:fragility}
For every integer $t$ with $t \geq 2$, for every graph $X$, there exists an integer $c$ such that for every $X$-minor-free graph $G$, for all integers $q$ with $q\geq 2$, if $\srtd_2(X)\leq t$, then
\begin{align*}
    \frate_q(G) &\leq c  \cdot q^{t-1},\\
    \frate_q(G) &\leq c  \cdot (\tw(G)+1) \cdot q^{t-2};\\
    \intertext{and if $\rtd_2(X)\leq t$, then}
    \frate_q(G) &\leq c  \cdot q^{t-1} \log q,\\
    \frate_q(G) &\leq c  \cdot (\tw(G)+1) \cdot q^{t-2} \log q.\\
\end{align*} 
\end{theorem}

Let us briefly discuss the tightness of the bounds in~\Cref{thm:centered,thm:wcol,thm:fragility}. 
Fix an integer $t$ and a graph $X$. 
All bounds involving $\tw(G)$ are tight up to a constant factor, while the remaining bounds are tight up to an $\mathcal{O}(q)$-factor.
There are two setups for the lower bound constructions: 
the first for $X$-minor-free graphs when $\rtd_2(X) = t$, and
the second for $X$-minor-free graphs when $\srtd_2(X) = t$.
Both of them were already present in the literature, though not stated in terms of our new parameters.
When $\rtd_2(X) = t$, this relates to a construction of graphs with bounded treewidth and large weak coloring numbers by Grohe et al.~\cite{Grohe15}.
When $\srtd_2(X) = t$, the lower bound relies on a construction of graphs with bounded simple treewidth and large centered chromatic numbers by Dębski et al.~\cite{Dbski2021}. 
We discuss the lower bounds in detail in~\Cref{thm:all_the_lower_bounds} and~\Cref{sec:lower-bounds}.

The next statement summarizes the growth rate of $\param_q$ for each $\param \in \{\wcol,\cen,\frate\}$ in proper minor-closed graph classes.

\begin{corollary}
\label{cor:general}
    Let $\param \in \{\cen, \wcol, \frate\}$,
    let $\mathcal{X}$ be a nonempty family of nonnull graphs, and
    let $\calC$ be the class of all graphs $G$ such that $X$ is not a minor of $G$ for all $X\in\calX$. 
    Let $t=\min\set{\rtd_2(X)\mid X\in\calX}$ and $s=\min\set{\srtd_2(X)\mid X\in\calX}$. 

    If $s\leq 2$, then 
    \[
        \max_{G\in\mathcal{C}} \param_q(G) = \Theta(1).
    \]
    If $(t,s)=(2,3)$, then 
    \[
        \max_{G\in\mathcal{C}} \param_q(G) = 
        \begin{cases}
            \Theta(q) & \textrm{if $\param=\cen$}, \\
            \Theta(\log q) & \textrm{if $\param \in \{\wcol, \frate\}$}.
        \end{cases}
    \]
    If $3\leq t=s$, then
    \[
    \max_{G\in\calC} \param_q(G) = \Omega(q^{t-2}) \ \text{ and } \ \max_{G\in\calC} \param_q(G) = \Oh(q^{t-1}).
    \]
    If $3\leq t<s$, then
    \[
    \max_{G\in\calC} \param_q(G) = \Omega(q^{t-2}\log q) \ \text{ and } \ \max_{G\in\calC} \param_q(G) = \Oh(q^{t-1}\log q).
    \]
    Moreover, if $\calX$ contains a planar graph, then the lower bounds are tight.
\end{corollary}

In particular, for every nonnull graph $X$ and each $\param\in\set{\cen,\wcol,\frate}$, 
the value of $f_{\param}$ is tied to $\rtd_2(X)$ as follows:
\begin{align*}
    f_{\param}(X) &\in \{ \rtd_2(X)-2, \rtd_2(X)-1\} &&\text{and}\\
    f_{\param}(X) &= \rtd_2(X)-2 &&\text{when $X$ is planar.}
\end{align*}
Several interesting corollaries follow from the general statement of~\cref{cor:general}. 
We present some of them below. 
See also~\Cref{table:centered,table:wcol}.

First of all, for every graph $X$ we have 
$\rtd_2(X)\leq \td(X)$, see~\Cref{lemma:rtd2_leq_td}. 
Thus,~\Cref{thm:wcol,thm:centered,thm:fragility} imply the following.
\begin{corollary}\label{cor:excluding_small_treedepth}
For every $\param\in\set{\wcol,\cen,\frate}$, 
for every integer $t$ with $t \geq 3$, 
for every graph $X$ with $\td(X)\leq t$, there exists an integer $c$ such that for every graph $G$,
if $G$ is $X$-minor-free, then for every integer $q$ with $q\geq2$,
\begin{align*}
\param_q(G) &\leq c\cdot q^{t-1}\log q, \\
\param_q(G) &\leq c\cdot (\tw(G)+1)\cdot q^{t-2} \log q.
\end{align*}
\end{corollary}
Since $\td(X)\leq\vc(X)+1$, 
this improves the result by van den Heuvel and Wood~\cite{vandenHeuvel2018} that $\wcol_q(G)=\Oh\left(q^{\vc(X)+1}\right)$ 
for all $X$-minor-free graphs $G$.
It also directly improves the result by Dujmović et al.~\cite{DHHJLMMRW24} that there exists an exponential function $g$ such that $\wcol_q(G)=\Oh\left(q^{g(\td(X))}\right)$ for all $X$-minor-free graphs $G$. 
For $\param\in\set{\cen,\frate}$, no bounds in terms of $\vc(X)$ or $\td(X)$ were known. 

The following corollary can be deduced from \Cref{cor:excluding_small_treedepth} and the fact that $\td(K_{s,t}) \leq s+1$ for all positive integers $s,t$.
\begin{corollary}\label{cor:excluding_Kst}
For every $\param\in\set{\wcol,\cen,\frate}$, 
for every positive integers $s,t$ with $s \geq 2$, 
there exists an integer $c$ such that for every graph $G$,
if $G$ is $K_{s,t}$-minor-free, then for every integer $q$ with $q\geq2$,
\begin{align*}
\param_q(G) &\leq c\cdot q^{s}\log q, \\
\param_q(G) &\leq c\cdot (\tw(G)+1)\cdot q^{s-1} \log q.
\end{align*}
\end{corollary}
Again, this improves the result by
van den Heuvel and Wood~\cite{vandenHeuvel2018} that $\wcol_q(G)=\Oh\left(q^{s+1}\right)$ for all $K_{s,t}$-minor-free graphs $G$.
No bounds like this were known for $\param\in\set{\cen,\frate}$.

Since for every nonnegative integer $g$, graphs of genus $g$ are $K_{3,2g+3}$-minor-free, we deduce the following.
\begin{cor}\label{cor:genus}
    For every $\param\in\set{\wcol,\cen,\frate}$,
    for every nonnegative integer $g$, 
    there exists an integer $c$ such that 
    for every graph $G$ of Euler genus at most $g$,
    for every integer $q$ with $q\geq 2$,
    \begin{align*}
        \param_q(G) &\leq c \cdot (\tw(G)+1) \cdot q^{2} \log q.
    \end{align*}  
\end{cor}
%It remains a challenge to nail down the asymptotics of $\param_q(G)$ when $G$ is planar. 
Determining the precise asymptotics of $\param_q(G)$ remains a challenge.
All we know is that it is in 
$\Omega(q^2\log q)$ and $\Oh(q^3\log q)$, and in case 
$(\param_q \mid q \in \posint) = (\wcol_q \mid q \in \posint)$ it is in $\Oh(q^3)$.
See Tables~\ref{table:centered} and \ref{table:wcol} again. 
\Cref{cor:genus} closes the gap when we have the additional assumption that $G$ has bounded treewidth. 

\begin{table}[tp]
    \centering
    \setlength{\tabcolsep}{1.25ex}
    \def\arraystretch{1.5}
    \begin{tabularx}{0.9\textwidth}{X@{\hspace{1.5ex}}c@{\hspace{1.5ex}}cc@{\hspace{1.5ex}}l}
        \toprule
        \bf Class $\mathcal{C}$ & {\bf lower bound}~\cite{Dbski2021} & \multicolumn{2}{c}{\bf upper bound} \\
        \midrule
        planar  & $\Omega(q^2\log q)$ & $\Oh(q^3 \log q)$ & \cite{Dbski2021} \\
        \midrule
        planar and $\tw \leq k$ & $\Omega(q^2 \log q)$ & $\Oh(q^2 \log q)$  &  \Cref{cor:genus} \\
        \midrule
        Euler genus $\leq g$ & $\Omega(q^2\log q)$ & $\Oh(q^3 \log q)$ &   \cite{Dbski2021} \\
        \midrule
        Euler genus $\leq g$ and $\tw \leq k$ & $\Omega(q^2 \log q)$ & $\Oh(q^2 \log q)$ & \Cref{cor:genus} \\
        \midrule
        outerplanar  & $\Omega(q \log q)$ & $\Oh(q \log q)$ & \cite{Dbski2021}  \\
        \midrule
        $K_{2,t}$-minor-free & $\Omega(q \log q)$ & $\Oh(q \log q)$ & \Cref{cor:excluding_Kst}  \\
        \midrule
        $\tw \leq k$ & $\binom{q+k}{k}$ & $\binom{q+k}{k}$ & \cite{PS19} \\
        \midrule
        $K_t$-minor-free & $\Omega(q^{t-2})$ & $\Oh(q^{t-1})$ & \Cref{thm:Kt_minor_free-centered} \\
        \midrule
        $K_{s,t}$-minor-free & $\Omega(q^{s-1} \log q)$ & $\Oh(q^{s} \log q)$ & \Cref{cor:excluding_Kst} \\
        \midrule
        $K_{s,t}$-minor-free and $\tw \leq k$ & $\Omega(q^{s-1} \log q)$ & $\Oh(q^{s-1} \log q)$ & \Cref{cor:excluding_Kst}\\
        \bottomrule
    \end{tabularx}

    \vspace{1mm}
    
    \caption{ 
    Lower and upper bounds on $\max_{G \in \mathcal{C}} \cen_q(G)$ for some minor-closed graph classes $\mathcal{C}$. 
    The variables $g,k,s,t$ are fixed positive integers with $s+3 \leq t \leq k$. 
    }
    \label{table:centered}
\end{table}

\begin{table}[tp]
    \centering
    \setlength{\tabcolsep}{1.25ex}
    \def\arraystretch{1.5}
    \begin{tabularx}{0.9\textwidth}{X@{\hspace{1.5ex}}c@{\hspace{1.5ex}}cc@{\hspace{1.5ex}}l}
        \toprule
        \bf Class $\mathcal{C}$ & {\bf lower bound}~\cite{JM22,Grohe15} & \multicolumn{2}{c}{\bf upper bound} \\
        \midrule
        planar  & $\Omega(q^2\log q)$ &  $\Oh(q^3)$ & \cite{vdHetal17} \\
        \midrule
        planar and $\tw \leq k$ & $\Omega(q^2 \log q)$ & $\Oh(q^2 \log q)$  &  \Cref{cor:genus} \\
        \midrule
        Euler genus $\leq g$ & $\Omega(q^2\log q)$ & $\Oh(q^3)$ & \cite{vdHetal17} \\
        \midrule
        Euler genus $\leq g$ and $\tw \leq k$ & $\Omega(q^2 \log q)$ & $\Oh(q^2 \log q)$ & \Cref{cor:genus} \\
        \midrule
        outerplanar  & $\Omega(q \log q)$ & $\Oh(q \log q)$ & \cite{JM22}  \\
        \midrule
        $K_{2,t}$-minor-free & $\Omega(q \log q)$ & $\Oh(q \log q)$ & \Cref{cor:excluding_Kst} \\
        \midrule
        $\tw \leq k$ & $\binom{q+k}{k}$ & $\binom{q+k}{k}$ & \cite{Grohe15} \\
        \midrule
        $K_t$-minor-free & $\Omega(q^{t-2})$ & $\Oh(q^{t-1})$ & \cite{vdHetal17} \\
        \midrule
        $K_{s,t}$-minor-free & $\Omega(q^{s-1} \log q)$ & $\Oh(q^{s} \log q)$ & \Cref{cor:excluding_Kst} \\
        \midrule
        $K_{s,t}$-minor-free and $\tw \leq k$ & $\Omega(q^{s-1} \log q)$ & $\Oh(q^{s-1} \log q)$ & \Cref{cor:excluding_Kst} \\
        \bottomrule
    \end{tabularx}
    
    \vspace{1mm}
    
    \caption{ 
    Lower and upper bounds on $\max_{G \in \mathcal{C}} \wcol_q(G)$ for some minor-closed graphs classes $\mathcal{C}$. 
    The variables $g,k,s,t$ are fixed positive integers with $s+3 \leq t \leq k$.
    }
    \label{table:wcol}
\end{table}

\subsubsection*{Organization of the paper}
In~\Cref{sec:cor}, we prove~\Cref{cor:general,cor:excluding_small_treedepth,cor:excluding_Kst}.
In \Cref{sec:preliminaries}, we introduce some notation and tools that will be used throughout the paper.
In \Cref{sec:outline-abstraction}, we introduce the key concepts behind our proofs, 
namely rich models, families of focused parameters, the coloring elimination property, and $(\param, \mathcal{X})$-bounding functions.
In \Cref{sec:nice_families}, we show that the families of parameters $(\cen_q \mid q \in \posint)$, $(\wcol_q \mid q \in \posint)$, and $(\frate_q \mid q \in \posint)$
can be described as families of focused graph parameters which satisfy a few key properties and fit into our framework as a result. 
In \Cref{sec:coloring-elimination}, we show that the classes $\Rt_t$ and $\SRt_t$ have the coloring elimination property for every nonnegative integer $t$. 
In \Cref{sec:proof_abstract_theorems}, we prove our main abstract theorems, which imply that passing from excluding a graph in $\Rt_t$ (resp. $\SRt_t$)
to excluding a graph in $\Rt_{t+1}$ (resp. $\SRt_{t+1}$) multiplies the growth rates of the considered families of graph parameters by $\bigO(q)$.
Next, we prove base cases, that is, we show bounds when a graph from $\SRt_2$, $\Rt_2$, or $\Rt_3$ is excluded (depending on the considered family of parameters).
In \Cref{sec:common_base}, we provide general techniques to prove these base cases.
In \Cref{sec:fragility,sec:wcol,sec:centered}, we apply this to each of the three considered families of parameters,
and we conclude the proofs of \Cref{thm:fragility,thm:wcol,thm:centered}. 
In~\Cref{sec:lower-bounds}, we discuss the lower bounds in detail.
In~\Cref{sec:td2_and_rtd2}, we give some basic facts on $2$-treedepth and its variants.
In~\Cref{ssec:pathwidth}, we recall a proof of an upper bound on weak coloring numbers of graphs with bounded pathwidth. 

\section{Proofs of \texorpdfstring{\Cref{cor:general,cor:excluding_small_treedepth}}{Corollaries 5 and 6}}
\label{sec:cor}

In this section, 
given our main theorems, i.e.\ \Cref{thm:centered,thm:wcol,thm:fragility}, 
and the lower bounds, see~\cref{thm:all_the_lower_bounds}, 
we wrap up the main corollaries from the introduction of the paper.
As mentioned before, the lower bound constructions already appeared in the literature, though without the context of the new parameters: $\srtd_2(\cdot)$ and $\rtd_2(\cdot)$.
For completeness, a uniform proof of \Cref{thm:all_the_lower_bounds} is given in \Cref{sec:lower-bounds}.

\begin{theorem}[\cite{Grohe15,JM22,Dbski2021,DS_2020}]\label{thm:all_the_lower_bounds}
    Let $\param \in \{\wcol,\cen,\frate\}$.
    For all positive integers $t,s$ with $s \geq 2$,
    \begin{align*}
        \max_{G \in \Rt_t} \param_q(G) &= \Omega(q^{t-1}), \\
        \max_{G \in \SRt_s} \param_q(G) &= \Omega(q^{s-2} \log q).
    \end{align*}
\end{theorem}

First, we show the following folklore observation.

\begin{obs}
    Let $\param\in\set{\wcol,\cen,\frate}$ 
    and let $q$ be a positive integer. 
    For every graph $G$,
    \[
        \param_q(G) \leq \td(G).
    \]
\end{obs}

Let $G$ be a graph.
If we take $\mathcal{Y}=\{\emptyset\}$ in the definition of $\frate_q(G)$, we immediately get $\frate_q(G)\leq \td(G)$ for all positive integers $q$.  
To argue the inequality for $\param\in\set{\wcol,\cen}$, 
it is convenient to work with the following definition of treedepth. 

Let $F$ be a rooted forest.
The \defin{vertex-height} of $F$ is the maximum number of vertices on a path from a root to a leaf in $F$, 
and the \defin{height} of a vertex $u$ in $F$ is the number of vertices in the path between $u$ and the root of its component in $F$.
For two vertices $u$ and $v$ in $F$, 
we say that $u$ is a \defin{descendant} of $v$ and $v$ is an \defin{ancestor} of $u$ in $F$ if $v$ lies on the path from a root to $u$ in $F$.
The \defin{closure} of $F$ is the graph with the vertex set $V(F)$ and the edge set $\{uv\mid  \textrm{$u$ is ancestor of $v$ in $F$}\}$. 
The treedepth of $G$ can be equivalently defined as the minimum vertex-height of a rooted forest whose closure contains $G$ as a subgraph.
An \defin{elimination ordering} of a rooted forest $F$ is an ordering 
$(x_1, \dots, x_{|V(F)|})$ of $V(F)$ such that for all $i,j\in[|V(F)|]$,  if $x_i$ is an ancestor of $x_j$ in $F$, then $i\leq j$.

Let $F$ be a rooted forest witnessing $\td(G)$ and let $\sigma$ be an elimination ordering of $F$. 
Let $q$ be a positive integer. 
Let $u$ be a vertex of $G$. 
Every vertex in $\WReach_q(G,\sigma,u)$  must be an ancestor of $u$ in $F$. 
Since every vertex in $G$ has at most $\td(G)$ ancestors in $F$, we conclude that  $\wcol_q(G)\leq \td(G)$.

Consider now a coloring $\phi$ of $V(G)$ so that $\phi(v)$ is the height of $v$ in $F$. Clearly, $\phi$ uses only colors in $[\td(G)]$.
For every connected subgraph $H$ of $G$, there is a vertex $v$ of $H$ such that $v$ has smaller height in $F$ than all other vertices in $H$. Therefore, $v$ has a unique color in $H$ under $\phi$. 
Thus, $\cen_q(G)\leq\td(G)$.

Another folklore observation is the following:
for every positive integer $\ell$,
if $G$ is a connected graph with no path on $\ell$ vertices,
then any DFS spanning tree has vertex-height at most $\ell-1$,
and so witnesses $\td(G) \leq \ell-1$.
Hence graphs excluding a fixed path (or more generally a fixed linear forest)
as a minor have bounded treedepth.

Now, we can prove \Cref{cor:general}.

\begin{proof}[Proof of \Cref{cor:general}]
    If $s \leq 2$, then
    $\mathcal{X}$ contains a linear forest. 
    In particular, $\mathcal{C}$ excludes a linear forest as a minor and therefore graphs in $\calC$ have treedepth bounded by a universal constant.
    It follows that $\max_{G \in \mathcal{C}} \param_q(G)$ is bounded by a constant independent of $q$.

    If $(t,s)=(2,3)$, then
    $\mathcal{X}$ contains a forest, and so by the Grid-Minor Theorem,
    $\mathcal{C}$ has bounded treewidth.
    First, suppose $\param \in \{\wcol, \frate\}$.
    Then, by \Cref{thm:wcol,thm:fragility},
    $\max_{G \in \mathcal{C}} \param_q(G) = \bigO(\log q)$.
    Moreover, since $s \geq 3$, every path is in $\mathcal{C}$.
    Paths are in $\SRt_2$ so by \Cref{thm:all_the_lower_bounds},
    $\max_{G \in \mathcal{C}} \param_q(G) = \Omega(\log q)$.
    Now, suppose $(\param_q \mid q \in \posint) = (\cen_q \mid q \in \posint)$.
    By \Cref{thm:centered} for $t=3$,
    $\max_{G \in \mathcal{C}} \cen_q(G) = \bigO(q)$.
    Moreover, $\mathcal{C}$ contains every path.
    Since paths have unbounded treedepth, and because $\cen_q(G) \geq \min\{\td(G),q+1\}$ for every graph $G$ and every positive integer $q$,
    we deduce that $\max_{G \in \mathcal{C}} \cen_q(G) \geq q+1$.

    If $3 \leq t = s$,
    then, by \Cref{thm:centered,thm:wcol,thm:fragility},
    $\max_{G \in \mathcal{C}} \param_q(G) = \bigO(q^{t-1})$.
    If $\mathcal{X}$ contains a planar graph, then
    $\mathcal{C}$ has bounded treewidth by the Grid-Minor Theorem,
    and so $\max_{G \in \mathcal{C}} \param_q(G) = \bigO(q^{t-2})$.
    For the lower bound,
    by the definition of $t$ and $s$, we have $\Rt_{t-1} \subseteq \mathcal{C}$.
    By \Cref{thm:all_the_lower_bounds}, 
    $\max_{G \in \mathcal{C}} \param_q(G) =\Omega(q^{t-2})$.

    If $3 \leq t = s-1$,
    then, by \Cref{thm:centered,thm:wcol,thm:fragility},
    $\max_{G \in \mathcal{C}} \param_q(G) = \bigO(q^{t-1} \log q)$.
    If moreover $\mathcal{X}$ contains a planar graph, then
    $\mathcal{C}$ has bounded treewidth by the Grid-Minor Theorem,
    and so $\max_{G \in \mathcal{C}} \param_q(G) = \bigO(q^{t-2} \log q)$.
    For the lower bound,
    by the definition of $t$ and $s$, we have $\SRt_{t} \subseteq \mathcal{C}$.
    By \Cref{thm:all_the_lower_bounds}, 
    $\max_{G \in \mathcal{C}} \param_q(G) =\Omega(q^{t-2}\log q)$.
\end{proof}

\Cref{cor:excluding_small_treedepth} follows from \Cref{thm:wcol,thm:centered,thm:fragility} and the following lemma.
For a class of graphs $\calX$, let \defin{$\Apex(\calX)$} be the class of all graphs $G$ such that $G \in \calX$ or $G - u \in \calX$ for some $u \in V(G)$. 

\begin{lemma}\label{lemma:rtd2_leq_td}
    For every graph $X$,
    \[
        \rtd_2(X) \leq \td(X).
    \]
\end{lemma}

\begin{proof}
    We prove by induction on $t$ that graphs of treedepth at most $t$ are in $\Rt_t$.
    This is true for $t=0$.
    Now suppose that $t \geq 1$, and let $X$ be a graph of treedepth at most $t$. 
    Let $C$ be a component of $X$. By the definition of treedepth, and by the induction hypothesis, $C \in \Apex(\{H \mid \td(H) \leq t-1\}) \subseteq \Apex(\Rt_{t-1})$.
    Since $\Apex(\Rt_{t-1}) \subseteq \Tree(\Rt_{t-1})$, 
    we deduce that every component of $X$ belongs to $\Tree(\Rt_{t-1})$,
    and since $\Tree(\Rt_{t-1}) = \Rt_t$ is closed under disjoint union, 
    it follows that $X \in \Rt_t$.
\end{proof}

\section{Preliminaries} \label{sec:preliminaries}
For each positive integer $k$, we write $\defin{\text{$[k]$}}=\{1,\ldots,k\}$ and $\defin{\text{$[0]$}} = \emptyset$.
All graphs considered in this paper are finite, simple, and undirected, unless explicitly stated otherwise.
The \defin{null graph} is the graph with no vertices.  
Connected graphs are nonnull.
A tree is defined as a connected forest, thus, trees and subtrees are also assumed to be nonnull.
A graph $G$ is \defin{edgeless} if it has no edges. 
We denote by \defin{$\edgeless$} the class of all the edgeless graphs.
For a positive integer $n$,
we denote by \defin{$K_n$} the \defin{complete graph} on the vertex set $[n]$; and we denote by \defin{$\overline{K_n}$} the edgeless graph on the vertex set $[n]$.

Let $S$ be a finite set and 
let $\sigma = (x_1, \dots, x_{|S|})$ be an ordering of $S$. 
We write $\defin{\text{$\min_{\sigma}S$}} = x_1$; $\defin{\text{$\max_{\sigma} S$}} = x_{|S|}$;  and $x_i \leq_\sigma x_j$ if and only if $i \leq j$, for all $i,j \in [|S|]$.
A collection $\mathcal{P}$ of subsets of a nonempty set $S$ is a \defin{partition} of $S$ 
if elements of $\mathcal{P}$ are nonempty, pairwise disjoint, and $\bigcup \mathcal{P} = S$.

Let $F$ be a rooted forest.
For a subtree $T$ of $F$, we define \defin{$\root(T)$} as the unique vertex in $T$ closest to any root in $F$.
For a vertex $x$ in $F$, 
we denote by \defin{$\subtree{F}{x}$} the subtree induced by all the descendants of $x$ in $F$.

Let $G$ be a graph.
For all $S,T \subset V(G)$, an \defin{\stpath{S}{T}} is a path in $G$ that is either 
a one-vertex path with a vertex in $S\cap T$, or 
a path with one endpoint in $S$ and the other endpoint in $T$ such that no internal vertices are in $S \cup T$. 
Sometimes we say that such a path is a path \defin{between} $S$ and $T$.
When $S$ and $T$ are singletons, say $S=\{u\}$ and $T=\{v\}$, we simplify the notation to \defin{\stpath{u}{v}}.

The \defin{distance} between two vertices $u$ and $v$ of $G$, denoted by \defin{$\dist_G(u,v)$}, is equal to $\infty$ if there is no \stpath{u}{v} in $G$ and otherwise, it is equal to the minimum number of edges in a \stpath{u}{v}.
The \defin{neighborhood} of a vertex $u$ in a graph $G$, denoted by \defin{$N_G(u)$}, is the set $\{v \in V(G) \mid uv \in E(G)\}$.
For every set $S$ of vertices of a graph $G$, let $\defin{\text{$N_G(S)$}}=\bigcup_{u \in S} N_G(u) \setminus S$.
A \defin{leaf addition} of a graph $G$ is a graph that is obtained from $G$ by adding one new vertex adjacent to at most one vertex.

Let $G_1, G_2$ be two graphs.
We denote by \defin{$G_1 \sqcup G_2$} the disjoint union of $G_1$ and $G_2$,
and by \defin{$G_1 \oplus G_2$} the graph obtained from $G_1 \sqcup G_2$ by adding all possible edges with one endpoint in $V(G_1)$ and the other in $V(G_2)$.
For every positive integer $k$, for every graph $G$, we write \defin{$k \cdot G$} for the union of $k$ disjoint copies of $G$.

\begin{obs}\label{obs:classes-closed}
    For every positive integer $t$, the classes $\Rt_t$ and $\SRt_t$ are closed under disjoint union.
    For every integer $t$, if $t \geq 2$, then $\Rt_t$ is closed under leaf addition, and if $t \geq 3$, then $\SRt_t$ is closed under leaf addition.
\end{obs}

For a graph $G$, a \defin{layering} of $G$ is a family $\calL=(L_i \mid i \in \NN)$ of pairwise disjoint subsets of $V(G)$ such that $\bigcup_{i \in \NN} L_i = V(G)$, and for every edge $uv$ of $G$, there exists $i \in \NN$ such that $u,v \in L_{i} \cup L_{i+1}$. 
A \defin{tree partition} of a graph $G$ is a pair $(T,\mathcal{P})$, where $T$ is a rooted tree and $\mathcal{P} = (P_{x} \mid x \in V(T))$ is a partition of $V(G)$ such that for every edge $uv$ in $G$ either there is $x \in V(T)$ with $u,v \in P_x$ or there is an edge $xy$ in $T$ with $u \in P_x$ and $v \in P_y$.
A \defin{tree partition} of $(G,S)$, where $G$ is a graph and $S \subset V(G)$ is a tree partition $\big(T,(P_{x} \mid x \in V(T))\big)$ of $G[S]$ such that for every component $C$ of $G - S$, there exists an edge $xy$ in $T$ such that $N_G(V(C)) \subset P_x \cup P_y$.
A \defin{path partition} of $G$ (resp. $(G,S)$) is a tree partition $\big(T,(P_{x} \mid x \in V(T))\big)$ of $G$ (resp. $(G,S)$)
where $T$ is a path.

Given a graph $G$ and a partition $\mathcal{P}$ of $V(G)$, the \defin{quotient graph} $G/\mathcal{P}$ is the graph with the vertex set $\mathcal{P}$ 
and two distinct $P,P' \in \mathcal{P}$ 
are adjacent in $G/\mathcal{P}$ if there are $u \in P$ and $u' \in P'$ such that $uu'$ is an edge in $G$.

For all graphs $A,B$, let $\defin{\text{$A \cup B$}} = (V(A) \cup V(B), E(A) \cup E(B))$ and let $\defin{\text{$A \cap B$}} = (V(A) \cap V(B), E(A) \cap E(B))$.

For a graph $H$, a \defin{model} of $H$ in a graph $G$ is a family $(B_x \mid x \in V(H))$ of pairwise disjoint subsets of $V(G)$ such that
\begin{enumerate}
    \item $G[B_x]$ is connected for every $x \in V(H)$, and
    \item there is an edge between $B_x$ and $B_y$ in $G$ for every $xy \in E(H)$.
\end{enumerate}
If $G$ has a model of $H$, then we say that $H$ is a \defin{minor} of $G$.
Otherwise, we say that $G$ is \defin{$H$-minor-free}.
Let $S \subseteq V(G)$.
A model $(B_x \mid x\in V(H))$ of $H$ in $G$ is \defin{$S$-rooted} if $B_x\cap S\neq\emptyset$ for each $x \in V(H)$.

Let $\calW = \big(T,(W_x\mid x\in V(T))\big)$ be a tree decomposition of a graph $G$.
For every subgraph $H$ of $G$, we denote by \defin{$\mathcal{W}\vert_H$} the tree decomposition $\big(T,(W_x \cap V(H) \mid x \in V(T))\big)$ of $H$.
An \defin{elimination ordering} of $\calW$ is an ordering $(u_1, \dots, u_{|V(G)|})$ of $V(G)$ such that for every $i \in [|V(G)|]$,
there exists $x \in V(T)$ such that 
\[\bigcup \{W_z \mid z \in V(T)\text{ and } u_i \in W_z \} \cap \{u_j \mid j \in [i-1]\} \subseteq W_x.\footnote{Equivalently, $(u_1, \dots, u_{|V(G)|})$ is an elimination ordering of $\mathcal{W}$ if and only if it is a perfect elimination ordering of the chordal graph obtained from $G$ by adding all possible edges between vertices in a same bag $W_x$ for every $x \in V(T)$.}\]

A tree decomposition
$\big(T,(W_x\mid x\in V(T))\big)$ of a graph $G$ is \defin{natural} if for every edge $e$ in $T$, 
for each component $T_0$ of $T \setminus e$, the graph $G\left[\bigcup\{W_z\mid z\in V(T_0) \}\right]\ \textrm{is connected.}$ 
The following lemma was first stated by Fraigniaud and Nisse~\cite{FN06}, see also~\cite{GJNW23}.
Later in the paper, we state and prove a variant of it: \Cref{lemma:making_a_td_natural}.
\begin{lemma}[{\cite[Theorem~1]{FN06}}]\label{lemma:natural_tree_decomposition}
    Let $G$ be a connected graph and let $\big(T,(W_x\mid x\in V(T))\big)$ be a tree decomposition of $G$. 
    There exists a natural tree decomposition $\big(T',(W'_{x'}\mid x'\in V(T'))\big)$ of $G$ such that 
    for every $x'\in V(T')$ there is $x\in V(T)$ with 
    $W'_{x'}\subseteq W_x$.
\end{lemma}

Graphs of bounded treewidth admit the following Helly-type property that is in some sense a basis of our techniques.

\begin{lemma}[{\cite[Statement (8.7)]{GM5}}]\label{lemma:helly_property_tree_decomposition}
    For every graph $G$, for every tree decomposition $\mathcal{W}$ of $G$, for every family $\mathcal{F}$ of connected subgraphs of $G$, for every nonnegative integer $d$, either
    \begin{enumerate}[label={\normalfont(\arabic*)}]
        \item there are $d+1$ pairwise vertex-disjoint subgraphs in $\mathcal{F}$, or
        \item there is a set $Z$ that is the union of at most $d$ bags of $\mathcal{W}$ such that $Z \cap V(F) \neq \emptyset$ for every $F \in \mathcal{F}$. \label{item:helly_property_tree_decomposition:hit}
    \end{enumerate}
\end{lemma}

Let $T$ be a tree.
For every $x,y \in V(T)$ with $xy \in E(T)$, we denote by \defin{$T_{x|y}$} the component of $x$ in $T \setminus \{xy\}$.
Observe that $T \setminus \{xy\}$ is the disjoint union of $T_{x|y}$ and $T_{y|x}$.
Suppose now that $T$ is rooted, and let $u$ and $v$ be two (not necessarily distinct) vertices of $T$.
The \defin{lowest common ancestor} of $u$ and $v$ in $T$, denoted by \defin{$\lca(T,u,v)$}, is the furthest vertex from the root that has $u$ and $v$ as descendants.
Let $Y \subset V(T)$.
The \defin{lowest common ancestor closure} of $Y$ in $T$ is the set $\defin{\text{$\LCA(T,Y)$}}=\{\lca(T,u,v)\mid u,v\in Y\}$.  
Observe that $\LCA(T,\LCA(T,Y)) = \LCA(T,Y)$.
The following lemma is folklore. 

\begin{lemma}[{\cite[Lemma~8]{DHHJLMMRW24}}]\label{lemma:increase_X_in_a_tree}
    Let $d$ be a positive integer,
    let $T$ be a rooted tree,
    and let $Y$ be a set of $d$ vertices of~$T$.
    If $X=\LCA(T,Y)$,
    then $|X|\leq 2d-1$, and for every component $C$ of $T-X$, $|N_T(V(C))| \leq 2$.
\end{lemma}

A similar lemma holds for tree decompositions.

\begin{lemma}[{\cite[Lemma~8]{DHHJLMMRW24}}]\label{lemma:increase_X_to_have_small_interfaces}
    Let $d$ be a positive integer.
    Let $G$ be a graph and let $\mathcal{W}=\big(T,(W_x\mid x\in V(T))\big)$ 
    be a rooted tree decomposition of $G$.
    Let $X$ be a set of $d$ vertices of $T$.
    Then $Z=\bigcup_{x\in\LCA(T,X)}W_x$ is the union of at most $2d-1$ bags of $\mathcal{W}$ 
    such that for every component $C$ of $G-Z$, $N_G(V(C))$
    is a subset of the union of at most two bags of $\mathcal{W}$.
    Moreover, if $\mathcal{W}$ is natural, then $N_G(V(C))$ intersects at most two components of $G-V(C)$.
\end{lemma}

We conclude this section by introducing the notion of
layered tree decompositions, and layered Robertson-Seymour decompositions.
A \defin{layered tree decomposition} of a graph $G$ is a pair $(\mathcal{W},\calL)$
where $\mathcal{W}= \big(T,(W_x \mid x \in V(T))\big)$ is a tree decomposition of $G$,
and $\calL=(L_i \mid i \in \NN)$
is a \defin{layering} of $G$.
The \defin{width} of $(\calW,\calL)$ is $\max |W_x \cap L_i|$ over all $x\in V(T)$ and $i \in \NN$.
The \defin{layered treewidth} of $G$, denoted by \defin{$\ltw(G)$}, is the minimum width of a layered tree decomposition of~$G$.

By definition, for every graph $G$, we have $\ltw(G) \leq \tw(G)$.
A typical example of a family of graphs of bounded layered treewidth but unbounded treewidth is the family of planar graphs.
Indeed, building upon ideas by Eppstein~\cite{Eppstein1999},
Dujmović, Morin, and Wood~\cite{Dujmovi2017} proved that planar graphs
have layered treewidth at
most three.

\begin{figure}[tp]
  \begin{center}
    \includegraphics{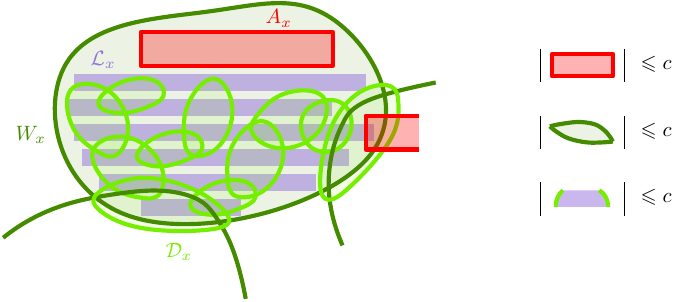}
  \end{center}
  \caption{
    A bag $W_x$ of $\calW$ in a layered RS-decomposition $(\calT, \calW, \calA, \calD, \calL)$ of width at most $c$. 
    The set $A_x$ (in red) is included in $W_x$. 
    The graph $\torso_{G,\mathcal{W}}(W_x)-A_x$ has a layering $\calL_x$ (in purple) and a tree decomposition $\calD_x$ (in green).
    Note that for every $y \in V(T) \setminus \{x\}$, $(W_x \cap W_y) \setminus A_x$ is a clique in $\torso_{G,\mathcal{W}}(W_x) - A_x$, 
    and so, it is contained in a single bag of $\calD_x$ and it at most two layers of $\calL_x$. 
  }
  \label{fig:lrs-outline}
\end{figure}

For a graph $G$,
a tree decomposition $\mathcal{W} = \big(T,(W_x \mid x \in V(T))\big)$,
and $x \in V(T)$, the \defin{torso} of $W_x$ in $G$, and $\mathcal{W}$, denoted by \defin{$\torso_{G,\mathcal{W}}(W_x)$}, 
is the graph with the vertex set $W_x$ and where two distinct vertices $u,v$ are adjacent
if $uv \in E(G)$, or there exists $y \in N_T(x)$ such that $u,v \in W_x \cap W_y$.
For every graph $G$ and every positive integer~$c$, a \defin{layered Robertson-Seymour decomposition} of $G$ (\defin{layered RS-decomposition} for short) of width at most $c$ 
is a tuple 
    \[(T,\calW,\calA,\calD,\calL)\]
where $T$ is a tree and
\begin{enumerate}[label=(lrs\arabic*)]
    \item $\calW = \big(T,(W_x \mid x \in V(T))\big)$ is a tree decomposition of $G$ of adhesion at most $c$; \label{LRS:adhesion}
    \item $\calA = \big(A_x \mid x \in V(T)\big)$ where $A_x \subset W_x$ and $|A_x| \leq c$ for every $x \in V(T)$; \label{LRS:apices}
    \item $\calD = \big(\calD_x \mid x \in V(T)\big)$ where $\calD_x = \big(T_x,(D_{x,z} \mid z \in V(T_x))\big)$ is a tree decomposition of $\torso_{G,\mathcal{W}}(W_x)-A_x$ for every $x \in V(T)$; \label{LRS:tree_decompositions}
    \item $\calL = \big(\calL_x \mid x \in V(T)\big)$ where $\calL_x = (L_{x,i}\mid i \in \mathbb{N})$ is a layering of $\torso_{G,\mathcal{W}}(W_x)-A_x$ for every $x\in V(T)$; and \label{LRS:layering}
    \item $|D_{x,z} \cap L_{x,i}| \leq c$ for all $x \in V(T)$, $z \in V(T_x)$, and $i \in \NN$.  \label{LRS:ltw}
\end{enumerate}
See~\Cref{fig:lrs-outline}.
Dujmović, Morin, and Wood~\cite{Dujmovi2017} proved that for every fixed positive integer $t$, $K_t$-minor-free graphs admits such decompositions with bounded width.

\begin{thm}[{\cite[Theorem~22 and Lemma~26]{Dujmovi2017}}]\label{theorem:Kt_free_product_structure_decomposition}
    For every positive integer $t$, there is a positive integer~$\cLRS(t)$ 
    such that every $K_t$-minor-free graph admits a layered RS-decomposition of width at most~$\cLRS(t)$. 
\end{thm}

We finish the preliminaries section with a simple observation.

\begin{obs}\label{obs:tw_leq_max_tw_of_a_torso}
    Let $G$ be a graph.
    For every tree decomposition $\big(T,(W_x \mid x \in V(T))\big)$ of $G$,
    we have
    \[
        \tw(G) \leq \max \{ \tw(\torso_{G,\mathcal{W}}(W_x)) \mid x \in V(T)\}.
    \]
\end{obs}

\section{The abstract framework}\label{sec:outline-abstraction}

The proofs of our main theorems (\Cref{thm:centered,thm:wcol,thm:fragility}) are by induction on rooted $2$-treedepth or simple rooted $2$-treedepth.
We present an abstract framework which will encapsulate the induction step for the three considered families of parameters: $(\cen_q \mid q \in \posint)$, $(\wcol_q \mid q \in \posint)$, and $(\frate_q \mid q \in \posint)$.

\subsection{Rich models} \label{sec:rich}
Let $G$ and $H$ be two graphs.
Let $\mathcal{F}$ be a family of connected subgraphs of $G$.
A model $(B_x\mid x\in V(H))$ of $H$ in $G$ is said to be \defin{$\mathcal{F}$-rich} if for every $x \in V(H)$, 
there exists $F \in \mathcal{F}$ such that $F$ is a subgraph of $G[B_x]$.

\begin{figure}[tp]
    \centering 
    \includegraphics[scale=1]{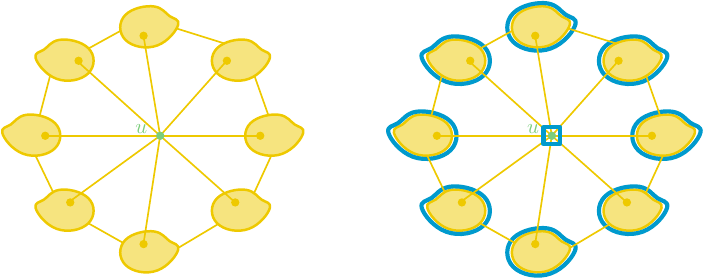} 
    \caption{
        On the left-hand side, we depict an $\mathcal{F}$-rich model of $X$, where $X$ is a cycle graph on $8$ vertices,
        and $\mathcal{F}$ is the family of all connected subgraphs of $G - \{u\}$ containing a neighbor of $u$ in $G$.
        On the right-hand side, we show how to construct, given an $\mathcal{F}$-rich model of $X$, a model of $K_1 \oplus X$.
    }
    \label{fig:F-rich-models}
\end{figure}

For example, if $H$ has $k$ vertices and no edges, then $G$ has an $\mathcal{F}$-rich model of $H$ if and only if 
$G$ contains $k$ pairwise disjoint members of $\mathcal{F}$.
An extreme case is when $\mathcal{F}$ contains every one-vertex subgraph of $G$.
Then every model of $H$ in $G$ is $\mathcal{F}$-rich.
See another example in~\Cref{fig:F-rich-models}.

It is convenient to define a notation to restrict families of connected subgraphs to a subgraph:
if $\mathcal{F}$ is a family of connected 
subgraphs of a graph $G$,
and if $H$ is a subgraph of $G$,
we denote by \defin{$\mathcal{F}\vert_{H}$} the family $\{F \in \mathcal{F} \mid F \subseteq H\}$.

\subsection{Families of focused parameters}

Let \defin{$\calZ$} be the family of all pairs $(G,S)$ where $G$ is a graph and $S \subset V(G)$.
We call a family
\[\param = (\param_{q} \colon \calZ \rightarrow [0,\infty) \mid q \in \posint )\]
a \defin{family of focused parameters}.\footnote{With a slight abuse of notation, we write $\param_{q}(G,S)$ instead of $\param_{q}((G,S))$.}

The following definition encapsulates the properties that a family of focused parameters has to satisfy so that the abstract induction step can be performed.

We say that a family of focused parameters is \defin{nice} if there exist positive integers $b$ and $b'$ such that
for every positive integer $q$, for every $(G,S) \in \calZ$, we have
\begin{enumerate}[label={\normalfont (n\arabic*)}]
    \item $\param_{q}(G,S) \leq \param_{q}(C,S \cap V(C))$ for some component $C$ of $G$;\label{item:nice:components} 
    \item $\param_{q}(G,S) \leq |S|$; \label{item:nice:at_most_size_of_S}
    \item $\param_{q}(G,S_0 \cup S_1) \leq \param_{q}(G,S_0) + \param_{q}(G-S_0,S_1)$ for all disjoint $S_0,S_1 \subseteq V(G)$; and \label{item:nice:subadditivity}
    \item for every tree partition $\big(T,(P_x \mid x \in V(T))\big)$ of $(G,S)$, \label{item:nice:tree-partitions}
        \[
            \param_{q}(G,S) \leq b'(q+1) \cdot \max_{x \in V(T)} \param_{bq} (G_x,P_x)
        \]
        where for every $x \in V(T)$, $G_x$ is the subgraph of $G$ induced by the union of $U_x = \bigcup \{P_z \mid z \in V(\subtree{T}{x})\}$ and all the vertex sets of the components of $G-S$ having a neighbor in~$U_x$.
\end{enumerate}
In~\Cref{sec:nice_families}, we define focused versions of $\cen_q(\cdot)$, $\wcol_q(\cdot)$, and $\frate_q(\cdot)$. 
We show that these focused versions satisfy \ref{item:nice:components}--\ref{item:nice:tree-partitions}. 
Therefore, the abstract theorems that we will state at the end of this section can be applied to prove the induction step of our main theorems.

\subsection{Coloring elimination property}

To bound the growth rate of our families of focused parameters in minor-closed graph classes, we rely on the structure of the excluded minors, namely, 
rooted $2$-treedepth and/or simple rooted $2$-treedepth.
The families with bounded $\rtd_2$ and $\srtd_2$ are built inductively in the same manner and they differ only in the base cases.
Therefore, we extract the main property of these two parameters to characterize a family $\calX$ of graphs for which the inductive step of our main theorems would apply.

Let $\mathcal{X}$ be a class of graphs.
We say that $\mathcal{X}$ has the \defin{coloring elimination property} if,
for every positive integer $k$, for every $X \in \mathcal{X}$,
there exists $Y \in \mathcal{X}$ such that
for all sets $S_1, \dots, S_k \subseteq V(Y)$ such that $\bigcup_{i\in[k]} S_i = V(Y)$, there exists $i \in [k]$ such that there is an $S_i$-rooted model of $X$ in~$Y$.
In this case, we say that $Y$ \defin{witnesses} the coloring elimination property of $\calX$ for $X$ and~$k$.

In \Cref{lemma:Rt_has_coloring_elimination_property}, in \Cref{sec:coloring-elimination}, we show that $\Rt_t$ and $\SRt_t$ have the coloring elimination property.
Here, we show that the property is preserved under \q{minor closure}.

\begin{lemma} \label{lemma:coloring_elimination_minor_closed}
    Let $\mathcal{X}$ be a class of graphs and let \[\mathcal{X}'=\{X'\mid \textrm{$X'$ is a minor of $X$ for some $X\in\calX$}\}.\]
    The class $\mathcal{X}$ has the coloring elimination property if and only if $\mathcal{X}'$ has the coloring elimination property.
\end{lemma}

\begin{proof}
    Suppose that $\calX'$ has the coloring elimination property.
    Let $k$ be a positive integer and let $X\in\calX$. 
    Since $\calX\subseteq\calX'$, 
    we can fix $Y' \in \calX'$ that witness the coloring elimination property of $\calX'$ for $X$ and $k$. 
    By the definition of $\calX'$, there is $Y\in\calX$ such that $Y$ contains a model $(B'_y\mid y\in V(Y'))$ of $Y'$.
    Fix $S_1,\dots,S_k \subset V(Y)$ whose union is $V(Y)$.
    For every $i\in [k]$, let $S'_i = \{y\in V(Y')\mid B'_y\cap S_i\neq \emptyset\}$.
    Since $\bigcup_{i\in[k]}S_i = V(Y)$, we have $\bigcup_{i\in[k]}S'_i = V(Y')$.
    Fix $j\in [k]$ and an $S'_j$-rooted model $(B_x\mid x\in V(X))$ of $X$ in $Y'$.
    For every $x\in V(X)$, let $D_x = \bigcup_{y\in B_x}B'_y$.
    By construction, $(D_x\mid x\in V(X))$ is an $S_j$-rooted model of $X$ in $Y$.
    Therefore, $\calX$ has the coloring elimination property.

    Next, suppose that $\calX$ has the coloring elimination property.
    Let $k$ be a positive integer and let $X'\in\calX'$.
    By the definition of $\calX'$, there is $X\in\calX$ such that $X$ contains a model $(B'_x\mid x\in V(X'))$ of $X'$.
    Let $Y \in \calX$ witness the coloring elimination property of $\calX$ for $X$ and $k$.
    Observe that $Y \in \calX'$.
    Fix $S_1,\dots,S_k \subset V(Y)$ whose union is $V(Y)$.
    Fix $j\in [k]$ and an $S_j$-rooted model $(B_x\mid x\in V(X))$ of $X$ in~$Y$.
    Finally, $(\bigcup_{y\in B'_x}B_y \mid x\in V(X'))$ is an $S_j$-rooted model of $X'$ in $Y$, which shows that $\calX'$ has the coloring elimination property.
\end{proof}

\subsection{The abstract induction step}\label{ssec:abstract-induction}
Let $g \colon \posint \rightarrow \posint$ be a function, let $\param = (\param_q \mid q \in \posint)$ be a family of focused parameters, let $X$ be a graph, and let $\alpha$ and $\beta$ be nonnegative integers.
We say that a graph $G$ is \defin{$(g,\param,X,\alpha,\beta)$-good}
if for every positive integer $q$, for every family $\calF$ of connected subgraphs of $G$, 
if $G$ has no $\calF$-rich model of $X$, 
then there exists $S \subset V(G)$ such that
\begin{enumerate}[label={\normalfont (g\arabic*)}]
    \item $S \cap V(F) \neq \emptyset$ for every $F \in \calF$;\label{item:par-bounding:hitting}
    \item for every component $C$ of $G-S$, $N_G(V(C))$ intersects at most $\alpha$ components of $G-V(C)$; and \label{item:par-bounding:nbrs}
    \item $\param_{q}(G,S) \leq \beta \cdot g(q)$.\label{item:par-bounding:bound}
\end{enumerate}
For a class of graphs $\calX$, we say that a function $g \colon \posint \rightarrow \posint$ is \defin{$(\param,\mathcal{X})$-bounding} if 
for every positive integer~$k$,
there exists a nonnegative integer $\alpha$ such that for every graph $X \in \mathcal{X}$,
there exists a nonnegative integer $\beta(X)$ such that
every $K_k$-minor-free graph $G$ is $(g,\param,X,\alpha,\beta(X))$-good.

In this case, when $k$ is fixed, we say that $\alpha$ and $\beta(\cdot)$ \defin{witness} $g$ being $(\param,\calX)$-bounding for $k$.
Essentially, 
we require that excluding an $\calF$-rich model of a graph $X$ in $G$ implies 
the existence of a set $S\subseteq V(G)$ hitting\footnote{A \defin{hitting set} of a family $\calF$ of connected subgraphs of a graph $G$ is a set $S \subset V(G)$ such that $S \cap V(F) \neq \emptyset$ for every $F \in \calF$.} $\calF$ such that the parameter focused on $S$ is bounded.
If for two classes of graphs $\calX$ and $\calX'$ we have $\calX \subset \calX'$, then the fact that a function $g$ is $(\param,\calX')$-bounding implies that $g$ is $(\param,\calX)$-bounding.

The following is the main abstract statement of this section.

\begin{theorem}\label{thm:abstract_induction_main}
    Let $(\param_q \mid q \in \posint)$ be a nice family of focused parameters,
    let $\mathcal{X}$ be a class of graphs having the coloring elimination property and closed under disjoint union and leaf addition.
    Let $g \colon \posint \to \posint$ be a function.
    If $g$ is $(\param,\mathcal{X})$-bounding, then there exists a positive integer $b$ such that
    $q \mapsto q \cdot g(bq)$ is $(\param,\Tree(\mathcal{X}))$-bounding.
\end{theorem}

In practice, we split this induction into the two steps as below.

\begin{theorem}\label{thm:main_X_to_A(X)}
    Let $(\param_q \mid q \in \posint)$ be a nice family of focused parameters,
    let $\mathcal{X}$ be a class of graphs closed under disjoint union and having the coloring elimination property,
    and let $g \colon \posint \to \posint$.
    If $g$ is $(\param,\mathcal{X})$-bounding, then
    there exists a positive integer $b$ such that $q \mapsto q \cdot g(bq)$
    is $(\param,\Apex(\mathcal{X}))$-bounding.
\end{theorem}

\begin{theorem}\label{thm:main_A(X)_to_T(X)}
    Let $(\param_q \mid q \in \posint)$ be a nice family of focused parameters,
    let $\mathcal{X}$ be a nonempty class of graphs closed under disjoint union and leaf addition,
    and let $g \colon \posint \to \posint$.
    If $g$ is $(\param,\Apex(\mathcal{X}))$-bounding, then
    $g$ is $(\param,\Tree(\mathcal{X}))$-bounding.
\end{theorem}

Clearly, \Cref{thm:main_X_to_A(X),thm:main_A(X)_to_T(X)} imply \Cref{thm:abstract_induction_main}.
We prove \Cref{thm:main_X_to_A(X),thm:main_A(X)_to_T(X)} in \Cref{sec:proof_abstract_theorems}.
We conclude the proofs of \Cref{thm:centered,thm:wcol,thm:fragility} by proving the base cases in \Cref{sec:common_base,sec:centered,sec:wcol,sec:fragility}.
We wrap up this section by a simple observation which will be reused several times to relate $(\param,\calX)$-bounding of certain functions with bounds such as in  \Cref{thm:centered,thm:wcol,thm:fragility}.

\begin{lemma}\label{lemma:par-bounding-to-bound}
    Let $\param = (\param_q \mid q \in \posint)$ be a family of focused parameters, let $\calX$ be a class of graphs, and let $g: \posint \rightarrow \posint$ be a function that is $(\param,\calX)$-bounding.
    Then, for every $X \in \calX$, there exists a constant $\beta(X)$ such that for every positive integer $q$ and for every $X$-minor-free graph $G$, we have 
        \[\param(G,V(G)) \leq \beta(X) \cdot g(q).\]
\end{lemma}
\begin{proof}
    Let $X \in \calX$ and $k = |V(X)|$.
    Since $g$ is $(\param,\calX)$-bounding, it follows that there exist nonnegative integers $\alpha$ and $\beta(X)$ such that every $K_k$-minor-free graphs $G$ are $(g,\param,X,\alpha,\beta(X))$-good.
    Let $G$ be an $X$-minor-free graph.
    It follows that $G$ is also $K_k$-minor-free.
    Let $\mathcal{F}$ be the family of all the one-vertex subgraphs of $G$.
    Since $G$ is $X$-minor-free,
    there is no $\mathcal{F}$-rich model of $X$ in $G$.
    Therefore,
    there exists $S \subseteq V(G)$ such that \ref{item:par-bounding:hitting}--\ref{item:par-bounding:bound} hold, and in particular,
    \begin{enumerate}
        \item $S \cap V(F) \neq \emptyset$ for every $F \in \mathcal{F}$ and \label{item:i-proof-}
        \item $\param_q(G,S) \leq \beta(X) \cdot g(q)$ for every positive integer $q$.
    \end{enumerate}
    Note that~\ref{item:i-proof-} implies that $S = V(G)$.
    Therefore, for every positive integer $q$,
    \[\param_q(G,V(G)) \leq \beta(X) \cdot g(q).\qedhere\]
\end{proof}

\section{Nice families of focused parameters}
\label{sec:nice_families}

\subsection{Weak coloring numbers}
\begin{figure}[tp]
    \centering 
    \includegraphics[scale=1]{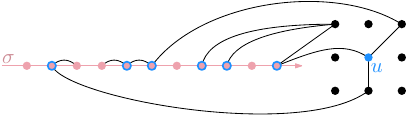} 
    \caption{The pink vertices correspond to the set $S$. The vertices in $S$ highlighted blue are in $\WReach_3[G,S,\sigma,u]$.} \label{fig:wcols}
\end{figure}

Let $G$ be a graph, let $q$ be a nonnegative integer, let $S \subseteq V(G)$, let $\sigma$ be an ordering of $S$, let $u \in V(G)$, and let $v \in S$.
We say that \defin{$v$ is weakly $q$-reachable from $u$ in $(G,S,\sigma)$}
if there is a \stpath{u}{v} $P$ in $G$ of length at most $q$ such that $\min_\sigma (V(P) \cap S) = v$.
We denote by \defin{$\WReach_q[G,S,\sigma,u]$} the set of all the weakly $q$-reachable vertices from $u$ in $(G,S,\sigma)$ and we write $\defin{\text{$\wcol_q(G,S,\sigma)$}} = \max_{u \in V(G)} |\WReach_q[G,S,\sigma,u]|$.
Finally, let $\wcol_q(G,S)$ be the minimum value of $\wcol_q(G,S,\sigma)$ among all orderings $\sigma$ of $S$.
For each of the defined objects, we drop $S$ when $S=V(G)$.
Namely, $v$ is weakly $q$-reachable from $u$ in $(G,\sigma)$ whenever $v$ is weakly $q$-reachable from $u$ in $(G,V(G),\sigma)$, $\WReach_q[G,\sigma,u] = \WReach_q[G,V(G),\sigma,u]$, $\wcol_q(G,\sigma) = \wcol_q(G,V(G),\sigma)$, and $\wcol_q(G) = \wcol_q(G,V(G))$.
This matches the definition given in~\Cref{sec:intro}.
See an illustration in \Cref{fig:wcols}.

Next, we state some simple observations concerning weak coloring numbers.
Ultimately, we will show that $(\wcol_q \mid q \in \posint)$ is a nice family of focused parameters.

\begin{obs}\label{obs:wcol_components}
    Let $G$ be a graph, let $S \subset V(G)$, and let $\mathcal{C}$ be the family of components of~$G$.
    For every positive integer $q$, we have
    \[
        \wcol_q(G,S) = \max_{C \in \mathcal{C}}\wcol_q(C,S \cap V(C)).
    \]
\end{obs}

\begin{obs}\label{obs:monotone}
    Let $G$ be a graph and let $S,U \subset V(G)$.
    For every nonnegative integer~$q$, we have
        \[\wcol_q(G - U,S \setminus U) \leq \wcol_q(G,S).\]
\end{obs}

\begin{obs}\label{obs:wcol_union}
    Let $G$ be a graph and let $S_1,S_2$ be disjoint subsets of $V(G)$.
    For every positive integer $q$, we have
    \[
        \wcol_q(G,S_1\cup S_2) \leq \wcol_q(G,S_1) + \wcol_q(G-S_1,S_2).
    \]
\end{obs}

\Cref{obs:wcol_components,obs:monotone} are clear from the definition and to see~\Cref{obs:wcol_union}, it suffices to take an ordering $\sigma_1$ of $S_1$ witnessing $\wcol_q(G,S_1)$ and an ordering $\sigma_2$ of $S_2$ witnessing $\wcol_q(G-S_1,S_2)$, and concatenate them so that all the vertices of $S_1$ lie before all the vertices of $S_2$.

\begin{figure}[tp]
    \centering 
    \includegraphics[scale=1]{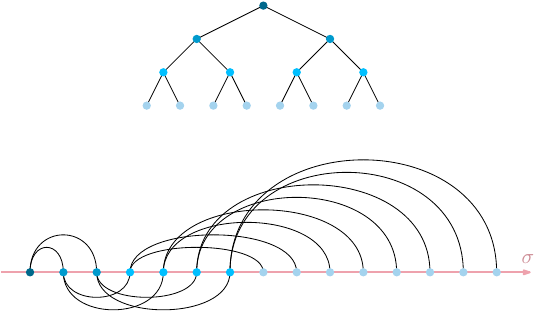} 
    \caption{An example of an elimination ordering of a complete binary tree of height $3$.} 
    \label{fig:tree}
\end{figure}

To derive an upper bound on weak coloring numbers of trees,
it suffices to root a given tree and take an elimination ordering of the vertices.
In such an ordering, for every positive integer~$q$, a vertex weakly $q$-reaches only its $q+1$ closest ancestors (including itself). 
See~\Cref{fig:tree}.

\begin{obs}\label{obs:wcol-trees}
Let $T$ be a rooted tree. 
Let $\sigma$ be an elimination ordering of $T$. 
Then
\[
\wcol_q(T,\sigma)\leq q+1.
\]
\end{obs}

We extend this idea to tree partitions of graphs.

\begin{lemma}\label{lemma:wcol_combining_tree_partitions}
    Let $G$ be a graph, let $S \subset V(G)$,
    let $\big(T,(P_x \mid x \in V(T))\big)$ be a tree partition of $(G,S)$.
    For every $x \in V(T)$, let $G_x$ be the subgraph of $G$ induced by the union of $U_x = \bigcup \{P_z \mid z \in V(\subtree{T}{x})\}$ and all the vertex sets of the components of $G-S$ having a neighbor in~$U_x$.
    For every positive integer $q$, we have
    \[
        \wcol_{q}(G,S) \leq (q+1) \cdot \max_{x \in V(T)} \wcol_{q} (G_x,P_x).
    \]
\end{lemma}

\begin{proof}
    Let $q$ be a positive integer and let $k = \max_{x \in V(T)} \wcol_{q} (G_x,P_x)$.
    Let $(x_1,\dots,x_{|V(T)|})$ be an elimination ordering of $T$.
    For every $x \in V(T)$, let $\sigma_x$ witness $\wcol_q(G_x,P_x) \leq k$.
    Let $\sigma$ be an ordering of $S = \bigcup_{x \in V(T)} P_x$ obtained by concatenating $\sigma_{x_1}, \dots, \sigma_{x_{|V(T)|}}$.
    
    It suffices to show that $\wcol_q(G,S,\sigma) \leq (q+1)k$.
    Let $u \in V(G)$ and consider $W = \WReach_q[G,S,\sigma,u]$.
    Assume that $W \neq \emptyset$.
    If $u \in S$, let $x(u) \in V(T)$ such that $u \in P_{x(u)}$.
    Otherwise, let
    $C$ be the component of $u$ in $G-S$ and 
    let $x(u)$ be the vertex of $T$ furthest from $\root(T)$ such that
    $N_G(V(C))$ intersects $P_{x(u)}$.  
    Observe that for every $x \in V(T)$, $W \cap P_x \subseteq \WReach_q[G_x,P_x,\sigma_x,u]$.
    Moreover, if $x$ is not among the closest $q+1$ ancestors of $x(u)$ (including $x(u)$ itself),
    then $W \cap P_x = \emptyset$.
    Altogether, this implies that $|W| \leq (q+1)k$ as desired.
\end{proof}

\begin{lemma}\label{lemma:wcol-nice}
    The family $(\wcol_q \mid q \in \posint)$ is nice.
\end{lemma}
\begin{proof}
    \Cref{obs:wcol_components} implies \ref{item:nice:components}; \ref{item:nice:at_most_size_of_S} is clear from the definition; \Cref{obs:wcol_union}
    implies \ref{item:nice:subadditivity}; and \ref{item:nice:tree-partitions} with $b'=b=1$
    follows from \Cref{lemma:wcol_combining_tree_partitions}.
\end{proof}

In~\Cref{thm:wcol}, we give variants of each upper bound by trading off a $q$-factor for a ($\tw(G)+1$)-factor.
In order to obtain the bound depending on treewidth, we apply our techniques to the following family of focused parameters.

\begin{lemma}\label{lemma:wcol-tw-nice}
For every graph $G$, for every $S \subseteq V(G)$, and for every positive integer $q$, let
    \[
       \param_q(G,S) = \frac{1}{\tw(G)+2} \wcol_q(G,S)
    \]    
    The family of focused parameters $(\param_q \mid q \in \posint)$ is nice.
\end{lemma}

\begin{proof}
    Since $\tw(G) \geq -1$ for every graph $G$, the parameters are well-defined.
    The family is nice by monotonicity of treewidth under taking subgraphs and because the family $(\wcol_q \mid q \in \posint)$ is nice by~\Cref{lemma:wcol-nice}.
\end{proof}

\subsection{Centered colorings}\label{ssec:centered-nice}

Let $S$ be a set. 
A \defin{coloring} of $S$ is a function $\phi\colon S \rightarrow C$ for some set $C$.
For each $u \in S$, we say that $\varphi(u)$ is the \defin{color} of $u$ 
and subsequently we say that $\varphi(S)$ is the \defin{set of colors used by $\varphi$}.
For two colorings $\phi_1$ and $\phi_2$ of $S$, we define the coloring \defin{$\phi_1 \times \phi_2$} of~$S$, called the \defin{product coloring} of $\phi_1$ and $\phi_2$, by $(\phi_1\times \phi_2)(u) = (\phi_1(u),\phi_2(u))$ for every $u \in S$.
Given $S' \subset S$, an element $u \in S'$ is a \defin{$\varphi$-center} of $S'$ if the color of $u$ is unique in $S'$ under $\varphi$, in other words, $\varphi(u) \notin \varphi(S' \setminus \{u\})$.
Let $G$ be a graph.
A \defin{coloring} of $G$ is a coloring of $V(G)$.
When $X \subset V(G)$ and $\varphi$ is a coloring of $G$, we denote by \defin{$\varphi\vert_{X}$}, the restriction of $\varphi$ to $X$.
Recall that a coloring $\varphi$ of $G$ is a $q$-centered coloring of $G$ for a positive integer $q$ if for every connected subgraph $H$ of $G$, either $|\varphi(V(H))| > q$ or $V(H)$ has a $\phi$-center.

In the case of centered colorings, we consider a focused family of parameters defined on precolored graphs.

Let $q$ be a positive integer, let $G$ be a graph, and let $\phi$ be a coloring of $G$.
For every $S \subseteq V(G)$, we define \defin{$\cen_q(G,\phi,S)$} to be the least nonnegative integer $k$ such that
there exists a coloring $\psi$ of $S$ using at most $k$ colors
such that for every connected subgraph $H$ of $G$ such that $V(H) \cap S \neq \emptyset$, one of the following is true:
\begin{enumerate}
    \item $|\phi(V(H))|>q$, or
    \item $|(\phi \times \psi)(V(H) \cap S)| > q$, or
    \item there is a $(\phi \times \psi)$-center of $V(H) \cap S$.
\end{enumerate}

The following statement is almost immediate from the definitions.

\begin{lemma}\label{lem:ordered-to-normal}
    Let $q$ be a positive integer, let $G$ be a graph, and let $\phi$ be a coloring of $G$ using $k$ colors.
    Then,
        \[\cen_q(G) \leq \cen_q(G,\phi,V(G)) \cdot k.\]
\end{lemma}
\begin{proof}
    It suffices to take the coloring $\phi \times \psi$ of $G$, where $\psi$ witnesses $\cen_q(G,\phi,V(G))$.
\end{proof}

We follow with a few simple observations on $\cen_q(\cdot,\cdot,\cdot)$.

\begin{obs}\label{obs:cen_components}
    Let $G$ be a graph, let $\varphi$ be a coloring of $G$, let $S \subset V(G)$, and let $\mathcal{C}$ be the family of components of~$G$.
    For every positive integer $q$, we have
    \[
        \cen_q(G,\varphi,S) = \max_{C \in \mathcal{C}}\cen_q(C,\varphi\vert_{V(C)},S \cap V(C)).
    \]
\end{obs}

\begin{obs}\label{obs:cen_union}
    Let $G$ be a graph, let $\varphi$ be a coloring of $G$, and let $S_1,S_2$ be disjoint subsets of $V(G)$.
    For every positive integer $q$, we have
    \[
        \cen_q(G,\varphi,S_1\cup S_2) \leq \cen_q(G,\varphi,S_1) + \cen_q(G-S_1,\varphi\vert_{V(G) \setminus S_1},S_2).
    \]
\end{obs}

\begin{obs}\label{obs:cen_at_most_S}
    Let $G$ be a graph, let $\varphi$ be a coloring of $G$, and let $S \subset V(G)$.
    \[
        \cen_q(G,\varphi,S) \leq |S|.
    \]
\end{obs}

\Cref{obs:cen_components} follows just from the definition of $\cen_q(G,\varphi,S)$, i.e.\ each connected subgraph $H$ of $G$ is a subgraph of a component of $G$.
~\Cref{obs:cen_union} holds as we can use two disjoint sets of 
$\cen_q(G,\varphi,S_1)$ and $\cen_q(G-S_1,\varphi\vert_{V(G) \setminus S_1},S_2)$ colors to color respectively $S_1$ and $S_2$.
Finally, for~\Cref{obs:cen_at_most_S} it suffices to set $\psi$ to be an injective coloring of $S$.

\begin{lemma}\label{lemma:centered_colorings_are_nice}
    Let $\mathcal{U}$ be an infinite graph,
    and for every positive integer $q$, let $\phi_q$ be a coloring of $\mathcal{U}$.
    For every positive integer $q$,
    for every graph $G$, and
    for every $S \subseteq V(G)$, let
    \[
        \param_q(G,S) =
        \begin{cases}
            \cen_q(G,\phi_q\vert_{V(G)},S) & \textrm{if $V(G) \subseteq V(\mathcal{U})$ and $E(G) \subseteq E(\mathcal{U})$,} \\
            0 & \textrm{otherwise.}
        \end{cases}
    \]
    The family of focused graph parameters $(\param_q \mid q \in \posint)$ is nice.
\end{lemma}

\begin{proof}
    Let $q$ be a positive integer
    and let $G$ be a graph.
    If it does not hold that $V(G) \subset V(\calU)$ and $E(G) \subset E(\calU)$, then $\param_q(G,S) = 0$ for every $S \subseteq V(G)$,
    and so \ref{item:nice:components}--\ref{item:nice:tree-partitions} hold trivially.
    Thus, we assume that $V(G) \subset V(\calU)$ and $E(G) \subset E(\calU)$

    Let $S \subseteq V(G)$.
    Note that for every component $C$ of $G$, we have $V(C) \subset V(G)$ and $E(C) \subset E(G)$.
    Also, $V(G - S) \subset V(G)$ and $E(G-S) \subset E(G)$.
    Therefore, \ref{item:nice:components} follows from \Cref{obs:cen_components},
    \ref{item:nice:at_most_size_of_S} follows from \Cref{obs:cen_at_most_S},
    and \ref{item:nice:subadditivity} follows from \Cref{obs:cen_union}.

    We now show \ref{item:nice:tree-partitions}.
    Let $\big(T,(P_x \mid x \in V(T))\big)$ be a tree partition of $(G,S)$ with $T$ rooted in $r \in V(T)$.
    For every $x \in V(T)$, let $G_x$ be the subgraph of $G$ induced by the union of $U_x = \bigcup \{P_z \mid z \in V(\subtree{T}{x})\}$ and all the vertex sets of the components of $G-S$ having a neighbor in~$U_x$.
    Let $k = \max\{\cen_q(G_x,\phi_q\vert_{V(G_x)}, P_x) \mid x \in V(T)\}$ and note that $V(G_x) \subset V(G)$ and $E(G_x) \subset E(G)$ for every $x \in V(T)$.
    For every $x \in V(T)$, let $\psi_x\colon P_x \to [k]$ witnesses $\cen_q(G_x,\phi_q\vert_{V(G_x)}, P_x) \leq k$.
    Let $\psi \colon S \to [k] \times \{0, \dots, q\}$ be defined by
    \[
        \psi(u) = \big(\psi_x(u), \dist_T(r,x) \bmod (q+1)\big)
    \]
    for every $x \in V(T)$ and for every $u \in P_x$.

    We claim that $\psi$ witnesses $\cen_q(G,\phi_q\vert_{V(G)},S) \leq (q+1)k$.
    Consider a connected subgraph $H$ of $G$ intersecting $S$.
    Let $x$ be the vertex of $T$ closest to the root such that $V(H) \cap P_x \neq \emptyset$.
    Observe that $H\subseteq G_x$.
    Since $\psi_x$ witnesses $\cen_q(G_x,\phi_q\vert_{V(G_x)},P_x)\leq k$, 
    one of the following is true: $|\phi_q(V(H))| >q$ or
    $|(\phi_q \times \psi_x)(V(H)\cap P_x)|>q$ or 
    there is a $(\phi_q\times\psi_x)$-center $u$ of $V(H)\cap P_x$.
    If one of the first two cases holds, the claim holds. 
    Thus, assume that the latter holds.
    If $u$ is a $(\phi_q\times\psi)$-center of $V(H) \cap S$, then again the claim holds.
    Otherwise, there exists $y \in V(T)$ distinct from $x$ such that 
    \[
        \dist_T(r,x)\equiv \dist_T(r,y)\mod (q+1),
    \]
    and $V(H) \cap P_y \neq \emptyset$.
    Since $H$ is connected, $x$ is an ancestor of $y$ by the choice of $x$.
    Altogether, this implies that $\dist_T(x,y) > q$ and $V(H)$ intersects $P_z$ for every $z$ in the \stpath{x}{y} in $T$.
    This implies that $|(\phi_q \times \psi)(V(H) \cap S)| \geq |\psi(V(H) \cap S)| > q$.
    Therefore, $\psi$ witnesses $\cen_q(G,\phi_q\vert_{V(G)},S) \leq (q+1)k$.
    This proves \ref{item:nice:tree-partitions} with $b'=b=1$,
    and concludes the proof of the lemma.
\end{proof}

\subsection{Fractional treedepth fragility rates}\label{ssec:ftdfr}
First, we define a focused version of treedepth.
Note that this is a different parameter from $\td(\cdot,\cdot)$ considered by us in another paper~\cite{quickly-excluding-apex-forest} and by Claus, Hodor, Joret, and Morin~\cite{CHMJ26}.
For every graph $G$ and $S \subset V(G)$, we define recursively
\[
    \defin{\text{$\ftd(G,S)$}} = 
    \begin{cases}
        0 & \textrm{if $S = \emptyset$,}\\
        \min_{u \in S} \ftd(G,S-\{u\}) + 1 &\textrm{if $G$ is connected, and}\\
        \max_{i \in [\ell]}\ftd(C_i, S \cap V(C_i))&\textrm{if $G$ consists of components $C_1,\dots,C_\ell$ and $\ell > 1$.}
    \end{cases}
\]
Straight from the definition, we have the following observation.

\begin{obs} \label{lemma:td_GS_components}
    Let $G$ be a graph, let $S\subseteq V(G)$, and let $\calC$ be the family of all the components of $G$.
    We have
    \[
        \ftd(G,S) = \max_{C\in\calC} \ftd(C,S\cap V(C)). 
    \]
\end{obs}

\begin{lemma} \label{lemma:td_GS_subadditivity}
    Let $G$ be graph and let $S_1$ and $S_2$ be disjoint subsets of $V(G)$.
    We have
    \[
        \ftd(G,S_1\cup S_2)\leq \ftd(G,S_1) + \ftd(G-S_1,S_2).
    \]
\end{lemma}

\begin{proof}
    We proceed by induction on $|V(G)|$.
    If $S_1 = \emptyset$, then the result is clear.
    Thus, assume that $|S_1| \geq 1$.
    If $G$ is not connected, then for every component $C$ of $G$, by the induction hypothesis,
    $\ftd(C,(S_1 \cup S_2) \cap V(C)) \leq \ftd(C,S_1\cap V(C)) + \ftd(C - (V(C) \cap S_1), S_2 \cap V(C)) \leq \ftd(G,S_1) + \ftd(G-S_1,S_2)$.
    This implies the result by \Cref{lemma:td_GS_components}.
    If $G$ is connected, then by the definition of $\ftd(G,S_1)$,
    there exists $u \in S_1$ such that $\ftd(G,S_1) = 1+\ftd(G-u,S_1 \setminus \{u\})$.
    By the induction hypothesis, we have 
    $\ftd(G-u, (S_1 \setminus \{u\}) \cup S_2) \leq \ftd(G - u,S_1 \setminus \{u\}) + \ftd(G-S_1,S_2) = \ftd(G,S_1)-1 + \ftd(G-S_1,S_2)$.
    Since $\ftd(G, S_1 \cup S_2) \leq 1 + \ftd(G-u, (S_1 \setminus \{u\}) \cup S_2)$, this implies the result.
\end{proof}

\begin{lemma} \label{lemma:td_GS_tree_partition}
    Let $G$ be a graph, let $S\subseteq V(G)$, and let $\big(T,(P_x\mid x\in V(T))\big)$ be a tree partition of $(G,S)$.
    We have
    \[
    \ftd(G,S) \leq h\cdot \max_{x\in V(T)} \ftd(G_x,P_x)
    \]
    where $h$ is the vertex-height of $T$ and for every $x\in V(T)$, $G_x$ is the subgraph of $G$ induced by the union of $U_x=\bigcup\{P_z\mid z\in V(\subtree{T}{x})\}$ and all the vertex sets of the components of $G-S$ having a neighbor in $U_x$.
\end{lemma}

\begin{proof}
    The statement is clear when $S = \emptyset$, hence, assume otherwise.
    We proceed by induction on $h$.
    The statement is clear again when $h = 1$, hence, assume $h \geq 2$.
    Let $r$ be the root of $T$.
    By possibly removing the vertex sets of all the components of $G$ disjoint from $S$,
    we suppose that every component of $G$ intersects $S$,
    and so $G=G_r$.
    Let $k = \max_{x\in V(T)} \ftd(G_x,P_x)$.
    For every child $x$ of $r$, $\big(\subtree{T}{x},(P_z\mid z\in V(T_x))\big)$ is a tree partition of $(G_x,S\cap V(G_x))$.
    Since $\subtree{T}{x}$ has vertex-height at most $h-1$, by induction, $\ftd(G_x,S\cap V(G_x))\leq (h-1)\max_{z\in V(\subtree{T}{x})}\ftd(G_z,P_z)\leq (h-1)k$.
    By \Cref{lemma:td_GS_components}, $\ftd(G-P_r,S-P_r)\leq \max\{\ftd(G_x,S\cap V(G_x))\mid x \in V(T-r), \parent(T,x)=r\} \leq (h-1)k$.
    By \Cref{lemma:td_GS_subadditivity}, $\ftd(G,S)\leq \ftd(G_r,P_r) + \ftd(G-P_r,S-P_r) \leq k + (h-1)k = h\cdot k$ since $G_r=G$.
    This concludes the proof.
\end{proof}

Let $q$ be a positive integer, $G$ be a graph, and $S \subset V(G)$.
The \defin{$q$th fractional treedepth fragility rate} of $(G,S)$,
denoted by \defin{$\frate_q(G,S)$},
is the minimum nonnegative integer $k$ such that there exists a family $\calY$ of subsets of $S$ such that
\begin{enumerate}
    \item $\ftd(G-Y,S\setminus Y) \leq k$ for every $Y \in \calY$, and
    \item there exists a probability distribution $\lambda$ on $\calY$ such that for every $u \in S$, we have $\sum_{u\in Y\in\calY}\lambda(Y) \leq \frac{1}{q}$. \label{item-frate-1-thin}
\end{enumerate}
We say that a probability distribution satisfying~\ref{item-frate-1-thin} is \defin{$q$-thin}.
In other words, the probability distribution of a random variable $\rvar{Y}$ taking subsets of $S$ as values is $q$-thin if and only if $\Pr[u \in \rvar{Y}] \leq \frac{1}{q}$ for every $u \in S$.
We will use the following abbreviation.
For a positive integer $k$, writing that \q{a random variable $\rvar{Y}$ witnesses $\frate(G,S) \leq k$}, 
we mean that $\rvar{Y}$ takes values in the set $\{Y \subseteq S \mid \ftd(G - Y, S \setminus Y) \leq k\}$ and the probability distribution of $\rvar{Y}$ is $q$-thin.
Furthermore, slightly abusing notation, we sometimes write that \q{$\rvar{Y}$ satisfies some statement} where $\rvar{Y}$ is a random variable, by which we mean that the statement is satisfied for every element that can be a value of $\rvar{Y}$. 

\begin{lemma}\label{ftdfr_are_nice}
    For every graph $G$, for every $S \subseteq V(G)$, and for every positive integer $q$, let
    \[
        \param_q(G,S) = \frac{1}{\tw(G)+2} \frate_q(G,S).
    \]
    The family of focused parameters
    $(\param_q \mid q \in \posint)$
    is nice.
\end{lemma}

\begin{proof}
    First, observe that, for every graph $G$, $1/(\tw(G)+2)$ is well-defined since $\tw(G)\geq -1$. 
    Let $G$ be a graph, let $S \subseteq V(G)$, and let $q$ be a positive integer.
    
    Let $\calC$ be the family of components of $G$.
    Let $k = \max_{C \in \calC} \frate_q(C,S\cap V(C))$.
    For every $C \in \calC$, let $\rvar{Y}_C$ witness $\frate_q(C,S\cap V(C))\leq k$.
    Let $\rvar{Y} = \bigcup_{C \in \calC} \rvar{Y}_C$.
    Using~\Cref{lemma:td_GS_components}, it is easy to verify that $\ftd(G-\rvar{Y},S-\rvar{Y})\leq k$.
    Let $u\in S$ and let $C \in \calC$ be such that $u \in V(C)$.
    Since the probability distribution of $\rvar{Y}_C$ is $q$-thin, we have
    $\Pr[u\in\rvar{Y}] = \Pr[u\in\rvar{Y}_C] \leq \frac{1}{q}$.
    Thus, $\rvar{Y}$ witnesses $\frate_q(G,S)\leq k$.
    By definition of $k$, there is $C\in\calC$ such that $\frate_q(C,S\cap V(C))=k$.
    Therefore, by monotonicity of treewidth, $\param_q(G,S) \leq \param_q(C,S \cap V(C))$, and so,~\ref{item:nice:components} holds.
    
    To show \ref{item:nice:at_most_size_of_S},
    consider a random variable $\rvar{Y}$ that is always equal to the empty set. 
    Then, since $\ftd(G,S)\leq |S|$, and because for every $u\in S$ we have $\Pr[u\in\rvar{Y}] = 0 \leq \frac{1}{q}$,
    we deduce that $\rvar{Y}$ witnesses $\frate_q(G,S) \leq |S|$, which gives~\ref{item:nice:at_most_size_of_S}.
    
    Let $S_1$ and $S_2$ be disjoint subsets of $S$ whose union is $S$.
    By monotonicity of treewidth, in order to obtain~\ref{item:nice:subadditivity}, 
    it suffices to prove that $\frate_q(G,S_1 \cup S_2) \leq \frate_q(G,S_1) + \frate_q(G-S_1,S_2)$.
    Let $\rvar{Y}_1$ and $\rvar{Y}_2$ witness $\frate_q(G,S_1) = k_1$ and $\frate_q(G-S_1,S_2) = k_2$ respectively.
    Let $\rvar{Y} = \rvar{Y}_1\cup\rvar{Y}_2$.
    By \Cref{lemma:td_GS_subadditivity},
    \[
        \ftd(G-\rvar{Y},S-\rvar{Y}) \leq \ftd(G-\rvar{Y}_1,S_1-\rvar{Y}_1) + \ftd\big((G-S_1)-\rvar{Y}_2,S_2-\rvar{Y}_2\big) \leq k_1 + k_2.
    \]
    Let $u\in S$ and assume that $u \in S_i$ for some $i \in [2]$.
    Since the probability distribution of $\rvar{Y}_i$ is $q$-thin, we have $\Pr[u\in\rvar{Y}] = \Pr[u\in\rvar{Y}_i] \leq \frac{1}{q}$.
    It follows that $\rvar{Y}$ witnesses $\frate_q(G,S) \leq k_1 + k_2$, which completes the proof of~\ref{item:nice:subadditivity}.
    
    To show \ref{item:nice:tree-partitions},
    consider a tree partition $\big(T,(P_x \mid x \in V(T))\big)$ of $(G,S)$ with $T$ rooted in $r \in V(T)$.
    For every $x \in V(T)$, let $G_x$ be the subgraph of $G$ induced by the union of $U_x = \bigcup \{P_z \mid z \in V(\subtree{T}{x})\}$ and all the vertex sets of the components of $G-S$ having a neighbor in~$U_x$.
    Let $k = \max_{x\in V(T)}\frate_{2q}(G_x,P_x)$.
    For every $x\in V(T)$, let $\rvar{Y}_x$ witness $\frate_{2q}(G_x,P_x)\leq k$.
    
    For each $i\in[2q]$, let $L_i=\{x\in V(T)\mid \dist_T(r,x) \equiv i \bmod 2q\}$.
    Let $\rvar{Y}'$ be a random variable with the uniform distribution on the set $\left\{\bigcup_{x\in L_i}P_x\mid i\in[2q]\right\}$, and let $\rvar{Y} = \rvar{Y'}\cup\bigcup_{x\in V(T)}\rvar{Y}_x$.
    For every $x\in V(T)$, let $\rvar{P}_x = P_x-\rvar{Y}$.
    Let $\rvar{F}$ be the forest obtained from $T$ by removing every vertex $x \in V(T)$ such that $\rvar{P}_x = \emptyset$.
    We root each component of $\rvar{F}$ in its vertex which is the closest to $r$ in $T$.
    Let $\rvar{C}$ be a component of $G-\rvar{Y}$, and let $\rvar{T}_\rvar{C}$ be the subtree of $\rvar{F}$ such that $V(\rvar{C})\cap \bigcup_{x\in V(\rvar{T}_\rvar{C})}\rvar{P}_x \neq \emptyset$.
    Observe that $\big(\rvar{T}_\rvar{C},(\rvar{P}_x \mid x \in V(\rvar{T}_\rvar{C}))\big)$ is a tree partition of $(\rvar{C},(S-\rvar{Y})\cap V(\rvar{C}))$.
    By the definition of $\rvar{Y}'$, all trees in $\rvar{F}$ have vertex-height at most $2q-1$.
    Therefore, by \Cref{lemma:td_GS_tree_partition}, 
    \begin{align*}
        \ftd(\rvar{C},(S-\rvar{Y})\cap V(\rvar{C}))&\leq (2q-1)\max_{x\in V(\rvar{T}_\rvar{C})}\ftd(\rvar{G}_x,\rvar{P}_x)\\
        &\leq (2q-1)\max_{x\in V(T)}\ftd(G_x-\rvar{Y}_x,P_x-\rvar{Y}_x)\leq (2q-1)k,
    \end{align*}
    where for each $x \in V(T)$, $\rvar{G}_x$ is the subgraph of $\rvar{C}$ induced by the union of $\rvar{U}_x = \bigcup \{\rvar{P}_z \mid z \in V((\rvar{T}_\rvar{C})_x)\}$ and all the vertex sets of the components of $\rvar{C}-(S \setminus \rvar{Y})$ having a neighbor in~$\rvar{U}_x$.
    
    Finally, let $u\in S$ and let $x\in V(T)$ such that $u\in P_x$.
    By the definition of $\rvar{Y}'$ and because the probability distribution of $\rvar{Y}_x$ is $q$-thin, the union bound yields
    \[
        \Pr[u\in\rvar{Y}] = \Pr[u\in\rvar{Y}'] + \Pr[u\in\rvar{Y}_x] \leq \frac{1}{2q} + \frac{1}{2q} = \frac{1}{q}.
    \]
    Therefore,~\ref{item:nice:tree-partitions} holds with $b'=b=2$.
    This concludes the proof of the lemma.
\end{proof}

\section{Coloring elimination property of \texorpdfstring{$\Rt_t$}{Rt} and \texorpdfstring{$\SRt_t$}{St}}
\label{sec:coloring-elimination}

In this section, we show that $\Rt_t$ and $\SRt_t$ have the coloring elimination property for every nonnegative integer $t$.
First, we show that the operators $\Apex$ and $\Tree$ preserve the property under some mild assumptions, see~\Cref{lemma:coloring_elimination_and_Tree}.
Therefore, to conclude, it suffices to check the base cases.

\begin{lemma}\label{lemma:coloring_elimination_and_Tree}
    Let $\mathcal{X}$ be a class of graphs closed under taking minors and under disjoint union.
    If $\mathcal{X}$ has the coloring elimination property, then
    $\Apex(\mathcal{X})$ and $\Tree(\mathcal{X})$ have the coloring elimination property.
\end{lemma}

\begin{proof}
    The statement is clear when $\mathcal{X}$ consists only of the null graph, thus, assume that $\mathcal{X}$ contains a nonnull graph.
    Let $k$ be a positive integer.

    For every $X \in \Tree(\mathcal{X})$,
    we say that a vertex $u \in V(X)$ is a \defin{root} of $X$
    if there exists a rooted forest decomposition $\big(F,(W_a \mid a \in V(F))\big)$ of $X$ of adhesion at most $1$ such that for every $a\in V(F)$ that is not a root, $G[W_a-W_{\parent(F,a)}]\in \calX$;
    for every $a$ that is a root, $|W_a|\leq 1$;
    and there exists a root $r$ of $F$ such that $W_r = \{u\}$.
    For every $X\in\Tree(\calX)$, $X$ has a root $u$ unless $X$ is the null graph because $\calX$ is closed under taking minors (in particular, for every $X \in \calX$ and $u \in V(X)$, we have $X - u \in \calX$). 
    This definition also applies to the members of $\Apex(\mathcal{X})$ since $\Apex(\mathcal{X}) \subseteq \Tree(\mathcal{X})$.
    For every $X\in \Apex(\calX)$, for every root $u$ of $X$, $X-u\in\calX$.

\begin{figure}[tp]
\begin{center}
\includegraphics{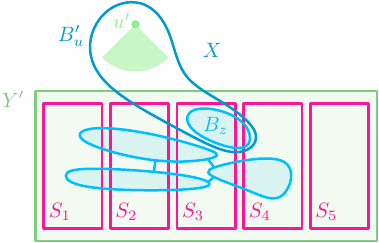}
\end{center}
\caption{
The vertex $u'$ corresponds to $K_1$ in $K_1 \oplus Y'$.
$S_1,\dots,S_5$ is a partition of $V(Y')$.
In the figure there is an $S_3$-rooted model of $K_1 \sqcup (X-u)$ in $Y'$ (light-blue).
We add $u'$ to the branch set $B_z$ where $z$ corresponds to $K_1$ in $K_1\sqcup (X-u)$, obtaining an $S_3$-rooted model of $X$ in $Y$.
}
\label{fig:coloring-property-apex}
\end{figure}

    \begin{claim}\label{claim:cep:apex}
        Let $X$ be a nonnull graph in $\Apex(\mathcal{X})$ and let $u$ be a root of $X$.
        There exists $Y \in \Apex(\mathcal{X})$ with a root $u'$ such that
        for every family of sets $S_1, \dots, S_k$ whose union is $V(Y)$,
        $Y$ contains an $S_i$-rooted model $(B_x \mid x \in V(X))$ of $X$
        for some $i \in [k]$ such that $u' \in B_u$.
    \end{claim}

    \begin{proofclaim}
        The proof is illustrated in~\Cref{fig:coloring-property-apex}.
        Let $X' = K_1 \sqcup (X-u)$.
        Since $\calX$ contains a nonnull graph and is closed under taking minors, it contains $K_1$.
        Since $u$ is a root of $X$ and $X \in \Apex(\calX)$, $X-u \in \calX$.
        In particular, $X' \in \calX$ as $\calX$ is closed under disjoint union.
        We denote by $z$ the vertex corresponding to $K_1$ in $X'$.
        Let $Y' \in \calX$ witness the coloring elimination property in $\calX$ for $X'$ and $k$.
        Furthermore, let $Y = K_1 \oplus Y'$, and let $u'$ be the vertex of $Y$ corresponding to $K_1$. 
        In particular, $Y \in \Apex(\calX)$ and $u'$ is a root of $Y$.
    
        Consider a family of sets $S_1, \dots, S_k \subset V(Y)$ whose union is $V(Y) = V(Y')\cup\{u'\}$.
        There exists $i \in [k]$ such that $Y'$ contains an $(S_i-\{u'\})$-rooted model $(B'_x \mid x \in V(X'))$ of $X'$.
        Let $B_u = B'_z \cup \{u'\}$ and $B_x = B'_x$ for every $x \in V(X') \setminus \{z\}$.
        We obtain that $(B_x \mid x \in V(X))$ is an $S_i$-rooted model of $X$ in $Y$.
        Since $u' \in B_u$, this completes the proof of the claim.
    \end{proofclaim}
    \Cref{claim:cep:apex} implies that $\Apex(\mathcal{X})$ has the coloring elimination property.
    A stronger assertion concerning the roots in~\Cref{claim:cep:apex} is needed to prove the next claim.
    \begin{claim}\label{claim:cep:tree}
        Let $X \in \Tree(\mathcal{X})$ be connected and let $u$ be a root of $X$.
        There exists $Y \in \Tree(\mathcal{X})$ and a root $u'$ of $Y$ such that
        for every family of sets $S_1, \dots, S_k$ whose union is $V(Y)$,
        there is an $S_i$-rooted model $(B_x \mid x \in V(X))$ of $X$ in $Y$ for some $i \in [k]$ such that $u' \in B_{u}$.
    \end{claim}

\begin{figure}[tp]
\begin{center}
\includegraphics{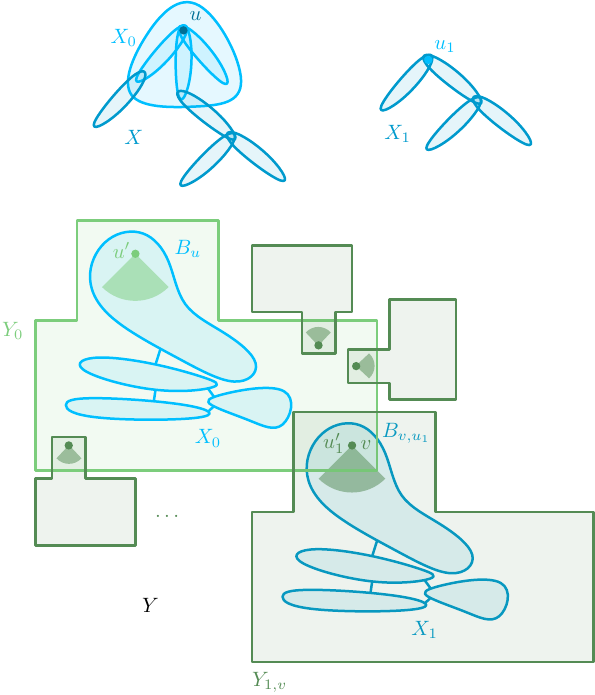}
\end{center}
\caption{
The graph $X$ is split into two graphs according to its rooted tree decomposition witnessing $X\in\Tree(\calX)$: 
$X_0\in \Apex(\calX)$ induced by all vertices of the bags containing the root $u$ of $X$
and $X_1\in \Tree(\calX)$ obtained by contracting $X_0$.
Since $X_0\in \Apex(\calX)$, \Cref{claim:cep:apex} gives $Y_0\in\Apex(\calX)$ which witnesses the coloring elimination property in $\Apex(\calX)$ for $X_0$. 
Since $X_1$ has a rooted tree decomposition witnessing $X_1\in\Tree(\calX)$ with smaller vertex-height, induction gives $Y_1\in\Tree(\calX)$ which witness the coloring elimination property in $\Tree(\calX)$ for $X_1$. 
We make a graph $Y\in\Tree(\calX)$ that witness the coloring elimination property in $\Tree(\calX)$ for $X$ by taking, for every $v \in V(Y_0)$,
a copy $Y_{1,v}$ of $Y_1$, and identifying $v$ and the vertex corresponding to a root $u'_1$ of $Y_1$.
This operation is allowed by the operator $\Tree$.
Moreover, for any family of $k$ sets covering $V(Y)$, a rooted model of $X$ can be found using a rooted model of $X_0$ in $Y_0$ and rooted models of $X_1$ in the copies $Y_{1,v}$ for each $v\in V(Y_0)$.
}
\label{fig:coloring-property-tree}
\end{figure}

    \begin{proofclaim}
        The objects defined in the proof of this claim are depicted in~\Cref{fig:coloring-property-tree}.
        There exists a rooted forest decomposition $\big(T,(W_a \mid a \in V(T))\big)$ of $X$ 
        witnessing the fact that $X \in \Tree(\mathcal{X})$ and with $\{u\}$ as a root bag.
        Without loss of generality, $W_a \neq \emptyset$ for every $a \in V(T)$.
        Since $X$ is connected, this implies that $T$ is a tree.
        Fix such a rooted forest decomposition.
        Let $r$ be the unique root of $T$.
        We proceed by induction on the vertex-height of $T$.
        When $T = K_1$, $X = K_1$ and the statement is clear.
        Now, suppose that $T \neq K_1$ and that the result holds for trees with smaller vertex-heights.

        Let $T_{0}$ be the subtree of $T$ induced by $\{a \in V(T) \mid u \in W_a\}$.
        Observe that 
        for every $a \in V(T_{0})$, $X[W_a \setminus \{u\}] \in \mathcal{X}$.
        Since $\mathcal{X}$ is closed under disjoint union,
        this implies that 
        \(
            X\left[\bigcup_{a \in V(T_{0})}\left( W_a\setminus \{u\}\right)\right] \in \calX
        \). 
        So, $X_0 = X\left[\bigcup_{a \in V(T_{0})}W_a\right]$ is in $\Apex(\mathcal{X})$ and $u$ is a root of $X_0$.
        By \Cref{claim:cep:apex},
        there exists $Y_0 \in \Apex(\mathcal{X})$ and a root $u'$ of $Y_0$ such that
        for every family of sets $S_1, \dots, S_k \subset V(Y_0)$ whose union is $V(Y_0)$,
        $Y_0$ contains an $S_i$-rooted model of $X_0$
        for some $i \in [k]$ such that the branch set corresponding to $u$ contains $u'$.
        
        Let $X_1$ be the graph obtained from $X$
        by identifying all the vertices in $V(X_0)$ into a single vertex $u_1$.
        Observe that $X_1\in\calX$ since $\calX$ is closed under taking minors.
        Let $T_1$ be the tree obtained from $T$ by identifying the vertices in $V(T_0)$ into a single vertex $r_1$.
        For every $a \in V(T_1)$, let 
        \[
            W_{1,a} = 
            \begin{cases}
                \{u_1\} &\textrm{if $a=r_1$,} \\
                W_{a} &\textrm{otherwise.}
            \end{cases}
        \]
        The rooted tree decomposition
        $\big(T_1, (W_{1,a} \mid a \in V(T_1))\big)$ of $X_1$ witnesses the fact that $X_1 \in \Tree(\mathcal{X})$.
        We claim that the vertex-height of $T_1$ is smaller than the vertex-height of $T$.
        Indeed, for every neighbor $s$ or $r$ in $T$, since $X$ is connected and $W_{s} \neq \emptyset$,
        we have $u \in W_s$, and so $s \in V(T_0)$. This proves that $N_T(r) \subseteq V(T_0)$,
        and so $T_1$ has smaller vertex-height than $T$.
        Therefore, by the induction hypothesis,
        fix a graph $Y_1 \in \Tree(\mathcal{X})$
        and a root $u'_1$ of $Y_1$ such that
        for every family of sets $S_1, \dots, S_k \subset V(Y_1)$ whose union is $V(Y_1)$,
        there is an $S_i$-rooted model of $X_1$ in $Y_1$ such that the branch set corresponding to $u_1$ contains $u_1'$.
        Since $X_1$ is connected, we assume without loss of generality that $Y_1$ is also connected.

        For every $v \in V(Y_0)$, let $Y_{1,v}$ be a copy of $Y_1$.
        We label vertices in $Y_{1,v}$ in the following way.
        The vertex corresponding to $u_1'$ is labeled as $v$ and a vertex corresponding to each $y \in V(Y_1) \setminus \{u_1'\}$ is labeled as $(y,v)$.
        Finally, let $Y$ be obtained from the disjoint union of $Y_0$ and $Y_{1,v}$ for each $v \in V(Y_0)$ by identifying vertices with the same labels.
        We claim that $Y$ witnesses the assertion of the claim.      

        First, we argue that $Y \in \Tree(\calX)$ and $u'$ is a root of $Y$.
        For each $v \in V(Y_0)$, let $\big(T'_{1,v}, (W'_{1,v,a} \mid a \in V(T'_{1,v}))\big)$ be a rooted tree decomposition of $Y_{1,v}$ witnessing $Y_1\in\Tree(\calX)$, and with $T'_{1,v}$ rooted in $r'_{1,v}$ such that $W'_{1,r'_{1,v}} = \{(u'_1,v)\}$.
        Let $T'$ be a tree obtained from a disjoint union of $T'_{1,v}$ for each $v \in V(Y_0)$ by identifying the roots of all copies into a vertex $r''$, and then adding the root $r'$ with a single neighbor $r''$.
        We define $U_{r'} = \{u'\}$, $U_{r''} = V(Y_0)$, and $U_a = W'_{1,v,a}$ for every $a \in V(T') \setminus \{r',r''\}$ where $a \in V(T'_{1,v})$.
        Since $Y_0\in \Apex(\calX)$, $Y_0-u'\in \calX$.
        Therefore, $\big(T', (U_{a} \mid a \in V(T'))\big)$ is a rooted tree decomposition of $Y$ witnessing that $Y \in \Tree(\calX)$.
        Additionally, $U_{r'} = \{u'\}$, hence, $u'$ is a root of $Y$.

        Let $S_1, \dots, S_k$ be a family of sets whose union is $V(Y)$.
        Let $v \in V(Y_0)$.
        Since $Y_{1,v}$ is a copy of $Y_1$ and $v$ is a root of $Y_{1,v}$,
        there exists $i(v) \in [k]$ and an $(S_{i(v)} \cap V(Y_{1,v}))$-rooted model $\calM_v = (B_{v,x} \mid x \in V(X_1))$ of 
        $X_1$ in $Y_{1,v}$ with $v \in B_{v,u_1}$.

        For each $j \in [k]$, let $S_j' = \{v \in V(Y_0) \mid i(v) = j\}$.
        Clearly, the union of $S_1',\dots,S_k'$ is equal to $V(Y_0)$.
        Therefore, fix $\ell \in [k]$ and an $S_\ell'$-rooted model $\calM_0 = (B_{x} \mid x \in V(X_0))$ of $X_0$ in $Y_0$ with $u' \in B_{u}$ since $u$ is a root of $X_0$ and $u'$ is a root of $Y_0$.
        For every $x \in V(X_0)$, let 
        \[D_x = B_x \cup \bigcup_{v\in B_{x}} B_{v,u_1}.\]
        In other words, we add to each branch set $B$ in $\calM_0$, the branch set of $\calM_v$ containing $v$ for every $v\in B$.
        Let $\calD = (D_x \mid x \in V(X_0))$.
        Consider 
        \[\calM = \calD \cup \bigcup_{v \in S_\ell'} (\calM_v\setminus \{B_{v,u_1}\}).\]
        Observe that $\calM$ is an $S_\ell$-rooted model in $Y$.
        Moreover, $\calM$ is a model of a supergraph of $X$ in $Y$.
        Finally, recall that $u' \in V(Y_0)$ is a root of $Y$ and $u' \in B_{u} \subset D_u$.
        This completes the proof of the claim.
    \end{proofclaim}
    Since for every $X \in \Tree(\mathcal{X})$, there exists $X' \in \Tree(\mathcal{X})$ connected such that $X$ is a subgraph of $X'$,
    \Cref{claim:cep:tree} implies that $\Tree(\mathcal{X})$ has the coloring elimination property, which concludes the proof of the lemma.
\end{proof}

The classes $\Rt_0$, $\SRt_0$, $\SRt_1$, and $\SRt_2$ are closed under disjoint union by definition.
If a class of graphs $\calX$ is closed under taking unions, then again by definition, the class $\Tree(\calX)$ are also closed under disjoint union.
In particular, we obtain that $\Rt_t$ and $\SRt_t$ are closed under disjoint union for every nonnegative integer $t$.
To apply~\Cref{lemma:coloring_elimination_and_Tree}, we also need to show that these classes are closed under taking minors.
This follows from~\Cref{lemma:T-is-minor-closed}.

\begin{lemma}
\label{lemma:T-is-minor-closed}
    Let $\calX$ be a class of graphs. 
    If $\mathcal{X}$ is closed under taking minors, then $\Tree(\mathcal{X})$ is closed under taking minors.
\end{lemma}

\begin{proof}
    Assume that $\calX$ is closed under taking minors.
    Let $H$ and $G$ be graphs such that $H$ is a minor of $G$. 
    Suppose that $G\in\Tree(\calX)$. 
    We want to show that $H \in \Tree(\mathcal{X})$.
    Let $(B_u \mid u \in V(H))$ be a model of $H$ in $G$.
    Consider a rooted forest decomposition $\big(F,(W_x \mid x \in V(F))\big)$ of $G$ witnessing the fact that $G \in \Tree(\mathcal{X})$. 
    Then, for every $x \in V(F)$, let $W'_x = \{u \in V(H) \mid B_u \cap W_x \neq \emptyset\}$.
    We claim that $\big(F,(W'_x \mid x \in V(F))\big)$ is a rooted forest decomposition of $H$ witnessing the fact that $H \in \Tree(\mathcal{X})$.
    Clearly, $\big(F,(W'_x \mid x \in V(F))\big)$ is a forest decomposition.
    For every root $x$ of $F$, $|W'_x| \leq |W_x| \leq 1$.
    Moreover, for every $xy \in E(F)$, $|W'_x \cap W'_y| \leq |W_x \cap W_y| \leq 1$.
    Hence, it only remains to prove that for every $x \in V(F)$, if $x$ is not a root of $F$, then the graph $H[W'_x \setminus W'_{\parent(F,x)}]$ lies in $\mathcal{X}$.
    Let $x$ be a non-root vertex of $F$. 
    Since the $(F,(W_x\mid x\in V(F)))$ has adhesion at most $1$, 
    the family $(B_u \cap W_x \mid u \in W'_x \setminus W'_{\parent(F,x)})$ is a model of $H[W'_x \setminus W'_{\parent(F,x)}]$ in $G[W_x \setminus W_{\parent(F,x)}]$. 
    Since $G[W_x \setminus W_{\parent(F,x)}] \in \mathcal{X}$,
    and as $\mathcal{X}$ is closed under taking minors, 
    we conclude that $H[W'_x \setminus W'_{\parent(F,x)}] \in \mathcal{X}$.
\end{proof}

\begin{lemma}\label{lemma:Rt_has_coloring_elimination_property}
    For every nonnegative integer $t$, the classes $\Rt_t$ and $\SRt_t$ have the coloring elimination property.
\end{lemma}

\begin{proof}
    For every nonnegative integer $t$, the classes $\Rt_t$ and $\SRt_t$ are closed under disjoint union.
    Additionally, by~\Cref{lemma:T-is-minor-closed}, they are closed under taking minors.
    Thus, by \Cref{lemma:coloring_elimination_and_Tree}, 
    it is enough to show that $\Rt_0$, $\SRt_0$, $\SRt_1$ and $\SRt_2$ have the coloring elimination property.
    
    Since $\Rt_0$ and $\SRt_0$ consist only of the null graph, the classes $\Rt_0$ and $\SRt_0$ have the coloring elimination property.
    Let $k$ be a positive integer.

    Consider $X \in \SRt_1$.
    Then $X$ consists of $|V(X)|$ isolated vertices.
    Let $Y$ consists of $k(|V(X)|-1) + 1$ isolated vertices.
    For every $S_1, \dots, S_k$ whose union is $V(Y)$,
    there exists $i \in [k]$ such that $|S_i| \geq |V(X)|$ by pigeonhole principle,
    and it follows that $Y$ contains an $S_i$-rooted model of $X$.
    This proves that $\SRt_1$ has the coloring elimination property.

    Consider $X \in \SRt_2$.
    Without loss of generality, $X$ is connected, and so, $X$ is a path, say on $\ell$ vertices.
    We set $Y$ to be a path on $k(\ell-1)+1$ vertices.
    Let $S_1, \dots, S_k$ be sets whose union is $V(Y)$.
    By the pigeonhole principle, there exists $i \in [k]$ such that $|S_i| \geq \ell$, and it follows that there is an $S_i$-rooted model of $X$ in $Y$.
    This shows that $\SRt_2$ has the coloring elimination property and ends the proof.
\end{proof}

\begin{figure}[tp]
    \centering
    \begin{tikzpicture}[scale=1.7]
        \node[myspacing] (param) at (0,1.15) {$\param_q \colon (G,S) \mapsto$};

        \node[myspacing] (frate) at (-3,0.66) {\color{green!50!black}$\frac{1}{\tw(G)+2}\frate_q(G,S)$};
        \node[myspacing] (wcol) at (0,0.66) {\color{blue!50!black}$\frac{1}{\tw(G)+2}\wcol_q(G,S)$};
        \node[myspacing] (cen) at (3,0.66) {\color{red!50!black}$\max_{\phi \text{ $(q,c)$-good}} \cen_q(G,\phi,S)$};
    
        \node[myspacing] (fratehelly) at (-3,0) {\color{green!50!black}\Cref{lem:edgeless-identity-bounding}};
        \node[myspacing] (wcolhelly) at (0,0) {\color{blue!50!black}\Cref{lem:edgeless-identity-bounding}};
        \node[myspacing] (cenhelly) at (3,0) {\color{red!50!black}Definition of $(q,c)$-good colorings};

        \node[box] (Ebounding) at (0,-1) {$q \mapsto 1$ is $(\param, \edgeless)$-bounding};
        \node at (2.5,-1) {\includegraphics{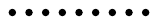}};
        \draw[implies] (fratehelly) -- (Ebounding);
        \draw[implies] (wcolhelly) -- (Ebounding);
        \draw[implies] (cenhelly) -- (Ebounding);

        \node[box] (S2bounding) at (0,-2) {$q \mapsto 1$ is $(\param, \SRt_2)$-bounding};
        \node at (2.5,-2) {\includegraphics{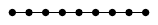}};
        \draw[implies] (Ebounding) -- node[right, myspacing] {\Cref{thm:R1_implies_S2}} (S2bounding);

        \node[box] (AS2bounding) at (0,-3) {$q \mapsto q$ is $(\param, \Apex(\SRt_2))$-bounding};
        \node at (2.5,-3) {\includegraphics{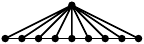}};
        \draw[implies] (S2bounding) -- node[right, myspacing] {\Cref{thm:main_X_to_A(X)}} (AS2bounding);
        
        \node[box] (S3bounding) at (0,-4) {$q \mapsto q$ is $(\param, \SRt_3)$-bounding};
        \node at (2.5,-4) {\includegraphics{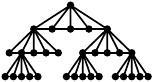}};
        \draw[implies] (AS2bounding) -- node[right, myspacing] {\Cref{thm:AS2_implies_S3}} (S3bounding);

        \node[box] (Stbounding) at (-2.5,-5.5) {$q \mapsto q^{t-2}$ is\\ $(\param,\SRt_t)$-bounding};
        \node[box] (AStbounding) at (2.5,-5.5) {$q \mapsto q^{t-2}$ is\\ $(\param,\Apex(\SRt_{t-1}))$-bounding};

        \draw[implies] (S3bounding) -- node[left=7.5mm, above=-2mm, myspacing] {$t \leftarrow 3$} (Stbounding);
        \draw[implies, transform canvas={yshift=3.3mm}] (Stbounding) -- node[above, myspacing] {\Cref{thm:main_X_to_A(X)}, $t \leftarrow t+1$} (AStbounding);
        \draw[implies, transform canvas={yshift=-3.3mm}] (AStbounding) -- node[below, myspacing] {\Cref{thm:main_A(X)_to_T(X)}} (Stbounding);
        
    \end{tikzpicture}
    \caption{A roadmap for the proofs of \Cref{thm:main_frate_bounded_tw_St,thm:main_wcol_bounded_tw_SRt,thm:main_cen_St} (the first parts of the statements of \Cref{thm:fragility,thm:centered,thm:wcol}). 
        We show that, for a carefully chosen nice family of focused parameters $(\param_q \mid q \in \posint)$, 
        $q \mapsto q^{t-2}$ is $(\param, \SRt_t)$-bounding for every integer $t$ with $t \geq 3$.
        The main part of the induction consists in applying \Cref{thm:main_X_to_A(X)} and \Cref{thm:main_A(X)_to_T(X)}.
        Since \Cref{thm:main_A(X)_to_T(X)} can be applied to $\mathcal{X} = \SRt_{t-1}$ only if $t \geq 4$ 
        (because $\mathcal{X}$ should be closed under leaf addition),
        the cases $t=2$ and $t=3$ are proved separately.
        Typical graphs in $\edgeless, \SRt_2, \Apex(\SRt_2), \SRt_3$ are depicted on the right side.}
    \label{fig:roadmap_wcol_bounded_tw_St}
\end{figure}

\begin{figure}[tp]
    \centering
    \vspace{10mm}
    \begin{tikzpicture}[scale=1.7]
        \node[myspacing] (param) at (0,1.15) {$\param_q(G,S) =$};

        \node[myspacing] (frate) at (-2.5,0.66) {\color{green!50!black}$\frac{1}{\tw(G)+2}\frate_q(G,S)$};
        \node[myspacing] (wcol) at (0,0.66) {\color{blue!50!black}$\frac{1}{\tw(G)+2}\wcol_q(G,S)$};
        \node[myspacing] (cen) at (2.5,0.66) {\color{red!50!black}$\max_{\phi \text{ $(q,c)$-good}} \cen_q(G,\phi,S)$};
    
        \node[myspacing] (frateforest) at (-2.5,0) {\color{green!50!black}\Cref{lemma:fragility_rate_base_case}};
        \node[myspacing] (wcolforest) at (0,0) {\color{blue!50!black}\Cref{lemma:bounded_tw_excluding_a_star_wcol}};
        \node[myspacing] (cenforest) at (2.5,0) {\color{red!50!black}\Cref{lemma:base_case_R2_cen_total}};
        
        \node[box] (R2bounding) at (-2.5,-1) {$q \mapsto \log (q+1)$ is $(\param, \Rt_2)$-bounding};
        \node at (0.25,-1) {\includegraphics{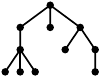}};
        \draw[implies] (frateforest) -- (R2bounding);
        \draw[implies] (wcolforest) -- (R2bounding);
        \node[box] (AR2bounding) at (2.5,-2.5) {$q \mapsto q \log (q+1)$ is $(\param, \Apex(\Rt_2))$-bounding};
        \node at (-0.25,-2.5) {\includegraphics{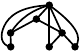}};
        \draw[implies] (cenforest) -- (AR2bounding);
        
        \node[box] (Rtbounding) at (-2.5,-4.5) {$q \mapsto q^{t-2} \log (q+1)$ is\\ $(\param,\Rt_t)$-bounding};
        \node[box] (ARtbounding) at (2.5,-4.5) {$q \mapsto q^{t-2} \log (q+1)$ is\\ $(\param,\Apex(\Rt_{t-1}))$-bounding};

        \draw[implies] (R2bounding) -- node[left=3.33mm, myspacing] {$t \leftarrow 2$} (Rtbounding);
        \draw[implies] (AR2bounding) -- node[left=3.33mm, myspacing] {$t \leftarrow 3$} (ARtbounding);
        \draw[implies, transform canvas={yshift=3.3mm}] (Rtbounding) -- node[above, myspacing] {\Cref{thm:main_X_to_A(X)}, $t \leftarrow t+1$} (ARtbounding);
        \draw[implies, transform canvas={yshift=-3.3mm}] (ARtbounding) -- node[below, myspacing] {\Cref{thm:main_A(X)_to_T(X)}} (Rtbounding);
    \end{tikzpicture}
    \caption{A roadmap for the proofs of \Cref{thm:main_frate_bounded_tw_Rt,thm:main_wcol_bounded_tw_Rt,thm:main_cen_Rt} (the second parts of the statements of \Cref{thm:fragility,thm:centered,thm:wcol}).}
    \label{fig:roadmap_wcol_bounded_tw_Rt}
\end{figure}

\section{Proof of the abstract induction step}
\label{sec:proof_abstract_theorems}
See~\Cref{fig:roadmap_wcol_bounded_tw_St,fig:roadmap_wcol_bounded_tw_Rt} for a roadmap of the proof.

\subsection{Preliminary reductions}

Let $(\param_q \mid q \in \posint)$ be a family of focused parameters,
let $g \colon \posint \to \posint$, and let $\mathcal{X}$ be a class of graphs.
In several cases, to show that $g$ is $(\param,\calX)$-bounding, we first show a weaker property that in fact (under some mild assumptions) implies the desired one, as shown in~\Cref{lemma:weakly_bounding_implies_bounding}.
Recall the definitions given in~\Cref{ssec:abstract-induction}.

The function $g$ is \defin{weakly $(\param,\mathcal{X})$-bounding} if for every positive integer $k$,
for every graph $X \in \mathcal{X}$,
there exist nonnegative integers $\alpha(X)$ and $\beta(X)$ such that every $K_k$-minor-free graph~$G$ is $(g,\param,X,\alpha(X),\beta(X))$-good.

It will also be convenient to assume that the graph $G$ in the definition above and the definition of $(\param,\mathcal{X})$-bounding is connected.
To this end, we show~\Cref{lem:reduce-to-G-connected}.
Similarly, when proving that some function is $(\param,\Tree(\calX))$-bounding, it is handy to assume that $X \in \Tree(\calX)$ is connected.
This is the purpose of~\Cref{lem:connected-graphs-in-T(X),cor:X-connected-in-Tree-X-bounding}.

\begin{lemma}\label{lemma:weakly_bounding_implies_bounding}
    Let $(\param_q \mid q \in \posint)$ be a nice family of focused parameters, let $\mathcal{X}$ be a class of graphs,
    and let $g\colon \posint \to \posint$.
    If $g$ is $(\param,\edgeless)$-bounding and weakly $(\param,\mathcal{X})$-bounding, then $g$ is $(\param,\mathcal{X})$-bounding.
\end{lemma}

\begin{figure}[tp]
    \begin{center}
    \includegraphics{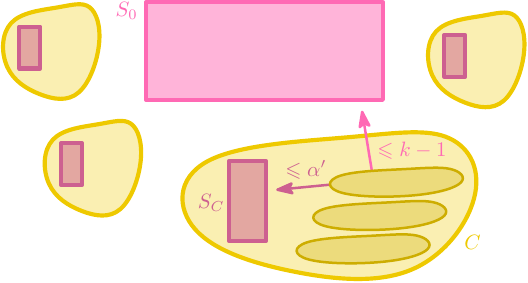}
    \end{center}
    \caption{
    An illustration of the proof of~\Cref{lemma:weakly_bounding_implies_bounding}.
    }
    \label{fig:weakly-to-bounding}
\end{figure}

\begin{proof}
    Let $k$ be a positive integer.
    Fix $\alpha_0(\cdot)$ and $\beta_0(\cdot)$ that witness $g$ being weakly $(\param,\mathcal{X})$-bounding for $k$ and fix $\alpha'$ and $\beta'(\cdot)$ witnessing that $g$ is $(\param,\edgeless)$-bounding for $k$.
    Since the only property distinguishing graphs in $\edgeless$ is the number of vertices, we treat $\beta'(\cdot)$ as a function taking nonnegative integers
    by setting $\beta'(d) = \beta'(\overline{K_d})$ for every positive integer $d$.
    Let
    \[\alpha = \alpha' + (k-1) \ \ \text{ and } \ \ \beta(X) = \beta_0(X) + \beta'\left(\textstyle\binom{\alpha_0(X)}{k}(k-1)+1\right) \text{ for every $X \in \calX$.}\]
    We show that $\alpha$ and $\beta(\cdot)$ witness $g$ being $(\param,\mathcal{X})$-bounding for $k$.
    Let $X \in \calX$ and let $G$ be a $K_k$-minor-free graph
    We show that $G$ is $(g, \param, X, \alpha, \beta(X))$-good.
    Let $q$ be a positive integer
    and let $\mathcal{F}$ be a family of connected subgraphs of $G$
    such that there is no $\mathcal{F}$-rich model of $X$ in $G$.
    Since $\alpha_0(\cdot)$ and $\beta_0(\cdot)$ witness $g$ being weakly $(\param,\calX)$-bounding for $k$,
    there exists $S_0 \subseteq V(G)$ such that
    \begin{enumerate}[label={\normalfont (g\arabic*-0)}]
        \item $S_0 \cap V(F) \neq \emptyset$ for every $F \in \mathcal{F}$; \label{item:gw1'}
        \item for every component $C$ of $G-S_0$, $N_G(V(C))$ intersects at most $\alpha_0(X)$ components of $G-V(C)$; and \label{item:gw2'}
        \item $\param_q(G,S_0) \leq \beta_0(X) \cdot g(q)$. \label{item:gw3'}
    \end{enumerate}
    Let $\mathcal{C}$ be the family of all the components of $G-S_0$,
    and let $C \in \mathcal{C}$.
    Let $\mathcal{U}$ be the family of the vertex sets of all the components of $G-V(C)$
    intersecting $N_G(V(C))$.
    By~\ref{item:gw2'}, we have $|\calU| \leq \alpha_0(X)$.
    Let $\mathcal{F}_C$ be the family of all the connected subgraphs $H$ of $C$
    such that $|\{U \in \mathcal{U} \mid N_G(V(H)) \cap U \neq \emptyset\}| \geq k$.
    We claim that there are no $N = \binom{\alpha_0(X)}{k}(k-1)+1$ pairwise disjoint
    members of $\mathcal{F}_C$.
    Otherwise, by the pigeonhole principle, there exist $U_1, \dots, U_k \in \mathcal{U}$ and $H_1, \dots, H_k \in \mathcal{F}_C$ pairwise disjoint
    such that there is an edge between $U_i$ and $V(H_j)$ in $G$ for all $i,j \in [k]$.
    However, then $(U_i \cup V(H_i) \mid i \in [k])$ is a model of $K_k$ in $G$, which is a contradiction proving that indeed there are no $N$ pairwise disjoint members of $\mathcal{F}_C$.
    
    Recall that $\alpha'$ and $\beta'(\cdot)$ witness $g$ being $(\param,\edgeless)$-bounding for $k$.
    Therefore, for every $C\in \mathcal{C}$, since $C$ has no $\mathcal{F}_C$-rich model of $\overline{K_N}$, i.e.\ $N$ pairwise disjoint members of $\mathcal{F}_C$, 
    there exists $S_{C} \subseteq V(C)$ such that
    \begin{enumerate}[label={\normalfont (g\arabic*-C)}]
        \item $S_{C} \cap V(F) \neq \emptyset$ for every $F \in \mathcal{F}_C$; \label{item:g1'}
        \item for every component $C'$ of $C-S_C$, $N_G(V(C'))$ intersects at most $\alpha'$ components of $C-V(C')$; and \label{item:g2'}
        \item $\param_q(C,S_C) \leq \beta'(N) \cdot g(q)$.\label{item:g3'}
    \end{enumerate}
    We set
    \[
        S = S_0 \cup \bigcup_{C \in \mathcal{C}} S_C.
    \]
    See~\Cref{fig:weakly-to-bounding}.
    It suffices to show that
    \begin{enumerate}[label={\normalfont (g\arabic*$^\star$)}]
        \item $S \cap V(F) \neq \emptyset$ for every $F \in \calF$; \label{item:par-bounding:hitting:g1'''}
        \item for every component $C$ of $G-S$, $N_G(V(C))$ intersects at most $\alpha$ components of $G-V(C)$; \label{item:par-bounding:nbrs:g2'''}
        \item $\param_{q}(G,S) \leq \beta(X) \cdot g(q)$.\label{item:par-bounding:bound:g3'''}
    \end{enumerate}
    It suffices to verify that $G$ is $(g,\param,X,\alpha,\beta)$-good.
    Since $S_0 \subset S$, by~\ref{item:gw1'}, $S \cap V(F) \neq \emptyset$ for every $F \in \mathcal{F}$, and so,~\ref{item:par-bounding:hitting:g1'''} holds.
    For every component $C'$ of $G-S$, there is a component $C \in \mathcal{C}$
    such that $C' \subseteq C$.
    Since $C' \not\in \mathcal{F}_C$ (by~\ref{item:g1'}), $N_G(V(C'))$ intersects at most $k-1$ components of $G-V(C)$. 
    Moreover, $N_G(V(C'))$ intersects at most $\alpha'$ components of $C-S_C$ (by~\ref{item:g2'}).
    Therefore, $N_G(V(C'))$ intersects at most $\alpha'+(k-1) = \alpha$ components
    of $G-V(C')$.
    Thus, \ref{item:par-bounding:nbrs:g2'''} holds.
    Finally,
    \begin{align*}
        \param_q(G,S)
        &\leq \param_q(G,S_0) + \param_q\left(G-S_0,\textstyle\bigcup_{C \in \mathcal{C}} S_C\right) &\textrm{by \ref{item:nice:subadditivity}} \\
        &\leq \param_q(G,S_0) + \textstyle\max_{C \in \mathcal{C}} \param_q(C,S_C) &\textrm{by \ref{item:nice:components}} \\
        &\leq \beta_0(X) \cdot g(q) + \beta'(N) \cdot g(q) &\textrm{by \ref{item:gw3'} and \ref{item:g3'}} \\
        &= \beta(X) \cdot g(q).
    \end{align*}
    Hence, \ref{item:par-bounding:bound:g3'''} holds.
    Altogether, we obtain that $g$ is $(\param,\mathcal{X})$-bounding,
    which concludes the proof of the lemma.
\end{proof}

\begin{lemma} \label{lem:reduce-to-G-connected}
    Let $\param$ be a nice family of focused parameters, 
    let $g\colon \posint \to \posint$ be a function.
    Let $X$ be a graph, let $\alpha$ and $\beta$ be nonnegative integers, and let $G$ be a graph.
    If every component of $G$ is $(g,\param,X,\alpha,\beta)$-good, then $G$ is $(g,\param,X,\alpha,\beta)$-good.
\end{lemma}
\begin{proof}
    Suppose that every component of $G$ is $(g,\param,X,\alpha,\beta)$-good.
    Let $q$ be a positive integer and let $\calF$ be a family of connected subgraphs of $G$ such that $G$ has no $\mathcal{F}$-rich model of $X$.
    Let $\mathcal{C}$ be the family of all the components of $G$.
    
    Let $C$ be a component of $G$.
    Since $C$ is $(g,\param,X,\alpha,\beta)$-good, there exists $S_C \subseteq V(C)$ such that
    \begin{enumerate}[label={\normalfont (g\arabic*-C)}]
        \item $S_C \cap V(F) \neq \emptyset$ for every $F \in \mathcal{F}\vert_C$; \label{item:reduce-to-strong-bounding:SC:hitting}
        \item for every component $C'$ of $C-S_C$, $N_C(V(C'))$ intersects at most $\alpha$ components of $C-V(C')$; and \label{item:reduce-to-strong-bounding:SC:connectivity}
        \item $\param_q(C,S_C) \leq \beta \cdot g(q)$. \label{item:reduce-to-strong-bounding:SC:param}
    \end{enumerate}
    We set $S = \bigcup_{C \in \mathcal{C}} S_C$ and we show that
    \begin{enumerate}[label={\normalfont (g\arabic*$^\star$)}]
        \item $S \cap V(F) \neq \emptyset$ for every $F \in \calF$;\label{item:par-reduce-to-strong-bounding:SC:hitting''}
        \item for every component $C'$ of $G-S$, $N_G(V(C'))$ intersects at most $\alpha$ components of $G-V(C')$; and \label{item:reduce-to-strong-bounding:SC:nbrs''}
        \item $\param_{q}(G,S) \leq \beta \cdot g(q)$.\label{item:reduce-to-strong-bounding:SC:bound''}
    \end{enumerate}
    For every $F \in \mathcal{F}$, there exists $C \in \mathcal{C}$
    such that $V(F) \subseteq V(C)$, and so by \ref{item:reduce-to-strong-bounding:SC:hitting}, $S_C \cap V(F) \neq \emptyset$, which implies $S \cap V(F) \neq \emptyset$.
    This gives~\ref{item:par-reduce-to-strong-bounding:SC:hitting''}.
    For every component $C'$ of $G-S$, there exists $C \in \mathcal{C}$
    such that $V(C') \subseteq V(C)$.
    Then, by \ref{item:reduce-to-strong-bounding:SC:connectivity}, $N_G(V(C')) = N_C(V(C'))$ intersects at most $\alpha$ components of $C-V(C')$, and so of $G-V(C')$.
    This shows~\ref{item:reduce-to-strong-bounding:SC:nbrs''}.
    Finally, by \ref{item:nice:components} and~\ref{item:reduce-to-strong-bounding:SC:param},
    \[
        \param_q(G,S) \leq \max_{C \in \mathcal{C}} \param_q(C,S_C) \leq \beta \cdot g(q).
    \]
    This yields~\ref{item:reduce-to-strong-bounding:SC:bound''}.
    Altogether, we obtain that $G$ is $(g,\param,X,\alpha,\beta)$-good.
\end{proof}

Let $\calX$ be a class of graphs.
We define \defin{$\Tree'(\calX)$} to be all the graphs $X \in \Tree(\calX)$ such that $X$ is connected and $X$ admits a rooted tree decomposition $\big(T,(W_x \mid x \in V(T))\big)$ witnessing $X \in \Tree(\mathcal{X})$ such that $W_a \neq \emptyset$ for every $a \in V(T)$.

\begin{lemma} \label{lem:connected-graphs-in-T(X)}
    Let $\calX$ be a class of graphs containing a nonnull graph.
    For every $X \in \Tree(\calX)$, there exists $X' \in \Tree'(\calX)$
    such that $X$ is a subgraph of $X'$.
\end{lemma}
\begin{proof}
    Let $Z$ be a nonnull graph in $\mathcal{X}$.
    Fix a vertex $z \in V(Z)$ and a graph $X \in \Tree(\mathcal{X})$.
    Assume that $X$ and $Z$ have disjoint vertex sets.
    Let $\mathcal{C}$ be the family of all components of $X$.
    Let $C\in\mathcal{C}$.
    We have $C \in \Tree(\mathcal{X})$.
    Since $C$ is connected, we fix $\big(T_C, (W_{C,x} \mid x \in V(T_C))\big)$ a rooted tree decomposition of $C$ witnessing the fact that $C \in \Tree(\mathcal{X})$, and let $r_C = \root(T_C)$.
    We can also assume without loss of generality that all bags except the one corresponding to the root are nonempty.
    Assume that the trees $T_C$ for $C \in \mathcal{C}$ have pairwise disjoint vertex sets.

\begin{figure}[tp]
\begin{center}
\includegraphics{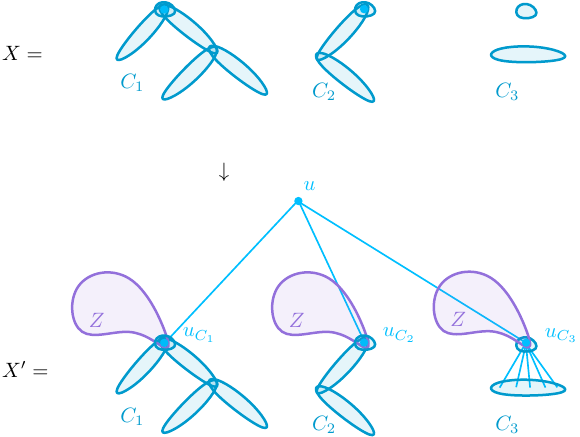}
\end{center}
\caption{
An illustration of the construction from the proof of~\Cref{lem:connected-graphs-in-T(X)}.
The graph $X$ has three components $C_1$, $C_2$, and $C_3$. In the top part of the figure, we depict a forest decomposition witnessing that $X \in \Tree(\calX)$.
Note that here $W_{C_3,r_{C_3}}$ is empty.
}
\label{fig:making-graph-connected}
\end{figure}

    We illustrate the following construction in~\Cref{fig:making-graph-connected}.
    Let $C \in \calC$.
    If $|W_{C,r_C}| = 1$, then we set $u_C$ to be the unique vertex of $W_{C,r_C}$ and we set $C' = C$.
    In this case, we also set $W_{C',x} = W_{C,x}$ for every $x \in V(T_C)$.
    Next, assume that $|W_{C,r_C}| = 0$.
    Let $C'$ be obtained from $C$ by adding a new vertex $u_C$ and making it adjacent to all vertices in $\bigcup \{W_{C,x} \mid \text{$x$ and $r_C$ are adjacent in $T_C$}\}$.
    We set $W_{C',x} = W_{C,x}$ for every $x \in V(T_C) \setminus \{r_C\}$ and $W_{C',r_C} = \{u_C\}$.
    In both cases $C'$ is connected and $\big(T_C, (W_{C',x} \mid x \in V(T_C))\big)$ is a rooted tree decomposition witnessing that $C' \in \Tree(\calX)$.

    For each $C \in \calC$, let $C''$ be obtained from the union of $C'$ and $Z$ by identifying $u_C$ and $z$.
    Finally, let $X'$ be obtained from the union of $C''$ for each $C \in \calC$ by adding a new vertex $u$ and making it adjacent to all vertices in the copy of $Z$ in each $C''$.
    By construction, $X'$ is connected.
    
    We shall prove that $X' \in \Tree(\calX)$.
    Let $T$ be obtained from the disjoint union of $T_C$ for $C \in \calC$ by adding a new vertex $r$ adjacent to each $r_C$ for $C \in \calC$.
    We set $\root(T) = r$.
    Let $x \in V(T)$.
    If $x = r$, we set $W_r = \{u\}$.
    Let $C \in \calC$ be such that $x \in V(T_C)$.
    If $x = r_C$, we set $W_x$ to be the set of vertices in the copy of $Z$ in $C''$.
    Finally, otherwise, we set $W_x = W_{C',x}$.
    It is easy to verify that $\big(T, (W_{x} \mid x \in V(T))\big)$ is a tree decomposition witnessing that $X' \in \Tree'(\calX)$.
    This ends the proof.
\end{proof}

\begin{corollary}\label{cor:X-connected-in-Tree-X-bounding}
    Let $(\param_q \mid q \in \posint)$ be a family of focused parameters, let $g \colon \posint \rightarrow \posint$, and let $\calX$ be a class of graphs containing a nonnull graph.
    Then, $g$ is weakly $(\param,\Tree(\calX))$-bounding if and only if $g$ is weakly $(\param,\Tree'(\calX))$-bounding.
\end{corollary}

\subsection{The proofs}

\begin{proof}[Proof of~\Cref{thm:main_X_to_A(X)}]
    Without loss of generality, assume that $\mathcal{X}$ contains a nonnull graph.
    Let $\mathcal{X}'=\{X'\mid \textrm{$X'$ is a minor of $X$ for some $X\in\calX$}\}$.
    Since $\calX$ is closed under disjoint union, $\calX'$ is also closed under disjoint union.
    Since $\calX$ has the coloring elimination property, by~\Cref{lemma:coloring_elimination_minor_closed}, $\calX'$ also has the coloring elimination property.
    Since $\Apex(\calX) \subset \Apex(\calX')$, if a function is $(\param,\Apex(\calX'))$-bounding, then this function is also $(\param,\Apex(\calX))$-bounding.
    Therefore, again without loss of generality, we will assume that $\calX = \calX'$, and so, $\calX$ is closed under taking minors.
    Since $\calX$ has a nonnull graph, $K_1 \in \calX$.
    
    Assume that $g$ is $(\param,\mathcal{X})$-bounding.
    Let $b$ and $b'$ be as in~\ref{item:nice:tree-partitions} for $(\param_q \mid q \in \posint)$.
    Let $k$ be a positive integer.
    Fix $\alpha'$ and $\beta'(\cdot)$ witnessing that $g$ is $(\param,\mathcal{X})$-bounding for $k$.
    Let $X \in \Apex(\mathcal{X})$ and $z \in V(X)$ such that $X - z \in \calX$.
    Note that $X' = K_1 \sqcup (X-z)$ is in $\mathcal{X}$.
    Let $Y'$ witness the coloring elimination property of $\calX$ for $X'$ and $\alpha'$.

    \begin{figure}[tp]
        \centering
        \includegraphics{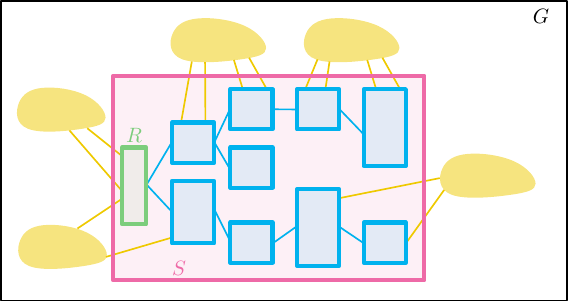}
        \caption{
            An illustration of the statement of~\Cref{claim:X-to-A(X)}.
        }
        \label{fig:excluding-tree-statement}
    \end{figure}

    \begin{claim}\label{claim:X-to-A(X)}
        Let $G$ be a connected $K_k$-minor-free graph,
        let $R$ be a nonempty set of at most $\alpha'$ vertices in $G$,
        let $\mathcal{F}$ be a family of connected subgraphs of $G$,
        and let $q$ be a positive integer.
        If $G$ has no $\mathcal{F}$-rich model of $X$, 
        then
        there exist $S \subseteq V(G)$, a tree $T$ rooted in $s \in V(T)$, and a tree partition $\big(T,(P_x \mid x \in V(T))\big)$
        of $(G,S)$
        with $P_s = R$
        such that
        \begin{enumerateOurAlph}
            \item $S \cap V(F)  \neq \emptyset$ for every $F \in \mathcal{F}$; \label{item:bounded_tw_claim_adding_apices_i}
            \item for every component $C$ of $G-S$, $N_G(V(C))$ intersects at most $2\alpha'$ components of $G-V(C)$; and \label{item:bounded_tw_claim_adding_apices_ii}
            \item for every $x \in V(T)$,
            \[
                \param_{q}\left(G_x, P_x\right) \leq \max\{\beta'(Y'), \alpha'\} \cdot g(q)
            \]
            where for every $x \in V(T)$, $G_x$ is the union of $U_x = \bigcup_{z \in V(\subtree{T}{x})} P_z$ and all the vertex sets of the components of $G-S$ having a neighbor in $U_x$. \label{item:bounded_tw_claim_adding_apices_vii}
        \end{enumerateOurAlph}
    \end{claim}

    \begin{proofclaim}
        The statement is illustrated in~\Cref{fig:excluding-tree-statement}.
        We proceed by induction on $|V(G)|$.
        If $\mathcal{F}\vert_{G-R} = \emptyset$, then the statement holds for $S=R$ and $T=K_1$.
        Thus, we suppose that $\mathcal{F}\vert_{G-R} \neq \emptyset$, and in particular $|V(G) \setminus R| > 0$.

        Let $\mathcal{F}'$ be the family of all connected subgraphs $H$ of $G-R$ such that $R \cap N_G(V(H)) \neq \emptyset$ and
        $F \subseteq H$ for some $F \in \mathcal{F}$.
        We argue that there is no $\mathcal{F}'$-rich model of $Y'$ in $G-R$.
        Suppose to the contrary that there is an $\mathcal{F}'$-rich model $(B_x \mid x \in V(Y'))$ of $Y'$ in $G-R$.
        For each $u \in R$, let $S_u = \{x \in V(Y') \mid u \in N_G(B_x)\}$.
        Since the model is $\mathcal{F}'$-rich, $\bigcup_{u \in R} S_u = V(Y')$.
        Since $Y'$ witnesses the coloring elimination property of $\mathcal{X}$ for $X'$ and $\alpha'$,
        there exists $u \in R$ such that $Y'$ contains an 
        $S_u$-rooted model of $X'$.
        As a consequence, there is an $\mathcal{F}$-rich model $(D_x \mid x \in V(X'))$ of $X'$ in $G-R$ 
        such that every branch set $D_x$ contains a neighbor of $u$ in $G$.
        Recall that $X' = K_1 \sqcup (X-z)$.
        Say that $y$ is the vertex of $K_1$ in $X'$.
        Let $D_z' = D_y \cup \{u\}$ and $D_x' = D_x$ for every $x \in V(X' \setminus\{y\})$.
        It follows that $(D_x' \mid x \in V(X))$ is an $\mathcal{F}$-rich model of $X$ in $G$, which is a contradiction.
        In turn, there is no $\mathcal{F}'$-rich model of $Y'$ in $G-R$.

        Since $\alpha'$ and $\beta'(\cdot)$ witness $g$ being $(\param,\calX)$-bounding for $k$, when applied to $Y'$, $G-R$ and $\mathcal{F}'$,
        we get a set $S_0 \subseteq V(G-R)$ such that   
        \begin{enumerate}[label={\normalfont (g\arabic*-0)}]
            \item $S_0 \cap V(F) \neq \emptyset$ for every $F \in \mathcal{F}'$;\label{thm:bounded_tw_main:item:hit'}
            \item for every component $C$ of $(G-R) - S_0$, $N_{G-R}(V(C))$ intersects at most $\alpha'$ components of $(G-R)-V(C)$; and \label{thm:bounded_tw_main:item:connected'}
            \item $\param_{q}(G-R,S) \leq \beta'(Y') \cdot g(q)$. \label{thm:bounded_tw_main:item:wcol'}
        \end{enumerate}
        Since $G$ is connected and $\mathcal{F}\vert_{G-R}$ is nonempty,
        there is a component of $G-R$ containing a member of $\mathcal{F}$, and so $\mathcal{F}'$ is nonempty.
        It follows by \ref{thm:bounded_tw_main:item:hit'} that $S_0 \neq \emptyset$.
    
        Let $\mathcal{C}_1$ be the family of all components of $(G-R)- S_0$ with no neighbors in $R$.
        Consider $C \in \mathcal{C}_1$.
        Let $A_C$ be the component of $G-R$ containing $C$, 
        and let $G_C$ be the graph obtained from $A_C$ by contracting each component of $A_C - V(C)$ into a single vertex.
        Let $R_{C}$ be the set of all vertices resulting from these contractions, that is $R_C = V(G_C) \setminus V(C)$.
        Observe that $R_C$ is nonempty since $G$ is connected, and $|R_C| \leq \alpha'$ by \ref{thm:bounded_tw_main:item:connected'}.
        Moreover, since $R \neq \emptyset$,  $|V(G_C)| < |V(G)|$.
        Since $G_C$ is a minor of $G$, $G_C$ has no $\mathcal{F}\vert_C$-rich model of $X$.
        By the induction hypothesis applied to $G_C$, $R_C$, and $\mathcal{F}\vert_C$, there exist $S_C \subseteq V(G_C)$, 
        a tree $T_C$ rooted in $s_C \in V(T_C)$, and a tree partition $\big(T_C,(P_{C,x} \mid x \in V(T_C))\big)$ of $(G_C,S_C)$ with $P_{C,s_C} = R_C$
        such that
        \begin{enumerateOurAlphPrim}
            \item $S_C \cap V(F) \neq \emptyset$ for every $F \in \mathcal{F}\vert_{C}$; \label{item:bounded_tw_claim_adding_apices_i-call-call}
            \item for every component $C'$ of $G_C-S_C$, $N_{G_C}(V(C'))$ intersects at most $2\alpha'$ components of $G_C-V(C')$; and \label{item:bounded_tw_claim_adding_apices_ii-call-call}
            \item for every $x \in V(T_C)$,
                \[
                    \param_{q}(G_{C,x}, P_{C,x}) \leq \max\{\beta'(Y'),\alpha'\} \cdot g(q)
                \]
                where for every $x \in V(T_C)$, $G_{C,x}$ is the union of $U_{C,x} = \bigcup_{z \in V(\subtree{T_{C}}{x})} P_z$ 
                and all the vertex sets of the components of $G_C-S_C$ having a neighbor in $U_{C,x}$. \label{item:bounded_tw_claim_adding_apices_vii-call-call}
        \end{enumerateOurAlphPrim}

        \begin{figure}[tp]
            \centering
            \includegraphics{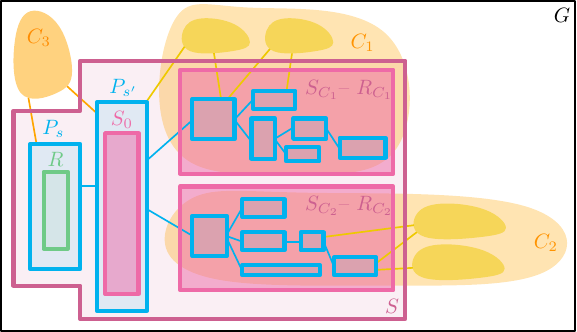}
            \caption{
                An illustration of the proof of~\Cref{claim:X-to-A(X)}.
                In the sketched case, $C_1,C_2 \in \calC_1$ and $C_3 \notin \calC_1$.
            }
            \label{fig:excluding-tree}
        \end{figure}
        
        We set
        \[
            S = R \cup S_0 \cup \bigcup_{C \in \mathcal{C}_1} (S_C \setminus R_C).
        \]
        Let $T$ be obtained from the disjoint union of $\{T_C \mid C\in\mathcal{C}_1\}$ by identifying the vertices 
        $\set{s_{C}\mid C\in\mathcal{C}_1}$ into a new vertex $s'$ and by adding the root $s$ adjacent only to $s'$ in $T$. 
        Let $P_s=R$, $P_{s'}=S_0$, and
        for each $C\in\mathcal{C}_1$, 
        $x \in V(T_{C} \setminus \{s_{C}\})$, 
        let $P_x = P_{C,x}$.
        See~\Cref{fig:excluding-tree}.

        In order to conclude the proof of the claim, we argue that $\big(T,(P_x \mid x \in V(T))\big)$ is a tree partition of $(G,S)$ 
        and that \ref{item:bounded_tw_claim_adding_apices_i}--\ref{item:bounded_tw_claim_adding_apices_vii} hold.

        Since for every $C \in \mathcal{C}_1$, 
        $R \cap N_G(V(C)) = \emptyset$, every edge in $G[S]$ containing a vertex in $R$ has another endpoint in $R \cup S_0 = P_s \cup P_{s'}$.
        Consider an edge $vw$ in $G[S]$ such that $v \in S_0$ and $w \in S_C$ for some $C \in \mathcal{C}_1$.
        Since $\big(T_C,(P_{C,x} \mid x \in V(T_C))\big)$ is a tree partition of $(G_C,S_C)$ with $P_{C,s_C} = R_C$ and $S_0 \subset V((G-R)-V(C))$, we conclude that $w \in P_x$ for some $x \in V(T_C)$ such that $s'x$ is an edge in $T$.
        Finally, for every edge $vw$ of $G[S]$ with $v,w \not\in R \cup S_0$, $vw$ is an edge of $G[S_C \setminus R_C]$ for some component $C \in \mathcal{C}_1$,
        and so $v \in P_{C,x}$ and $w \in P_{C,y}$ for adjacent or identical vertices $x,y$ of $T_C$.
        Then $v \in P_x$ and $w \in P_y$.
        For every component $C'$ of $G-S$, either $N_G(V(C')) \subseteq R \cup S_0 = P_s \cup P_{s'}$, or $C' \subseteq C$ for some $C \in \mathcal{C}_1$.
        In the latter case, $C'$ is a component of $G_C - S_C$, and $N_G(R) \cap V(C) = \emptyset$.
        Since $\big(T_C,(P_{C,x} \mid x \in V(T_C))\big)$ is a tree partition of $(G_C,S_C)$,
        there is $x,y \in V(T_C)$ such that $N_{G_C}(V(C')) \subseteq P_{C,x} \cup P_{C,y}$.
        If $x=s_C$, then $N_G(V(C')) \subseteq P_s \cup P_y$ and $sy \in E(T)$.
        If $y=s_C$, then $N_G(V(C')) \subseteq P_s \cup P_x$ and $sx \in E(T)$.
        Finally, if $x,y \neq s_C$, then $N_G(V(C')) \subseteq P_s \cup P_x$ and $xy \in E(T)$.
        This proves that $\big(T,(P_x \mid x \in V(T))\big)$ is a tree partition of $(G,S)$.

        Let $F \in \mathcal{F}$.
        If $(R \cup S_0) \cap V(F) \neq \emptyset$, then $S \cap V(F) \neq \emptyset$.
        Otherwise, $F \subseteq (G - R) - S_0$ and therefore, by \ref{thm:bounded_tw_main:item:hit'}, we have $F \notin \mathcal{F}'$.
        In this case, there is a component $C \in \mathcal{C}_1$ such that $F \in \mathcal{F}\vert_C$.
        By \ref{item:bounded_tw_claim_adding_apices_i-call-call}, $(S_C \setminus R_C) \cap V(F) \neq \emptyset$ and so $S \cap V(F) \neq \emptyset$.
        This proves~\ref{item:bounded_tw_claim_adding_apices_i}.

        Consider a component $C'$ of $G - S$.
        If $C' \subseteq C$ for some $C \in \mathcal{C}_1$, then by~\ref{item:bounded_tw_claim_adding_apices_ii-call-call}, it follows that $N_{G_C}(V(C'))$ intersects at most $2\alpha'$ components of $G_C-V(C')$, 
        and so $N_{G}(V(C'))$ intersects at most $2\alpha'$ components of $G-V(C')$.
        Otherwise, $C'$ is a component of $(G - R) - S_0$ such that $N_G(R) \cap V(C') \neq \emptyset$.
        By~\ref{thm:bounded_tw_main:item:connected'}, $N_{G-R}(V(C'))$ intersects at most $\alpha'$ components of $(G-R)-V(C')$, and therefore, $N_G(V(C'))$ intersects at most $\alpha'+|R| \leq 2\alpha'$ components of $G-V(C')$.
        This proves~\ref{item:bounded_tw_claim_adding_apices_ii}.

        Finally, we argue~\ref{item:bounded_tw_claim_adding_apices_vii}.
        Let $x \in V(T)$.
        If $x = s$, $P_{x} = P_s = R$, which has size at most $\alpha'$, and thus the assertion holds by \ref{item:nice:at_most_size_of_S}.
        If $x = s'$, then $G_{s'}$ is a union of components of $G-R$,
        and $P_x = P_{s'} = S_0$, and so, by~\ref{item:nice:components} and~\ref{thm:bounded_tw_main:item:wcol'}, 
        \[
             \param_{q}(G_{s'}, S_0) \leq \param_{q}(G-R, S_0) \leq \beta'(Y') g(q).
        \]
        If $x \in V(T_C - \{s_C\})$ for some $C \in \mathcal{C}_1$, we have $\subtree{T}{x} = \subtree{T_C}{x}$ and $G_x = G_{C,x}$.
        Thus, the asserted inequality follows from~\ref{item:bounded_tw_claim_adding_apices_vii-call-call}.
        This ends the proof of the claim.
    \end{proofclaim}
    Let 
    \[\alpha = 2\alpha' \ \ \text{ and } \ \ \beta(X) = 2b'\max\{\beta'(Y'), \alpha'\}.\]
    
    In order to conclude the proof, we argue that every $K_k$-minor-free graph $G$ is $(q \mapsto q g(bq),\param,X,\alpha,\beta(X))$-good.
    Since $X$ was any graph in $\Apex(\calX)$, this will show that $q \mapsto q g(bq)$ is $(\param,\Apex(\calX))$-bounding.

    By~\cref{lem:reduce-to-G-connected}, it suffices to show that every connected $K_k$-minor-free graph $G$ is $(q \mapsto q g(bq),\param,X,\alpha,\beta(X))$-good.
    Let $G$ be a connected $K_k$-minor-free graph, let $q$ be a positive integer, and let $\mathcal{F}$ a family of connected subgraphs of $G$ such that
    there is no $\mathcal{F}$-rich model of $X$ in $G$.
    We show that
    \begin{enumerate}[label={\normalfont (g\arabic*$^\star$)}]
        \item $S \cap V(F) \neq \emptyset$ for every $F \in \calF$;\label{item:par-bounding:hitting-16}
        \item for every component $C$ of $G-S$, $N_G(V(C))$ intersects at most $\alpha$ components of $G-V(C)$; and \label{item:par-bounding:nbrs-16}
        \item $\param_{q}(G,S) \leq \beta \cdot g(q)$.\label{item:par-bounding:bound-16}
    \end{enumerate}
    
    By~\Cref{claim:X-to-A(X)}, applied for $R$ being an arbitrary singleton in $V(G)$,
    there exist $S \subseteq V(G)$, a tree $T$ rooted in $s \in V(T)$, and a tree partition $\big(T,(P_x \mid x \in V(T))\big)$
    of $(G,S)$
    with $P_s = R$ such that \ref{item:bounded_tw_claim_adding_apices_i}, \ref{item:bounded_tw_claim_adding_apices_ii}, and~\ref{item:bounded_tw_claim_adding_apices_vii} are satisfied. 
   
    Items \ref{item:bounded_tw_claim_adding_apices_i} and \ref{item:bounded_tw_claim_adding_apices_ii} imply~\ref{item:par-bounding:hitting-16} and~\ref{item:par-bounding:nbrs-16}, respectively. 
    Additionally,~\ref{item:nice:tree-partitions} and~\ref{item:bounded_tw_claim_adding_apices_vii} imply~\ref{item:par-bounding:bound-16}, since
    \begin{align*}
        \param_{q}(G,S) \leq b' (q+1) \cdot \max_{x \in V(T)}\param_{bq}(G_x,P_x) \leq \beta(X) \cdot qg(bq).
    \end{align*}
    Altogether, $q \mapsto qg(bq)$ is $(\param,\Apex(\mathcal{X}))$-bounding,
    which concludes the proof.    
\end{proof}

\begin{proof}[Proof of \Cref{thm:main_A(X)_to_T(X)}]
    The statement is clear when $\calX$ contains only the null graph, thus, we assume that $\calX$ has a nonnull graph.
    Suppose that $g$ is $(\param, \Apex(\mathcal{X}))$-bounding.
    Since $\mathcal{X}$ is nonempty and closed under leaf addition, every edgeless graph is a subgraph of a graph in $\mathcal{X} \subseteq \Apex(\mathcal{X})$.
    Hence, since $g$ is $(\param, \Apex(\mathcal{X}))$-bounding, we deduce that $g$ is $(\param, \edgeless)$-bounding.
    Therefore, by \Cref{lemma:weakly_bounding_implies_bounding}
    it suffices to show that $g$ is weakly $(\param, \Tree(\mathcal{X}))$-bounding.
    Furthermore, by~\Cref{cor:X-connected-in-Tree-X-bounding}, it suffices to show that $g$ is weakly $(\param,\Tree'(\calX))$-bounding.
    Suppose that the assertion does not hold, i.e.\ $g$ is not weakly $(\param,\Tree'(\calX))$-bounding.
    This is witnessed by some graph $X \in \Tree'(\calX)$ and a positive integer $k$.
    There is a rooted tree decomposition $\big(T,(W_a \mid a \in V(T))\big)$ of $X$ witnessing the fact that $X \in \Tree'(\mathcal{X})$ (i.e.\ such that $W_a \neq \emptyset$ for every $a \in V(T)$).
    We consider such a triplet $\big(k,X,(T,(W_a \mid a \in V(T)))\big)$
    with $T$ of minimum vertex-height.

    It is easy to check that $T \neq K_1$, as otherwise $X$ is the null graph or a single vertex graph.
    In both cases $X$ can not be a counterexample.
    For every $Y \in \Tree'(\mathcal{X})$ which is not a counterexample,
    we denote by $\alpha_0(Y)$ and $\beta_0(Y)$ integers such that
    every $K_k$-minor-free graph $G$ is $(g,(\param_q\mid q \in \posint),Y,\alpha_0(Y),\beta_0(Y))$-good.
    Let $\alpha'$ and $\beta'(\cdot)$ witness $g$ being $(\param, \Apex(\mathcal{X}))$-bounding for $k$.

    In~\Cref{fig:graphs-X-0-1-2-3}, we illustrate the object constructed in the proof.
    We denote by $r$ the root of $T$, and let $z \in V(X)$ be such that $W_r = \{z\}$.
    Let $T_0$ be the subtree of $T$ induced by $\{a \in V(T) \mid z \in W_a\}$.
    Observe that for every $a \in V(T_0)$, $X[W_a \setminus \{z\}] \in \mathcal{X}$.
    Since $\mathcal{X}$ is closed under disjoint union,
    this implies that $X[\bigcup_{a \in V(T_0)} W_a \setminus \{z\}] \in \mathcal{X}$.
    Let
    \[
        X_0 = X\left[\textstyle\bigcup_{a \in V(T_0)} W_a\right].
    \]
    Let $X_1$ be obtained from $X_0$ by adding for every $x \in V(X_0)$
    a new vertex $\ell_x$ adjacent to $x$.
    Then, since $\mathcal{X}$ is closed under leaf addition and
    $X_1 \setminus \{z\} \in \mathcal{X}$, we have
    $X_1 \in \Apex(\mathcal{X})$.

    For every $x \in V(X_0)$,
    let $\mathcal{R}_x$ be family of all the components $R$ of $T - V(T_0)$
    such that $x \in \bigcup_{a \in V(R)} W_a$,
    and let 
    \[
        X_x = 
        \begin{cases}
            X\left[\bigcup_{R \in \mathcal{R}_x, a \in V(R)} W_a\right] &\text{if $\mathcal{R}_x \neq \emptyset$,} \\
            X[\{x\}] &\text{if $\mathcal{R}_x = \emptyset$.}
        \end{cases}
    \]

    Now, let $X_2$ be obtained from $X$ 
    by identifying all the vertices in $V(X_0)$ into a single vertex~$z_2$.
    By the definition of the forest decomposition witnessing that $X \in \Tree'(\calX)$, for every $x \in V(X_0)$, the graph $X_2$ contains a copy of $X_x$ where $x$ plays the role of $z_2$.
    Then, let $X_3$ be obtained from the union of two copies $X_{2,1}$ and $X_{2,2}$ by identifying the vertices $z_2$ of both copies.

    \begin{figure}[tp]
        \centering
        \includegraphics{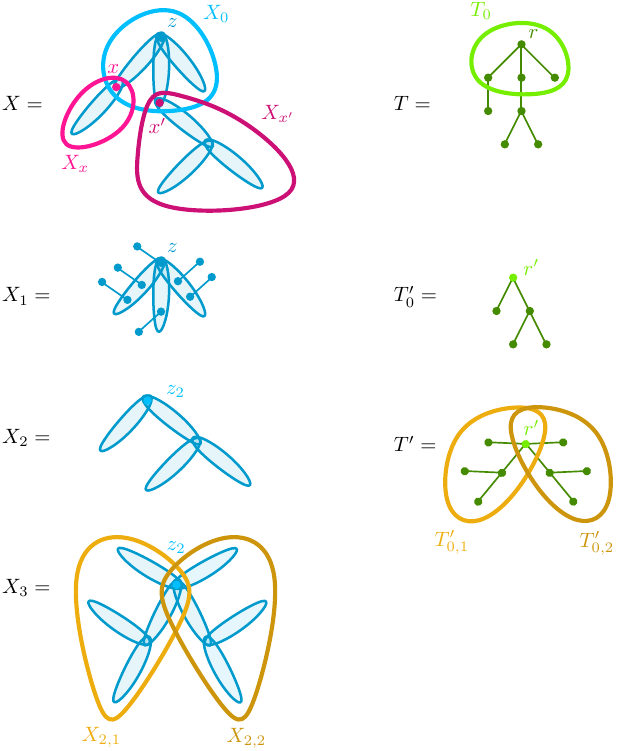}
        \caption{
            An illustration of graphs in the proof of~\Cref{thm:main_A(X)_to_T(X)}.
        }
        \label{fig:graphs-X-0-1-2-3}
    \end{figure}
    
    We now build a rooted tree decomposition $\big(T',(W'_a \mid a \in V(T'))\big)$
    of $X_3$ witnessing the fact that $X_3 \in \Tree'(X)$.
    Also, $T'$ has smaller vertex-height than $T$.
    Let $T'_0$ be obtained from $T$ by contracting $T_0$ into a single vertex $r'$,
    and then let $T'$ be obtained from the union of two copies $T'_{0,1}$ and $T'_{0,2}$
    of $T'_0$ by identifying the vertices $r'$ of both copies.
    Then, for every $a \in V(T'_{0,i})$ for each $i \in \{1,2\}$,
    let $W'_a$ be obtained by adding $z_2$ of $X_3$ to the set of vertices playing roles of $W_a \setminus V(X_0)$ in $X_{2,i}$.
    Then $\big(T',(W_a \mid a \in V(T'))\big)$ is a rooted tree decomposition of $X_3$
    witnessing the fact that $X_3 \in \Tree'(\mathcal{X})$.
    We claim that the vertex-height of $T'$ is smaller than the vertex-height of $T$.
    Indeed, for every neighbor $s$ of $r$ in $T$, since $X$ is connected and $W_s \neq \emptyset$,
    we have $z \in W_s$, and so $s \in V(T_0)$. This proves that $N_T(r) \subseteq V(T_0)$,
    and so $T'$ has smaller vertex-height than $T$.
    It follows that $X_3$ is not a counterexample and we can define
    \[
        \alpha_0(X) = \alpha' + \alpha_0(X_3)
        \quad \text {and} \quad
        \beta_0(X) = \beta'(X_1) + \beta_0(X_3).
    \]

    Let $q$ be a positive integer, 
    let $G$ be a $K_k$-minor-free graph,
    and let $\mathcal{F}$ be a family of connected subgraphs of $G$
    such that there is no $\mathcal{F}$-rich model of $X$ in $G$.
    Let $\mathcal{F}'$ be the family of all the connected subgraphs $H$
    of $G$ such that there is an $\mathcal{F}\vert_{H}$-rich model of $X_3$
    in $H$.

    \begin{figure}[tp]
        \centering
        \includegraphics{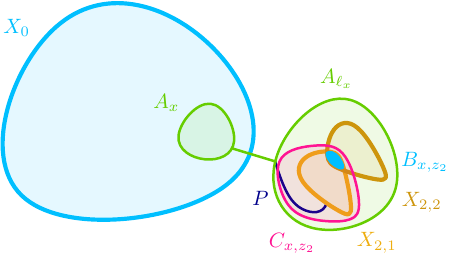}
        \caption{
            An illustration of the construction of the $\mathcal{F}$-rich model $(E_y \mid y \in V(X))$ of $X$ in $G$ in the proof of~\Cref{thm:main_A(X)_to_T(X)}.
        }
        \label{fig:building-model-of-X}
    \end{figure}
    
    We claim that there is no $\mathcal{F}'$-rich model of $X_1$ in $G$.
    Suppose for contradiction that there is an $\mathcal{F}'$-rich model $(A_x \mid x \in V(X_1))$ of $X_1$ in $G$.
    Let $x \in V(X_0)$.
    The following construction is illustrated in~\Cref{fig:building-model-of-X}.
    The graph $G[A_{\ell_x}]$ contains a member of $\mathcal{F}'$, and so, an $\mathcal{F}$-rich model $(B_{x,y} \mid y \in V(X_3))$ of $X_3$.
    Since $G[A_{\ell_x}]$ is connected, there is an \stpath{N(A_x)}{(\bigcup_{y \in V(X_3)} B_{x,y})} $P$ in $G[A_{\ell_x}]$.
    Let $y \in V(X_3)$ be such that $P$ is an \stpath{N(A_x)}{B_{x,y}} path, let $i \in [2]$ be such that $y \in X_{2,i}$, and let $\{i,j\} = [2]$.
    For every $y \in V(X_{2,j})$, let
    \[
        C_{x,y} = 
        \begin{cases}
            V(P) \cup \bigcup_{y' \in V(X_{2,i})} B_{x,y'} &\textrm{if $y =z_2$,} \\
            C_{x,y} &\textrm{otherwise}.
        \end{cases}
    \]
    It follows that $(C_{x,y} \mid y \in V(X_{2,j}))$ is a model of $X_{2,j}$ (which is a copy of $X_2$) in $G[A_{\ell_x}]$ such that $C_{x,z_2}$ contains a vertex adjacent to $A_x$ in $G$.
    By restricting this model to $(V(X_x) \setminus \{x\}) \cup \{z_2\}$,
    and relabeling $z_2$ with $x$,
    we obtain a model $(D_{x,y} \mid y \in V(X_x))$ of $X_x$ in $G[A_{\ell_x}]$
    such that $D_{x,x}$ is adjacent to $A_x$ in $G$.
    Finally, for every $y \in V(X)$, let
    \[
        E_y = 
        \begin{cases}
            A_y \cup D_{y,y} &\textrm{if $y \in V(X_0)$,} \\
            D_{x,y} &\textrm{if $y \not\in V(X_0)$, for $x \in V(X_0)$ such that $y \in V(X_x)$.}
        \end{cases}
    \]
    Then, $(E_y \mid y \in V(X))$ is an $\mathcal{F}$-rich model of $X$ in $G$,
    a contradiction.
    This proves that $G$ has no $\mathcal{F}'$-rich model of $X_1$ in $G$.

    Recall that $\alpha'$ and $\beta'$ witness $g$ being $(\param,\Apex(\calX))$-bounding for $k$.
    Since $G$ has no $\calF'$-rich model of $X_1$ in $G$, there exists $S_0 \subseteq V(G)$ such that
    \begin{enumerate}[label={\normalfont(g\arabic*-0)}]
        \item $S_0 \cap V(F) \neq \emptyset$ for every $F \in \mathcal{F}'$;
        \label{claim:bounded_tw_adding_apices_less_technical_J-final:hitting}
        \item for every component $C$ of $G-S_0$, $N_G(V(C))$ intersects 
        at most $\alpha'$ components of $G-V(C)$; and
        \label{claim:bounded_tw_adding_apices_less_technical_J-final:connected}
        \item $\param_{q}(G,S_0) \leq \beta'(X_1) \cdot g(q)$ \label{claim:bounded_tw_adding_apices_less_technical_J-final:wcol}
    \end{enumerate}
    Since $S_0 \cap V(F) \neq \emptyset$ for every $F \in \mathcal{F}'$, $G-S_0$ has no $\mathcal{F}\vert_{G-S_0}$-rich model of $X_3$.
    Therefore, as $X_3$ is not a counterexample, 
    there exists $S_1 \subset V(G-S_0)$ such that
    \begin{enumerate}[label=(g\arabic*-1)]
        \item $S_1 \cap V(F) \neq \emptyset$ for every $F \in \mathcal{F}\vert_{G-S_0}$; \label{item:claim:TX:hitting'}
        \item for every component $C$ of $(G-S_0)-S_1$, $N_{G-S_0}(V(C))$ intersects at most $\alpha_0(X_3)$ components of $(G-S_0)-V(C)$; and \label{item:claim:TX:component'}
        \item $\param_{q}(G-S_0,S_1) \leq \beta_0(X_3) \cdot g(q)$. \label{item:claim:TX:param'}
    \end{enumerate}
    Finally, let
    \[
        S = S_0 \cup S_1.
    \]
    We now argue 
    \begin{enumerate}[label=(g\arabic*$^\star$)]
        \item $S \cap V(F) \neq \emptyset$ for every $F \in \mathcal{F}\vert_{G}$; \label{item:claim:TX:hitting:star}
        \item for every component $C$ of $G-S$, $N_{G}(V(C))$ intersects at most $\alpha_0(X)$ components of $G-V(C)$; and \label{item:claim:TX:component:star}
        \item $\param_{q}(G,S) \leq \beta_0(X) \cdot g(q)$. \label{item:claim:TX:param:star}
    \end{enumerate}
    First, observe that for every $F \in \mathcal{F}$,
    either $S_0 \cap V(F) \neq \emptyset$ by \ref{claim:bounded_tw_adding_apices_less_technical_J-final:hitting},
    or $F \in \mathcal{F}\vert_{G-S_0}$ and so $S_1 \cap V(F)\neq \emptyset$ by
    \ref{item:claim:TX:hitting'}.
    This proves \ref{item:claim:TX:hitting:star}.

    For every component $C$ of $G-S$,
    $N_G(V(C))$ intersects at most $\alpha_0(X_3)$ components of $(G-S_0)-V(C)$
    by \ref{item:claim:TX:component'},
    and at most $\alpha'$ components of $G-S_0$.
    Therefore, $N_G(V(C))$ intersects at most $\alpha' + \alpha_0(X_3) = \alpha_0(X)$
    components of $G-V(C)$, which proves \ref{item:claim:TX:component:star}.

    Finally,
    \begin{align*}
        \param_q(G,S)
        &\leq \param_q(G,S_0) + \param_q(G-S_0,S_1) &&\textrm{by \ref{item:nice:subadditivity}} \\
        &\leq \beta'(X_1) \cdot g(q) + \beta_0(X_3) \cdot g(q)
        &&\textrm{by \ref{claim:bounded_tw_adding_apices_less_technical_J-final:wcol} and \ref{item:claim:TX:param'}} \\
        &\leq \beta_0(X) \cdot g(q).
    \end{align*}
    This proves \ref{item:claim:TX:component:star}.
    It follows that $X$ is not a counterexample to $g$ being weakly $(\param,\Tree'(\calX))$-bounding, which is a contradiction that concludes the proof.
\end{proof}

\section{Common arguments for base cases}\label{sec:common_base}

Recall that $\edgeless$ is the class of all edgeless graphs.
We begin this section by showing that for certain families of focused parameters $\param$, that the function $q \mapsto 1$ is $(\param,\edgeless)$-bounding (\Cref{lem:edgeless-identity-bounding}).
In~\Cref{sec:base:srtd2}, we study graphs with no $\calF$-rich models of a fixed linear forest.
Namely, we show that for each family of focused parameters $\param$, if a function is $(\param,\edgeless)$-bounding, then this function is also $(\param,\SRt_2)$-bounding (\Cref{thm:R1_implies_S2}), and if a function is $(\param,\Apex(\SRt_2))$-bounding, then it is also $(\param,\SRt_3)$-bounding (\Cref{thm:AS2_implies_S3}).
In~\cref{sec:base:rtd2}, we study graphs with no $\calF$-rich models of a fixed forest.
This will be used in \Cref{sec:fragility,sec:wcol} to show that $q \mapsto \log q$ is $(\param,\mathcal{R}_2)$-bounding,
and in \Cref{sec:centered} to show that $q \mapsto q \log q$ is $(\param,\mathcal{R}_3)$-bounding,
for some families of focused parameters $(\param_q \mid q \in \posint)$.

\begin{lemma}\label{lem:edgeless-identity-bounding}
    Let $(\param_q \mid q \in \posint)$ be a family of focused parameters.
    For every graph $G$, for every $S \subset V(G)$, and for every positive integer $q$, assume that $\param_q(G,S) \leq |S|$ and define
    \[
        \param'_q(G,S) = \frac{1}{\tw(G)+2}  \param_q(G,S).
    \]
    Then, the function $q \mapsto 1$ is $(\param',\edgeless)$-bounding where $\param' = (\param'_q \mid q \in \posint)$.
\end{lemma}

\begin{proof}
    We set $\alpha = 2$ and for every $X \in \edgeless$, we set $\beta(X) = \max\{0,2|V(X)|-3\}$.
    We claim that for every $X \in \edgeless$, every graph $G$ is $(q \mapsto 1, \param',X,\alpha,\beta(X))$-good.
    Let $X \in \edgeless$.
    Let $q$ be a positive integer and let $\calF$ be a family of connected subgraphs of $G$ such that $G$ has no $\calF$-rich model of $X$, or in other words, $G$ does not contain $|V(X)|$ pairwise disjoint members of $\calF$.
    We have to prove that there exists $S \subset V(G)$ such that
    \begin{enumerate}[label={\normalfont (g\arabic*$^\star$)}]
        \item $S \cap V(F)  \neq \emptyset$ for every $F \in \mathcal{F}$; \label{item:edgeless-identity-bounding:hitting}
        \item for every component $C$ of $G-S$, $N_G(V(C))$ intersects at most two components of $G-V(C)$; and \label{item:edgeless-identity-bounding:components}
        \item $\param'_q(G,S) \leq \max\{0,2|V(X)|-3\}$. \label{item:edgeless-identity-bounding:param}
    \end{enumerate}
    If $|V(X)| \leq 1$, then $\mathcal{F} = \emptyset$,
    and so $S=\emptyset$ satisfies \ref{item:edgeless-identity-bounding:hitting}--\ref{item:edgeless-identity-bounding:param}.
    Now assume $|V(X)|\geq 2$.
    Let $\calW$ be a tree decomposition of $G$ of width $\tw(G)$.
    By~\Cref{lemma:natural_tree_decomposition}, there exists a natural tree decomposition of $\calW' = \big(T,(W_x \mid x \in V(T)\big)$ of $G$ of the same width as $\calW$.
    Since $G$ has no $|V(X)|$ pairwise disjoint members of $\calF$, by~\Cref{lemma:helly_property_tree_decomposition}, there exists $Z \subset V(T)$ of size at most $|V(X)| - 1$ such that if $S' = \bigcup_{x \in Z} W_x$, then $S' \cap V(F) \neq \emptyset$ for every $F \in \calF$.
    We root $T$ arbitrarily and define
    \[
        S = \bigcup_{x \in \LCA(T,Z)} W_x.
    \]
    Since $S' \subset S$, we have~\ref{item:edgeless-identity-bounding:hitting}.
    By~\Cref{lemma:increase_X_to_have_small_interfaces}, we have $|\LCA(T,Z)| \leq 2|V(X)| - 3$ and for every component $C$ of $G - Z$, $N_G(V(C))$ intersects at most two components of $G - V(C)$.
    Thus,~\ref{item:edgeless-identity-bounding:components} follows.
    Finally, 
    \[
        \param'_q(G,S) 
        = \frac{1}{\tw(G)+2} \param_q(G,S) 
        \leq \frac{|S|}{\tw(G) + 2} 
        \leq \frac{(2|V(X)| - 3)(\tw(G) + 1)}{\tw(G) +2} 
        \leq 2|V(X)| - 3. 
    \]
    Hence, we also have~\ref{item:edgeless-identity-bounding:param}.
    This completes the proof.
\end{proof}

\subsection{Excluding a linear forest}
\label{sec:base:srtd2}
Recall that $\edgeless$ is the class of all edgeless graphs,
and $\SRt_2$ is the class of all linear forests.
In the following proof, we denote by \defin{$P_\ell$}
the path on the vertex set $[\ell]$ such that
$i$ and $i+1$ are adjacent for every $i \in [\ell-1]$.
This proof of the following theorem is a modification of an argument due to Paul, Protopapas, and Thilikos~\cite[Theorem~3.22]{paul2023universal}.

\begin{theorem}\label{thm:R1_implies_S2}
    Let $(\param_q \mid q \in \posint)$ be a nice family of focused parameters,
    and let $g\colon \posint \to \posint$.
    If $g$ is $(\param,\edgeless)$-bounding,
    then $g$ is $(\param, \SRt_2)$-bounding.
\end{theorem}
\begin{proof}
    Suppose that $g$ is $(\param,\calE)$-bounding.
    By~\Cref{lemma:weakly_bounding_implies_bounding}, it suffices to show that $g$ is weakly $(\param, \SRt_2)$-bounding.
    
    Let $k$ be a positive integer.
    Let $\alpha'$ and $\beta'(\cdot)$ witness that $g$ is $(\param,\edgeless)$-bounding for $k$.
    We can assume that $\beta'$ takes a positive integer $d$ as an argument when $X = \overline{K_d}\in\edgeless$.
    For every $X \in \SRt_2$, let
    \[\alpha(X) = (|V(X)|-1)\alpha' \ \text{ and } \ \beta(X) = (|V(X)| -1)\beta'(2).\]
    We prove that $\alpha(\cdot)$ and $\beta(\cdot)$ witness that $g$ is weakly $(\param,\calS_2)$-bounding for $k$.

    Since every $X \in \SRt_2$ is a subgraph of a path on $|V(X)|$ vertices, it suffices to prove that for every positive integer $\ell$, every $K_k$-minor-free graph $G$ is $(g,\param,P_\ell,(\ell-1)\alpha',(\ell-1)\beta'(2))$-good.
    By~\Cref{lem:reduce-to-G-connected}, we may also assume that $G$ is connected.
    The proof is by induction on $\ell$.
    Let $G$ be a connected $K_k$-minor-free graph, let $q$ be a positive integer, and let $\mathcal{F}$
    be a family of connected subgraphs of $G$ such that $G$ has
    no $\mathcal{F}$-rich model of $P_\ell$.
    We claim that there exists $S \subseteq V(G)$ such that
    \begin{enumerate}[label={\normalfont (g\arabic*$^\star$)}]
        \item $S \cap V(F)  \neq \emptyset$ for every $F \in \mathcal{F}$; \label{item:local_claim_td:hitting}
        \item for every component $C$ of $G-S$, $N_G(V(C))$ intersects at most $(\ell-1)\alpha'$ components of $G-V(C)$; and \label{item:local_claim_td:components}
        \item $\param_q(G,S) \leq (\ell-1) \beta'(2) \cdot g(q)$. \label{item:local_claim_td:param}
    \end{enumerate}
    If $\ell=1$, then $\mathcal{F}$ is empty and the assertion holds with $S=\emptyset$ by \ref{item:nice:at_most_size_of_S}.
    Thus, we suppose that $\ell>1$.

    Let $\mathcal{F}'$ be the family of all the connected subgraphs $H$ of $G$ such that $H$ contains an $\mathcal{F}$-rich model of $P_{\ell-1}$.
    We claim that there are no two disjoint members of $\mathcal{F}'$ in $G$.
    Indeed, suppose to the contrary that $(A_1, \dots, A_{\ell-1})$ and $(B_1, \dots ,B_{\ell-1})$ are two disjoint $\mathcal{F}$-rich models of $P_{\ell-1}$ in $G$.
    Since $G$ is connected, there is a \stpath{(\bigcup_{k\in[\ell-1]} A_k)}{(\bigcup_{k \in [\ell-1]} B_k)} $Q$ in $G$.
    Let $i,j \in [\ell-1]$ be such that $Q$ is a \stpath{A_i}{B_j}. 
    By possibly reversing the orderings $(A_1, \dots, A_{\ell-1})$ and $(B_1, \dots, B_{\ell-1})$, 
    we assume that $j \leq \frac{\ell-1}{2} \leq i$.
    Note that $(A_1, \dots, A_i, (V(Q) \setminus A_i) \cup B_j, \dots, B_{\ell-1})$
    is an $\mathcal{F}$-rich model of $P_{\ell'}$ in $C$, where $\ell' = i + \ell-1-(j-1) \geq \frac{\ell-1}{2} + \ell - \frac{\ell-1}{2} = \ell$.
    This contradicts the fact that $G$ has no $\mathcal{F}$-rich model of $P_\ell$.
    Therefore, there are no two disjoint members of $\mathcal{F}'$.
    
    Since $\alpha'$ and $\beta'(\cdot)$ witness $g$ being $(\param,\edgeless)$-bounding for $k$, there exists $S_{0} \subseteq V(C)$ such that
    \begin{enumerate}[label={\normalfont (g\arabic*-0)}]
        \item $S_{0} \cap V(F) \neq \emptyset$ for every $F \in \mathcal{F}'$; \label{item:base_case_td:S0C:hit} 
        \item for every component $C$ of $G-S_{0}$, $N_G(V(C))$ intersects at most $\alpha'$ components of $G-V(C)$; and \label{item:base_case_td:S0C:components} 
        \item $\param_q(G,S_{0}) \leq \beta'(2) \cdot g(q)$. \label{item:base_case_td:S0C:param} 
    \end{enumerate}
    Now, by \ref{item:base_case_td:S0C:hit},
    there is no $\mathcal{F}$-rich model of $P_{\ell-1}$ in $G - S_{0}$.
    Therefore, by the induction hypothesis applied to $\ell-1$, $q$, $G-S_{0}$, and $\mathcal{F}\vert_{G- S_{0}}$,
    there exists $S_{1} \subseteq V(G)$ such that
    \begin{enumerate}[label={\normalfont (g\arabic*-1)}]
        \item $S_{1} \cap V(F) \neq \emptyset$ for every $F \in \mathcal{F}\vert_{G - S_{0}}$; \label{item:R1_implies_S2:hit'}
        \item for every component $C$ of $(G - S_{0}) - S_{1}$, $N_{G - S_{0}}(V(C))$ intersects at most $(\ell-2)\alpha'$ components of $(G - S_{0})-V(C)$; and \label{item:R1_implies_S2:components'}
        \item $\param_q(G-S_{0}, S_{1}) \leq (\ell-2) \beta'(2) \cdot g(q)$. \label{item:R1_implies_S2:param'}
    \end{enumerate}

    Now, let
    \[
        S = S_0 \cup S_1.
    \]
    We now argue \ref{item:local_claim_td:hitting}, \ref{item:local_claim_td:components}, and \ref{item:local_claim_td:param}.
    First, for every $F \in \mathcal{F}$,
    either $S_0 \cap V(F) \neq \emptyset$ or $F \in \calF\vert_{G - S_0}$, and so, $S_1 \cap F \neq \emptyset$ by~\ref{item:R1_implies_S2:hit'}.
    This proves~\ref{item:local_claim_td:hitting}.
    By \ref{item:base_case_td:S0C:components} and \ref{item:R1_implies_S2:components'}, for every component $C$ of $G - S$, $N_G(V(C))$ intersects at most $\alpha' + (\ell-2)\alpha' = (\ell-1) \alpha'$ components of $G - V(C)$.
    This gives~\ref{item:local_claim_td:components}.
    Finally, 
    \begin{align*}
        \param_q(G,S)&\leq \param_q(G, S_{0}) + \param_q(G-S_{0}, S_{1}) &&\textrm{by \ref{item:nice:subadditivity}} \\
        &\leq (\ell-1)\beta'(2) \cdot g(q) &&\textrm{by \ref{item:base_case_td:S0C:param} and \ref{item:R1_implies_S2:param'}.}
    \end{align*}
    This proves \ref{item:local_claim_td:param}.
\end{proof}

\begin{theorem}\label{thm:AS2_implies_S3}
    Let $(\param_q \mid q \in \posint)$ be a nice family of focused parameters,
    and let $g\colon \posint \to \posint$.
    If $g$ is $(\param, \Apex(\SRt_2))$-bounding,
    then $g$ is $(\param, \SRt_3)$-bounding.
\end{theorem}

\begin{proof}
    Suppose that $g$ is $(\param,\Apex(\SRt_2))$-bounding.
    Since $\edgeless \subseteq \Apex(\SRt_2)$, $g$ is also $(\param, \edgeless)$-bounding.
    Therefore, by~\Cref{lemma:weakly_bounding_implies_bounding}, it suffices to show that $g$ is weakly $(\param, \SRt_3)$-bounding.
    Recall that $\SRt_3 = \Tree(\SRt_2)$.
    By~\Cref{cor:X-connected-in-Tree-X-bounding}, it suffices to show that $g$ is weakly $(\param,\Tree'(\SRt_2))$-bounding.
    Suppose that the assertion does not hold, i.e.\ $g$ is not weakly $(\param,\Tree'(\SRt_2))$-bounding.
    This is witnessed by some graph $X \in \Tree'(\SRt_2)$ and a positive integer $k$.
    There is a rooted tree decomposition $\big(T,(W_a \mid a \in V(T))\big)$ of $X$ witnessing the fact that $X \in \Tree'(\SRt_2)$ (i.e.\ such that $W_a \neq \emptyset$ for every $a \in V(T)$).
    We consider such a triplet $\big(k,X,(T,(W_a \mid a \in V(T)))\big)$ with $T$ of minimum vertex-height.

    We will define $\alpha_0(X)$ and $\beta_0(X)$ such that every $K_k$-minor-free graph $G$ is $(g,(\param_q\mid q \in \posint),X,\alpha_0(X),\beta_0(X))$-good obtaining a contradiction and so completing the proof.
    
        For every $Y \in \Tree'(\mathcal{X})$ which is not a counterexample,
        we denote by $\alpha_0(Y)$ and $\beta_0(Y)$ integers such that
        every $K_k$-minor-free graph $G$ is $(g,(\param_q\mid q \in \posint),Y,\alpha_0(Y),\beta_0(Y))$-good.
    Let $\alpha'$ and $\beta'(\cdot)$ witness $g$ being $(\param,\Apex(\SRt_2))$-bounding for $k$.

    If $T=K_1$, then $X \in \Apex(\SRt_2)$, hence, it suffices to set $\alpha_0(X) = \alpha'$ and $\beta_0(X) = \beta'(X)$.
    Thus, we may suppose that the vertex-height $h$ of $T$ is at least $2$.

        Let $T_0$ be the subtree of $T$ induced by the vertices of $T$ at distance less than $h$ from the root.
        Let $X_0 = X[\bigcup_{a \in V(T_0)} W_a]$.
        For every $x \in V(X_0)$, 
        let 
        \[
            X_x = X \left[\textstyle\bigcup_{a \in V(T) \setminus V(T_0), x \in W_a} W_a\right].
        \]
        See the top part of~\Cref{fig:AS2-to-S3}.
        For every $x \in V(X_0)$, 
        $X_x \in \Apex(\SRt_2)$
        and if $V(X_x) \neq \emptyset$,
        then $x$ has degree at most $3$ in $X_0$.
        Let $\ell$ be the maximum of $|V(X_x)|-1$ over all $x \in V(X_0)$.
        
        Let 
        \[
            \alpha_0(X) = \alpha_0(X_0) + \alpha' 
            \quad \text{and} \quad
            \beta_0(X) = \beta_0(X_0) + \beta'(K_1 \oplus P_{4\ell}).
        \]

        Let $q$ be a positive integer,
        let $G$ be a $K_k$-minor-free graph,
        let $\mathcal{F}$ be a family of connected subgraphs
        of $G$ such that $G$ has no $\mathcal{F}$-rich model of $X$.

        \begin{figure}[tp]
            \centering
            \includegraphics{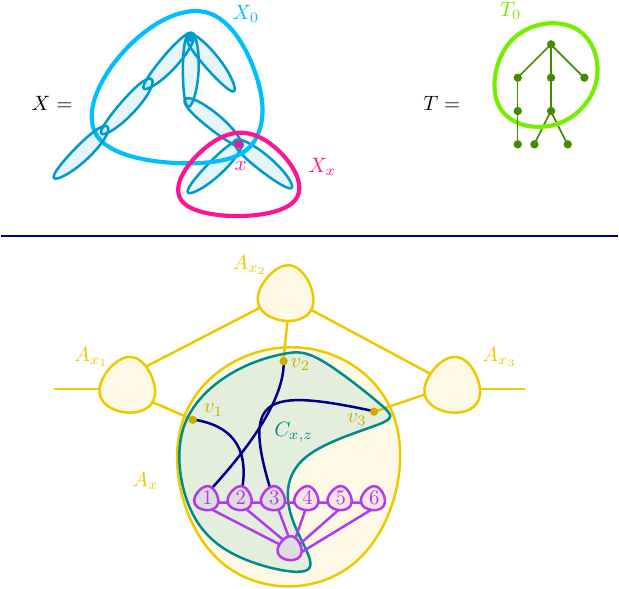}
            \caption{
                An illustration of the proof of~\Cref{thm:AS2_implies_S3}.
                In the top part of the figure, we depict the general setup.
                In the bottom part of the figure, we illustrate the construction of the $\mathcal{F}$-rich model $(E_x \mid x \in V(X))$ of $X$ in $G$.
                In this case, $y_1 = 2$, $y_2 = 1$, and $y_3 = 3$. Moreover, $I = \{4,5,6\}$.
                We obtain a model of $K_1 \oplus P_3$ in $G[A_x]$ such that the branch set corresponding to the vertex of $K_1$ is adjacent to $A_{x_1}$, $A_{x_2}$, and $A_{x_3}$.
            }
            \label{fig:AS2-to-S3}
        \end{figure}

        Let $\mathcal{F}'$ be the family of all the connected subgraphs $H$ of $G$
        such that $H$ contains an $\mathcal{F}$-rich model of $K_1 \oplus P_{4\ell}$.
        We claim that there is no $\mathcal{F}'$-rich model of $X_0$ in $G$.
        Suppose to the contrary otherwise, i.e.\ let $(A_x \mid x \in V(X_0))$ be an $\mathcal{F}'$-rich model of $X_0$ in $G$.
        The following constuction is illustrated in~\Cref{fig:AS2-to-S3}.
        Let $x \in V(X_0)$ such that $V(X_x) \neq \emptyset$.
        The subgraph $G[A_x]$ contains
        an $\mathcal{F}$-rich model $(B_{x,y} \mid y \in V(K_1 \oplus P_{4\ell}))$
        of $K_1 \oplus P_{4\ell}$.
        Let $x_1, \dots, x_d$ be the neighbors of $x$ in $X_0$.
        Recall that $d \leq 3$.
        For every $i \in [d]$,
        $A_x$ contains a neighbor $v_i$ of $A_{x_i}$.
        Since $G[A_x]$ is connected, for every $i \in [d]$, there exists $y_i \in V(K_1 \oplus P_{4\ell})$ and a \stpath{v_i}{y_i} path $Q_i$ in $G[A_x]$ internally disjoint from $\bigcup_{y \in V(K_1 \oplus P_{4\ell})} B_{x,y}$.
        Since $d \leq 3$, by the pigeonhole principle,
        there is an interval $I \subseteq [4\ell]$ of size $\ell$
        such that $I \cap \{y_1, \dots, y_d\} = \emptyset$.
        Let $z$ be the vertex of $K_1$ in $K_1\oplus P_\ell$, where we label the vertices of $P_\ell$ by the elements of $I$.
        For every $y \in I \cup \{z\}$, let
        \[
            C_{x,y} = 
            \begin{cases}
                \bigcup_{i \in [d]} V(Q_i) \cup \bigcup_{y \in V(K_1 \oplus P_{4\ell}) \setminus I} B_{x,y}&\textrm{if $y =z$,} \\
                B_{x,y} &\textrm{otherwise}.
            \end{cases}
        \]
        We obtain an $\mathcal{F}$-rich model $(C_{x,y}\mid y\in V(K_1\oplus P_{\ell}))$ of $K_1\oplus P_{\ell}$ in $G[A_x]$ such that $C_{x,z}$ contains a neighbor of $A_{x_i}$ for every $i\in [d]$.
        Now, since $\ell \geq |V(X_x)|-1$ and $X_x \in \Apex(\SRt_2)$, $K_1 \oplus P_\ell$ contains
        $X_x$ as a subgraph.
        Therefore, we obtain an $\mathcal{F}$-rich
        model $(D_{x,y} \mid y \in V(X_x))$ of $X_x$ in $G[A_x]$ such that $D_{x,x}$
        contains a neighbor of $A_y$ for every $y \in N_X(x)$.
        Finally, for every $x \in V(X)$, let
        \[
            E_x =
            \begin{cases}
                D_{x',x} & \textrm{if $x \not\in V(X_0)$ and there exists $x' \in V(X_0)$ is such that $x \in V(X_{x'})$,} \\
                D_{x,x} & \textrm{if $x \in V(X_0)$ and $V(X_x) \neq \emptyset$, } \\
                A_x & \textrm{if $x \in V(X_0)$ and $V(X_x) = \emptyset$.}
            \end{cases}
        \]
        Then, $(E_x \mid x \in V(X))$ is an $\mathcal{F}$-rich model of $X$ in $G$, which is a contradiction with $G$ having an $\mathcal{F}'$-rich model of $X_0$.
        This proves that there is no $\mathcal{F}'$-rich model of $X_0$ in $G$.

        Since $\big(T_0, (W_a \mid a \in V(T_0))\big)$ is a tree decomposition of $X_0$
        witnessing the fact that $X_0 \in \Tree'(\SRt_2)$,
        and because $T_0$ has smaller vertex-height than $T$,
        we can apply the induction.
        Therefore,
        there exists $S_0 \subseteq V(G)$ such that
        \begin{enumerate}[label={\normalfont (g\arabic*-0)}]
            \item $S_0 \cap V(F) \neq \emptyset$ for every $F \in \mathcal{F}'$; \label{item:claim:AS2_implies_S3:hitting:call}
            \item for every component $C$ of $G-S$, $N_G(V(C))$ intersects at most $\alpha_0(X_0)$ components of $G-V(C)$; and \label{item:claim:AS2_implies_S3:components:call}
            \item $\param_q(G,S) \leq \beta_0(X_0) \cdot g(q)$. \label{item:claim:AS2_implies_S3:param:call}
        \end{enumerate}

        By \ref{item:claim:AS2_implies_S3:hitting:call},
        there is no $\mathcal{F}$-rich model of $K_1 \oplus P_{4\ell}$ in $G-S_0$.
        Since $\alpha'$ and $\beta'(\cdot)$ witness $g$ being $(\param,\Apex(\SRt_2))$-bounding for $k$ and $K_1 \oplus P_{4\ell} \in \Apex(\SRt_2)$,
        there exists $S_1 \subseteq V(G-S_0)$ such that
        \begin{enumerate}[label={\normalfont (g\arabic*-1)}]
            \item $S_1 \cap V(F) \neq \emptyset$ for every $F \in \mathcal{F} \vert_{G-S_0}$; \label{item:claim:AS2_implies_S3:g'':hit}
            \item for every component $C$ of $(G-S_0)-S_1$, $N_{G-S_0}(V(C))$ intersects at most $\alpha'$ components of $(G-S_0)-V(C)$; and \label{item:claim:AS2_implies_S3:g'':components}
            \item $\param_q(G-S_0,S_1) \leq \beta'(K_1 \oplus P_{4\ell}) \cdot g(q)$. \label{item:claim:AS2_implies_S3:g'':param}
        \end{enumerate}

        Let
        \[
            S = S_0 \cup S_1.
        \]
        We now argue that
        \begin{enumerate}[label={\normalfont (g\arabic*$^\star$)}]
            \item $S \cap V(F) \neq \emptyset$ for every $F \in \mathcal{F} \vert_{G}$; \label{item:claim:AS2_implies_S3:gstar:hit}
            \item for every component $C$ of $G-S$, $N_{G}(V(C))$ intersects at most $\alpha_0(X)$ components of $G-V(C)$; and \label{item:claim:AS2_implies_S3:gstar:components}
            \item $\param_q(G,S) \leq \beta_0(X) \cdot g(q)$. \label{item:claim:AS2_implies_S3:gstar:param}
        \end{enumerate}
        First, for every $F \in \mathcal{F}$, either $S_0 \cap V(F) \neq \emptyset$,
        or $F \in \mathcal{F}\vert_{G-S_0}$ and so $S_1 \cap V(F) \neq \emptyset$ by \ref{item:claim:AS2_implies_S3:g'':hit}.
        In both cases, $S \cap V(F) \neq \emptyset$.
        This proves \ref{item:claim:AS2_implies_S3:gstar:hit}.

        Second, for every component $G$ of $G-S$,
        there is a component $C'$ of $G-S_0$ containing $C$.
        Since $N_G(V(C'))$ intersects at most $\alpha_0(X_0)$ components of $G-V(C')$ by \ref{item:claim:AS2_implies_S3:components:call}, 
        and $N_{G-S_0}(V(C))$ intersects at most $\alpha'$ components of $(G-S_0)-V(C)$ by \ref{item:claim:AS2_implies_S3:g'':components},
        we conclude that $N_G(V(C))$ intersects at most $\alpha_0(X_0) + \alpha' = \alpha_0(X)$ components of $G-V(C)$.
        This proves \ref{item:claim:AS2_implies_S3:gstar:components}.

        Finally,
        \begin{align*}
            \param_q(G,S) 
            &\leq \param_q(G,S_0) + \param_q(G-S_0,S_1) &&\textrm{by \ref{item:nice:subadditivity}} \\
            &\leq (\beta_0(X_0) + \beta'(K_1 \oplus P_{4\ell})) \cdot g(q) &&\textrm{by \ref{item:claim:AS2_implies_S3:param:call} and \ref{item:claim:AS2_implies_S3:g'':param}} \\
            &= \beta_0(X) \cdot g(q).
        \end{align*}
        This shows \ref{item:claim:AS2_implies_S3:gstar:param},
        and concludes the proof.
\end{proof}

With \Cref{thm:R1_implies_S2,thm:AS2_implies_S3}, we may reduce proving $(\param,\SRt_t)$-bounding of certain function to showing only $(\param,\edgeless)$-bounding.

\begin{corollary}\label{cor:SRt-edgeless-to-t}
    Let $\param = (\param_q \mid q \in \posint)$ be a nice family of focused parameters and let $g\colon \posint \rightarrow \posint$ be a polynomial function.
    If $g$ is $(\param,\edgeless)$-bounding, then for every positive integer $t$, the function $q \mapsto g(q)\cdot q^{t-2}$ is $(\param,\SRt_t)$-bounding.
\end{corollary}
\begin{proof}
    Assume that $g$ is $(\param,\edgeless)$-bounding.
    The proof is by induction on $t$.
    Recall that $\SRt_1 = \edgeless$ and assume that $t\geq 2$.
    Since $g$ is $(\param,\edgeless)$-bounding, by~\Cref{thm:R1_implies_S2}, $g$ is also $(\param, \SRt_2)$-bounding.
    This shows the inductive claim for $t = 2$.
    Next, suppose that $t \geq 3$.
    By~\Cref{obs:classes-closed}, $\SRt_{t-1}$ is closed under disjoint union, and by \Cref{lemma:Rt_has_coloring_elimination_property}, $\SRt_{t-1}$ has the coloring elimination property.
    \Cref{thm:main_X_to_A(X)} implies that $q\mapsto g(q)\cdot q$ is $(\param, \Apex(\SRt_2))$-bounding.
    Then, \Cref{thm:AS2_implies_S3} implies that  $q \mapsto g(q)\cdot q$ is $(\param, \SRt_3)$-bounding.
    This shows the inductive claim for $t = 3$.
    Finally, suppose that $t \geq 4$.
    In particular, $\SRt_{t-1}$ is closed under leaf addition (see~\Cref{obs:classes-closed} again).
    By induction, $q \mapsto g(q) \cdot q^{t-3}$ is $(\param, \SRt_{t-1})$-bounding.
    Therefore, by \Cref{thm:abstract_induction_main}, 
    $q \mapsto g(q) \cdot q^{t-2}$ is $(\param, \SRt_{t})$-bounding, as desired.
\end{proof}

\subsection{Excluding a forest}
\label{sec:base:rtd2}

In this section, we show structural properties for graphs with no $\mathcal{F}$-rich model of a given forest.
The main source of inspiration for the proofs in this subsection is the paper
of Dujmović, Hickingbotham, Joret, Micek, Morin, and Wood~\cite{Dujmovi2023}.

For all positive integers $h$ and $d$, we denote by \defin{$F_{h,d}$} the (rooted) complete $d$-ary tree of vertex-height $h$. 
In particular, $F_{2,d}$ is the star with $d$ leaves.
For every forest $X$, there exist positive integers $h$ and $d$ such that $X$ is a subgraph of $F_{h,d}$.

\begin{lemma}\label{lemma:bounded_tw_excluding_a_star}
    Let $d$ be a positive integer.
    Let $G$ be a connected graph, let $\mathcal{D} = \big(T,(W_x \mid x \in V(T))\big)$ be a tree decomposition of $G$, and
    let $\mathcal{F}$ be a family of connected subgraphs of $G$
    such that
    $G$ has no $\mathcal{F}$-rich model of $F_{2,d}$.
    For every $u \in V(G)$, 
    there is a set $S \subseteq V(G)$ and a path partition $(P_0, \dots, P_\ell)$ of $(G,S)$ with $P_0 = \{u\}$ such that
    \begin{enumerateOurAlph}
        \item $S \cap V(F) \neq \emptyset$ for every $F \in \mathcal{F}$ and \label{item:bounded_tw_excluding_a_star_i}
        \item for every $i \in [\ell]$, $P_i$ is contained in the union of at most $d$ bags of $\mathcal{D}$. \label{item:bounded_tw_excluding_a_star_iv}
    \end{enumerateOurAlph}
\end{lemma}

\begin{proof}
    We illustrate some objects appearing in the proof in \Cref{fig:excluding-star-tw}.
    We proceed by induction on $|V(G)|$.
    Let $u\in V(G)$.
    If $\mathcal{F}\vert_{G-\{u\}}$ is empty, then it suffices to take $\ell=0$ and $P_0 = \{u\}$. 
    In particular, this is the case for $|V(G)| = 1$.
    Therefore, assume $\mathcal{F}\vert_{G-\{u\}} \neq \emptyset$ and $|V(G)| > 1$.
    Let $\mathcal{F}_0$ be the family of all the connected subgraphs $A$ of $G-\{u\}$ such that $A$ contains some member of $\mathcal{F}$ and
    $V(A) \cap N_G(u) \neq \emptyset$.
    Since $\mathcal{F}\vert_{G-\{u\}} \neq \emptyset$ and $G$ is connected, $\mathcal{F}_0$ is nonempty.
    
    Observe that any collection of $d+1$ pairwise disjoint $A_1, \dots, A_{d+1} \in \mathcal{F}_0$ yields an $\mathcal{F}$-rich model of $F_{2,d}$ in $G$.
    Indeed, it suffices to take $\{u\} \cup A_{d+1}$ as the branch set corresponding to the root of $F_{2,d}$ and $A_1,\ldots,A_d$ as the branch sets of the remaining $d$ vertices of $F_{2,d}$. 
    Therefore, there is at most $d$ pairwise disjoint members of $\mathcal{F}_0$.
    Thus, by \Cref{lemma:helly_property_tree_decomposition} applied to $G - \{u\}$ and $\mathcal{F}_0$, 
    there exists a set  $Z \subseteq V(G-\{u\})$ included in the union of at most $d$ bags of $\mathcal{D}$ 
    such that $Z \cap V(F) \neq \emptyset$ for every $F \in \mathcal{F}_0$.
    Suppose that $Z$ is inclusion-wise minimal for this property.
    Note that since $\mathcal{F}_0 \neq \emptyset$, $Z$ is nonempty.

    \begin{figure}[tp]
        \centering 
        \includegraphics[scale=1]{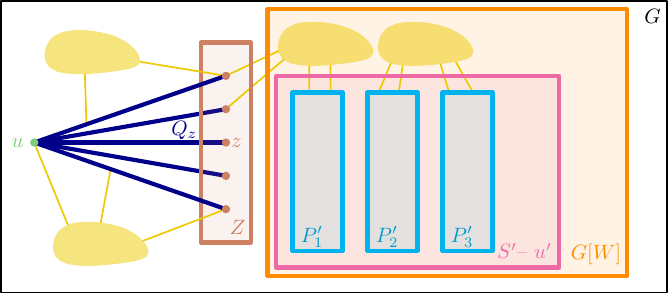} 
        \caption{
        An illustration of the proof of \Cref{lemma:bounded_tw_excluding_a_star}.
        }
        \label{fig:excluding-star-tw}
    \end{figure} 

    Let $\mathcal{C}_0$ be the family of all the components $C$ of $G-(\{u\}\cup Z)$ such that $N_G(u)\cap V(C)=\emptyset$. 
    Let $W = \bigcup_{C \in \mathcal{C}_0}V(C)$.
    Let $z \in Z$.
    By the minimality of $Z$, there exists $A_z \in \mathcal{F}_0$ with $V(A_z) \cap (Z \setminus z) = \emptyset$.
    We have $A_z - \{z\} \notin \mathcal{F}_0$, thus, $z \in V(A_z)$.
    Since $A_z \in \mathcal{F}_0$, $V(A_z) \cap N_G(u) \neq \emptyset$, and so, there is a \stpath{u}{z} $Q_z$ in $G[\{u\} \cup A_z] \subset G - (Z-z)$.
    For every component $C \in \mathcal{C}_0$, we have $V(C) \cap V(Q_z) = \emptyset$, hence, $W \cap V(Q_z) = \emptyset$.

    Let $Q = \bigcup_{z \in Z} V(Q_z)$ and let $G'$ be the graph obtained from $G[W \cup Q]$ by contracting $Q$ into a single vertex $u'$.
    Note that $G'$ is a minor of $G$, and $V(G') = \{u'\} \cup W$.
    Moreover, since $Z \neq \emptyset$, we have $|V(G')| < |V(G)|$.

    For every $x \in V(T)$, let 
    \[
    W'_x =
    \begin{cases}
        W_x & \textrm{if $W_x \subseteq W$} \\
        (W_x \cap W)  \cup \{u'\} & \textrm{otherwise.}
    \end{cases}
    \]
    It follows that $\mathcal{D}' = \big(T,(W'_x \mid x \in V(T)) \big)$ is a tree decomposition of $G'$.
    By the induction hypothesis applied on $G', \mathcal{D}',u',\mathcal{F}\vert_{G[W]}$, there is a set $S' \subseteq V(G')$ and a path partition $(P'_0, \dots, P'_{\ell'})$ of $(G',S')$ with $P'_0=\{u'\}$ such that
    \begin{enumerateOurAlphPrim}
        \item $S' \cap V(F) \neq \emptyset$ for every $F \in \mathcal{F}\vert_{G[W]}$ and \label{item:bounded_tw_excluding_a_star_i'}
        \item for every $i \in [\ell']$, $P'_i$ is contained in the union of at most $d$ bags of $\mathcal{D}'$. \label{item:bounded_tw_excluding_a_star_iv'}
    \end{enumerateOurAlphPrim}

    Let $S = \{u\} \cup Z \cup (S'\setminus \{u'\})$, let $P_0 = \{u\}$, let $P_1 = Z$, and let $P_i = P'_{i-1}$ for every $i \in \{2, \dots, \ell'+1\}$ -- we set $\ell = \ell'+1$.
    We claim that $(P_0, \dots, P_\ell)$ is a path partition of $G[S]$ satisfying \ref{item:bounded_tw_excluding_a_star_i}--\ref{item:bounded_tw_excluding_a_star_iv}, which will complete the proof of the lemma.

    Let $i,j \in \{0, \dots, \ell\}$ with $i<j$ and assume that there is an edge incident to a vertex in $P_i$ and a vertex in $P_j$ in $G$.
    If $i\geq 2$, then $P_i \subseteq P'_{i-1}$ and $P_j = P'_{j-1}$, which implies $|i-j|\leq 1$ since $(P_0',\dots,P_{\ell'}')$ is a path partition of $G'[S']$.
    Otherwise, $i\in \{0,1\}$.
    If $i = 0$, then $j = 1$ since $u$ has no neighbors in $W$.
    If $i=1$, then $j=2$ since $N_{G'}(u') \subset P_1'$.
    For every component $C$ of $G-S$, either $C \cap W = \emptyset$ and so $N(V(C)) \subseteq \{u\} \cup Z = P_0 \cup P_1$,
    or $V(C) \subseteq W$, and so $C$ is a component of $G'-S'$.
    It follows that there exists $i,j \in \{0, \dots, \ell'\}$ with $|i-j| \leq 1$ such that $N_{G'}(V(C)) \subseteq P'_i \cup P'_j$
    and so $N_G(V(C)) \subseteq P_{i+1} \cup P_{j+1}$.
    It follows that $(P_0,\dots,P_\ell)$ is a path partition of $(G,S)$.

    Let $F \in \mathcal{F}$.
    If $(\{u\} \cup Z) \cap V(F) \neq \emptyset$, then $S \cap V(F) \neq \emptyset$.
    Otherwise, $F \subset C$ for some component of $G - (\{u\} \cup Z)$, and in particular, $C \notin \mathcal{F}_0$.
    In this case, $N_G(u) \cap V(C) \neq \emptyset$, hence, $C \in \mathcal{C}_0$, thus, $F \in \mathcal{F}\vert_{G[W]}$, and finally, $S' \cap V(F) \neq \emptyset$ by \ref{item:bounded_tw_excluding_a_star_i'}.
    This proves~\ref{item:bounded_tw_excluding_a_star_i}.

    Finally, $P_0$ is contained in one bag of $\mathcal{D}$, $P_1$ is contained in at most $d$ bags of $\mathcal{D}$ by the definition of $Z$, 
    and for every $i \in \{2, \dots, \ell\}$, $P_i$ is contained in at most $d$ bags of $\mathcal{D}$ by \ref{item:bounded_tw_excluding_a_star_iv'}.
    Therefore,~\ref{item:bounded_tw_excluding_a_star_iv} holds, which concludes the proof of the lemma.
\end{proof}

To generalize \Cref{lemma:bounded_tw_excluding_a_star} to graphs with no $\mathcal{F}$-rich model of $F_{h,d}$ for $h > 2$, we need the following straightforward property.

\begin{lemma}\label{lemma:rooting_a_model_of_Fhd}
    Let $h$ and $d$ be positive integers.
    Let $G$ be a connected graph.
    If there is a model $(B_x \mid x \in V(F_{h,d+1}))$ of $F_{h,d+1}$ in $G$,
    then for every $u \in V(G)$,
    there is a model $(B'_x \mid x \in V(F_{h,d}))$ of $F_{h,d}$ in $G$ such that
    \begin{enumerateOurAlph}
        \item $u \in B'_s$, where $s$ is the root of $F_{h,d}$ and \label{lemma:rooting_a_model_of_Fhd:i}
        \item for every $x \in V(F_{h,d})$, $B_y \subseteq B'_x$ for some $y \in V(F_{h,d+1})$. \label{lemma:rooting_a_model_of_Fhd:ii}
    \end{enumerateOurAlph}
\end{lemma}

\begin{figure}[tp]
    \centering 
    \includegraphics[scale=1]{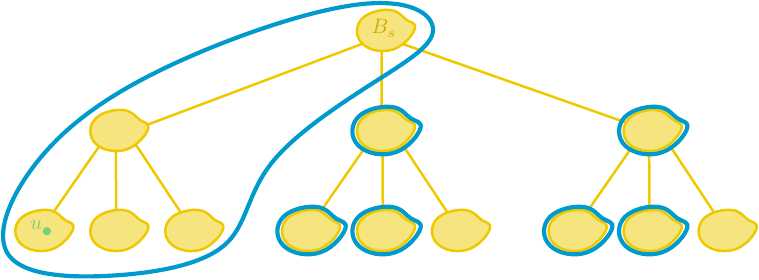} 
    \caption{
    An illustration of the proof of \Cref{lemma:rooting_a_model_of_Fhd}.
    In yellow, we depict a model of $F_{3,3}$.
    In blue, we depict a model of $F_{3,2}$ such that $u$ is in the branch set corresponding to the root.
    }
    \label{fig:augmenting-models-of-trees}
\end{figure}

\begin{proof}
    See~\Cref{fig:augmenting-models-of-trees}.
    In this proof, we treat $F_{h,d}$ as a subgraph of $F_{h,d+1}$ with the same root $s$.
    Suppose that there is a model $(B_x \mid x \in V(F_{h,d+1}))$ of $F_{h,d+1}$ in $G$. 
    Since $G$ is connected, we can assume that $\bigcup_{x \in V(F_{h,d+1})} B_x = V(G)$.
    Let $u\in V(G)$.
    There is a subtree $T'$ of $F_{h,d+1}$ rooted in a child of $s$ such that $u \in \bigcup_{x \in V(T') \cup \{s\}} B_x$.
    Define $B_s' = \bigcup_{x \in V(T') \cup \{s\}} B_x$ and $B_x' = B_x$ for every $x \in V(F_{h,d+1}) \setminus (\{s\} \cup V(T'))$.
    The collection $(B'_x \mid x \in V(F_{h,d}))$ is a model of $F_{h,d}$ in $G$ satisfying \ref{lemma:rooting_a_model_of_Fhd:i} and \ref{lemma:rooting_a_model_of_Fhd:ii}.
\end{proof}

\begin{lemma}\label{lemma:excluding_a_forest_general}
    Let $h$ and $d$ be positive integers.
    Let $G$ be a graph,
    let $\mathcal{D} = \big(T,(W_x \mid x \in V(T))\big)$ be a tree decomposition of $G$,
    and let $\mathcal{F}$ be a family of connected subgraphs of $G$
    such that $G$ has no $\mathcal{F}$-rich model of $F_{h+1,d}$.
    There exist pairwise disjoint $S_1, \dots, S_h \subseteq V(G)$,
    and for every $a \in [h]$, there exists a path partition $(P_{a,0}, \dots, P_{a,\ell_a})$ of $(G-(S_1\cup\dots\cup S_{a-1}), S_a)$ such that
    for $S = \bigcup_{a \in [h]} S_a$,
    \begin{enumerateOurAlph}
        \item $S \cap V(F) \neq \emptyset$ for every $F \in \mathcal{F}$ and \label{item:bounded_tw_excluding_a_tree_i}
        \item for every $a \in [h]$ and 
            for every $i \in [\ell_a]$, $P_{a,i}$ is contained in the union of at most $d + h-1$ bags of $\mathcal{D}$. \label{item:bounded_tw_excluding_a_tree_iv}
    \end{enumerateOurAlph}
\end{lemma}

\begin{proof}
    We proceed by induction on $h$.
    For $h=1$, the result follows from \Cref{lemma:bounded_tw_excluding_a_star}.
    Next, suppose that $h>1$ and that the result holds for $h-1$.

    Let $\mathcal{F}'$ be the family of all the connected subgraphs $H$
    of $G$ containing an $\mathcal{F}$-rich model of $F_{h,d+1}$.
    We claim that there is no $\mathcal{F}'$-rich model of $F_{2,d}$ in $G$.

    Suppose for a contradiction that there is an $\mathcal{F}'$-rich model
    $\big(B_x \mid x \in V(F_{2,d})\big)$ of $F_{2,d}$ in $G$.
    We will find an $\mathcal{F}$-rich model of $F_{h,d}$ in $G$.
    We denote by $c$ the center of the star $F_{2,d}$
    and by $x_1, \dots, x_d$ the leaves of $F_{2,d}$.
    We denote by $r$ the root of $F_{h,d+1}$.
    Let $i \in [d]$.
    Since $(B_x \mid x \in V(F_{2,d}))$ is $\mathcal{F}'$-rich,
    $B_{x_i}$ contains a neighbor $u_i$ of $B_c$,
    and $G[B_{x_i}]$ has an $\mathcal{F}\vert_{G[B_{x_i}]}$-rich model
    $(B^i_x \mid x \in V(F_{h,d+1}))$ of $F_{h,d+1}$.
    By \Cref{lemma:rooting_a_model_of_Fhd},
    there is an $\mathcal{F}$-rich model $\mathcal{M}_i$
    of $F_{h,d}$ in $G[B_{x_i}]$ 
    whose branch set of the root contains $u_i$.
    It follows that the union of the models $\mathcal{M}_1, \dots, \mathcal{M}_d$,
    together with $B_c$ for the branch set of the root,
    yields an $\mathcal{F}$-rich model of $F_{h+1,d}$ in $G$, a contradiction.
    This proves that there is no $\mathcal{F}'$-rich model of $F_{2,d}$ in $G$.

    Therefore, by \Cref{lemma:bounded_tw_excluding_a_star} applied for an arbitrary $u$,
    there is a set $S_1 \subseteq V(G)$ and a path partition $(P_{0}, \dots, P_{\ell})$ of $G[S_1]$ such that
    \begin{enumerate}[label=\ref{lemma:bounded_tw_excluding_a_star}.(\alph*)]
        \item $S_1 \cap V(F) \neq \emptyset$ for every $F \in \mathcal{F}'$ and \label{item:bounded_tw_excluding_a_star_i:call}
        \item for every $i \in [\ell]$, $P_{i}$ is contained in the union of at most $d$ bags of $\mathcal{D}$. \label{item:bounded_tw_excluding_a_star_iv:call}
    \end{enumerate}

    By \ref{item:bounded_tw_excluding_a_star_i:call}, $G - S_1$ has no $\mathcal{F}$-rich model of $F_{h,d+1}$.
    Therefore, by the induction hypothesis applied to $G-S_1$ with the tree decomposition $\mathcal{D}\vert_{G-S_1}$,
    there exist pairwise disjoint $S'_{1}, \dots, S'_{h-1} \subseteq V(G-S_1)$,
    and for every $a \in [h-1]$ there is a path partition $(P'_{a,0}, \dots, P'_{a,\ell_{a}})$ of $((G - S_1)-(S'_{1}\cup\dots\cup S'_{a-1}), S'_{a})$ such that
    for $S' = \bigcup_{a \in [h-1]} S'_{a}$,
    \begin{enumerateOurAlphPrim}
        \item $S' \cap V(F) \neq \emptyset$ for every $F \in \mathcal{F}\vert_{G - S_1}$ and \label{item:bounded_tw_excluding_a_tree_i:induction_call}
        \item for every $a \in [h-1]$,
            for every $j \in [\ell_a]$, $P_{a,j}$ is contained in the union of at most $(d + 1) + (h-1) - 1 = d + h - 1$ bags of $\mathcal{D}\vert_{G-S_1}$. \label{item:bounded_tw_excluding_a_tree_iv:induction_call}
    \end{enumerateOurAlphPrim}

    Now, for every $a \in [h-1]$, let
    $S_{a+1} = S'_{a}$,
    and let $(P_{a+1,0}, \dots, P_{a+1,\ell_{a+1}}) = (P'_{a,0}, \dots, P'_{a,\ell_{a}})$.
    Finally, let $(P_{1,0}, \dots, P_{1,\ell_1}) = (P_0, \dots, P_\ell)$.
    Then, $S_1, \dots, S_{h}$ is as claimed,
    which concludes the proof of the lemma.
\end{proof}

\begin{lemma}\label{lemma:base_case:path_decomposition_partition}
    Let $h$ and $d$ be positive integers.
    Let $G$ be a graph,
    let $\mathcal{D} = \big(T,(W_x \mid x \in V(T))\big)$ be a tree decomposition of $G$,
    and let $\mathcal{F}$ be a family of connected subgraphs of $G$
    such that $G$ has no $\mathcal{F}$-rich model of $F_{h+1,d}$.
    There exist pairwise disjoint $S_1, \dots, S_{h+1} \subseteq V(G)$,
    and for every $a \in [h+1]$ there exists a path partition $(P_{a,0}, \dots, P_{a,\ell_a})$ of $(G-(S_1\cup\dots\cup S_{a-1}), S_a)$ such that
    for $S = \bigcup_{a \in [h+1]} S_a$,
    \begin{enumerateOurAlph}
        \item $S \cap V(F) \neq \emptyset$ for every $F \in \mathcal{F}$; \label{item:base_case:path_decomposition_partition:i}
        \item for every component $C$ of $G-S$, 
            $N_G(V(C))$ intersects at most two components of $G-S$; and \label{item:base_case:path_decomposition_partition:ii}
        \item for every $a \in [h+1]$ and
            for every $i \in [\ell_a]$, $P_{a,i}$ is contained in the union of at most $4h(d + h-1)$ bags of $\mathcal{D}$. \label{item:base_case:path_decomposition_partition:iii}
    \end{enumerateOurAlph}
\end{lemma}

\begin{proof}
    We proceed by induction on $|V(G)|$.

    First, suppose that $G$ is not connected, let $\mathcal{C}$ be the family of all the components of $G$,
    and let $C \in \mathcal{C}$.
    By the induction hypothesis,
    there exist pairwise disjoint $S_{C,1}, \dots, S_{C,h+1} \subseteq V(C)$,
    and for every $a \in [h+1]$ there exists a path partition $(P_{C,a,0}, \dots, P_{C,a,\ell_{C,a}})$ of $(C-(S_{C,1}\cup \dots\cup S_{C,a-1}), S_{C,a})$ such that
    for $S_C = \bigcup_{a \in [h+1]} S_{C,a}$,
    \begin{enumerate}[label={\normalfont (\alph*-C)}]
        \item $S_C \cap V(F) \neq \emptyset$ for every $F \in \mathcal{F}\vert_C$; 
        \item for every component $C'$ of $C-S_C$, 
            $N_C(V(C'))$ intersects at most two components of $C-S_C$; and
        \item for every $a \in [h+1]$ and
            for every $i \in [\ell_{C,a}]$, $P_{C,a,i}$ is contained in the union of at most $4h(d + h-1)$ bags of $\mathcal{D}\vert_C$. 
    \end{enumerate}
    
    Then, for every $a \in [h+1]$, let
    \[
        S_a = \bigcup_{C \in \mathcal{C}} S_{C,a},
    \]
    and let $(P_{a,0}, \dots, P_{a,\ell_a})$ be a concatenation of
    $(P_{C,a,0}, \dots, P_{C,a,\ell_{C,a}})$ over all $C \in \mathcal{C}$.
    These objects satisfy the outcome of the lemma.

    Next, suppose that $G$ is connected.
    By \Cref{lemma:natural_tree_decomposition}, we can assume that $\mathcal{D}$ is a natural tree decomposition.
    By \Cref{lemma:excluding_a_forest_general},
    there exist pairwise disjoint $S_1, \dots, S_h \subseteq V(G)$,
    and for every $a \in [h]$ there exists a path partition $(P_{a,0}, \dots, P_{a,\ell_a})$ of $(G-(S_1\cup \dots\cup S_{a-1}), S_a)$ such that
    for $S' = \bigcup_{a \in [h]} S_a$,
    \begin{enumerate}[label=\ref{lemma:excluding_a_forest_general}.(\alph*)]
        \item $S' \cap V(F) \neq \emptyset$ for every $F \in \mathcal{F}$ and \label{item:bounded_tw_excluding_a_tree_i'}
        \item for every $a \in [h]$,
            for every $i \in [\ell_a]$, $P_{a,i}$ is contained in the union of at most $d + h-1$ bags of $\mathcal{D}$. \label{item:bounded_tw_excluding_a_tree_iv'}
    \end{enumerate}

    Let $C'$ be a component of $G-S'$.
    By \ref{item:bounded_tw_excluding_a_tree_iv'}, $N_G(V(C'))$ is included in the union of at most $2h(d+h-1)$ bags of $\mathcal{D}$.
    By \Cref{lemma:increase_X_to_have_small_interfaces}, 
    there is a set $Z_{C'} \subseteq V(G)$ containing $N_{G}(V(C'))$ which is included in the union of at most $4h(d+h-1)$ bags of $\mathcal{D}$
    such that for every component $C$ of $G-Z_{C'}$, $N_G(V(C))$ intersects at most two components of $G-V(C)$.

    Let
    \[
    S_{h+1} = \bigcup_{C' \in \mathcal{C}'} (Z_{C'} \cap V(C'))
    \]
    where $\mathcal{C}'$ is the family of all the component of $G-S'$
    and let $S = \bigcup_{a \in [h+1]} S_a$.

    First, since $S' \subseteq S$, $S \cap V(F) \neq \emptyset$ for every $F \in \mathcal{F}$, which proves \ref{item:base_case:path_decomposition_partition:i}.
    Second, for every component $C$ of $G-S$, there exists $C' \in \mathcal{C}'$ such that $C \subseteq C'$,
    and so $C$ is a component of $G-Z_{C'}$, which implies that $N_G(V(C))$ intersects at most two components of $G-V(C)$, and so \ref{item:base_case:path_decomposition_partition:ii} holds.
    Finally, for an arbitrary ordering $C_1, \dots, C_{\ell_{h+1}}$ of $\mathcal{C}'$,
    $(Z_{C_1} \cap V(C_1), \dots, Z_{C_{\ell_{h+1}}} \cap V(C_{\ell_{h+1}}))$ 
    is a path partition of $(G-S',S_{h+1})$,
    and for every $i \in [\ell_{h+1}]$, $Z_{C_i} \cap V(C_i)$ is included in the union of at most $4h(d+h-1)$ bags of $\mathcal{D}$.
    This proves \ref{item:base_case:path_decomposition_partition:iii} and concludes the proof of the lemma.
\end{proof}

\section{Fractional treedepth fragility rates}
\label{sec:fragility}

In this section, we prove~\Cref{thm:fragility}.
In~\Cref{ssec:fragility-bd-tw}, we prove the bounds containing treewidth, and in~\Cref{ssec:fragility-general}, we prove the general bounds.
The proofs follow the abstract framework that we have developed.
In particular, the proofs are by induction and the induction step will easily follow from~\Cref{thm:abstract_induction_main}.
Recall the roadmaps to the proofs  in~\Cref{fig:roadmap_wcol_bounded_tw_St,fig:roadmap_wcol_bounded_tw_Rt}.
Also note that we will work with an alternative probabilistic definition of fractional treedepth fragility rates that we introduced in~\Cref{ssec:ftdfr}.

As mentioned in the introduction, known results imply a bound of $\Oh(q^{t-1})$ on $q$th fractional treedepth fragility rates of $K_t$-minor-free graphs.
In~\Cref{sec:kt-minor-free-ftdfr-proof}, we give a proof of this result, which also serves as an illustration of our methods.

\begin{theorem} \label{thm:Kt_minor_free-fragility}
    Let $t$ be an integer with $t\geq 2$.
    There exists an integer $c$ such that for every $K_t$-minor-free graph $G$ and for every positive integer $q$,
    \[
    \frate_q(G)\leq c\cdot q^{t-1}.
    \]
\end{theorem}

\subsection{\texorpdfstring{$K_t$}{Kt}-minor-free graphs} 
\label{sec:kt-minor-free-ftdfr-proof}

As indicated, in this subsection, we prove~\Cref{thm:Kt_minor_free-fragility}.
We begin by stating the first two ingredients of the proof.

\begin{theorem}[{\cite[Theorem~28]{DS_2020}}]\label{thm:fftdr-bd-tw}
    For every positive integer $w$, there exists a constant $c_{\ref{thm:fftdr-bd-tw}}(w)$ such that 
    for every graph $G$ with $\tw(G) \leq w$ and for every positive integer $q$ 
    \[
        \frate_q(G) \leq c_{\ref{thm:fftdr-bd-tw}}(w) \cdot  q^{w}.
    \]
\end{theorem}

\begin{theorem}[{\cite[Theorem~4]{ISW22}}]\label{thm:ISW}
    For all positive integers $t,w$ with $t \geq 2$,
    every $K_t$-minor-free graph $G$ with $\tw(G) \leq w$ admits a vertex-partition $\calP$ such that $G / \calP$ has treewidth at most $t-2$ 
    and each part of $\calP$ has at most $w+1$ elements.\footnote{In other words, there exists a graph $A$ of treewidth at most $t-2$ such that 
    $G \subseteq A \boxtimes K_{w+1}$, where $\boxtimes$ denotes the strong product.} 
\end{theorem}

Theorems~\ref{thm:fftdr-bd-tw} and~\ref{thm:ISW} imply that for every $K_t$-minor-free graph $G$ of bounded treewidth, 
we have $\frate_q(G) = \bigO(q^{t-2})$ as we show below.
Let $t$ and $w$ be positive integers, let $G$ be a $K_t$-minor-free with $\tw(G) \leq w$ and let 
$\calP$ be a partition of $V(G)$ as in~\Cref{thm:ISW}.
Since $\tw(G / \calP) \leq t-2$, we have $\frate_{q}(G/\calP) = \Oh(q^{t-2})$ by~\Cref{thm:fftdr-bd-tw}.
In particular, there exists a random variable $\rvar{Z}'$ over subsets of $\calP$ such that $\td(G/\calP - \rvar{Z}') = \Oh(q^{t-2})$,
and $\rvar{Z}'$ is equipped with a $q$-thin probability distribution.
Let $\rvar{Z} = \bigcup \rvar{Z}'$.
Since every element of $\calP$ has at most $w+1$ vertices of $G$, $ \td(G - \rvar{Z}) \leq (w+1) \cdot \td(G/\calP - \rvar{Z}')$.
It follows that $\rvar{Z}$ witnesses that $\frate_q(G) = \Oh( (w+1) \cdot q^{t-2})$.

Lifting this reasoning to the general $K_t$-minor-free case requires the following technical lemma, which will be useful also in the proof of~\Cref{thm:fragility}.

\begin{lemma}\label{lemma:fftdr_reduction_to_bd_tw}
    For every positive integer $t$, 
    there exists a positive integer $c_{\ref{lemma:fftdr_reduction_to_bd_tw}}(t)$ such that 
    for every $K_t$-minor-free graph $G$ and every positive integer $q$, 
    there exists a random variable $\rvar{Y}$ over subsets of $V(G)$ such that 
    $\tw(G - \rvar{Y}) \leq c_{\ref{lemma:fftdr_reduction_to_bd_tw}}(t) \cdot q$,
    and $\rvar{Y}$ is equipped with a $q$-thin probability distribution.
\end{lemma}

\begin{proof}
    Let $t$ be a positive integer.
    By \Cref{theorem:Kt_free_product_structure_decomposition}, we get an integer $\cLRS(t)$ such that every $K_t$-minor-free graph admits a layered RS-decomposition of width at most $\cLRS(t)$.
    We set $c_{\ref{lemma:fftdr_reduction_to_bd_tw}}(t)=2\cLRS(t)$.
    Let $G$ be a $K_t$-minor-free graph.
    Therefore, $G$ admits a layered RS-decomposition $(T,\calW,\calA,\calD,\calL)$ of width at most $\cLRS(t)$.
    Let $q$ be a positive integer.
    We root $T$ at an arbitrary vertex $s$.
    Let $\calW = \big(T,(W_x \mid x \in V(T))\big)$,
    let $\calA = \big(A_x \mid x \in V(T)\big)$,
    let $\calD = \big(\calD_x \mid x \in V(T)\big)$, and
    let $\calL = \big(\calL_x \mid x \in V(T)\big)$,
    where $\calL_x = (L_{x,i}\mid i \in \mathbb{N})$
    and $\mathcal{D}_x = \big(T_x, (D_{x,z} \mid z \in V(T_x))\big)$ for every $x\in V(T)$.
    For every vertex $u \in V(G)$, let $s_u$ be the root of the subtree of $T$ induced by $\{x \in V(T) \mid u \in W_x\}$.
    For every $x \in V(T)$, let $\rvar{i}(x)$ be a random variable over $\{0, \dots, q-1\}$ with the uniform distribution.
    Let $\rvar{Y}$ be a the random variable defined by
    \[
    \rvar{Y} = \big\{u \in V(G) \mid u \in \textstyle\bigcup_{k \in \NN} L_{s_u,\rvar{i}(s_u) +kq}\big\}.
    \]
    First, observe that for every $u \in V(G)$,
    \begin{align*}
        \Pr[u \in \rvar{Y}] 
        &= 
        \begin{cases}
            0 &\textrm{if $u \in A_{s_u}$,} \\
            \frac{1}{q} &\textrm{otherwise,}
        \end{cases}
    \end{align*}
    and so the distribution of $\rvar{Y}$ is $q$-thin.

    Next, we bound $\tw(G-\rvar{Y})$.
    Let $x \in V(T)$.
    Let $U_x$ be defined as $W_x \cap W_{\parent(T,x)}$ if $x \neq s$, and as $\emptyset$ otherwise.
    Note that $|U_x| \leq \cLRS(t)$.
    For every component $C$ of $\torso_{G,\mathcal{W}}(W_x)-A_x-U_x-\rvar{Y}$,
    for every $z \in V(T_x)$,
    $V(C)$ intersects at most $q-1$ layers of $\mathcal{L}_x$ and so
    \[
        |D_{x,z} \cap V(C)| \leq \cLRS(t) \cdot (q-1).
    \]
    It follows that $\tw(C) \leq \cLRS(t) \cdot (q-1) -1$ as witnessed by $\big(T_x,(D_{x,z} \cap V(C) \mid z \in V(T_x))\big)$.
    Therefore, $\tw(\torso_{G,\mathcal{W}}(W_x)-A_x-U_x-\rvar{Y}) \leq \cLRS(t) \cdot (q-1) -1$,
    and so
    \begin{align*}
        \tw(\torso_{G,\mathcal{W}}(W_x) - \rvar{Y}) 
        &\leq |A_x| + |U_x| + \cLRS(T) \cdot (q-1) - 1 \\
        &\leq \cLRS(t) + \cLRS(t) +  \cLRS(T) \cdot (q-1) - 1\\
        &\leq \cLRS(t) \cdot (q+1) - 1 \leq 2 \cLRS(t) \cdot q.
    \end{align*}
    By \Cref{obs:tw_leq_max_tw_of_a_torso} applied to $G-\rvar{Y}$ and its tree decomposition $\big(T,(W_x\setminus \rvar{Y} \mid x \in V(T))\big)$,
    this implies that $\tw(G-\rvar{Y}) \leq c_{\ref{lemma:fftdr_reduction_to_bd_tw}}(t) \cdot q$ and completes the proof.
\end{proof}

\begin{proof}[Proof of~\Cref{thm:Kt_minor_free-fragility}]
    Let $t$ be an integer with $t \geq 2$, let $c_{\ref{lemma:fftdr_reduction_to_bd_tw}}(\cdot)$ and $c_{\ref{thm:fftdr-bd-tw}}(\cdot)$ be as in~\Cref{lemma:fftdr_reduction_to_bd_tw,thm:fftdr-bd-tw}, respectively.
    We set $c = 2^t \cdot c_{\ref{lemma:fftdr_reduction_to_bd_tw}}(t) \cdot c_{\ref{thm:fftdr-bd-tw}}(t-2)$.
    Let $G$ be a $K_t$-minor-free graph and let $q$ be a positive integer.
    By~\Cref{lemma:fftdr_reduction_to_bd_tw}, there exists a random variable $\rvar{Y}$ over subsets of $V(G)$ such that $\tw(G - \rvar{Y}) \leq c_{\ref{lemma:fftdr_reduction_to_bd_tw}}(t) \cdot 2q$, and $\rvar{Y}$ is equipped with a $2q$-thin probability distribution.
    By~\Cref{thm:ISW}, $G-\rvar{Y}$ admits a vertex-partition $\calP_\rvar{Y}$ such that $(G-\rvar{Y}) / \calP_\rvar{Y}$ has treewidth at most $t-2$ and each part of $\calP_\rvar{Y}$ has at most $c_{\ref{lemma:fftdr_reduction_to_bd_tw}}(t) \cdot 2q+1$ elements.
    Furthermore, by~\Cref{thm:fftdr-bd-tw}, $\frate_{2q}((G-\rvar{Y}) / \calP_\rvar{Y}) \leq c_{\ref{thm:fftdr-bd-tw}}(t-2) \cdot (2q)^{t-2}$, and so, there exists a random variable $\rvar{Z}'$ over subsets of $\calP_\rvar{Y}$ such that $\td((G-\rvar{Y}) / \calP_\rvar{Y} - \rvar{Z}') \leq c_{\ref{thm:fftdr-bd-tw}}(t-2) \cdot (2q)^{t-2}$, and $\rvar{Z}'$ is equipped with a $2q$-thin probability distribution.
    Let 
    \[
        \textstyle \rvar{Z} = \rvar{Y} \cup \bigcup \rvar{Z}'.
    \]
    For every $u \in V(G)$,
    if $\rvar{P}_u \in \mathcal{P}_\rvar{Y}$ is such that $u \in \rvar{P}_u$,
    then
    \begin{align*}
        \Pr[u \in \rvar{Z}]
        \leq \Pr[u \in \rvar{Y}] + \Pr[\rvar{P}_u \in \rvar{Z}']\leq \frac{1}{2q} + \frac{1}{2q} = \frac{1}{q}.
    \end{align*}
    Therefore, the probability distribution of $\rvar{Z}$ is $q$-thin.
    Finally, 
    \[
        \td(G - \rvar{Z}) 
        \leq (c_{\ref{lemma:fftdr_reduction_to_bd_tw}}(t) \cdot 2q+1) \cdot (c_{\ref{thm:fftdr-bd-tw}}(t-2) \cdot (2q)^{t-2}) 
        \leq c \cdot q^{t-1}.
    \]
    This proves that $\frate_q(G) \leq c \cdot q^{t-1}$.
\end{proof}

\subsection{The bounded treewidth case} \label{ssec:fragility-bd-tw}

\Cref{thm:main_frate_bounded_tw_St,thm:main_frate_bounded_tw_Rt} imply the bounds involving treewidth in~\Cref{thm:fragility}.
In~\Cref{lemma:fragility_rate_base_case}, we consider graphs of bounded treewidth excluding $\calF$-rich models of a forest, which gives a necessary base case in \Cref{thm:main_frate_bounded_tw_Rt}.

\begin{theorem}\label{thm:main_frate_bounded_tw_St}
    Let $t$ be an integer with $t \geq 2$.
    For every $X \in \SRt_t$, there exists an integer $c$ such that
    for every positive integer $q$ and
    for every $X$-minor-free graph $G$,
    \[
        \frate_q(G) \leq c \cdot(\tw(G)+1) \cdot q^{t-2}.
    \]
\end{theorem}

\begin{proof}
    For every graph $G$, for every $S \subset V(G)$, and for every positive integer $q$, let
    \[
        \param_q(G,S) = \frac{1}{\tw(G)+2} \frate_q(G,S)
    \]
    By \Cref{ftdfr_are_nice}, the family of focused parameters $\param = (\param_q \mid q \in \posint)$ is nice.
    By~\Cref{lem:edgeless-identity-bounding}, $q \mapsto 1$ is $(\param,\edgeless)$-bounding.
    Therefore, by \Cref{cor:SRt-edgeless-to-t}, 
    $q \mapsto q^{t-2}$ is $(\param, \SRt_{t})$-bounding.
    By~\Cref{lemma:par-bounding-to-bound}, this implies that for every $X \in \SRt_t$, there exists a constant $\beta = \beta(X)$ such that for every positive integer $q$, we have
    \[\frate_q(G) = \frate_q(G,V(G)) \leq \beta \cdot (\tw(G)+2) \cdot q^{t-2} \leq 2\beta \cdot (\tw(G)+1) \cdot q^{t-2}.\qedhere\]
\end{proof}

\begin{lemma}\label{lemma:td_path_partition}
    Let $G$ be a graph and let $S \subseteq V(G)$.
    For every path partition $(P_1, \dots, P_\ell)$ of $(G,S)$,
    \[
        \ftd(G,S) \leq \lceil\log(\ell+1)\rceil \cdot \max_{i \in [\ell]} |P_i|.
    \]
\end{lemma}

\begin{proof}
    Let $(P_1, \dots, P_\ell)$ be a path partition of $(G,S)$,
    and let $N = \max_{i \in [\ell]} |P_i|$.
    It is enough to show that $\ftd(G,S) \leq k \cdot N$ if $\ell = 2^k - 1$.
    We proceed by induction on $k$.
    When $k=1$, the statement is clear.
    Next, suppose that $k>1$.
    Let $S_0 = P_{2^{k-1}}$.
    Note that $S_0$ intersects every path in $G$ between $\bigcup_{i \in [2^{k-1}-1]} P_i$ and $\bigcup_{i \in \{2^{k-1}+1, \dots, 2^k\}} P_i$. Thus, there exists a partition $(V_1,V_2)$ of $V(G - S_0)$ such that $\bigcup_{i \in [2^{k-1}-1]} P_i \subseteq V_1$,
    $\bigcup_{i \in \{2^{k-1}+1, \dots, 2^k-1\}} P_i \subseteq V_2$, and there is no edge between $V_1$ and $V_2$ in $G$.
    We obtain,
    \begin{align*}
        \ftd(G,S) 
        &\leq \ftd(G,S_0) + \ftd(G-S_0,S-S_0) &&\text{by~\Cref{lemma:td_GS_subadditivity}}\\
        & = \max\{\ftd((G-S_0)[V_1], S \cap V_1), \ftd((G-S_0)[V_2], S \cap V_2)\} &&\text{by~\Cref{lemma:td_GS_components}}\\
        &\leq |S_0| + (k-1) \cdot N &&\text{by induction}\\
        &\leq k \cdot N. && \pushright{\mbox{\qedhere}}
    \end{align*}
\end{proof}

\begin{lemma}\label{lemma:fragility_rate_base_case}
    For every positive integer $q$, for every graph $G$, and for every $S \subseteq V(G)$, let
    \[
        \param_q(G,S) = \frac{1}{\tw(G)+2} \frate_q(G,S).
    \]
    The function $q\mapsto \log (q+1)$ is $(\param,\Rt_2)$-bounding for $\param =(\param_q \mid q \in \posint)$.
\end{lemma}

\begin{proof}
    Recall that $\Rt_2$ is the class of all forests, and let $X \in \Rt_2$.
    Let $h$ and $d$ be positive integers such that $X \subseteq F_{h+1,d}$ and let
    \[
        \beta(X) = 8(h+1)h(d+h-1).
    \]
    We will prove that every graph is $(q\mapsto \log (q+1),\param,X,2,\beta(X))$-good.
    Let $G$ be a graph and let $\calF$ be a family of connected subgraphs of $G$ such that $G$ has no $\mathcal{F}$-rich model of $X$.
    Let $\mathcal{D}$ be a tree decomposition of $G$ of width $\tw(G)$.
    Since $G$ does not contain an $\mathcal{F}$-rich model of $X$,
    by \Cref{lemma:base_case:path_decomposition_partition},
    there exist pairwise disjoint sets $S_1, \dots, S_{h+1} \subseteq V(G)$,
    such that for every $a \in [h+1]$,
    there exists a path partition $(P_{a,0}, \dots, P_{a,\ell_a})$ of $(G-(S_1\cup\dots\cup S_{a-1}), S_a)$ such that
    for $S = \bigcup_{a \in [h+1]} S_a$,
    \begin{enumerate}[label=\ref{lemma:base_case:path_decomposition_partition}.(\alph*)]
        \item $S \cap V(F) \neq \emptyset$ for every $F \in \mathcal{F}$; \label{lemma:base_case:path_decomposition_partition-i}
        \item for every component $C$ of $G-S$, 
            $N_G(V(C))$ intersects at most two components of $G-S$; and \label{lemma:base_case:path_decomposition_partition-ii}
        \item for every $a \in [h+1]$, 
            for every $i \in [\ell_a]$, $P_{a,i}$ is contained in the union of at most $4h(d + h-1)$ bags of $\mathcal{D}$. \label{lemma:base_case:path_decomposition_partition-iii}
    \end{enumerate}

    Let $q$ be a positive integer.
    We shall prove that
            \begin{enumerate}[label={\normalfont (g\arabic*$^\star$)}]
            \item $S \cap V(F) \neq \emptyset$ for every $F \in \mathcal{F}$; \label{item:lemma:fragility_rate_base_case:hit}
            \item for every component $C$ of $G-S$, $N_{G}(V(C))$ intersects at most two components of $G-V(C)$; and \label{item:lemma:fragility_rate_base_case:components}
            \item $\param_q(G,S) \leq \beta(X) \cdot \log (q+1)$. \label{item:lemma:fragility_rate_base_case:param}
        \end{enumerate}
    Clearly,~\ref{lemma:base_case:path_decomposition_partition-i} and~\ref{lemma:base_case:path_decomposition_partition-ii} imply~\ref{item:lemma:fragility_rate_base_case:hit} and~\ref{item:lemma:fragility_rate_base_case:components}, respectively.
    It remains to prove~\ref{item:lemma:fragility_rate_base_case:param}.
    
    For every $a \in [h+1]$, let $\rvar{i}(a)$ be a random variable over $\{0, \dots, q-1\}$ with the uniform distribution.
    Now, let
    \[
        \rvar{Y} = \bigcup \{P_{a,i} \mid a \in [h+1], i \in \{0, \dots, \ell_a\}, i \equiv \rvar{i}(a) \bmod q \}.
    \]
    First, observe that for every $u \in S$, $\Pr[u \in \rvar{Y}] = \frac{1}{q}$,
    and so the distribution of $\rvar{Y}$ is $q$-thin.
    Moreover, for every $a \in [h+1]$,
    since $G-(S_1 \cup \dots \cup S_{a-1}) - \rvar{Y}$ has a path partition with $q-1$ parts,
    by \Cref{lemma:td_path_partition} and~\ref{lemma:base_case:path_decomposition_partition-iii}.
    \begin{align*}
        \ftd(G- (S_1 \cup \dots \cup S_{a-1}) - \rvar{Y}, S_a - \rvar{Y})
        &\leq \lceil\log(q-1) + 1\rceil\cdot 4h(d+h-1) \cdot (\tw(G)+1) \\
        &\leq 8h(d+h-1) \cdot (\tw(G)+1) \cdot \log (q+1)     
    \end{align*}
    By repeatedly applying \Cref{lemma:td_GS_subadditivity}, we obtain
    \[
        \ftd(G-\rvar{Y}, S \setminus \rvar{Y})
        \leq 8(h+1)h(d+h-1) \cdot (\tw(G)+1) \cdot \log (q+1) = c \cdot (\tw(G)+1) \cdot \log (q+1).
    \]
    Hence,
    \[
        \frate_q(G,S) \leq c \cdot (\tw(G)+1) \cdot \log (q+1),
    \]
    which concludes the proof the lemma.   
\end{proof}

\begin{theorem}\label{thm:main_frate_bounded_tw_Rt}
    Let $t$ be an integer with $t \geq 2$.
    For every $X \in \Rt_t$, there exists an integer $c$ such that
    for every integer $q$ with $q \geq 2$ and
    for every $X$-minor-free graph $G$,
    \[
        \frate_q(G) \leq c \cdot(\tw(G)+1) \cdot q^{t-2} \log q.
    \]
\end{theorem}

\begin{proof}
    For every graph $G$, for every $S \subset V(G)$, and for every positive integer $q$, let
    \[
        \param_q(G,S) = \frac{1}{\tw(G)+2} \frate_q(G,S)
    \]
    By \Cref{ftdfr_are_nice}, the family of focused parameters $\param = (\param_q \mid q \in \posint)$ is nice.
    We show by induction on $t$ (integers with $t \geq 2$) that $q \mapsto q^{t-2}\log (q+1)$ is $(\param, \Rt_t)$-bounding.
    When $t=2$, this is given by \Cref{lemma:fragility_rate_base_case}.
    Next, suppose that $t \geq 3$.
    By~\Cref{obs:classes-closed}, $\Rt_{t-1}$ is closed under disjoint union and leaf addition, and by \Cref{lemma:Rt_has_coloring_elimination_property}, $\Rt_{t-1}$ has the coloring elimination property.
    By induction, $q \mapsto q^{t-3}\log (q+1)$ is $(\param, \Rt_{t-1})$-bounding.
    Therefore, by \Cref{thm:abstract_induction_main}, 
    $q \mapsto q^{t-2} \log (q+1)$ is $(\param, \Rt_{t})$-bounding, as desired.
    By~\Cref{lemma:par-bounding-to-bound}, this implies that for every $X \in \Rt_t$, there exists a constant $\beta = \beta(X)$ such that for every integer $q$ with $q \geq 2$, we have
    \[\frate_q(G) = \frate_q(G,V(G)) \leq \beta \cdot (\tw(G)+2) \cdot q^{t-2}\log (q+1) \leq 4\beta \cdot (\tw(G)+1) \cdot q^{t-2} \log q.\qedhere\]
\end{proof}

\subsection{The general case}\label{ssec:fragility-general}
In this subsection, we conclude discussing fractional treedepth fragility rates by proving the remaining bounds in~\Cref{thm:fragility}.
In fact, these bounds (\Cref{thm:main_frate_St,thm:main_frate_Rt}) follow easily from the bounds in terms of treewidth (\Cref{thm:main_frate_bounded_tw_St,thm:main_frate_bounded_tw_Rt}, respectively) and~\Cref{lemma:fftdr_reduction_to_bd_tw}.
Note that the proofs of~\Cref{thm:main_frate_St,thm:main_frate_Rt} are almost the same.

\begin{theorem}\label{thm:main_frate_St}
    Let $t$ be an integer with $t \geq 2$
    and let $X \in \SRt_t$.
    There exists an integer $c$ such that
    for every positive integer $q$ and
    for every $X$-minor-free graph $G$,
    \[
        \frate_q(G) \leq c \cdot q^{t-1}.
    \]
\end{theorem}

\begin{proof}
    Let $G$ be an $X$-minor-free graph
    and let $q$ be a positive integer.
    By \Cref{lemma:fftdr_reduction_to_bd_tw} applied with $2q$,
    there exists a random variable $\rvar{Y}_1$ over subsets of $V(G)$
    such that $\tw(G-\rvar{Y}_1) \leq c_{\ref{lemma:fftdr_reduction_to_bd_tw}}(|V(X)|)\cdot 2q$,
    and $\rvar{Y}_1$ is equipped with a $2q$-thin probability distribution.
    By \Cref{thm:main_frate_bounded_tw_St}, there exists
    a positive integer $c_0$ such that
    \[
        \frate_{2q}(G-\rvar{Y}_1) \leq c_0 \cdot (\tw(G-\rvar{Y}_1)+1) \cdot (2q)^{t-2}.
    \]
    Let $\rvar{Y}_2$ be a random variable over subsets of $V(G - \rvar{Y}_1)$ witnessing this inequality and let $\rvar{Y} = \rvar{Y}_1 \cup \rvar{Y}_2$.
    We have,
    \begin{align*}
        \td(G-\rvar{Y}) &\leq  \frate_{2q}(G-\rvar{Y}_1)\\
        &\leq 2^{t-2}c_0 \cdot (\tw(G-\rvar{Y}_1)+1) \cdot q^{t-2}  \\
        &\leq 2^{t-2}c_0 \cdot (c_{\ref{lemma:fftdr_reduction_to_bd_tw}}(|V(X)|)\cdot 2q+1) \cdot q^{t-2} 
        \leq c \cdot q^{t-1} 
    \end{align*}   
    for $c = 2^{t-1} c_0 \cdot (c_{\ref{lemma:fftdr_reduction_to_bd_tw}}(|V(X)|)+1)$.
    Additionally, for every $u \in V(G)$,
    \[\Pr[u \in \rvar{Y}] \leq \Pr[u \in \rvar{Y}_1] + \Pr[u \in \rvar{Y}_2] \leq \frac{1}{2q} + \frac{1}{2q} = \frac{1}{q}\] 
    since $\rvar{Y}_2$ is equipped with a $2q$-thin probability distribution.
    Therefore, $\frate_q(G) \leq c \cdot q^{t-1}$.
\end{proof}

\begin{theorem}\label{thm:main_frate_Rt}
    Let $t$ be a positive integer,
    and let $X \in \Rt_t$.
    There exists an integer $c$ such that
    for every integer $q$ with $q \geq 2$ and
    for every $X$-minor-free graph $G$,
    \[
        \frate_q(G) \leq c \cdot q^{t-1} \log q.
    \]
\end{theorem}

\begin{proof}
    Let $G$ be an $X$-minor-free graph
    and let $q$ be a positive integer.
    By \Cref{lemma:fftdr_reduction_to_bd_tw} applied with $2q$,
    there exists a random variable $\rvar{Y}_1$ over subsets of $V(G)$
    such that $\tw(G-\rvar{Y}_1) \leq c_{\ref{lemma:fftdr_reduction_to_bd_tw}}(|V(X)|)\cdot 2q$,
    and $\rvar{Y}_1$ is equipped with a $2q$-thin probability distribution.
    By \Cref{thm:main_frate_bounded_tw_Rt}, there exists
    a positive integer $c_0$ such that
    \begin{align*}
        \frate_{2q}(G-\rvar{Y}_1) 
        &\leq c_0 \cdot (\tw(G-\rvar{Y}_1)+1) \cdot (2q)^{t-2} \log (2q) \\
        &\leq 2^{t-1}c_0 \cdot (\tw(G-\rvar{Y}_1)+1) \cdot q^{t-2} \log q 
    \end{align*}
    Let $\rvar{Y}_2$ be a random variable over subsets of $V(G - \rvar{Y}_1)$ witnessing this inequality and let $\rvar{Y} = \rvar{Y}_1 \cup \rvar{Y}_2$.
    We have,
    \begin{align*}
        \td(G-\rvar{Y}) &\leq  \frate_{2q}(G-\rvar{Y}_1)\\
        &\leq 2^{t-1}c_0 \cdot (\tw(G-\rvar{Y}_1)+1) \cdot q^{t-2} \log q \\
        &\leq 2^{t-1}c_0 \cdot (c_{\ref{lemma:fftdr_reduction_to_bd_tw}}(|V(X)|)\cdot 2q+1) \cdot q^{t-2} \log q
        \leq c \cdot q^{t-1} \log q
    \end{align*}   
    for $c = 2^{t}c_0 \cdot (c_{\ref{lemma:fftdr_reduction_to_bd_tw}}(|V(X)|)+1)$,
    Additionally, for every $u \in V(G)$,
    \[\Pr[u \in \rvar{Y}] \leq \Pr[u \in \rvar{Y}_1] + \Pr[u \in \rvar{Y}_2] \leq \frac{1}{2q} + \frac{1}{2q} = \frac{1}{q}\] 
    since $\rvar{Y}_2$ is equipped with a $2q$-thin probability distribution.
    Therefore, $\frate_q(G) \leq c \cdot q^{t-1} \log q$.
\end{proof}

\section{Weak coloring numbers}\label{sec:wcol}

In this section, we prove~\Cref{thm:wcol}.
In~\Cref{sec:wcol-tw}, we prove the bounds involving treewidth and in~\Cref{sec:wcol-general}, we prove the remaining bounds.
Recall the roadmaps in~\Cref{fig:roadmap_wcol_bounded_tw_St,fig:roadmap_wcol_bounded_tw_Rt}.

Before we start the proof, we state one more well-known result on weak coloring numbers.
Trees admit linear weak coloring number, see~\Cref{obs:wcol-trees}.
For paths, Joret and Micek~\cite{JM22} proved a logarithmic bound.
We provide a sketch of the proof later in the caption of~\Cref{fig:path}.

\begin{lemma}[{\cite[Theorem~1]{JM22}}]\label{lemma:wcol_paths}
   For every positive integer $q$ and for every path $P$,
   $\wcol_q(P) \leq 2+\lceil \log q \rceil$.
\end{lemma}

\subsection{The bounded treewidth case}\label{sec:wcol-tw}

\Cref{thm:main_wcol_bounded_tw_SRt,thm:main_wcol_bounded_tw_Rt} imply the bounds involving treewidth in~\Cref{thm:wcol}.
In~\Cref{lemma:bounded_tw_excluding_a_star_wcol}, we consider graphs of bounded treewidth excluding $\calF$-rich models of a forest, which gives a necessary base case in \Cref{thm:main_wcol_bounded_tw_Rt}.

\begin{theorem}\label{thm:main_wcol_bounded_tw_SRt}
    Let $t$ be a positive integer.
    For every $X \in \SRt_t$, there exists an integer $c$ such that for every integer $q$ with $q \geq 2$ and
    for every $X$-minor-free graph $G$,
    \[
        \wcol_q(G) \leq c \cdot (\tw(G)+1) \cdot q^{t-2}.
    \]
\end{theorem}

\begin{proof}
    For every graph $G$, for every $S \subset V(G)$, and for every positive integer $q$, let
    \[
        \param_q(G,S) = \frac{1}{\tw(G)+2} \wcol_q(G,S)
    \]
    By \Cref{lemma:wcol-tw-nice}, the family of focused parameters $\param = (\param_q \mid q \in \posint)$ is nice.
    By~\Cref{lem:edgeless-identity-bounding}, $q \mapsto 1$ is $(\param,\edgeless)$-bounding.
    Therefore, by \Cref{cor:SRt-edgeless-to-t}, 
    $q \mapsto q^{t-2}$ is $(\param, \SRt_{t})$-bounding.
    By~\Cref{lemma:par-bounding-to-bound}, this implies that for every $X \in \SRt_t$, there exists a constant $\beta = \beta(X)$ such that for every positive integer $q$, we have
    \[\wcol_q(G) = \wcol_q(G,V(G)) \leq \beta \cdot (\tw(G)+2) \cdot q^{t-2} \leq 2\beta \cdot (\tw(G)+1) \cdot q^{t-2}.\qedhere\]
\end{proof}

\begin{lemma}\label{lemma:bounded_tw_excluding_a_star_wcol}
    For every positive integer $q$, for every graph $G$, and for every $S \subseteq V(G)$, let
    \[
        \param_q(G,S) = \frac{1}{\tw(G)+2} \wcol_q(G,S).
    \]
    The function $q\mapsto \log (q+1)$ is $(\param,\Rt_2)$-bounding for $\param =(\param_q \mid q \in \posint)$.
\end{lemma}

\begin{proof}
    Recall that $\Rt_2$ is the class of all forests, and let $X \in \Rt_2$.
    Let $h$ and $d$ be positive integers such that $X \subseteq F_{h+1,d}$ and let
    \[
        \beta(X) = 24(h+1)h(d+h-1).
    \]
    We will prove that every graph is $(q\mapsto \log (q+1),\param,X,2,\beta(X))$-good.
    Let $G$ be a graph and let $\calF$ be a family of connected subgraphs of $G$ such that $G$ has no $\mathcal{F}$-rich model of $X$.
    Let $\mathcal{D}$ be a tree decomposition of $G$ of width $\tw(G)$.
    Since $G$ does not contain an $\mathcal{F}$-rich model of $X$,
    by \Cref{lemma:base_case:path_decomposition_partition},
    there exist pairwise disjoint sets $S_1, \dots, S_{h+1} \subseteq V(G)$,
    such that for every $a \in [h+1]$,
    there exists a path partition $(P_{a,0}, \dots, P_{a,\ell_a})$ of $(G-(S_1\cup\dots\cup S_{a-1}), S_a)$ such that
    for $S = \bigcup_{a \in [h+1]} S_a$,
    \begin{enumerate}[label=\ref{lemma:base_case:path_decomposition_partition}.(\alph*)]
        \item $S \cap V(F) \neq \emptyset$ for every $F \in \mathcal{F}$; \label{lemma:base_case:path_decomposition_partition-i-wcol}
        \item for every component $C$ of $G-S$, 
            $N_G(V(C))$ intersects at most two components of $G-S$; and \label{lemma:base_case:path_decomposition_partition-ii-wcol}
        \item for every $a \in [h+1]$, 
            for every $i \in [\ell_a]$, $P_{a,i}$ is contained in the union of at most $4h(d + h-1)$ bags of $\mathcal{D}$. \label{lemma:base_case:path_decomposition_partition-iii-wcol}
    \end{enumerate}

    Let $q$ be a positive integer.
    We shall prove that
            \begin{enumerate}[label={\normalfont (g\arabic*$^\star$)}]
            \item $S \cap V(F) \neq \emptyset$ for every $F \in \mathcal{F}$; \label{item:lemma:bounded_tw_excluding_a_star_wcol:hit}
            \item for every component $C$ of $G-S$, $N_{G}(V(C))$ intersects at most two components of $G-V(C)$; and \label{item:lemma:bounded_tw_excluding_a_star_wcol:components}
            \item $\param_q(G,S) \leq \beta(X) \cdot \log (q+1)$. \label{item:lemma:bounded_tw_excluding_a_star_wcol:param}
        \end{enumerate}
    Clearly,~\ref{lemma:base_case:path_decomposition_partition-i-wcol} and~\ref{lemma:base_case:path_decomposition_partition-ii-wcol} imply~\ref{item:lemma:bounded_tw_excluding_a_star_wcol:hit} and~\ref{item:lemma:bounded_tw_excluding_a_star_wcol:components}, respectively.
    It remains to prove~\ref{item:lemma:bounded_tw_excluding_a_star_wcol:param}.
    
    Let $a \in [h+1]$.
    For convenience, we set $P_{a,\ell_a+1} = \emptyset$.
    Consider the path $Q$ with $V(Q) = \{0,\dots,\ell_a+1\}$ where two numbers are connected by an edge whenever they are consecutive.
    Let $\sigma' = (i_0,\dots, i_{\ell_a})$ be an ordering of $\{0,\dots,\ell_a+1\}$ given by \Cref{lemma:wcol_paths}, that is such that $\wcol_q(Q,\sigma') \leq 2+\lceil \log q \rceil \leq 3 \log (q+1)$.
    For each $i \in \{0,\dots,\ell_a+1\}$, let $\sigma_i$ be an arbitrary ordering of $P_i$.
    Let $\sigma$ be the concatenation of $\sigma_{i_0}\dots \sigma_{i_{\ell_a+1}}$ in this order.

    Let $u \in V(G - (S_1 \cup \dots \cup S_{a-1}))$,
    and let $W = \WReach_r[G- (S_1 \cup \dots \cup S_{a-1}),S_a,\sigma,u]$.
    We argue that
    \[
        |W| \leq 24h(d + h-1) \cdot (\tw(G)+1) \cdot \log (q+1).
    \]
    Suppose $W \neq \emptyset$.
    Let $j_u \in V(T)$ be such that if $u \in S_a$, then $u \in P_{j_u}$, and otherwise, $j_u \in \{0,\dots,\ell\}$ is the least value such that $P_{j_u}$ intersects $N_G(C)$, where $C$ is the component of $G-S$ containing $u$.
    Let $A = \WReach_q[Q,\sigma',j_u] \cup \WReach_q[Q,\sigma',j_u+1]$.
    In particular, $|A| \leq 2\cdot \wcol_q(Q,\sigma') \leq 2 \cdot 3 \log (q+1)$.
    Since $(P_{a,0}, \dots, P_{a,\ell_a+1})$ is a path partition of $(G- (S_1 \cup \dots \cup S_{a-1}), S_a)$, we have 
    \(
        W \subset \bigcup_{j \in A} P_j.
    \)
    By~\ref{lemma:base_case:path_decomposition_partition-iii-wcol},
    for every $j \in \{0,\dots,\ell+1\}$, $P_j$ is contained in the union of at most $d$ bags of $\mathcal{D}$, hence $|P_j| \leq 2h(d + h-1)(\tw(G)+1)$.
    It follows that
    \[
        |W| \leq |A| \cdot 4h(d + h-1)(\tw(G)+1) \leq 24h(d + h-1)(\tw(G)+1)\log (q+1).
    \]
    In particular,
    \[
        \wcol_q(G-(S_1\cup\dots\cup S_{a-1}),S_a) \leq  24h(d + h-1) (\tw(G)+1)\log q.
    \]
    Finally, applying \Cref{obs:wcol_union} repeadetly, we obtain
    \begin{align*}
        \wcol_q(G,S) 
        &\leq \sum_{a \in [h+1]} \wcol_q(G-(S_1\cup\dots\cup S_{a-1}),S_a) \\
        &\leq 24(h+1)h(d + h-1)  (\tw(G)+1)\log (q+1)
        = \beta(X) \cdot (\tw(G)+1) \cdot \log (q+1).
    \end{align*}
    This shows \ref{item:lemma:bounded_tw_excluding_a_star_wcol:param} and concludes
    the proof of the lemma.
\end{proof}

\begin{theorem}\label{thm:main_wcol_bounded_tw_Rt}
    Let $t$ be a positive integer.
    For every $X \in \Rt_t$, there exists an integer $c$ such that for every integer $q$ with $q \geq 2$ and
    for every $X$-minor-free graph $G$,
    \[
        \wcol_q(G) \leq c \cdot (\tw(G)+1) \cdot q^{t-2} \log q.
    \]
\end{theorem}

\begin{proof}
    For every graph $G$, for every $S \subset V(G)$, and for every positive integer $q$, let
    \[
        \param_q(G,S) = \frac{1}{\tw(G)+2} \wcol_q(G,S)
    \]
    By \Cref{ftdfr_are_nice}, the family of focused parameters $\param = (\param_q \mid q \in \posint)$ is nice.
    We show by induction on $t$ (integers with $t \geq 2$) that $q \mapsto q^{t-2}\log (q+1)$ is $(\param, \Rt_t)$-bounding.
    When $t=2$, this is given by \Cref{lemma:bounded_tw_excluding_a_star_wcol}.
    Next, suppose that $t \geq 3$.
    By~\Cref{obs:classes-closed}, $\Rt_{t-1}$ is closed under disjoint union and leaf addition, and by \Cref{lemma:Rt_has_coloring_elimination_property}, $\Rt_{t-1}$ has the coloring elimination property.
    By induction, $q \mapsto q^{t-3}\log (q+1)$ is $(\param, \Rt_{t-1})$-bounding.
    Therefore, by \Cref{thm:abstract_induction_main}, 
    $q \mapsto q^{t-2} \log (q+1)$ is $(\param, \Rt_{t})$-bounding, as desired.
    By~\Cref{lemma:par-bounding-to-bound}, this implies that for every $X \in \Rt_t$, there exists a constant $\beta = \beta(X)$ such that for every integer $q$ with $q\geq 2$, we have
    \[\wcol_q(G) = \wcol_q(G,V(G)) \leq \beta \cdot (\tw(G)+2) \cdot q^{t-2}\log (q+1) \leq 4\beta \cdot (\tw(G)+1) \cdot q^{t-2} \log q.\qedhere\]
\end{proof}

\subsection{The general case}\label{sec:wcol-general}

In this subsection, we prove the remaining bounds in~\Cref{thm:wcol}, i.e.\ \Cref{thm:main_wcol_SRt,thm:main_wcol_Rt}.
Since the proofs follow the general framework that we developed, in particular, since the induction step follows from~\Cref{thm:abstract_induction_main}, most of this subsection is devoted to the base cases.
In~\Cref{lemma:excluding_a_star_patch,lemma:excluding_a_star_wcol} we study graphs excluding $\calF$-rich models of a star, and in~\Cref{lem:new_base_case_wcol} we study graphs excluding $\calF$-rich models of general trees.
Note that these involved base cases are necessary for the proof of~\Cref{thm:main_wcol_Rt} in contrast to~\Cref{thm:main_wcol_SRt}.

A key ingredient in the proofs of the bounds in~\Cref{thm:wcol} involving treewidth is the Helly property of graphs of bounded treewidth:~\Cref{lemma:helly_property_tree_decomposition}.
Such a Helly property is false for $K_k$-minor-free graphs.\footnote{Consider the $n \times n$ planar grid $G$ and the family $\calF$ consisting of all the unions of one row and one column in the grid. The packing number of $\calF$ is $1$ but there is no hitting set of $\calF$ with less than $n$ elements. Additionally, each planar grid is $K_5$-minor-free.}
Therefore, we consider a relaxed statement, where instead of bounding the size of a hitting set (see~\Cref{lemma:helly_property_tree_decomposition}.\ref{item:helly_property_tree_decomposition:hit}), we restrict its structure.
Namely, we require it to be the union of a small number of geodesics, see~\Cref{lemma:E-P_in_Kk_minor_free_graphs}.

Before stating this generalized Helly property, let us pause to discuss geodesics in the context of weak coloring numbers.
Let $G$ be a graph.
For all $u,v \in V(G)$, a \defin{\stgeodesic{u}{v}} in $G$ is an \stpath{u}{v} in $G$ of minimum length.
A \defin{geodesic} in $G$ is a path in $G$ that is a \stgeodesic{u}{v} for some $u,v \in V(G)$.
The following folklore statement makes geodesics a fundamental tool in bounding weak coloring numbers. 

\begin{lemma}[{\cite[Lemma 23]{DHHJLMMRW24}}]
\label{lemma:intersection_ball_with_geodesic}
    Let $G$ be a graph and let $q$ be a nonnegative integer. 
    For every geodesic $Q$ in $G$ and for every vertex $v \in V(G)$,
    \[
    |\{u \in V(Q) \mid \dist_G(u,v) \leq q\}| \leq 2q+1.
    \]
\end{lemma}

Combining~\Cref{obs:wcol_union} and \Cref{lemma:intersection_ball_with_geodesic}, we obtain the following.

\begin{obs}\label{obs:geodesics}
    Let $G$ be a graph, let $S \subset V(G)$, let $\ell$ be a positive integer, and let $Q_1,\dots,Q_\ell$ be geodesics in $G$.
    For every nonnegative integer $q$, we have
        \[\wcol_q(G,S \cup V(Q_1) \cup \dots \cup V(Q_\ell)) \leq \wcol_q(G,S) + \ell(2q+1).\]
\end{obs}

The next lemma was proved by Dujmovi\'c et al.~\cite[Lemma~21]{DHHJLMMRW24}.
Its proof relies on the Graph Minor Structure Theorem by Robertson and Seymour.
We use the lemma to obtain the advertised generalized Helly property for $K_k$-minor-free-graphs:~\Cref{lemma:E-P_in_Kk_minor_free_graphs}.

\begin{lemma}[{\cite[Lemma~21]{DHHJLMMRW24}}]\label{lemma:E-P_Kk-free_original_statement}
    For every positive integer $k$, there exists a positive integer $c_{\ref{lemma:E-P_Kk-free_original_statement}}(k)$ such that for every positive integer $d$, for every $K_k$-minor-free graph $G$, for every family $\cgF$ of connected subgraphs of $G$ either
    \begin{enumerate}[label=\normalfont(\arabic*)]
        \item there are $d$ pairwise vertex-disjoint subgraphs in $\mathcal{F}$; or \label{lemma:E-P_Kk-free_original_statement:item:i}
        \item there exists $A \subseteq V(G)$ with $|A| \leq (d-1) c_{\ref{lemma:E-P_Kk-free_original_statement}}(k)$,
            and there exists a subgraph $X$ of $G$ which is the union of at most $(d-1)^2c_{\ref{lemma:E-P_Kk-free_original_statement}}(k)$ geodesics in $G-A$,
            such that for every $F \in \cgF$ we have $(V(X) \cup A) \cap V(F) \neq \emptyset$.  \label{lemma:E-P_Kk-free_original_statement:item:ii}
    \end{enumerate}
\end{lemma}

\begin{corollary}\label{lemma:E-P_in_Kk_minor_free_graphs}
    For all positive integers $k$ and $d$, there exists a positive integer $c_{\ref{lemma:E-P_in_Kk_minor_free_graphs}}(k,d)$ such that for every connected $K_k$-minor-free graph $G$, for every family $\cgF$ of connected subgraphs of $G$, 
    either there are $d$ pairwise vertex-disjoint subgraphs in $\mathcal{F}$, 
    or there exists $S \subseteq V(G)$ such that 
    \begin{enumerateOurAlph}
        \item $S \cap V(F) \neq \emptyset$ for every $F \in \mathcal{F}$;\label{lem:EP:item:hit}
        \item $G[S]$ is connected; and \label{lem:EP:item:connected}
        \item $\wcol_q(G,S) \leq c_{\ref{lemma:E-P_in_Kk_minor_free_graphs}}(k,d) \cdot q$ for every positive integer $q$. \label{lem:EP:item:wcol}
    \end{enumerateOurAlph}
\end{corollary}

\begin{proof}
    Let $c_{\ref{lemma:E-P_in_Kk_minor_free_graphs}}(k,d) = 12(d-1)^2c_{\ref{lemma:E-P_Kk-free_original_statement}}(k)$.
    Let $G$ be a $K_k$-minor-free graph and let $\mathcal{F}$ be a family of connected subgraphs of $G$.
    Suppose that there are no $d$ pairwise disjoint members of $\mathcal{F}$,
    and hence,~\Cref{lemma:E-P_Kk-free_original_statement}.\ref{lemma:E-P_Kk-free_original_statement:item:ii} holds, 
    yielding $A \subset V(G)$ and a subgraph $X$ of $G$ such that $|A| \leq (d-1) c_{\ref{lemma:E-P_Kk-free_original_statement}}(k)$ and     
    $X$ is the union of at most $(d-1)^2 c_{\ref{lemma:E-P_Kk-free_original_statement}}(k)$ geodesics in $G-A$.
    Note that $G[A \cup V(X)]$ has at most $|A| + (d-1)^2 c_{\ref{lemma:E-P_Kk-free_original_statement}}(k)$ components.
    Let $Q_1, \dots, Q_\ell$ be a family of at most 
    $(d-1) c_{\ref{lemma:E-P_Kk-free_original_statement}}(k) + (d-1)^2 c_{\ref{lemma:E-P_Kk-free_original_statement}}(k) - 1$ 
    geodesics in $G$ 
    such that the set 
    \[S = A \cup V(X) \cup \bigcup_{i \in [\ell]} V(Q_i)\] 
    induces a connected subgraph in $G$.
    In particular, $\ell \leq 2(d-1)^2c_{\ref{lemma:E-P_Kk-free_original_statement}}(k)$.
    For every positive integer $q$,
    by~\Cref{obs:geodesics} and since $\wcol_q(G,A \cup V(X)) \leq |A \cup V(X)|$, 
    \begin{align*}
        \wcol_q(G,S) 
        &\leq \wcol_q(G,A \cup V(X)) + \ell \cdot (2q+1)\\
        &\leq |A| + (d-1)^2c_{\ref{lemma:E-P_Kk-free_original_statement}}(k) \cdot (2q+1) + \ell \cdot (2q+1) \\
        &\leq (d-1)c_{\ref{lemma:E-P_Kk-free_original_statement}}(k) + (d-1)^2c_{\ref{lemma:E-P_Kk-free_original_statement}}(k)(2r+1) + 2(d-1)^2c_{\ref{lemma:E-P_Kk-free_original_statement}}(k) \cdot(2q+1)\\
        &\leq 4(d-1)^2c_{\ref{lemma:E-P_Kk-free_original_statement}}(k) \cdot (2q+1) \leq c_{\ref{lemma:E-P_in_Kk_minor_free_graphs}}(k,d) \cdot q. \qedhere
    \end{align*}    
\end{proof}

\begin{theorem}\label{thm:main_wcol_SRt}
    Let $t$ be a positive integer.
    For every $X \in \SRt_t$, there exists an integer $c$ such that for every integer $q$ with $q \geq 2$ and
    for every $X$-minor-free graph $G$,
    \[
        \wcol_q(G) \leq c \cdot q^{t-1}.
    \]
\end{theorem}

\begin{proof}
    By \Cref{lemma:wcol-nice}, the family of focused parameters 
    $\wcol = (\wcol_q \mid q \in \posint)$ is nice.
    By~\Cref{lemma:E-P_in_Kk_minor_free_graphs,obs:wcol_components}, $q \mapsto q$ is $(\wcol,\edgeless)$-bounding.
    Therefore, by \Cref{cor:SRt-edgeless-to-t}, 
    $q \mapsto q^{t-1}$ is $(\wcol, \SRt_{t})$-bounding.
    By~\Cref{lemma:par-bounding-to-bound}, this implies that for every $X \in \SRt_t$, there exists a constant $\beta = \beta(X)$ such that for every positive integer $q$, we have
    \[\wcol_q(G) = \wcol_q(G,V(G)) \leq \beta \cdot q^{t-1}.\qedhere\]
\end{proof}

\begin{lemma}\label{lemma:excluding_a_star_patch}
    Let $k$ and $d$ be positive integers.
    Let $G$ be a connected $K_k$-minor-free graph, and
    let $\mathcal{F}$ be a family of connected subgraphs of $G$
    such that
    $G$ has no $\mathcal{F}$-rich model of $F_{2,d}$.
    For every nonempty $U \subset V(G)$ such that $G[U]$ is connected, 
    there is a path decomposition $(W_0, \dots, W_\ell)$ of $G$ with $\ell \geq 1$
    and sets $R_2, \dots, R_\ell \subseteq V(G)$
    such that
    for $S = U \cup \bigcup_{i\in\{2, \dots,\ell\}} (W_{i-1} \cap W_i)$,
    \begin{enumerateOurAlph}
        \item $W_0 = U$; \label{item:excluding_a_star_patch_i}
        \item $S \cap V(F) \neq \emptyset$ for every $F \in \mathcal{F}$; \label{item:excluding_a_star_patch_ii}
        \item $G[S]$ is connected; \label{item:excluding_a_star_patch_iii}
        \item $G[R_i]$ is connected for every $i \in \{2, \dots, \ell\}$; \label{item:excluding_a_star_patch_iv}
        \item $W_{i-1} \cap W_i \subseteq R_{i} \subseteq \bigcup_{j \in \{0, \dots, i-1\}} W_j$ for every $i \in \{2, \dots, \ell\}$; \label{item:excluding_a_star_patch_v}
        \item $W_{i}$ and $W_{i+2}$ are disjoint for every $i \in \{0, \dots, \ell-2\}$; and \label{item:excluding_a_star_patch_vi}
        \item $\wcol_q(G,R_i) \leq (c_{\ref{lemma:E-P_in_Kk_minor_free_graphs}}(k,d+1)+3) \cdot q$ for every $i \in \{2, \dots, \ell\}$ and for every positive integer~$q$. \label{item:excluding_a_star_patch_vii}
    \end{enumerateOurAlph}
\end{lemma}

\begin{proof}
    The statement of the lemma is visualized in~\Cref{fig:excluding-star-statement}.
    Some of the objects defined in the proof are depicted in \Cref{fig:excluding-star}.
    We proceed by induction on $|V(G)|-|U|$.
    Let $U \subseteq V(G)$ be nonempty such that $G[U]$ is connected.
    If $\mathcal{F}\vert_{G-U} = \emptyset$, then it suffices to take $W_0=W_1=U$, $\ell=1$.
    In particular, this is the case for $U = V(G)$.
    Therefore, assume $|U|<|V(G)|$ and $\mathcal{F}\vert_{G-U} \neq \emptyset$.
    Let $\mathcal{F}_0$ be the family of all the connected subgraphs $A$ of $G-U$
    such that $A$ contains a member of $\mathcal{F}$ and $V(A) \cap N_G(U) \neq \emptyset$.
    We argue that $\mathcal{F}_0 \neq \emptyset$.
    Since $\mathcal{F}\vert_{G-U} \neq \emptyset$,
    there is a component $C$ of $G-U$ containing a member of $\mathcal{F}$.
    Since $G$ is connected, $V(C) \cap N_G(U) \neq \emptyset$ and so $C \in \mathcal{F}_0$.

    \begin{figure}[tp]
        \centering 
        \includegraphics[scale=1]{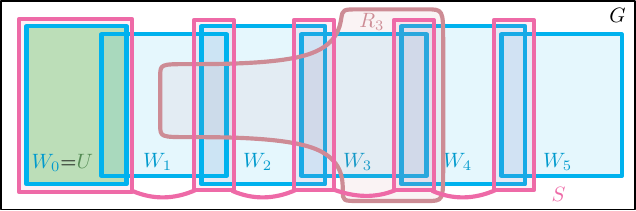} 
        \caption{
        We illustrate the statement of \Cref{lemma:excluding_a_star_patch}.
        The green area is a given set $U$, which should become $W_0$.
        The set $R_3$ has to satisfy $W_3 \cap W_4 \subset R_3 \subset W_0 \cup W_1 \cup W_2 \cup W_3 \cup W_4$ and we want $\wcol_r(G,R_3)$ to be low.
        The set $S$ is connected.
        We do not depict all the sets $R_i$ for readability.
        }
        \label{fig:excluding-star-statement}
    \end{figure} 

    Observe that any collection of $d+1$ pairwise disjoint $A_1, \dots, A_{d+1} \in \mathcal{F}_0$ yields
    an $\mathcal{F}$-rich model of $F_{2,d}$. Indeed, it suffices to take $U \cup A_{d+1}$ as the branch set 
    corresponding to the root of $F_{2,d}$ and $A_1, \dots, A_d$ as the branch sets of the remaining $d$ vertices
    of $F_{2,d}$.
    Therefore, there are no $d+1$ pairwise disjoint members of $\mathcal{F}_0$,
    and thus, by Lemma~\ref{lemma:E-P_in_Kk_minor_free_graphs}
    applied to $G$ and $\mathcal{F}_0$, there exists a set $S_0 \subseteq V(G)$ such that
    \begin{enumerate}[label={\normalfont\ref{lemma:E-P_in_Kk_minor_free_graphs}.(\makebox[\mywidth]{\alph*})}]
        \item $S_0 \cap V(F) \neq \emptyset$ for every $F \in \mathcal{F}_0$; \label{lem:EP:item:hit'_patch}
        \item $G[S_0]$ is connected; and \label{lem:EP:item:connected'_patch}
        \item $\wcol_q(G,S_0) \leq c_{\ref{lemma:E-P_in_Kk_minor_free_graphs}}(k,d+1) \cdot q$ for every positive integer $q$. \label{lem:EP:item:wcol'_patch} 
    \end{enumerate}
    Since $\mathcal{F}_0 \neq \emptyset$, we have $S_0 \setminus U \neq \emptyset$.
    Let $Q$ be a \stgeodesic{U}{S_0} in $G$ (possibly just a one-vertex path), and let $S_1 = S_0 \cup V(Q)$.
    Note that by~\ref{lem:EP:item:connected'_patch}, $G[S_1]$ is connected.

    \begin{figure}[tp]
        \centering 
        \includegraphics[scale=1]{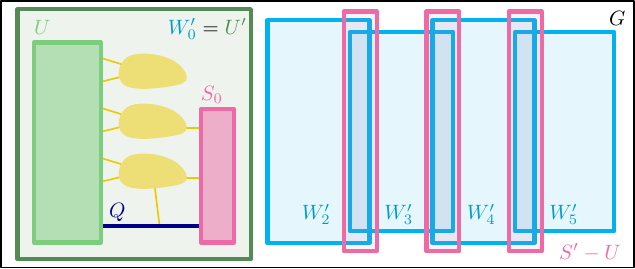} 
        \caption{
        An illustration of the objects considered in \Cref{lemma:excluding_a_star_patch}.
        Note that $S_0$ may intersect $U$.
        }
        \label{fig:excluding-star}
    \end{figure} 

    Let $\mathcal{C}_0$ be the family of all the components $C$ of $G-U-S_1$ such that $N_G(U) \cap V(C) = \emptyset$.
    Let $U' = V(G) \setminus \bigcup_{C \in \mathcal{C}_0} V(C)$.
    Observe that $|U'| > |U|$ since $S_0 \setminus U \neq \emptyset$ and $U'$ contains $U \cup S_0$.
    Let $\mathcal{F}' = \{F \in \mathcal{F} \mid V(F) \cap U' = \emptyset\}$.
    By the induction hypothesis applied to $G$, $\mathcal{F}'$ and $U'$,
    there is a path decomposition $(W'_0, \dots, W'_{\ell'})$ of $G$
    and sets $R'_2, \dots, R'_{\ell'} \subseteq V(G)$
    such that
    for $S' = U' \cup \bigcup_{i\in\{2, \dots,\ell'\}} (W'_{i-1} \cap W'_i)$,
    \begin{enumerateOurAlphPrim}
        \item $W'_0 = U'$; \label{item:excluding_a_star_patch_i'}
        \item $S' \cap V(F) \neq \emptyset$ for every $F \in \mathcal{F}'$; \label{item:excluding_a_star_patch_ii'}
        \item $G[S']$ is connected; \label{item:excluding_a_star_patch_iii'}
        \item $G[R'_i]$ is connected for every $i \in \{2, \dots, \ell'\}$; \label{item:excluding_a_star_patch_iv'}
        \item $W'_{i-1} \cap W'_i \subseteq R'_i \subseteq \bigcup_{j \in \{0, \dots, i-1\}} W'_j$ for every $i \in \{2, \dots, \ell'\}$; \label{item:excluding_a_star_patch_v'}
        \item $W'_{i}$ and $W'_{i+2}$ are disjoint for every $i \in \{0, \dots, \ell'-2\}$; and \label{item:excluding_a_star_patch_vi'}
        \item $\wcol_q(G,R'_i) \leq (c_{\ref{lemma:E-P_in_Kk_minor_free_graphs}}(k,d+1)+3) \cdot q$ for every $i \in \{2, \dots, \ell\}$ and for every positive integer~$q$. \label{item:excluding_a_star_patch_vii'}
    \end{enumerateOurAlphPrim}
    
    Let $\ell = \ell'+1$, $W_0 = U$, $W_1 = U'$, $W_2 = (W'_1 \setminus U') \cup (S_1 \setminus U)$, $W_i = W'_{i-1}$ for every $i \in \{3, \dots, \ell\}$,
    $R_2 = S_1$, and $R_i = R'_{i-1}$ for every $i \in \{3, \dots, \ell\}$.
    Note that~\ref{item:excluding_a_star_patch_i} holds by construction.
    We claim that $(W_0, \dots, W_\ell)$ is a path decomposition of $G$
    and \ref{item:excluding_a_star_patch_ii}--\ref{item:excluding_a_star_patch_vii} hold, which completes the proof of the lemma.

    Let $u \in V(G)$.
    We claim that $I = \{i \in \{0, \dots, \ell\} \mid u \in W_i\}$ is an interval.
    Since $(W'_0, \dots, W'_{\ell'})$ is a path decomposition of $G$,
    $I' = \{i \in \{0, \dots, \ell'\} \mid u \in W'_i\}$ is an interval.
    If $u \not \in U' = W'_0$, then $I =  \{i \in \{2, \dots, \ell\} \mid u \in W'_{i-1}\} = \{i+1 \mid i \in I'\}$,
    which is an interval too.
    Next, suppose that $u \in U'$, and so $0 \in I'$.
    If $u \not\in S_1 \setminus U$, then $u \not \in W_2$ and $u \not\in W'_i$ for every $i \in \{2, \dots, \ell'\}$ by \ref{item:excluding_a_star_patch_i'} and \ref{item:excluding_a_star_patch_vi'}.
    Hence
    $I = \{0,1\}$ if $u \in U$, and $I = \{1\}$ otherwise,
    which is an interval in both cases.
    If $u \in S_1 \setminus U$, then $u \not\in U = W_0$,
    and so $I = \{1\} \cup \{i+1 \mid i \in I' \setminus \{0\}\}$, which is an interval.
    This proves that $I$ is an interval.

    Let $uv$ be an edge of $G$.
    We claim that there exists $i \in \{0, \dots, \ell\}$ such that $u,v \in W_i$.
    If there exists $i'\in \{2,\dots,\ell'\}$ such that $u,v \in W'_{i'}$, then $u,v \in W'_{i'} = W_{i'-1}$ and we are done.
    Next, suppose that $u$ and $v$ are not both in $W'_i$ for every $i \in \{2, \dots, \ell\}$.
    Since $(W'_0, \dots, W'_{\ell'})$ is a path decomposition of $G$,
    there exists $i' \in \{0, 1\}$ such that $u,v \in W'_{i'}$.
    If $i' = 0$, then $u,v \in W'_0 = U' = W_1$.
    Now suppose that $u$ and $v$ are not both in $W'_0=U'$; in particular, $i' = 1$.
    Without loss of generality, assume that $v \notin U'$.
    It follows that $v \in W_1' \setminus U' \subset W_2$.
    Let $C$ be the component of $v$ in $G-U-S_1$.
    Since $v \not\in U'$, $C$ belongs to $\mathcal{C}_0$, and so $N_G(V(C)) \cap U = \emptyset$.
    It follows that $N_G(V(C)) \cap U = \emptyset$, and so, $u \in S_1 \setminus U$.
    Therefore, $u \in W_2$, which concludes the claim.
    Furthermore, we obtained that $(W_0, \dots, W_\ell)$ is a path decomposition of $G$.

    We now prove \ref{item:excluding_a_star_patch_ii}.
    Consider $F \in \mathcal{F}$.
    If $F$ intersects $U$, then $S \cap V(F) \neq \emptyset$ since $U \subseteq S$.
    If $F$ intersects $S_1 \setminus U$, then $F$ intersects $W_1 \cap W_2 \subseteq S$.
    Now suppose that $F$ is disjoint from $U \cup S_1$.
    Let $C$ be the component containing $F$ in $G-U-S_1$.
    Since $C$ is disjoint from $S_0 \subseteq S_1$, by \ref{lem:EP:item:hit'_patch}, $C$ is not a member of $\mathcal{F}_0$.
    This implies that $N_G(U) \cap V(C) = \emptyset$, and thus, $C \in \mathcal{C}_0$.
    In particular, $F$ is disjoint from $U'$, and so, $F \in \mathcal{F}'$.
    By \ref{item:excluding_a_star_patch_ii'}, $S' \cap V(F) \neq \emptyset$,
    hence, there exists $i \in \{2, \dots, \ell'\}$ such that $V(F)$ intersects $W'_{i-1} \cap W'_i$.
    It follows that $W'_{i-1} \cap W'_i = W_{i} \cap W_{i+1}$ and so $S \cap V(F) \neq \emptyset$.
    This proves \ref{item:excluding_a_star_patch_ii}.

Let us pause to underline a simple observation that follows directly from the construction,~\ref{item:excluding_a_star_patch_i'}, and~\ref{item:excluding_a_star_patch_vi'}:
    \begin{enumerate}[label=($\star$)]
        \setcounter{enumi}{6}
        \item for every $i \in \{3,\dots,\ell\}$, we have $W_{i-1} \cap W_i = W_{i-2}' \cap W_{i-1}'$. \label{star:intersections}
    \end{enumerate}

    By \ref{item:excluding_a_star_patch_iii'}, $S' = U' \cup \bigcup_{i \in \{2, \dots, \ell'\}} (W'_{i-1} \cap W'_i)$ induces a connected subgraph of $G$.
    In particular, every component of $G[S']-U'$ has a neighbor in $N_G(V(G)\setminus U') \subseteq S_1$.
    Since $G[U \cup S_1]$ is connected, it follows that $(S' \setminus U') \cup U \cup S_1$ induces a connected subgraph of $G$.
    However, $S = (S' \setminus U') \cup U \cup S_1$ by~\ref{star:intersections}, which yields~\ref{item:excluding_a_star_patch_iii}.

    For every $i \in \{3, \dots, \ell\}$,
    $R_i = R'_{i-1}$ induces a connected subgraph of $G$ by \ref{item:excluding_a_star_patch_iv'}, and $R_2 = S_1$ induces a connected subgraph of $G$ by definition, hence~\ref{item:excluding_a_star_patch_iv} follows.

    For the proof of~\ref{item:excluding_a_star_patch_v}, first, observe that by construction,
    $\bigcup_{j \in \{0, \dots, i-1\}} W_j = \bigcup_{j \in \{0, \dots, i-2\}} W'_j$ for every $i \in \{2, \dots, \ell\}$.
    In particular, it follows that
    $R_i = R'_{i-1} \subseteq \bigcup_{j \in \{0, \dots, i-2\}} W'_j =\bigcup_{j \in \{0, \dots, i-1\}} W_j$ for every $i \in \{3,\dots,\ell\}$ by \ref{item:excluding_a_star_patch_v'}.
    Moreover, $R_2 = S_1 \subseteq W_1$.
    It remains to show that $W_{i-1} \cap W_i \subseteq R_i$ for every $i \in \{2, \dots, \ell'\}$.
    For $i=2$,
    $W_1 \cap W_2 = S_1 \setminus U \subseteq R_2$.
    For $i \in \{3, \dots, \ell\}$, $W_{i-1} \cap W_i = W'_{i-2} \cap W'_{i-1} \subseteq R'_{i-1} = R'_i$ by~\ref{star:intersections} and~\ref{item:excluding_a_star_patch_v'}.
    This gives \ref{item:excluding_a_star_patch_v}.

    For every $i \in \{3, \dots, \ell-2\}$,
    $W_i \cap W_{i+2} = W'_{i-1} \cap W'_{i+1} = \emptyset$ by \ref{item:excluding_a_star_patch_vi'}.
    Moreover, $W'_3=W_4$ is disjoint from $S_1 \setminus U \subseteq W'_0$ by \ref{item:excluding_a_star_patch_vi'}.
    Hence $W_2 \cap W_4 = W'_1 \cap W'_3 = \emptyset$ by \ref{item:excluding_a_star_patch_vi'}.
    Similarly, $W'_2=W_3$ is disjoint from $U' = W'_0$ by \ref{item:excluding_a_star_patch_vi'}.
    Hence $W_1 \cap W_3 = W'_0 \cap W'_2 = \emptyset$.
    Finally, $W_0 \cap W_2 = \emptyset$ by construction, and so, \ref{item:excluding_a_star_patch_vi} holds.

    It remains to show \ref{item:excluding_a_star_patch_vii}.
    First, for every $i \in \{3, \dots, \ell\}$,
    $R_i = R'_{i-1}$ and so $\wcol_q(G,R_i) \leq (c_{\ref{lemma:E-P_in_Kk_minor_free_graphs}}(k,d+1)+3)\cdot q$ for every positive integer $q$ by \ref{item:excluding_a_star_patch_vii'}.
    Moreover, $R_2 = S_1 = V(Q) \cup S_0$.
    Hence by \ref{lem:EP:item:wcol'_patch} and \Cref{obs:geodesics},
    \[
    \wcol_q(G,R_2) \leq c_{\ref{lemma:E-P_in_Kk_minor_free_graphs}}(k,d+1) \cdot q + (2q+1) \leq (c_{\ref{lemma:E-P_in_Kk_minor_free_graphs}}(k,d+1)+3) \cdot q
    \]
    for every positive integer $q$.
    This shows that \ref{item:excluding_a_star_patch_vii} holds, which concludes the proof of the lemma.
\end{proof}

\begin{figure}[tp]
    \centering 
    \includegraphics[scale=1]{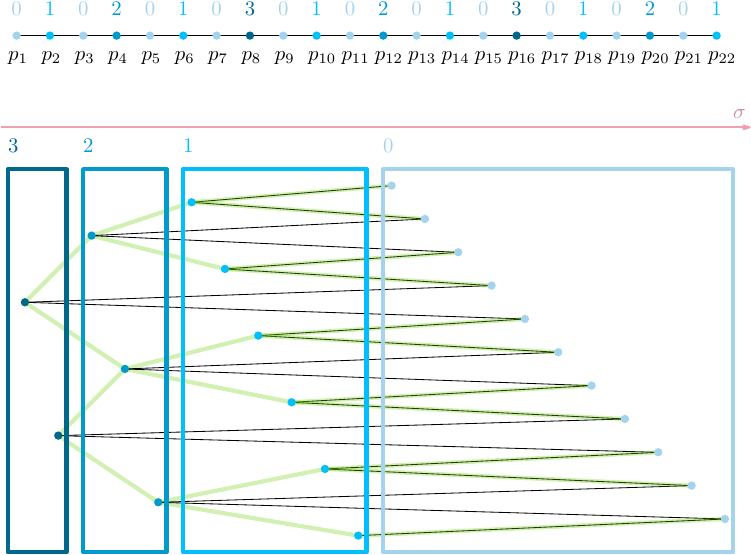} 
    \caption{
        Consider the path $P = p_1\cdots p_{22}$.
        Let $r = 7$ and $s = 3 = \left\lceil \log q \right\rceil$.
        We claim that $\wcol_q(P) \leq 2+s$.
        We mark every eighth vertex with $3$ (this corresponds to the set $I_3$ in the proof of~\Cref{lemma:excluding_a_star_wcol}).
        Then, we mark each fourth unmarked vertex with $2$ (the set $I_2 \setminus I_3$), every second unmarked vertex with $1$ (the set $I_1 \setminus I_2$), and finally all remaining vertices with $0$ (the set $I_0 \setminus I_1$).
        Next, we construct an ordering on the vertices $p_1,\dots,p_{22}$.
        First, preserving the original ordering in the path, we put the vertices marked with $3$, then similarly the ones marked with $2$, with $1$, and with $0$.
        It suffices to argue that for every $u \in V(P)$, we have $|\WReach_q[P,\sigma,u]| \leq 5$.
        We build an auxiliary tree (depicted in green in the figure), where for every $i \in \{3,2,1\}$, we connect every vertex marked with $i$ with the closest vertices in the path marked with $i-1$.
        It is not hard to see that $\WReach_q[P,\sigma,u]$ consists only of the ancestors of $u$ in the auxiliary tree.
    }
    \label{fig:path}
\end{figure} 

The proof of the next lemma follows ideas used to show that $\wcol_q(P)\leq 2+\lceil \log q \rceil$ 
for every path $P$ and every positive integer $q$~\cite[Theorem~1]{JM22}.
See \Cref{fig:path}.

\begin{lemma}\label{lemma:excluding_a_star_wcol}
    Let $k$ and $d$ be positive integers.
    For every connected $K_k$-minor-free graph $G$, 
    for every family $\mathcal{F}$ of connected subgraphs of $G$, 
    if $G$ has no $\mathcal{F}$-rich model of $F_{2,d}$, then
    there is a set $S \subseteq V(G)$ such that
    \begin{enumerate}[label={\normalfont(\makebox[\mywidth]{\alph*})}]
        \item $S \cap V(F) \neq \emptyset$ for every $F \in \mathcal{F}$;\label{lemma:excluding_a_star_patch_wcol:assertion:hitting}
        \item $G[S]$ is connected; and \label{lemma:excluding_a_star_patch_wcol:assertion:connected}
        \item $\wcol_q(G,S) \leq 5(c_{\ref{lemma:E-P_in_Kk_minor_free_graphs}}(k,d+1)+3) \cdot q \log q$ for every integer $q$ with $q \geq 2$. \label{lemma:excluding_a_star_patch_wcol:assertion:wcol}
    \end{enumerate}
\end{lemma}

\begin{proof}   
    Let $q$ be an integer with $q \geq 2$,
    let $G$ be a connected $K_k$-minor-free graph, 
    let $\mathcal{F}$ be a family of connected subgraphs of $G$, 
    and suppose that $G$ has no $\mathcal{F}$-rich model of $F_{2,d}$.
    Let $U$ be an arbitrary singleton of a vertex in $G$.
    \Cref{lemma:excluding_a_star_patch} applied to $G$, $\mathcal{F}$, and $U$ gives
    a path decomposition $(W_0, \dots, W_\ell)$
    and sets $R_2, \dots, R_\ell \subseteq V(G)$
    such that for $S' = U \cup \bigcup_{i \in \{2, \dots, \ell\}} (W_{i-1} \cap W_i)$,
    \begin{enumerate}[label={\normalfont\ref{lemma:excluding_a_star_patch}.(\makebox[\mywidth]{\alph*})}]
        \item $W_0 = U$; \label{item:excluding_a_star_recall_i}
        \item $S' \cap V(F) \neq \emptyset$ for every $F \in \mathcal{F}$; \label{item:excluding_a_star_recall_ii}
        \item $G[S']$ is connected; \label{item:excluding_a_star_recall_iii}
        \item $G[R_i]$ is connected for every $i \in \{1, \dots, \ell\}$; \label{item:excluding_a_star_recall_iv}
        \item $W_{i-1} \cap W_i \subseteq R_{i} \subseteq \bigcup_{j \in \{0, \dots, i-1\}} W_j$ for every $i \in \{2, \dots, \ell\}$; \label{item:excluding_a_star_recall_v}
        \item $W_{i}$ and $W_{i+2}$ are disjoint for every $i \in \{0, \dots, \ell-2\}$; and \label{item:excluding_a_star_recall_vi}
        \item $\wcol_q(G,R_i) \leq (c_{\ref{lemma:E-P_in_Kk_minor_free_graphs}}(k,d+1)+3) \cdot q$ for every $i \in \{2, \dots, \ell\}$. \label{item:excluding_a_star_recall_vii}
    \end{enumerate}

    For convenience, we set $R_1 = U$.

    Let $s = \lceil \log(q+1) \rceil$.
    For every $i \in \{0, \dots, s\}$, let $I_i = \{i \in \{1, \dots, \ell\} \mid j = 0 \mod 2^i\}$.
    We construct recursively families $\set{R'_j}_{j \in \{1,\dots,\ell\}}$ and 
    $\set{S_i}_{i \in \{0,\dots,s\}}$ of subsets of $V(G)$ and a family $\{\sigma_j\}_{j \in \{1,\dots,\ell\}}$ such that 
    $\sigma_j$ is an ordering of $R_j'$ for every $j \in \{1, \dots, \ell\}$.
    For every $j \in I_s$, let 
    \[ 
        R_j' = R_j \setminus \textstyle\bigcup_{a \in  \{0,\dots , j-2^s-1\}} W_{a}
    \]
    and let $S_s = \bigcup_{ j \in I_s} R'_j$.
    Let $j \in I_s$.
    If $j < 2 \cdot 2^s$, then $j = 2^s$ and $R_j' = R_j$, and so by~\ref{item:excluding_a_star_recall_vii}, $\wcol_q(G,R_j') \leq (c_{\ref{lemma:E-P_in_Kk_minor_free_graphs}}(k,d+1)+3) \cdot q$.
    Now assume that $j \geq 2 \cdot 2^s$.
    Since $(W_0, \dots, W_\ell)$ is a path decomposition of $G$, $W_{j-2^s-1} \cap W_{j-2^s}$ separates $\bigcup_{a \in \{0, \dots, j-2^s-1\}} W_a$ and $\bigcup_{a \in \{j-2^s, \dots, \ell\}} W_a$ in $G$.
    Since $W_{j-2^s-1} \cap W_{j-2^s} \subseteq R_{j-2^s}'$ (by \ref{item:excluding_a_star_recall_v}), by \Cref{obs:wcol_components}, we obtain
    \begin{align*}
        \wcol_q\left(G-\textstyle\bigcup_{a \in \{1, \dots, j-2^s\} \cap I_s} R'_a, R'_j\right) 
        &= \wcol_q\left(G - (W_{j-2^s-1} \cap W_{j-2^s}), R'_j\right) \\
        &= \wcol_q\left(G - \textstyle\bigcup_{a \in \{0, \dots ,j-2^s-1\}} W_a, R'_j\right).
    \end{align*}
    Finally,
    \begin{align*}
        \wcol_q\left(G - \textstyle\bigcup_{a \in \{0, \dots ,j-2^s-1\}} W_a, R'_j\right)
        &\leq \wcol_q(G,R_j) && \textrm{by \Cref{obs:monotone}} \\
        &\leq (c_{\ref{lemma:E-P_in_Kk_minor_free_graphs}}(k,d+1)+3)\cdot q && \textrm{by \ref{item:excluding_a_star_recall_vii}}.
    \end{align*}
    Let $\sigma_j$ be an ordering of $R'_j$ such that
    \[
        \wcol_q\left(G - \textstyle\bigcup_{a \in \{1, \dots, j-2^s\} \cap I_s} R'_a, R'_j, \sigma_j\right) \leq (c_{\ref{lemma:E-P_in_Kk_minor_free_graphs}}(k,d+1)+3)\cdot q.
    \]

    Next, let $i \in \{0, \dots, s-1\}$ and assume that $S_{i+1}$ is defined.
    Now, for every $j \in I_i \setminus I_{i+1}$, let
    \[
    R'_j = \left(R_j \setminus \textstyle\bigcup_{a \in \{0, \dots, j-2^i-1\}} W_a \right) \setminus S_{i+1},
    \]
    and let $S_i = \bigcup_{j \in I_i} R'_j$.
    Note that $S_{i+1} \subset S_i$.
    Also note that for every $j \in I_i$, $W_{j-1} \cap W_j \subset R_j'$ by~\ref{item:excluding_a_star_recall_vi}.
    Let $j \in I_i \setminus I_{i+1}$.
    We have $j-2^i \in I_{i+1}$, and therefore,
    $W_{j-2^i-1} \cap W_{j-2^i} \subseteq R_{j-2^i}' \subset S_{i+1}$.
    Since $(W_0 ,\dots, W_\ell)$ is a path decomposition of~$G$,
    $W_{j-2^i-1} \cap W_{j-2^i}$ separates $\bigcup_{a \in \{0, \dots, j-2^i-1\}} W_a$ and $\bigcup_{a \in \{j-2^i, \dots, \ell\}} W_a$ in $G$.
    It follows by \Cref{obs:wcol_components} that
    \begin{align*}
    \wcol_q(G-S_{i+1}, R'_j)
    &= \wcol_q\left(G-(W_{j-2^i-1} \cap W_{j-2^i}), R'_j\right) \\
    &= \wcol_q\left(G-\textstyle\bigcup_{a \in \{0, \dots, j-2^i-1\}} W_a,R'_j\right).
    \end{align*}
    Furthermore,
    \begin{align*}
        \wcol_q\left(G-\textstyle\bigcup_{a \in \{0, \dots, j-2^i-1\}} W_a,R'_j\right) 
        &\leq \wcol_q(G,R_j) && \textrm{by \Cref{obs:monotone}} \\
        &\leq (c_{\ref{lemma:E-P_in_Kk_minor_free_graphs}}(k,d+1)+3) \cdot q && \textrm{by \ref{item:excluding_a_star_recall_vii}}.
    \end{align*}
    Let $\sigma_j$ be an ordering of $R'_j$ such that
    \[
    \wcol_q(G-S_{i+1}, R'_j, \sigma_j) \leq (c_{\ref{lemma:E-P_in_Kk_minor_free_graphs}}(k,d+1)+3) \cdot q.
    \]

    We define $S = S_0$.
    Now, it suffices to show that~\ref{lemma:excluding_a_star_patch_wcol:assertion:hitting}--\ref{lemma:excluding_a_star_patch_wcol:assertion:wcol} hold.
    Since $S' \subseteq S$, \ref{lemma:excluding_a_star_patch_wcol:assertion:hitting} holds by~\ref{item:excluding_a_star_recall_ii}.
    
    Recall that $G[S']$ is connected by \ref{item:excluding_a_star_recall_iii}.
    Next, let $C$ be a component of $G[R_j']$ for some fixed $j \in I_s$.
    If $V(C) \cap (W_{j-1} \cap W_j) \neq \emptyset$, then $V(C) \cap S' \neq \emptyset$, and so, $G[S' \cup V(C)]$ is connected.
    Thus, assume that $V(C) \cap (W_{j-1} \cap W_j) = \emptyset$.
    However, $W_{j-1} \cap W_j \subset R_j'$, hence, $C$ has a neighbor in $W_{j-2^s-1}$, in particular, in $W_{j-2^s-1} \cap W_{j-2^s} \subset S'$.
    Hence, again $G[S' \cup V(C)]$ is connected.
    In particular, we have just proved that $G[S' \cup S_s]$ is connected.
    Next, suppose that $G[S' \cup S_{i+1}]$ is connected for some $i \in \{0,\dots,s-1\}$.
    Let $C$ be a component of $G[R_j']$ for some fixed $j \in I_i$.
    If $V(C) \cap (W_{j-1} \cap W_j) \neq \emptyset$, then $V(C) \cap S' \neq \emptyset$, and so, $G[S' \cup V(C)]$ is connected.
    Thus, assume that $V(C) \cap (W_{j-1} \cap W_j) = \emptyset$.
    However, $W_{j-1} \cap W_j \subset R_j'$, hence, $C$ has a neighbor in $W_{j-2^i-1} \cup S_{i+1}$, in particular, in $(W_{j-2^i-1} \cap W_{j-2^i}) \cup S_{i+1} \subset S' \cup S_{i+1}$.
    Hence, $G[S' \cup S_{i+1} \cup V(C)]$ is connected.
    Finally, $G[S'\cup S_0] = G[S]$ is connected, which yields~\ref{lemma:excluding_a_star_patch_wcol:assertion:connected}.
    
    The sets $\{R'_j\}_{j \in \{2, \dots, \ell\}}$ are pairwise disjoint, and they partition $S$.
    Let $\sigma$ be an ordering of $S$ such that
    \begin{enumerate}
        \item $\sigma$ extends $\sigma_j$, for every $j \in \{1, \dots, \ell\}$;
        \item for every $j,j' \in I_s$ with $j<j'$, for all $u \in R'_j$ and $v \in R'_{j'}$, $u<_\sigma v$; and
        \item for every $i \in \{0, \dots, s-1\}$, for all $u \in S_{i+1}$ and $v \in S_i \setminus S_{i+1}$, $u<_\sigma v$.
    \end{enumerate}
    For convenience, let $R'_0 = \emptyset$ and $W_j=R'_j=\emptyset$ for every integer $j$ with $j > \ell$.

    We now show \ref{lemma:excluding_a_star_patch_wcol:assertion:wcol}.
    Let $u \in V(G)$.
    We will show that $|\WReach_q[G,S,\sigma,u]| \leq 5(c_{\ref{lemma:E-P_in_Kk_minor_free_graphs}}(k,d+1)+3)\cdot q \log q$.
    Let $j_u \in \{0, \dots, \ell\}$ be minimum such that $u \in W_{j_u}$.
    We claim that
    \[
    |\WReach_q[G,S,\sigma,u] \cap S_s| \leq 2(c_{\ref{lemma:E-P_in_Kk_minor_free_graphs}}(k,d+1)+3) \cdot q.
    \]
    Let $\alpha = \max\{0\} \cup \{a \in I_s \mid a \leq j_u\}$,
    and let $\beta = \alpha + 2^s$.
    Thus, if $\beta \leq \ell$, then $\beta \in I_s$.
    Next, we argue that
    \[
    \WReach_q[G,S,\sigma,u] \cap S_s \subseteq R'_{\alpha} \cup R'_{\beta}.
    \]
    Suppose to the contrary that there is a vertex $v \in \WReach_q[G,S,\sigma,u] \cap S_s$ with
    $v \not\in R'_{\alpha} \cup R'_{\beta}$.
    Let $\gamma \in I_s \setminus \{\alpha, \beta\}$ be such that $v \in R'_\gamma$.
    Then either $\gamma < \alpha$, or $\gamma > \beta$.
    First assume that $\gamma < \alpha$.
    Since $R'_\gamma \subseteq \bigcup_{a \in \{0, \dots, \gamma-1\}} W_a$ and because $(W_0, \dots, W_\ell)$ is a path decomposition of $G$,
    every \stpath{u}{v} in $G$ intersects $W_{a-1} \cap W_{a}$ for each $a \in \{\gamma, \dots, j_u\}$.
    Since $(W_{a-1} \cap W_a)_{a \in \{1, \dots, \ell\}}$ are pairwise disjoint,
    we deduce that $\dist_G(u,v) \geq j_u - \gamma \geq \alpha - \gamma \geq 2^s > q$,
    which contradicts the fact that $v \in \WReach_r[G,S,\sigma,u]$.
    Finally, assume $\gamma > \beta$.
    Note that $\gamma \leq \ell$ since $R'_\gamma \neq \emptyset$ as $v \in R'_\gamma$.
    Since $R'_\gamma \subseteq \bigcup_{a \in \{\gamma-2^s, \dots, \gamma-1\}} W_a$,
    and because $(W_0, \dots, W_\ell)$ is a path decomposition of $G$,
    every \stpath{u}{v} in $G$ intersects $W_{\beta-1} \cap W_\beta$.
    However, for every $w \in W_{\beta-1} \cap W_\beta$, we have $w<_\sigma v$, thus, $v \not\in \WReach_q[G,S,\sigma,u]$, which is a contradiction.
    We obtain that $\WReach_q[G,S,\sigma,u] \cap S_s \subseteq R'_{\alpha} \cup R'_{\beta}$.

    For every $\epsilon \in \{\alpha,\beta\}$, by definition of $\sigma$, we have
    \begin{align*}
        \WReach_q[G,S,\sigma,u] \cap R'_\epsilon
        &\subseteq \WReach_q\left[G - \textstyle\bigcup_{a \in I_s \cap \{1,\dots,\alpha-1\}}  R'_a, R'_\epsilon, \sigma_\epsilon, u\right] \\
        &\subseteq \WReach_q\left[G - \textstyle\bigcup_{a \in \{1, \dots, j_u-2^s\} \cap I_s} R'_a, R'_\epsilon, \sigma_\epsilon, u\right]
    \end{align*}
    and therefore,
    \[
    |\WReach_q[G,S,\sigma,u] \cap R'_\epsilon| \leq (c_{\ref{lemma:E-P_in_Kk_minor_free_graphs}}(k,d+1)+3) \cdot q.
    \]
    In particular,
    \begin{align*}
    |\WReach_q[G,S,\sigma,u] \cap S_s| 
    &\leq |\WReach_q[G,S,\sigma,u] \cap (R'_\alpha \cup R'_\beta)| \\
    &\leq 2(c_{\ref{lemma:E-P_in_Kk_minor_free_graphs}}(k,d+1)+3) \cdot q.
    \end{align*}
    
    Next, let $i \in \{0, \dots, s-1\}$. We claim that
    \[
    |\WReach_q[G,S,\sigma,u] \cap (S_i \setminus S_{i+1})| \leq (c_{\ref{lemma:E-P_in_Kk_minor_free_graphs}}(k,d+1)+3) \cdot q.
    \]
    Since each vertex of $S_{i+1}$ precedes each vertex of $S_i$ in $\sigma$, we have
    \[
    \WReach_q[G,S,\sigma,u] \cap (S_i \setminus S_{i+1}) \subseteq \WReach_q[G-S_{i+1}, S-S_{i+1}, \sigma \vert_{S \setminus S_{i+1}}, u] \cap (S_i \setminus S_{i+1}).
    \]
    Let $\alpha = \max\{a \in I_{i+1} \mid a \leq j_u\}$ and $\beta = \alpha + 2^{i+1}$.
    Let $C$ be the component of $u$ in $G-S_{i+1}$.
    Since $W_{\alpha-1} \cap W_\alpha, W_{\beta-1} \cap W_\beta \subseteq S_{i+1}$,
    and because $(W_0, \dots, W_\ell)$ is a path decomposition of $G$,
    $V(C) \cap S \subseteq \bigcup_{a \in \{\alpha, \dots, \beta-1\}} W_a$.
    We deduce that
    \[
        \WReach_q[G-S_{i+1},S\setminus S_{i+1}, \sigma\vert_{S \setminus S_{i+1}}, u] \cap (S_i \setminus S_{i+1})
        \subseteq \textstyle\bigcup_{a \in \{\alpha, \dots, \beta-1\}} W_a.
    \]
    Since the only members of $I_i \setminus I_{i+1}$ in $\{\alpha+1, \dots, \beta-1\}$ is $\gamma = \alpha+2^i$,
    we in fact have
    \[
    \WReach_q[G-S_{i+1},S\setminus S_{i+1}, \sigma\vert_{S \setminus S_{i+1}}, u] \cap (S_i \setminus S_{i+1})
    \subseteq R'_\gamma,
    \]
    and we deduce that
    \begin{align*}
        |\WReach_q[G-S_{i+1},S\setminus S_{i+1}, \sigma\vert_{S \setminus S_{i+1}}, u] \cap (S_i \setminus S_{i+1})| &\leq \wcol_q(G - S_{i+1}, R_\gamma', \sigma_\gamma) \\
        &\leq (c_{\ref{lemma:E-P_in_Kk_minor_free_graphs}}(k,d+1)+3) \cdot q.
    \end{align*}

    For convenience let $S_{s+1} = \emptyset$.
    Since $S = S_0$, it follows that
    \begin{align*}
        |\WReach_q[G,S,\sigma,u]|
        &\leq \sum_{i \in \{0, \dots, s\}} |\WReach_q[G,S,\sigma,u] \cap (S_i \setminus S_{i+1})| \\
        &\leq (s+2) \cdot (c_{\ref{lemma:E-P_in_Kk_minor_free_graphs}}(k,d+1)+3) \cdot q \\
        &\leq 5(c_{\ref{lemma:E-P_in_Kk_minor_free_graphs}}(k,d+1)+3)\cdot q \log q. \qedhere
    \end{align*}
\end{proof}

\begin{lemma}\label{lem:new_base_case_wcol}
    Let $k$, $h$, and $d$ be positive integers with $h \geq 2$.
    There is an integer $c_{\ref{lem:new_base_case_wcol}}(h,d,k)$ such that
    for every connected $K_k$-minor-free graph $G$,
    for every family $\mathcal{F}$ of connected subgraphs of $G$,
    if $G$ has no $\mathcal{F}$-rich model of $F_{h,d}$, then
    there is a set $S \subseteq V(G)$
    such that
    \begin{enumerateOurAlph} 
        \item $S \cap V(F) \neq \emptyset$ for every $F \in \mathcal{F}$;\label{item:new_base_case_wcol_i}
        \item $G[S]$ is connected; and \label{item:new_base_case_wcol_ii}
        \item $\wcol_q(G,S) \leq c_{\ref{lem:new_base_case_wcol}}(h,d,k) \cdot q \log q$ for every integer $q$ with $q \geq 2$.\label{item:new_base_case_wcol_iii}
    \end{enumerateOurAlph}
\end{lemma}

\begin{proof}
    We proceed by induction on $h$.
    For $h=2$, the result is given by \Cref{lemma:excluding_a_star_wcol} setting $c_{\ref{lem:new_base_case_wcol}}(1,d,k) = 5(c_{\ref{lemma:E-P_in_Kk_minor_free_graphs}}(k,d+1)+3)$.
    Next, assume $h>2$ and that $c_{\ref{lem:new_base_case_wcol}}(h-1,d,k)$ witnesses the assertion for $h-1$.
    Let $c_{\ref{lem:new_base_case_wcol}}(h,d,k) = 5(c_{\ref{lemma:E-P_in_Kk_minor_free_graphs}}(k,d+1)+3) + 3 + c_{\ref{lem:new_base_case_wcol}}(h-1,d+1,k)$.

    Let $q$ be an integer with $q \geq 2$,
    let $G$ be a connected $K_k$-minor-free graph,
    and let $\mathcal{F}$ be a family of connected subgraphs of $G$.
    Suppose that $G$ has no $\mathcal{F}$-rich model of $F_{h,d}$.
    Let $\mathcal{F}'$ be the family of all the connected subgraphs $H$ of $G$ 
    such that $H$ contains an $\mathcal{F}\vert_H$-rich model of $F_{h-1,d+1}$.
    We claim that there is no $\mathcal{F}'$-rich model of $F_{2,d}$ in $G$.
    Suppose to the contrary that $\big(B_x \mid x \in V(F_{2,d})\big)$ is such a model.
    Let $s$ be the root of $F_{2,d}$ and let $s'$ be the root of $F_{h-1,d}$.
    For every $x \in V(F_{2,d}) \setminus \{s\}$,
    by \Cref{lemma:rooting_a_model_of_Fhd},
    there is an $\mathcal{F}$-rich model $\big(C_y \mid y \in V(F_{h-1,d})\big)$ of $F_{h-1,d}$ in $G[B_x]$ such that $C_{s'}$ contains a vertex of $N_G(B_s) \cap B_x$.
    The union of these models together with $B_s$ yields an $\mathcal{F}$-rich model of $F_{h,d}$ in $G$, which is a contradiction.
    See \Cref{fig:forest-models}.
    
    \begin{figure}[tp]
        \centering 
        \includegraphics[scale=1]{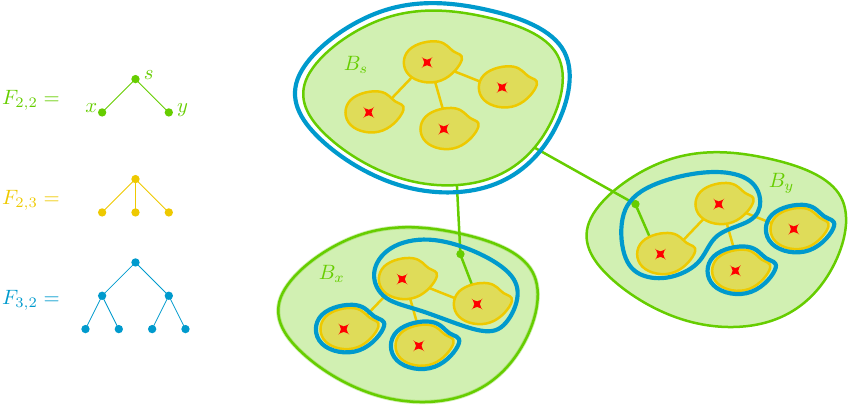} 
        \caption{
            We provide an example of the construction of an $\mathcal{F}$-rich model of $F_{h,d}$ in $G$ assuming that there is an $\mathcal{F}'$-rich model of $F_{2,d}$ in $G$ in the case where $h=3$ and $d=2$.
            In green, we depict an $\mathcal{F}'$-rich model of $F_{2,d} = F_{2,2}$ in the graph.
            Each branch set contains an $\mathcal{F}$-rich model of $F_{h-1,d+1} = F_{2,3}$.
            We depict these models in yellow and the red stars are the elements of $\mathcal{F}$.
            The obtained model of $F_{h,d} = F_{3,2}$ we depict in blue.
            Note that this model is $\mathcal{F}$-rich.
        }
        \label{fig:forest-models}
    \end{figure} 

    Since $G$ has no $\mathcal{F}'$-rich model of $F_{2,d}$, 
    by \Cref{lemma:excluding_a_star_wcol}, 
    there is a set $S_0 \subseteq V(G)$ such that
    \begin{enumerate}[label={\normalfont\ref{lemma:excluding_a_star_wcol}.(\makebox[\mywidth]{\alph*})}]
        \item for every $F \in \mathcal{F}'$, $S_0 \cap V(F) \neq \emptyset$; \label{lemma:star:item:hitting}
        \item $G[S_0]$ is connected; and \label{lemma:star:item:connected}
        \item $\wcol_q(G,S_0) \leq 5(c_{\ref{lemma:E-P_in_Kk_minor_free_graphs}}(k,d+1)+3) \cdot q \log q$. \label{lemma:star:item:wcol}
    \end{enumerate}
    Let $C$ be a component of $G - S_0$. 
    By~\ref{lemma:star:item:hitting}, $C \notin \mathcal{F}'$, and so, $C$ has no $\mathcal{F}\vert_C$-rich model of $F_{h-1,d+1}$.
    Therefore, by the induction hypothesis, there is a set $S_C \subseteq V(C)$ such that
    \begin{enumerateOurAlphPrim}
        \item $S_C \cap V(F) \neq \emptyset$ for every $F \in \mathcal{F}\vert_C$; \label{item:new_base_case_wcol_i_induction}
        \item $C[S_C]$ is connected; and \label{item:new_base_case_wcol_iii_induction}
        \item $\wcol_q(C,S_C) \leq c_{\ref{lem:new_base_case_wcol}}(h-1,d+1,k) \cdot q \log q$.\label{item:new_base_case_wcol_ii_induction}
    \end{enumerateOurAlphPrim}
    Let $Q_C$ be an \stgeodesic{S_C}{N_G(S_0)} in $G$.
    In particular, $Q_C$ is a geodesic in $C$.
    Let $\mathcal{C}$ be the family of components of $G-S_0$ and let 
    \[
    S = S_0 \cup \bigcup_{C \in \mathcal{C}} (S_C \cup V(Q_C)).
    \]
    See~\Cref{fig:excluding-forest-combine} for an illustration.
    We claim that \ref{item:new_base_case_wcol_i}--\ref{item:new_base_case_wcol_iii} hold.

    \begin{figure}[tp]
        \centering 
        \includegraphics[scale=1]{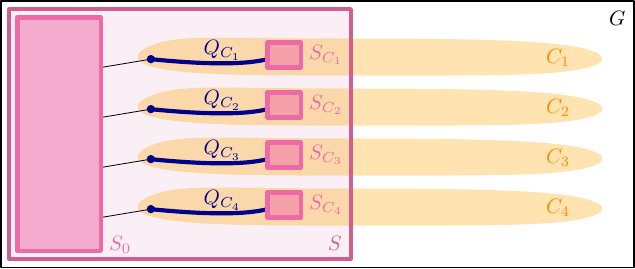} 
        \caption{
            An illustration of the construction of the set $S$ in the proof of \Cref{lem:new_base_case_wcol}.
        }
        \label{fig:excluding-forest-combine}
    \end{figure} 

    Let $F \in \mathcal{F}$.
    If $S_0 \cap V(F) = \emptyset$, then $F \subset C$ for some component $C$ of $G - S_0$.
    In particular, $F \in \mathcal{F}\vert_C$, and thus, by~\ref{item:new_base_case_wcol_i_induction}, $S_C \cap V(F) \neq \emptyset$, which proves~\ref{item:new_base_case_wcol_i}.
    The graph $G[S]$ is connected by construction, \ref{item:new_base_case_wcol_iii_induction} and \ref{lemma:star:item:connected}, which yields~\ref{item:new_base_case_wcol_ii}.
    The following sequence of inequalities concludes the proof of~\ref{item:new_base_case_wcol_iii} and the lemma:
    \begin{align*}
        \wcol_q(G,S) &\leq \wcol_q(G,S_0) + \wcol_q\left(G-S_0,\bigcup_{C \in \mathcal{C}} \big(S_C \cup V(Q_C)\big)\right)&&\textrm{by~\Cref{obs:wcol_union}}\\ 
        &\leq \wcol_q(G,S_0) + \max_{C \in \mathcal{C}} \wcol_q(C, S_C \cup V(Q_C))&&\textrm{by~\Cref{obs:wcol_components}}\\ 
        &\leq \wcol_q(G,S_0) + \max_{C \in \mathcal{C}} \wcol_q(C, S_C) + (2q+1)&&\textrm{by~\Cref{obs:geodesics}}\\ 
        &\leq (5(c_{\ref{lemma:E-P_in_Kk_minor_free_graphs}}(k,d+1)+3) + c_{\ref{lem:new_base_case_wcol}}(h-1,d+1,k) + 3) \cdot q\log q \hspace{-4pt}&&\textrm{by~\ref{lemma:star:item:wcol} and \ref{item:new_base_case_wcol_ii_induction}}\\ 
        &= c_{\ref{lem:new_base_case_wcol}}(h,d,k) \cdot q \log q. && \pushright{\mbox{\qedhere}}
    \end{align*}
\end{proof}

\begin{corollary}\label{cor:base_case_wcol}
    The function $q\mapsto q\log (q+1)$ is $(\wcol,\Rt_2)$-bounding for the family of focused parameters $\wcol =(\wcol_q \mid q \in \posint)$.
\end{corollary}

\begin{proof}
    Recall that $\Rt_2$ is the class of all forests.
    Let $k$ be a positive integer and let $X \in \Rt_2$.
    Let $k'$ be an integer such that all $K_k$-minor-free graphs $G$ are $k'$-degenerate or in other words, $\wcol_1(G) \leq k'$. 
    We pick integers $h$ and $d$ such that $X \subset F_{h,d}$ and we set 
        \[\beta(k,X) = \max\{k',c_{\ref{lem:new_base_case_wcol}}(h,d,k)\}.\]
    We show that every $K_k$-minor-free graph is $(q\mapsto q\log (q+1),\wcol,X,1,\beta(k,X))$-good.
    Let $G$ be a $K_k$-minor-free graph and let $\calF$ be a family of connected subgraphs of $G$ such that $G$ has no $\mathcal{F}$-rich model of $X$ (and so of $F_{h,d}$).

    Let $C$ be a component of $G$.
    By~\Cref{lem:new_base_case_wcol}, there exists $S_C \subset V(C)$ such that
    \begin{enumerate}[label={\normalfont\ref{lem:new_base_case_wcol}.(\makebox[\mywidth]{\alph*})}]
        \item $S_C \cap V(F) \neq \emptyset$ for every $F \in \mathcal{F}\vert_C$;\label{item:new_base_case_wcol_i-cor}
        \item $C[S_C]$ is connected; and \label{item:new_base_case_wcol_ii-cor}
        \item $\wcol_q(C,S_C) \leq c_{\ref{lem:new_base_case_wcol}}(h,d,k) \cdot q \log q$ for every integer $q$ with $q \geq 2$.\label{item:new_base_case_wcol_iii-cor}
    \end{enumerate}
    For $\calC$ being the set of the components of $G$, we set
        \[S = \bigcup_{C \in \calC} S_C.\]
    
    Let $q$ be a positive integer.
    We shall prove that
            \begin{enumerate}[label={\normalfont (g\arabic*$^\star$)}]
            \item $S \cap V(F) \neq \emptyset$ for every $F \in \mathcal{F}$; \label{item:lemma:fragility_rate_base_case:hit-cor}
            \item for every component $C'$ of $G-S$, $N_{G}(V(C'))$ intersects at most one component of $G-V(C')$; and \label{item:lemma:fragility_rate_base_case:components-cor}
            \item $\param_q(G,S) \leq \beta(k,X) \cdot \log (q+1)$. \label{item:lemma:fragility_rate_base_case:param-cor}
        \end{enumerate}
    Clearly,~\ref{item:new_base_case_wcol_i-cor} and~\ref{item:new_base_case_wcol_ii-cor} imply~\ref{item:lemma:fragility_rate_base_case:hit-cor} and~\ref{item:lemma:fragility_rate_base_case:components-cor}, respectively.
    Finally,~\ref{item:new_base_case_wcol_iii-cor} and~\Cref{obs:wcol_components} imply~\ref{item:lemma:fragility_rate_base_case:param}. 
\end{proof}

\begin{theorem}\label{thm:main_wcol_Rt}
    Let $t$ be a positive integer.
    For every $X \in \Rt_t$, there exists an integer $c$ such that for every integer $q$ with $q \geq 2$ and
    for every $X$-minor-free graph $G$,
    \[
        \wcol_q(G) \leq c \cdot q^{t-1} \log q.
    \]
\end{theorem}

\begin{proof}
    By \Cref{lemma:wcol-nice}, the family of focused parameters 
    $\wcol = (\wcol_q \mid q \in \posint)$ is nice.
    We show by induction on $t$ (integers with $t \geq 2$) that $q \mapsto q^{t-1}\log (q+1)$ is $(\wcol, \Rt_t)$-bounding.
    When $t=2$, this is given by \Cref{cor:base_case_wcol}.
    Next, suppose that $t \geq 3$.
    By~\Cref{obs:classes-closed}, $\Rt_{t-1}$ is closed under disjoint union and leaf addition, and by \Cref{lemma:Rt_has_coloring_elimination_property}, $\Rt_{t-1}$ has the coloring elimination property.
    By induction, $q \mapsto q^{t-2}\log (q+1)$ is $(\wcol, \Rt_{t-1})$-bounding.
    Therefore, by \Cref{thm:abstract_induction_main}, 
    $q \mapsto q^{t-1} \log (q+1)$ is $(\wcol, \Rt_{t})$-bounding.
    This completes the inductive proof that for every integer $t$ with $t\geq 2$, $q \mapsto q^{t-1}\log (q+1)$ is $(\wcol, \Rt_{t})$-bounding, as desired.
    By~\Cref{lemma:par-bounding-to-bound}, this implies that for every $X \in \Rt_t$, there exists a constant $\beta = \beta(X)$ such that for every integer $q$ with $q\geq 2$, we have
    \[\wcol_q(G) = \wcol_q(G,V(G)) \leq \beta \cdot q^{t-1}\log (q+1) \leq 2\beta  \cdot q^{t-1} \log q.\qedhere\]
\end{proof}

\section{Centered chromatic numbers}
\label{sec:centered}

In this section, we prove~\Cref{thm:centered}.
% The plan for this section is as follows. 
First, we show the existence of suitable precolorings of $K_t$-minor-free graphs in \Cref{sec:kt-minor-free-centered-proof-good}.
Then, in \Cref{sec:cen:the_bounds},
we follow the roadmaps depicted in~\Cref{fig:roadmap_wcol_bounded_tw_St,fig:roadmap_wcol_bounded_tw_Rt}
to prove \Cref{thm:centered}. Most of the work lies in the base case
corresponding to excluding a graph in $\Apex(\Rt_2)$,
and which is inspired by the proof of D\k{e}bski, Felsner, Micek, and Schröder~\cite{Dbski2021} that outerplanar graphs
have $q$-centered chromatic numbers in $\bigO(q \log q)$.
We set our base case to excluding a graph in $\Apex(\mathcal{R}_2)$ 
instead of $\mathcal{R}_2$ (as we did in \Cref{sec:wcol,sec:fragility})
because graphs excluding a forest as a minor do not have
$q$th centered chromatic numbers in $\bigO(\log q)$, but in $\Theta(q)$.
This particularity of centered chromatic numbers 
is highlighted in the case $(t,s)=(2,3)$ of \Cref{cor:general}.

We remark that one of the conference versions of this paper~\cite{centered_paper_soda} contains a direct proof of the special case of $K_t$-minor-free graphs 
(i.e.\ \Cref{thm:Kt_minor_free-centered}).
This proof avoids the abstract framework.
Its main building blocks are analogous to~\cref{thm:Kt_minor_free-fragility}, where we bound the fractional treedepth fragility rates for $K_t$-minor-free graphs.
Namely, a bound for centered chromatic numbers of bounded treewidth due to Pilipczuk and Siebertz~\cite{PS19}, and extending the idea of Illingworth, Scott, and Wood~\cite[Theorem~4]{ISW22}.
Additionally, we need the material of~\Cref{sec:kt-minor-free-centered-proof-good}.

We start with some simple combinatorial statements, whose proofs are given for completeness.
Let $T$ be a tree and let $\calQ$ be a collection of connected subgraphs of $T$ whose vertex sets partition $V(T)$.
For every $X \subset V(T)$, let $\defin{\text{$\calQ(X)$}} = \{Q \in \calQ \mid V(Q) \cap X \neq \emptyset\}$.

\begin{lemma}\label{lemma:LCA-stuff}
    Let $T$ be a tree and let $\calQ$ be a collection of connected subgraphs of $T$ whose vertex sets partition $V(T)$.
    Let $X, Y \subseteq V(T)$ with $X\subseteq Y$ and $\LCA(T,X) = X$.
    If $\calQ(X) = \calQ(Y)$,
    then $\calQ(Y)=\calQ(\LCA(T,Y))$.
\end{lemma}

\begin{proof}
    Suppose that $\calQ(X)=\calQ(Y)$.
    Since $Y \subset \LCA(T,Y)$, we have $\calQ(Y) \subset \calQ(\LCA(T,Y))$.
    Thus, it suffices to prove that $\calQ(\LCA(T,Y)) \subset \calQ(Y)$.
    Consider $Q \in \calQ(\LCA(T,Y))$.
    There exist $y_1,y_2 \in Y$ such that $\lca(T,y_1,y_2) \in Q$.
    Let $Q_1,Q_2 \in \mathcal{Q}$ be such that
    $y_1 \in V(Q_1)$, $y_2 \in V(Q_2)$.
    In particular, $Q_1,Q_2 \in \calQ(Y) = \calQ(X)$.
    Hence, there exist $x_1 \in V(Q_1) \cap X$ and $x_2 \in V(Q_2) \cap X$.
    
    If the roots of $Q_1$ and $Q_2$ are not in an ancestor-descendant relation in $T$, then for all $q_1\in V(Q_1)$ and $q_2 \in V(Q_2)$, we have $\lca(T,q_1,q_2) = \lca(T, \root(Q_1),\root(Q_2))$.
    In particular, 
        \[\lca(T,x_1,x_2) = \lca(T, \root(Q_1),\root(Q_2)) = \lca(T,y_1,y_2) \in V(Q).\]
    Observe that $\lca(T,x_1,x_2) \in \LCA(T,X) = X$, thus, $V(Q) \cap X \neq \emptyset$, and so, $Q \in \calQ(X) = \calQ(Y)$. 
    
    Therefore, we can assume that one of the roots say $\root(Q_1)$, is an ancestor of the other $\root(Q_2)$ in $S$.
    It follows that for all $q_1\in V(Q_1)$ and $q_2 \in V(Q_2)$, $\lca(T,q_1,q_2)$ is either equal to $\root(Q_1)$ or is a descendant of $\root(Q_1)$.
    In both cases, $\lca(T,q_1,q_2)$ lies in the path from $q_1$ to $\root(Q_1)$ in $S$.
    Since $Q_1$ is connected, we obtain that $\lca(T,q_1,q_2) \in V(Q_1)$.
    In particular, $\lca(T,x_1,x_2) \in V(Q_1)$, and so, $Q_1 = Q$ implying $Q \in \calQ(Y)$, which ends the proof.
\end{proof}

\begin{lemma}\label{lemma:projections_on_a_torso_stay_connected}
    Let $G$ be a graph, 
    let $\mathcal{W} = \big(T,(W_x \mid x \in V(T))\big)$ be a tree decomposition of $G$,
    and let $x \in V(T)$.
    If $H$ is a connected subgraph of $G$ intersecting $W_x$, then
    $V(H) \cap W_x$ induces a connected subgraph of $\torso_{G,\mathcal{W}}(W_x)$.
\end{lemma}
\begin{proof}
    Let $H$ be a connected subgraph of $G$ intersecting $W_x$, and suppose to the contrary that the subgraph $H'$ of $\torso_{G,\mathcal{W}}(W_x)$ induced by $V(H) \cap W_x$ is not connected.
    Let $u,v \in V(H')$ be vertices in distinct components of $H$ such that the distance between $u$ and $v$ is minimal in $H$.
    Note that the internal vertices of a shortest path between $u$ and $v$ in $H$ do not lie in $W_x$ as otherwise we obtain a pair of vertices of $H'$ in distinct components that are closer in $H$.
    By the properties of tree decompositions, such a shortest path has all its vertices in $\bigcup_{z \in V(T_{y \mid x})} W_z$ for some $y \in N_T(x)$,
    and so $u,v \in W_x \cap W_y$.
    It follows that $u$ and $v$ are adjacent in $\torso_{G,\mathcal{W}}(W_x)$, and so, in $H'$, which contradicts the assumption and completes the proof.
\end{proof}

We also need a more explicit version of \Cref{lemma:natural_tree_decomposition}.
Indeed, at one point in the proof, we would like to assume that a tree decomposition that we consider is natural but still preserves certain properties.
For this reason, we introduce the following definitions and we prove~\Cref{lemma:making_a_td_natural}.

\begin{figure}[tp]
    \centering
    \includegraphics{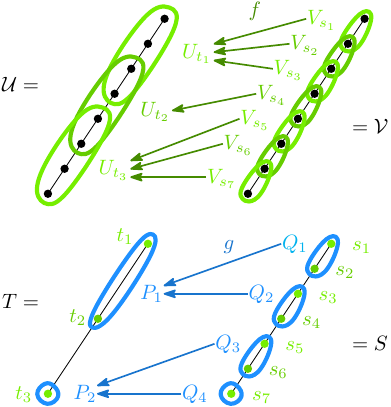}
    \caption{
        In this example, $G$ is an $8$-vertex path.
        $\calP = \{P_1,P_2\}$ and $\calQ = \{Q_1,Q_2,Q_3,Q_4\}$.
        One can check that $(\calV,\calQ)$ refines $(\calU,\calP)$ 
        which is witnessed by $f$ and $g$.
        Note that $f$ is a function on $V(S)$ but for the readability reasons we depict it as it acted on bags of $\calV$.}
    \label{fig:refinement}
\end{figure}

Let $G$ be a graph. 
A pair $(\calU,\calP)$ is a \defin{normal pair} of 
$G$ if $\calU=\big(T,(U_t\mid t\in V(T))\big)$ is a tree decomposition of $G$, 
and $\calP$ is a collection of connected subgraphs 
of $T$ whose vertex sets partition $V(T)$.
Consider two normal pairs 
$(\calU,\calP)$ and $(\calV,\calQ)$ of $G$ with $\calU=\big(T,(U_t\mid t\in V(T))\big)$ and $\calV=\big(S,(V_s\mid s\in V(S))\big)$. 
We say that $(\calV,\calQ)$ \defin{refines}, or is a \defin{refinement} of $(\calU,\calP)$ if there exist $f\colon V(S) \rightarrow V(T)$ and $g\colon \calQ \rightarrow \calP$ such that
\begin{enumerate}[label=(r\arabic*)]
    \item for every $s \in V(S)$, 
        \[V_s \subset U_{f(s)};\] \label{item:making_a_td_natural:i}
    \item for every $Q \in \calQ$, 
        \[f(V(Q)) \subset V(g(Q)).\]\label{item:making_a_td_natural:ii}
\end{enumerate}
See an example in~\cref{fig:refinement}.

\begin{lemma}\label{lemma:making_a_td_natural}
    Let $G$ be a connected graph and let $(\calU,\calP)$ be a normal pair of $G$. 
    There exists $(\calV,\calQ)$ a refinement of $(\calU,\calP)$ such that $\calV$ is natural.
\end{lemma}

\begin{proof}
For every positive integer $i$, for every tree decomposition $\mathcal{W}$ of $G$,
we denote by $n_i(\mathcal{W})$ the number of bags of $\mathcal{W}$ of size $i$.
Then let $n(\mathcal{W}) = (n_{|V(G)|}(\mathcal{W}), \dots, n_0(\mathcal{W}))$. 

Since the refinement relation is reflexive, every normal pair of $G$ has a refinement.
Consider the refinement $(\calV,\calQ)$ of $(\calU,\calP)$ with $n(\calV)$ minimal in the lexicographic order.
We claim that $\calV$ is natural.

Let $\calV = \big(S, (V_s \mid s \in V(S))\big)$.
Suppose to the contrary that there exist $x,y \in V(S)$ such that $xy \in E(S)$
and $G\left[\bigcup_{s \in V(S_{x \mid y})} V_s\right]$ is not connected, and let $\mathcal{C}$ be the family of its components.

\begin{figure}[tp]
  \begin{center}
    \includegraphics{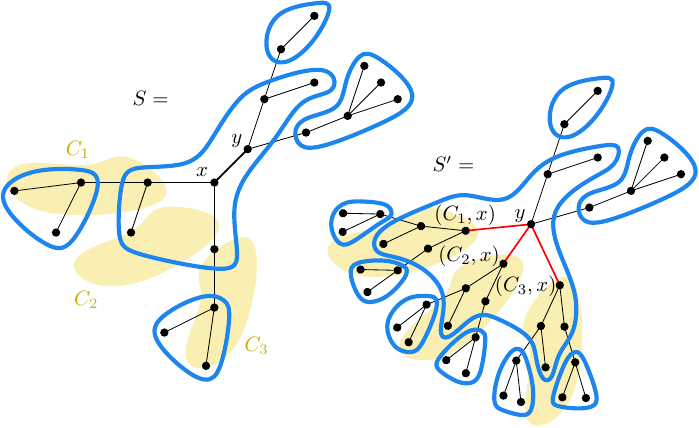}
  \end{center}
  \caption{
  Proof of~\Cref{lemma:making_a_td_natural}: construction of $S'$ from $S$ and of $\calQ'$ from $\calQ$. The blue partition of $S$ is $\calQ$ and the blue partition of $S'$ is $\calQ'$.
  }
  \label{fig:finidng-natural-td}
\end{figure}

The plan for obtaining a contradiction is finding a refinement $(\calV',\calQ')$ of $(\calV,\calQ)$ such that $n(\calV') < n(\calV)$ in the lexicographic order.
For every $C \in \mathcal{C}$, let $S_C$ be a copy of $S_{x \mid y}$
with vertex set $\{(C,s) \mid s \in V(S_{x \mid y})\}$ and edge set
$\{(C,s_1)(C,s_2) \mid s_1s_2 \in E(S_{x \mid y})\}$.
Let $S'$ be the tree obtained from the disjoint union of $S_{y|x}$ with all the $S_C$ over $C \in \mathcal{C}$ by adding the edges $(C,x)y$ for each $C \in \mathcal{C}$.
See~\cref{fig:finidng-natural-td}.
Let $f\colon V(S') \to V(S)$ be such that for every $z \in V(S')$,
\[
f(z) = \begin{cases} 
             s &\textrm{if $z=s$ for some $s\in V(S_{y \mid x})$,}\\
             s &\textrm{if $z=(C,s)$ for some $s\in V(S_{x \mid y})$ and $C\in\calC$.}
        \end{cases}
\]
Also, for every $z \in V(S')$, let
\[
V'_{z} = \begin{cases} 
V_{s}&\textrm{if $z=s$ for some $s\in V(S_{y|x})$,}\\
V_{s} \cap V(C)&\textrm{if $z=(C,s)$ for some $s\in V(S_{x \mid y})$ and $C\in\calC$.}
\end{cases}
\]
Note that by construction, for every $z\in V(S')$, $V'_z \subseteq V_{f(z)}$.
Moreover, $\calV'=\big(S',(V'_{z}\mid z\in V(S'))\big)$ is a tree decomposition of $G$.

To complete the construction, we specify a collection $\calQ'$ of subtrees of $S'$ whose vertex sets partition $V(S')$.
Simultaneously, we define $g\colon \calQ' \rightarrow \calQ$.
First, for every $Q \in \calQ$ with $Q \cap V(S_{y\mid x}) \neq \emptyset$, let 
    \[Q'=S'\left[(V(Q) \cap V(S_{y \mid x})) \cup \bigcup_{C \in \mathcal{C}} \{(C,z) \mid z \in V(Q) \cap V(S_{x \mid y})\}\right]\textrm{ and }g(Q')=Q.\]
Next, for every $Q \in \calQ$ with $Q \subseteq V(S_{x\mid y})$, for every $C \in \calC$, let
    \[Q_C = S'\left[\{(C,z) \mid z \in Q\}\right]\textrm{ and }g\big(Q_C\big)=Q.\]
Finally, let 
\begin{align*}
    \calQ' &= \left\{ Q' \mid Q \in \calQ \text{ with } V(Q) \cap V(S_{y\mid x}) \neq \emptyset  \right\} \cup \bigcup_{C \in \calC} \left\{ Q_C \mid V(Q) \subseteq V(S_{x\mid y})  \right\}.
\end{align*}
It follows that $\calQ'$ is a collection of subtrees of $S'$ whose sets of vertices partition $V(S')$ and 
and $f(V(Q')) = Q = g(Q)$ for every $Q' \in \calQ'$.
In other words,~\ref{item:making_a_td_natural:ii} holds.
Note that~\ref{item:making_a_td_natural:i} also holds, and thus, $(\calV',\calQ')$ refines $(\calV,\calQ)$.
Moreover, by transitivity of the refinement relation, $(\calV',\calQ')$ refines $(\calU,\calP)$.

To conclude the proof, we show that $n(\mathcal{V}')$ is less than $n(\mathcal{V})$ in the lexicographic order, which leads to a contradiction.
For every $s \in V(S_{x \mid y})$, we have $\sum_{C \in \calC} |V_{(C,s)}'| = |V_s|$. 
Thus, either $|V_{(C,s)}'| < |V_s|$ for every $C \in \calC$ or there is exactly one $C \in \calC$ with $|V_{(C,s)}'| = |V_s|$ and $|V_{(D,s)}'| = 0$ for all $D \in \calC\setminus\{C\}$.
Moreover, for every $s \in V(S_{y \mid x})$, we have $V_s = V_s'$.
Since $G$ is connected, $V_x$ intersects every component in $\mathcal{C}$.
Also, note that we supposed $|\calC| \geq 2$.
Consider the maximum integer $i_0$ such that there exists $s \in V(S_{x \mid y})$ 
with $|V_s|=i_0$ and $V_s$ intersects at least two components in $\mathcal{C}$.
From the construction, we obtain that $n_i(\mathcal{V}') = n_i(\mathcal{V})$ for every integer $i$ with $i > i_0$ and $n_{i_0}(\mathcal{V}')< n_{i_0}(\mathcal{V})$.
We conclude that $n(\mathcal{V}')$ is less than $n(\mathcal{V})$ in the lexicographic order.
\end{proof}

\subsection{Good colorings: a generalized Helly property} \label{sec:kt-minor-free-centered-proof-good}

Similarly as in the case of weak coloring numbers, we need a variant of the Helly property as in~\Cref{lemma:helly_property_tree_decomposition} for $K_t$-minor-free graphs.
Again, instead of bounding the size of a hitting set (as in~\Cref{lemma:helly_property_tree_decomposition}.\ref{item:helly_property_tree_decomposition:hit}), we require it to be \q{well-structured}.
In this case, expressing this structure is more involved.

Let $c$ be a positive integer.
A coloring $\phi$ of $G$ is \defin{$(q,c)$-good} if
for every subgraph $G_0$ of $G$,
for every family $\mathcal{F}$ of connected subgraphs of $G_0$,
for every positive integer $d$,
if there are no $d+1$ pairwise disjoint members of $\mathcal{F}$, then
there exists $Z \subseteq V(G_0)$ such that
\begin{enumerate}[label=(qc\arabic*)]
    \item $Z \cap V(F) \neq \emptyset$ for every $F \in \calF$; \label{centered-Helly-hitting}
    \item for every component $C$ of $G_0-Z$, $N_{G_0}(V(C))$ intersects at most two components of $G_0-V(C)$; and \label{centered-Helly-components}
    \item for every $P \subseteq Z$, there exists $R \subseteq Z$ with
        \begin{enumerate}
            \item $P \subseteq R$,
            \item $\cen_p(G_0,\phi,R) \leq c \cdot d$, and
            \item $|R \setminus P| \leq c \cdot d$.
        \end{enumerate} \label{centered-Helly-restriction}
\end{enumerate}

Note that \ref{centered-Helly-restriction} applied for $P=Z$
yields
\begin{enumerate}[resume*]
    \item $\cen_q(G_0,\phi,Z)\leq c\cdot d$.\label{centered-Helly-coloring}
\end{enumerate}
Moreover, for all graphs $G$ and $H$ where $H$ is a subgraph of $G$,
if $\phi$ is a $(q,c)$-good of $G$, then $\phi\vert_{V(H)}$ is a $(q,c)$-good coloring of $H$.
This implies the following straightforward observation.
\begin{obs}\label{obs:qc-good-on-components}
    Let $q$ and $c$ be positive integers and let $G_1,\dots,G_\ell$ be pairwise vertex-disjoint graphs.
    Assume that for each $i \in [\ell]$, there exists a $(q,c)$-good coloring $\phi_i$ of $G_i$.
    Let $G$ be the disjoint union of $G_1,\dots,G_\ell$ and let $\phi$ be the coloring such that $\phi\vert_{V(G_i)} = \phi_i$ for every $i \in [\ell]$.
    Then, $\phi$ is $(q,c)$-good.
\end{obs}

In~\Cref{ssec:centered-nice}, we considered an abstract framework of building families of focused parameters related to centered chromatic numbers based on precolored graphs.
We will apply the framework to graphs where the precoloring is a $(q,c)$-good coloring.

\begin{lemma}\label{lem:cen-good-const-edgeless-bounding}
    Let $c$ be a positive integer.
    Let $\calU$ be an infinite graph, and for every positive integer $q$, let $\phi_q$ be a coloring of $\calU$ such that for every finite subgraph $G$ of $\calU$, $\phi_q\vert_{V(G)}$ is a $(q,c)$-good coloring of $G$.
    For every positive integer $q$,
    for every graph $G$, and
    for every $S \subseteq V(G)$,
    let
    \[
        \param_q(G,S) =
        \begin{cases}
            \cen_q(G,\phi_q\vert_{V(G)},S) & \textrm{if $V(G) \subseteq V(\mathcal{U})$ and $E(G) \subseteq E(\mathcal{U})$,} \\
            0 & \textrm{otherwise.}
        \end{cases}
    \]
    Then, $q \mapsto 1$ is $(\param,\edgeless)$-bounding, where $\param = (\param_q \mid q \in \posint)$ is a family of focused parameters.
\end{lemma}
\begin{proof}
    Let $X = \overline{K_d} \in \edgeless$.
    For each graph $G$ for which it does not hold that $V(G) \subset V(\calU)$ and $E(G) \subset E(\calU)$, choosing the empty set as $S$ witnesses that $G$ is $(q \mapsto 1, \param, X, 1, 0)$-good.
    Let $G$ be a graph for which $V(G) \subset V(\calU)$ and $E(G) \subset E(\calU)$.
    We claim that $G$ is $(q \mapsto 1,\param,X,2,cd)$-good. 
    Let $q$ be a positive integer and let $\calF$ be a family of connected subgraphs of $G$ such that $G$ has no $\calF$-rich model of $X$, i.e.\ there are no $d$ pairwise disjoint members of $\calF$.
    By assumption, $\varphi_q\vert_{V(G)}$ is $(q,c)$-good.
    It follows that there exists $S \subset V(G)$ such that 
    \begin{enumerate}[label=(qc\arabic*')]
        \item $S \cap V(F) \neq \emptyset$ for every $F \in \calF$; \label{centered-Helly-hitting-app}
        \item for every component $C$ of $G-S$, $N_{G}(V(C))$ intersects at most two components of $G-V(C)$; and \label{centered-Helly-components-app} \addtocounter{enumi}{1}
    \item $\cen_q(G,\varphi_q\vert_{V(G)},S)\leq c\cdot d$.\label{centered-Helly-coloring-app}
    \end{enumerate}
    Item~\ref{centered-Helly-coloring-app} implies that $\param_q(G,S) \leq cd$, which shows that indeed $G$ is $(q \mapsto 1,\param,X,2,cd)$-good, and so, $q \mapsto 1$ is $(\param,\edgeless)$-bounding.
\end{proof}

A simple example of $(q,2)$-good coloring of a graph $G$ using $\tw(G)+1$ colors is constructed in~\Cref{lemma:p1_good_colorings_when_bounded_tw} below.
We will use it in the proof of bounds in~\Cref{thm:centered} involving treewidth.
Next, relying on layered RS-decompositions (see \Cref{theorem:Kt_free_product_structure_decomposition}), 
we construct in~\Cref{lemma:Kt_free_graphs_have_centered_Helly_colorings} a $(q,c)$-good coloring using $q+1$ colors for each $K_t$-minor-free.
This is an essential step in the proof of the bounds in~\Cref{thm:centered}.

\begin{lemma}\label{lemma:p1_good_colorings_when_bounded_tw}
    Let $q$ be a positive integer and let $G$ be a graph.
    There is a $(q,2)$-good coloring of $G$ using $\tw(G)+1$ colors.
\end{lemma}

\begin{proof}
    By \Cref{lemma:natural_tree_decomposition}, we may fix a natural tree decomposition $\big(T,(W_x \mid x \in V(T))\big)$ of $G$ of width $\tw(G)$.
    Let $\phi\colon V(G) \to [\tw(G)+1]$ be a proper coloring  of the graph $\big(V(G), \bigcup_{x \in V(T)} \binom{W_x}{2}\big)$, which is possible since this graph has treewidth $\tw(G)$ and chromatic number at most treewidth plus $1$.
    Observe that for every $x \in V(T)$, $\phi\vert_{W_x}$ is injective.

    We now show that $\phi$ is a $(q,2)$-good coloring of $G$.
    Let $G_0$ be a subgraph of $G$,
    let $d$ be a positive integer,
    and let $\mathcal{F}$ be a family of connected subgraphs of $G_0$ with no $d+1$ pairwise disjoint members.
    By \Cref{lemma:helly_property_tree_decomposition} and \Cref{lemma:increase_X_to_have_small_interfaces},
    there is a set $X \subseteq V(T)$ of size at most $2d-1$ such that for $Z = \bigcup_{x \in X} V(T)$,
    \begin{enumerate}
        \item $Z \cap V(F) \neq \emptyset$ for every $F \in \mathcal{F}$, and so \ref{centered-Helly-hitting} holds and
        \item for every component $C$ of $G_0-Z$, $N_{G_0}(V(C))$ intersects at most two components of $G_0-V(C)$, and so \ref{centered-Helly-components} holds.
    \end{enumerate}
    Now, let $\psi \colon Z \to X$ such that $z \in W_{\psi(x)}$ for every $z \in Z$.
    Note that $\phi\vert_{Z} \times \psi$ is injective.
    Hence, for every connected subgraph $H$ of $G$ intersecting $Z$,
    any vertex in $V(H) \cap Z$ is a $(\phi\vert_Z \times \psi)$-center of $V(H) \cap Z$.

    Finally, for every $P \subseteq Z$,
    $\cen_q(G_0,\phi,P) \leq |\psi(P)| \leq |X| \leq 2 d$.
    This shows \ref{centered-Helly-restriction} and proves that $\phi$ is a $(q,2)$-good coloring of $G$.
\end{proof}

The remainder of this subsection is devoted to proving the following statement.

\begin{lemma}\label{lemma:Kt_free_graphs_have_centered_Helly_colorings}
    For every positive integer $t$, there exists a positive integer $c_{\ref{lemma:Kt_free_graphs_have_centered_Helly_colorings}}(t)$ 
    such that
    for every $K_t$-minor-free graph $G$ and every positive integer $q$,
    $G$ admits a $(q,c_{\ref{lemma:Kt_free_graphs_have_centered_Helly_colorings}}(t))$-good coloring using $q+1$ colors.
\end{lemma}

The proof of~\Cref{lemma:Kt_free_graphs_have_centered_Helly_colorings} is rather technical, however, the general idea of the proof is relatively simple to describe.
For purposes of this short discussion, consider a simpler variant of good colorings 
where we drop~\ref{centered-Helly-components}, which is required in the real definition for very technical reasons.
Observe that in the definition of good colorings, a subgraph $G_0$ occurs only in~\ref{centered-Helly-components}, 
hence, if we ignore this item, we can also assume that $G_0 = G$.
We say that a coloring is a \defin{simple} $(q,c)$-good coloring if it admits the relaxed definition above.

For graphs of bounded treewidth, there is a straightforward way of constructing (simple) good colorings, as shown in~\Cref{lemma:p1_good_colorings_when_bounded_tw}.
A true inspiration for the definition of good colorings comes from the case of graphs of bounded layered treewidth.
Let $G$ be a graph, let $q$ be a positive integer, 
and let $\ltw(G)$ be witnessed by $(\calW,\calL)$, where $\calW=\big(T,(W_x \mid x\in V(T))\big)$ and $\calL=(L_i \mid i \in \NN)$. 
For all $i \in \NN$ and $v\in L_i$, we define 
$\phi(v) = i \bmod (q+1)$. 
We claim that $\phi$ is a simple $(q,\ltw(G))$-good coloring of $G$.
Let $\mathcal{F}$ be a family of connected subgraphs of $G$ with no $d+1$ pairwise disjoint members. 
Let $Z$ be the union of at most $d$ bags of $\calW$ that is a hitting set of $\mathcal{F}$ in $G$ (given by~\Cref{lemma:helly_property_tree_decomposition}). 
In particular, $|Z \cap L_i| \leq \ltw(G) \cdot d$ for every $i\in \NN$.
Let $\psi \colon Z \to [\ltw(G) \cdot d]$ be any coloring of $Z$ which is injective on $Z \cap L_i$ for every $i \in \NN$.
It remains to verify~\ref{centered-Helly-coloring}.
Let $H$ be a connected subgraph of $G$ with $V(H)\cap Z\neq\emptyset$. 
Since $H$ is connected and $\calL$ is a layering of $G$, $V(H)$ intersects a set of consecutive layers in $\calL$.
If this set has at least $q+1$ elements, then $|\varphi(V(H))| > q$, and otherwise any element of $V(H) \cap Z$ is a $(\phi\times\psi)$-center of $V(H) \cap Z$.

To lift this idea to the class of general $K_t$-minor-free graphs, we use a layered RS-decomposition (see~\Cref{sec:preliminaries}). 

Given a layered RS-decomposition $(T,\calW,\calA,\calD,\calL)$ of a graph $G$ of bounded width, the strategy to find a simple good coloring of $G$ is the following.
Let $q$ be a positive integer.
The intention is to mimic the proof for the bounded layered treewidth case.
That is, ignoring the overlapping between the bags of $\calW$, we set $\varphi(v) \equiv i \bmod (q+1)$ for every $i \in \NN$ and every $v \in L_{x,i}$.
We claim that $\varphi$ is a simple $(q,c)$-good coloring, where $c$ depends only on the width of the given decomposition.
Let $\mathcal{F}$ be a family of connected subgraphs of $G$ with no $d+1$ pairwise disjoint members.
We glue all tree decompositions in $\calD$ in a natural way to obtain a tree decomposition $\calU$ of $G$.
The vertices of $A_x$ do not occur in bags of tree decompositions of $\calD$, but we can add them to all bags of $\calU$ corresponding to $W_x$.
The important property of $\calU$ is that~\ref{LRS:ltw} is preserved.
Next, we apply~\Cref{lemma:helly_property_tree_decomposition} to $\calF$ and $\calU$ obtaining a hitting set $Z$ of $\mathcal{F}$ in $G$, which is the union of at most $d$ bags of $\calU$.
Note that there are also at most $d$ bags of $\calW$, whose union contains $Z$.
Say that these bags correspond to the vertices in $X \subset V(T)$.
The next step is to ``disconnect'' elements of $Z$ in different bags of $\calW$.
To this end, we root $T$ arbitrarily and we define $B$ as the union of $\bigcup_{x \in X}A_x$ and the union of adhesions between $W_x$ and the bag corresponding to the parent of $x$ in $\calW$ for all $x \in X$.
In particular, the size of $B$ is bounded.
Finally, we define a coloring $\psi$ of $Z$ so that $\psi$ is injective on $B$ and $\psi$ is injective on $L_{x,i} \cap Z$ but the colors used are disjoint from $\psi(B)$.
Checking that $\psi$ witnesses $\varphi$ being a simple $(q,c)$-good coloring of $G$ is very similar to the bounded layered treewidth case.

\begin{proof}[Proof of \Cref{lemma:Kt_free_graphs_have_centered_Helly_colorings}]
    Let $t$ be a positive integer. 
    Let $\cLRS(t)$ be the constant from~\Cref{theorem:Kt_free_product_structure_decomposition}. 
    We set
    \[
    c_{\ref{lemma:Kt_free_graphs_have_centered_Helly_colorings}}(t) = \cLRS(t) \left(12\cLRS(t) + 10\right).
    \]
    
    Let $G$ be a $K_t$-minor-free graph, 
    and $q$ be a positive integer.
    By \Cref{theorem:Kt_free_product_structure_decomposition}, $G$ admits a layered RS-decomposition 
    $(T, \mathcal{W}, \mathcal{A}, \mathcal{D}, \mathcal{L})$ 
    of width at most $\cLRS(t)$.
    Let $\mathcal{W} = (W_x \mid x \in V(T))$,
    $\mathcal{A} = (A_x \mid x \in V(T))$,
    $\mathcal{D} = \big(\big(T_x,(D_{x,z} \mid z \in V(T_x))\big) \mid x \in V(T)\big)$, and
    $\mathcal{L} = ((L_{x,i} \mid i \in \NN) \mid x \in V(T))$.
    We root $T$ in an arbitrary vertex $r \in V(T)$.

    We define a coloring $\phi \colon V(G) \rightarrow \{0,\dots,q\}$ as follows. For every $v\in V(G)$ let
    \[
        \phi(v) = 
        \begin{cases}
            0&\textrm{if $v\in A_r$,}\\
            i\bmod (q+1)&\textrm{if $v\in W_r \setminus A_r$ and $v \in L_{r,i}$,}\\
            0&\textrm{if $v\in A_x \setminus W_{\parent(T,x)}$ for $x\in V(T)\setminus \set{r}$,}\\
            i\bmod (q+1)&\textrm{if $v\in W_x \setminus (A_x\cup W_{\parent(T,x)})$ and $v \in L_{x,i}$ for $x\in V(T) \setminus \set{r}$.}\\
        \end{cases}
    \]

    It remains to  show that $\phi$ is a $(q,c_{\ref{lemma:Kt_free_graphs_have_centered_Helly_colorings}}(t))$-good coloring of $G$.
    Fix a subgraph $G_0$ of $G$. 

    Now, we reduce to the case where $G_0$ is connected. 
    Indeed, if $G_0$ is not connected, consider the family $\mathcal{C}$ of all the components of $G_0$.
    Let $d$ be a nonnegative integer and let $\mathcal{F}$ be a family of connected subgraphs of $G_0$ such that there are no $d+1$ disjoint members of $\mathcal{F}$.
    For every $C \in \mathcal{C}$, let $d_C$ be the smallest integer such that there are no $d_C+1$ disjoint members of $\mathcal{F}_C = \{F \in \mathcal{F} \mid F \subseteq C\}$.
    Clearly $d_C \leq d$. 
    Assuming we can prove the result when $G_0$ is connected,
    we apply it for $C$ and $\mathcal{F}_C$ and obtain a set $Z_C \subseteq V(C)$ and $\psi_C \colon Z_C \to [c \cdot d]$ such that \ref{centered-Helly-hitting}--\ref{centered-Helly-components} holds.
    Now take $Z = \bigcup_{C \in \mathcal{C}} Z_C$ and let $\psi$ be defined by $\psi(u) = \psi_C(u)$ for every $C \in \mathcal{C}$ and every $u \in Z_C$.
    Recall that each $F\in\calF$ is connected, so $F\subseteq C$ for some $C\in\calC$, and therefore $Z_C\cap V(F)\neq\emptyset$. 
    Thus,~\ref{centered-Helly-hitting} holds.
    Since every connected subgraph $H$ of $G_0$
    is a subgraph of $C$ for some $C \in \mathcal{C}$, 
    \ref{centered-Helly-coloring} holds. 
    Finally,~\ref{centered-Helly-components} holds because for every component $C'$ of $G-Z$, there exists $C \in \mathcal{C}$ such that $C' \subseteq C$.

    Therefore, from now on we assume that $G_0$ is connected.
    We start by building a normal pair $(\calU,\calP)$ where $\calU$ is a tree decomposition of $G_0$ obtained from ``gluing'' tree decompositions in $\calD$ along the edges of $T$.
    This construction is depicted in~\Cref{fig:new-td}.

\begin{figure}[tp]
    \centering
    \includegraphics{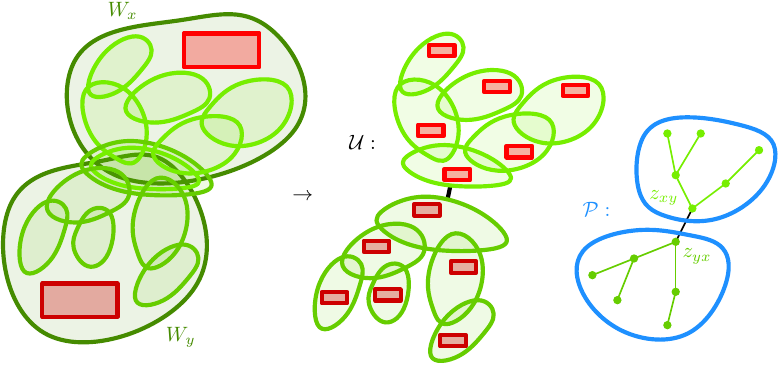}
    \caption{
        Tree decompositions $\calD_x$ and $\calD_y$ are ``glued'' in a natural way. 
        Additionally, we add all corresponding apices (members of $A_x$ and $A_y$) to each bag of the corresponding part in the tree decomposition $\calU$.
        The partition $\calP$ indicates from which of $\{T_x \mid x \in V(T)\}$ a given vertex comes from.
    }
    \label{fig:new-td}
\end{figure}

    For all $x,y \in V(T)$ with $xy \in E(T)$, $W_x \cap W_y$ induces a clique in both $\torso_{G,\mathcal{W}}(W_x)$ and $\torso_{G,\mathcal{W}}(W_y)$, hence, there exist $z_{xy} \in V(T_x)$ and $z_{yx} \in V(T_y)$ such that
    $W_x \cap W_y = D_{x,z_{xy}} \cap D_{y,z_{yx}}$.
    Let $T_0$ be the tree defined by
    \begin{align*}
        V(T_0) &= \{(x,z) \mid x \in V(T), z \in V(T_x)\} \text{ and }\\
        E(T_0) &= \{(x,z_1)(x,z_2) \mid x \in V(T), z_1z_2 \in E(T_x)\} \cup \{(x, z_{xy}) (y,z_{yx}) \mid xy \in E(T)\}.
    \end{align*}
    In other words, $T_0$ is obtained from the disjoint union of $T_x$ over $x \in V(T)$ by adding the edges between $z_{xy}$ and $z_{yx}$ for adjacent vertices $x$ and $y$ in $T$.
    Next, let $U_{(x,z)} = (D_{x,z} \cup A_x) \cap V(G_0)$ for every $(x,z) \in V(T_0)$.
    We claim that $\mathcal{U} = \big(T_0, (U_{(x,z)} \mid (x,z) \in V(T_0))\big)$ is a tree decomposition of $G_0$.
    For every $v\in V(G_0)$, $T_0[\{(x,z)\in V(T_0)\mid v\in U_{(x,z)}\}]$ is isomorphic to the connected subgraph of $T_0$ formed by replacing every vertex $x$ of $T(v)=T[\{x\in V(T)\mid v\in W_x\}]$ with $T_x[\{z\in V(T_x)\mid v\in D_{x,z}\}]$ when $v\notin A_x$ or with $T_x$ when $v\in A_x$, and every edge $xy$ of $T(v)$ with $z_{xy}z_{yx}$.
    For every edge $uv\in E(G_0)$, there exists $x\in V(T)$ such that $u,v\in W_x$.
    If $u, v \in W_x \setminus A_x$, then there exists $z\in V(T_x)$ such that $u,v\in D_{x,z}$ 
    and thus $u,v\in U_{(x,z)}$.
    If $u\in A_x$ and $v\in W_x \setminus A_x$, then for any $z\in V(T_x)$ with $v\in D_{x,z}$ we have $v\in U_{(x,z)}$ and $u\in A_x\cap V(G_0)\subseteq U_{(x,z)}$.
    If $u,v\in A_x$, then for every $z \in V(T_x)$, we have $u,v\in U_{(x,z)}$.
    Thus, indeed $\mathcal{U}$ is a tree decomposition of $G_0$.
    
    Observe that $\mathcal{P} = \{T_0[\{x\} \times V(T_x)] \mid x \in V(T)\}$ is a collection of subtrees of $T_0$ whose vertex sets partition $V(T_0)$. 
    Therefore, $\big(\mathcal{U},\mathcal{P}\big)$ is a normal pair of $G_0$.
    
    By \Cref{lemma:making_a_td_natural}, 
    there exists $\mathcal{V} = \big(S,(V_y \mid y \in V(S))\big)$
    and $\mathcal{Q}$ such that $\mathcal{V}$ is natural and $(\mathcal{V},\mathcal{Q})$ is a normal pair of $G$ refining $(\mathcal{U},\mathcal{P})$. 
    Among all such pairs $(\mathcal{V},\mathcal{Q})$, we take one with $|\mathcal{Q}|$ minimum.
    Let $f\colon V(S) \to V(T_0)$ and $g \colon \mathcal{Q} \to \mathcal{P}$ witness the refinement relation.
    For every $Q \in \mathcal{Q}$, let $x(Q) \in V(T)$ 
    be such that $g(Q) = T_0[\{x(Q)\} \times V(T_{x(Q)})]$.
    See~\Cref{fig:VQ}.

    In the next two claims, we show that the new tree decomposition $\calV$ in some sense preserves small adhesions (only ones coming from $\calW$) and small intersections with layers of the layerings in $\calL$.
    In the first claim, we exploit the minimality of $(\calV,\calQ)$. 
    \begin{claim}\label{claim:small_adhesions_between_Q}
        Let $y,y' \in V(S)$ and $Q,Q' \in \mathcal{Q}$ be such that
        $y \in V(Q)$, $y' \in V(Q')$ and $Q \neq Q'$.
        Then
        \[
        |V_y \cap V_{y'}| \leq \cLRS(t).
        \]
    \end{claim}

    \begin{proofclaim}
        By properties of tree decompositions,
        there is an edge $zz'$ in $S$ such that $V_y \cap V_{y'} \subseteq V_{z} \cap V_{z'}$.
        Therefore, without loss of generality, we assume that $yy'$ is an edge in $S$.
        We argue that $g(Q) \neq g(Q')$.
        Suppose to the contrary that $g(Q)=g(Q')$. 
        Since $yy'$ is an edge in $S$, the subgraphs $Q$ and $Q'$ are adjacent in $S$, and thus, $Q'' = S[V(Q) \cup V(Q')]$ is a connected subgraph of $S$.
        Observe that $(\mathcal{V}, \mathcal{Q} \setminus \{Q,Q'\} \cup \{Q''\})$ is a normal pair which refines $(\mathcal{U},\mathcal{P})$
        as witnessed by $f$ and $g'$ defined by $g'(P) = g(P)$ for every $P \in \mathcal{Q} \setminus \{Q,Q'\}$ and
        $g'(Q'') = g(Q) = g(Q')$.
        However, this contradicts the minimality of $|\mathcal{Q}|$.
        We obtain that indeed $g(Q) \neq g(Q')$, and so, $x(Q) \neq x(Q')$.
        In particular,
        \[
            |V_y \cap V_z|
            \leq |U_{(x(Q),f(y))} \cap U_{(x(Q'),f(y'))}| 
            \leq |W_{x(Q)} \cap W_{x(Q')}|  
            \leq \cLRS(t),
        \]
        where the first inequality follows from~\ref{item:making_a_td_natural:i}, the second from the properties of tree decompositions, and the last from~\ref{LRS:adhesion}.
    \end{proofclaim}

\begin{figure}[tp]
    \centering
    \includegraphics{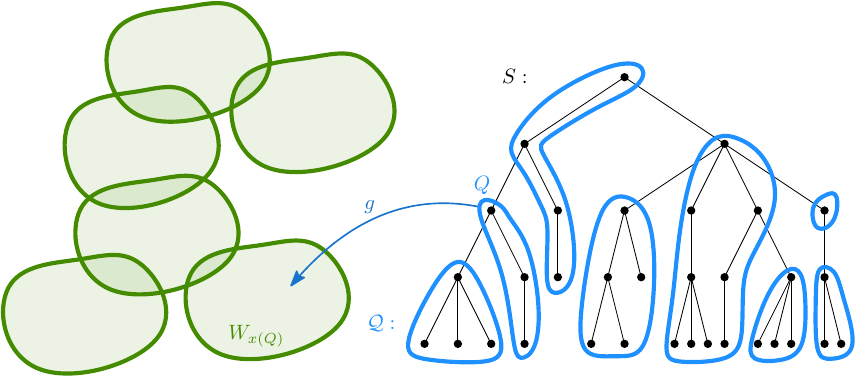}
    \caption{
        Intuitively, the part of the tree decomposition $\calV$ corresponding to the vertices in $Q$ comes from $\calD_{x(Q)}$. 
    }
    \label{fig:VQ}
\end{figure}
    
    \begin{claim}\label{claim:Vy_intersection_Lxi_small}
    Let $x\in V(T)$, $y\in V(S)$, and let $i\in \NN$. Then 
    \[
        |V_y\cap L_{x,i}|\leq  \cLRS(t).
    \]
    \end{claim}
    \begin{proofclaim}
    Since $(\calV,\calQ)$ refines $(\calU,\calP)$, we have
    $V_y\subseteq U_{f(y)}$. Let $f(y)=(x',z)$ where $x'\in V(T)$ and $z\in V(T_{x'})$. Thus, 
    $U_{f(y)}=U_{(x',z)}=(D_{x',z}\cup A_{x'})\cap V(G_0)$.
    Recall that $A_x\cap L_{x,i}=\emptyset$. Thus in the case of $x=x'$,
    \[
        |V_y\cap L_{x,i}| 
        \leq |U_{(x,z)} \cap L_{x,i}| 
        \leq |(D_{x,z}\cup A_x)\cap L_{x,i}|
        \leq |D_{x,z}\cap L_{x,i}|
        \leq \cLRS(t),
    \]
    where the last inequality follows by~\ref{LRS:ltw}.
    Finally, in the case of $x\neq x'$, we have 
    $U_{(x',z)}\subseteq W_{x'}$ and $L_{x,i}\subseteq W_{x}$, and so,
    \[
        |V_y\cap L_{x,i}| 
        \leq |U_{(x',z)} \cap L_{x,i}| 
        \leq |W_{x'}\cap W_x|
        \leq \cLRS(t),
    \]
    where the last inequality follows from~\ref{LRS:adhesion}.
    \end{proofclaim}
    
    We root $S$ in an arbitrary vertex $\root(S)$.
    For every $X \subset V(S)$, we define
        \[\calQ(X) = \set{Q\in\calQ\mid\ V(Q) \cap X \neq \emptyset}.\]

    Let $d$ be a positive integer and let $\mathcal{F}$ be a family of connected subgraphs of $G_0$
    such that there are no $d+1$ pairwise disjoint members of $\mathcal{F}$.
    In the remainder of the proof, we construct $Z$ satisfying \ref{centered-Helly-hitting}--\ref{centered-Helly-restriction}.
    
    By \Cref{lemma:helly_property_tree_decomposition}, there exists $X_0 \subseteq V(S)$ of size at most $d$ such that $\bigcup_{y \in X_0} V_y$ intersects every member of $\mathcal{F}$.
    Let
        \[X_1 = \LCA(S,X_0).\]
    By~\Cref{lemma:increase_X_in_a_tree},  
    \begin{equation}\label{eq:X1-and-Q1}
    |\calQ(X_1)|\leq |X_1| \leq 2d-1.
    \end{equation}

    The set $\bigcup_{y \in X_1} V_y$ is a hitting set of $\calF$, thus, it is a good candidate to be $Z$.
    However, we still need to define $\psi$ and the ultimate goal is to repeat the idea of coloring that we applied in the bounded layered treewidth case.
    Since $\calL_x$ is a layering of $\torso_{G,\mathcal{W}}(W_x)-A_x$, we need to take into account the vertices in $A_x$ when building the final $Z$.
    Moreover, the coloring $\varphi$ may not be compatible with layerings on the adhesions of $\calW$, we also need to consider some of them.
    To this end, we will add some more vertices to $Z$, which we ultimately color injectively with a separate palette of colors.
    Because we work with $\calV$ instead of $\calU$ we need the following notions of projections. 
    See also~\Cref{fig:projections} for an illustration.

\begin{figure}[tp]
    \centering
    \includegraphics{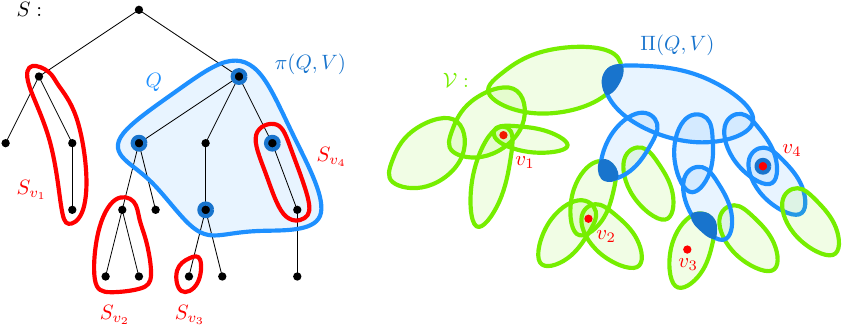}
    \caption{
        Here, $V = \{v_1,v_2,v_3,v_4\}$.
        Note that $\bigcup_{y \in V(Q)} V_y$ is separated from $V$ by $\Pi(Q,V)$ in $G_0$ as we show later in~\Cref{claim:Pi-separates-stuff}. 
    }
    \label{fig:projections}
\end{figure}

    Let $Q \in \mathcal{Q}$
    and $v \in V(G_0)$.
    Let $S_v = S[\set{z\in V(S)\mid v\in V_z}]$.
    Thus, $Q$ and $S_v$ are two subtrees of $S$. 
    We define projections $\pi(Q,v)$ and $\Pi(Q,v)$ as follows.
    If $Q$ and $S_v$ share a vertex, then choose such a vertex $z$ arbitrarily and set $\pi(Q,v) = z$ and $\Pi(Q,v)=\set{v}$. 
    If $V(Q)$ and $V(S_v)$ are disjoint, then 
    consider the shortest \stpath{V(Q)}{V(S_v)} path in $S$. 
    Let $z$ be the endpoint of that path in $Q$ and 
    let $y$ be the vertex adjacent to $z$ in that path. 
    Then set $\pi(Q,v) = z$ and $\Pi(Q,v) = V_y\cap V_{z}$.
    The definition can be naturally extended to subsets of $V(G_0)$.
    For each $V\subseteq V(G_0)$, let 
    \[\text{$\pi(Q,V)=\set{\pi(Q,v)\mid v\in V}$ \ and \ 
    $\Pi(Q,V)=\bigcup_{v \in V} \Pi(Q,v)$.}\]
    Note that by~\Cref{claim:small_adhesions_between_Q},
    \begin{equation}\label{eq:size_Pi}
        |\Pi(Q,V)| \leq \cLRS(t) \cdot |V|.
    \end{equation}
    The key property of these objects is the following. 
    See~\Cref{fig:projections} again.
    \begin{claim}\label{claim:Pi-separates-stuff}
        Let $Q \in \calQ$ and let $V \subset V(G_0)$.
        Every connected subgraph $H$ of $G_0 - \Pi(Q,V)$ that intersects $\bigcup_{y \in V(Q)} V_y$ is disjoint from $V$.
    \end{claim}
    \begin{proof}
        It suffices to observe that by the properties of tree decompositions and construction, $\Pi(Q,V)$ intersects every path between $\bigcup_{y \in V(Q)} V_y$ and $V$ in $G_0$.
        Thus, if a connected subgraph $H$ of $G_0$ has a vertex is both $\bigcup_{y \in V(Q)} V_y$ and $V$, then it has a vertex in $\Pi(Q,V)$.
    \end{proof}

    Let
    \[
    Y_Q =
    \begin{cases}
        \pi\left(Q, A_{x(Q)}\right) &\textrm{if $x(Q)=r$,}\\
        \pi\left(Q, A_{x(Q)} \cup \left(W_{x(Q)} \cap W_{\parent(T,x(Q))}\right)\right) & \textrm{if $x(Q) \neq r$.}
    \end{cases}
    \]
    Recall that $|A_{x(Q)}|\leq \cLRS(t)$ by~\ref{LRS:apices},
    and if $x(Q)\neq r$, then $|W_{x(Q)} \cap W_{\parent(T,x(Q))}|\leq \cLRS(t)$ by~\ref{LRS:adhesion}.
    Thus, we have 
    \begin{equation}
    \begin{aligned}
    |Y_Q|
    & \leq |A_{x(Q)}| 
    &&\leq \cLRS(t) && \textrm{if $x(Q)=r$,}\\
    |Y_Q|
    & 
    \leq |A_{x(Q)}|+|W_{x(Q)}\cap W_{\parent(T,x(Q))}| 
    &&\leq 2\cLRS(t) && \textrm{if $x(Q)\neq r$.}
    \end{aligned}
    \label{eq:YQ}
    \end{equation}
    Moreover, by definition $Y_Q \subset V(Q)$.

    Let
    \[
        X_2 = X_1 \cup \{\root(Q) \mid Q \in \mathcal{Q}(X_1)\} \cup \bigcup \{ Y_Q \mid Q \in \mathcal{Q}(X_1)\}
    \]
    and
    \[
        X_3 = \LCA(S,X_2).
    \]
    Observe that
    \begin{equation}
    \begin{aligned}
    |X_3|
    &\leq 2|X_2| - 1 
    &&\textrm{by~\Cref{lemma:increase_X_in_a_tree}}\\
    &\leq 2(|X_1| + |\calQ(X_1)| + |\calQ(X_1)|\cdot 2\cLRS(t))-1
    &&\textrm{by~\eqref{eq:YQ}}\\
    &\leq 2((2d-1) + (2d-1) + (2d-1)\cdot 2\cLRS(t))-1 &&\textrm{by~\eqref{eq:X1-and-Q1}}\\
    &\leq (8\cLRS(t)+8)d.
    \end{aligned}\label{eq:size-X_3}
    \end{equation}

    Recall that by~\Cref{lemma:LCA-stuff}, for all $X, Y \subseteq V(S)$ with $X\subseteq Y$ and $\LCA(S,X) = X$, if $\calQ(X) = \calQ(Y)$, then $\calQ(Y)=\calQ(\LCA(S,Y))$.
    Since by construction $\calQ(X_1) = \calQ(X_2)$, we can apply the above with $X = X_1$ and $Y = X_2$ to obtain that
    \[
        \mathcal{Q}(X_1) = \calQ(X_2) = \mathcal{Q}(X_3).
    \]
    
    Finally, the set $Z$ that witnesses the assertion of the claim is given by
    \[
        Z = \bigcup_{y \in X_3} V_y.
    \]
    
    Next, consider $Q \in \mathcal{Q}(X_1)$
    and let 
    \[
        B_Q =
        \begin{cases}
            \Pi\left(Q, A_{x(Q)}\right) &\textrm{if $x(Q)=r$,} \\
            \Pi\left(Q, A_{x(Q)} \cup \left(W_{x(Q)} \cap W_{\parent(T,x(Q))}\right)\right) & \textrm{if $x(Q) \neq r$.}
        \end{cases}
    \]
    Since for every $v \in V(G_0)$, $\Pi(Q,v) \subseteq V_{\pi(Q,v)}$ and $Y_Q \subseteq X_2 \subseteq X_3$, we have 
    \[
        B_Q \subseteq \bigcup_{y \in Y_Q} V_y \subseteq \bigcup_{y \in X_2} V_y \subseteq Z.
    \]
    Recall that $|A_{x(Q)}|\leq \cLRS(t)$ by~\ref{LRS:apices},
    and if $x(Q)\neq r$, then $|W_{x(Q)} \cap W_{\parent(T,x(Q))}|\leq \cLRS(t)$ by~\ref{LRS:adhesion}.
    Thus, by \eqref{eq:size_Pi}, we have
    \begin{equation}
        \begin{aligned}
            |B_Q|
            &\leq |A_{x(Q)}| \cdot \cLRS(t)
            &&\leq (\cLRS(t))^2
            && \textrm{if $x(Q)=r$,}\\
            |B_Q|
            &\leq \left|A_{x(Q)}\cup(W_{x(Q)}\cap W_{\parent(T,x(Q))}))\right| \cdot \cLRS(t) 
            &&\leq 2(\cLRS(t))^2
            &&\textrm{if $x(Q) \neq r$.}
        \end{aligned}
        \label{eq:A'Q}
    \end{equation}
    
    Let 
    \[
        B = \left(\bigcup_{Q \in \mathcal{Q}(X_1), \root(S) \not\in V(Q)} V_{\root(Q)} \cap V_{\parent(S,\root(Q))}\right) \cup \left(\bigcup_{Q \in \mathcal{Q}(X_1)} B_Q\right).
    \]
    Recall that for every $Q \in \mathcal{Q}(X_1)$ with $\root(S) \not\in V(Q)$,
    $\root(Q) \in X_2$ and so $V_{\root(Q)} \cap V_{\parent(S,\root(Q))} \subseteq V_{\root(Q)} \subseteq Z$. 
    Moreover, for every $Q \in \mathcal{Q}(X_1)$, $B_Q \subseteq Z$.
    Therefore,
    \[
        B\subseteq Z. 
    \]
    Observe also that
    \begin{equation}
        \begin{aligned}
            |B| 
            &\leq |\mathcal{Q}(X_1)|\cdot \cLRS(t)  + |\mathcal{Q}(X_1)|\cdot 2(\cLRS(t))^2 && \textrm{by \Cref{claim:small_adhesions_between_Q} and \eqref{eq:A'Q}} \\
            &\leq (2d-1)\cdot \cLRS(t) 
            + (2d-1)\cdot 2(\cLRS(t))^2
            &&\textrm{by \eqref{eq:X1-and-Q1}}\\
            & \leq \cLRS(t)(4\cLRS(t)+2)d.
        \end{aligned}
        \label{eq:size-of-B}
    \end{equation}
    
    \begin{claim}\label{claim:B_separates_the_Qs}
        For all distinct $Q_1,Q_2 \in \mathcal{Q}(X_1)$,
        there is no path between $\bigcup_{y \in V(Q_1)} V_y \setminus B$ and
        $\bigcup_{y \in V(Q_2)} V_y \setminus B$ in $G_0 - B$.
        In particular, $\bigcup_{y \in V(Q_1)} V_y \setminus B$ and $\bigcup_{y \in V(Q_2)} V_y \setminus B$ are disjoint.
    \end{claim}

    \begin{proofclaim}
        Let $Q_1,Q_2$ be distinct members of $\mathcal{Q}(X_1)$.
        Recall that $Q_1,Q_2$ are disjoint subtrees of $S$.
        Note that there exists $i \in \{1,2\}$
        such that $\root(Q_i) \neq \root(S)$ and every path in $S$ between $Q_1$ and $Q_2$
        goes through the edge $\root(Q_i)\parent(S,\root(Q_i))$ of $S$.
        Hence by properties of tree decompositions, 
        $V_{\root(Q_i)} \cap V_{\parent(S,\root(Q_i))}$ intersects 
        every path between $\bigcup_{y \in V(Q_1)} V_y$ and $\bigcup_{y \in V(Q_2)} V_y$ in $G_0$.
        Since $V_{\root(Q_i)} \cap V_{\parent(S,\root(Q_i))} \subseteq B$, this proves the claim.
    \end{proofclaim}

    As a consequence of \Cref{claim:B_separates_the_Qs},
    $\big\{\big(\bigcup_{y \in V(Q)} V_y\big) \cap (Z \setminus B) \mid Q \in \mathcal{Q}(X_1)\big\}$ is a family of pairwise disjoint sets covering $Z \setminus B$ (because $\mathcal{Q}(X_1) = \mathcal{Q}(X_3)$).
    We want to refine this family 
    with the layerings of the torsos.
    To this end, we need the following fact where we substantially use the fact that $(\calV,\calQ)$ refines $(\calU,\calP)$.
    \begin{claim}\label{claim:V(Q)-stays-in-one-W}
        For every $Q \in \calQ(X_1)$ and $V \subset V(G_0)$ such that $\Pi(Q,V) \subset B_Q$, we have
            \[\bigcup_{y \in V(Q)} V_y \setminus B \subseteq W_{x(Q)} \setminus V.\]
    \end{claim}
    \begin{proofclaim}
        Let $Q \in \calQ(X_1)$ and let $V \subset V(G_0)$ be such that $\Pi(Q,V) \subset B_Q$.
        Since $\Pi(Q,V) \subset B_Q \subset B$, then
            \[\bigcup_{y \in V(Q)} V_y \setminus B \subset \bigcup_{y \in V(Q)} V_y \setminus B_Q \subset \bigcup_{y \in V(Q)} V_y \setminus \Pi(Q,V).\]
        By~\cref{claim:Pi-separates-stuff}, $\bigcup_{y \in V(Q)} V_y \setminus \Pi(Q,V)$ is disjoint from $V$.
        Thus, to conclude the claim, it suffices to show that $\bigcup_{y \in V(Q)} V_y \subset W_{x(Q)}$.
        Recall that $g(Q) = T_0[\{x(Q)\}\times T_{x(Q)}]$.
        We have
            \[\bigcup_{y \in V(Q)} V_y \subset \bigcup_{y \in V(Q)} U_{f(y)} \subset \bigcup_{z \in V(g(Q))} U_z = \bigcup_{z \in \{x(Q)\} \times T_{x(Q)}} U_z \subset W_{x(Q)}\]
        where the first inclusion follows from~\ref{item:making_a_td_natural:i}, the second from~\ref{item:making_a_td_natural:ii}, and the last one from the construction of $\calU$.
    \end{proofclaim}
    
    Since  $\Pi(Q,A_{x(Q)}) \subset B_Q$ for every $Q \in \calQ(X_1)$, by~\cref{claim:V(Q)-stays-in-one-W},
        \[\textstyle \big(\bigcup_{y \in V(Q)} V_y \setminus B\big) \cap A_{x(Q)} = \emptyset. \]
    Recall that for every $Q \in \calQ(X_1)$, $(L_{x(Q),i} \mid  i \in \NN)$ is a layering of
    $\torso_{G,\mathcal{W}}(W_{x(Q)}) - A_{x(Q)}$.
    It follows that the family
    \[\textstyle
    \calB = \left\{\big(\bigcup_{y \in V(Q)} V_y\big) \cap L_{x(Q),i} \cap (Z \setminus B) \mid Q \in \mathcal{Q}(X_1), i \in \NN \right\}
    \]
    is a family of pairwise disjoint sets covering $Z \setminus B$.
    Moreover, the members of $\calB$ are of reasonable size, namely, for every $Q \in \calQ(X_1)$ and $i \in \NN$,
    \begin{equation} 
    \begin{aligned}
    \textstyle
        \left|\big(\bigcup_{y \in V(Q)} V_y\big) \cap L_{x(Q),i} \cap (Z \setminus B)\right| &\leq |L_{x(Q),i} \cap Z| \leq \sum_{y \in X_3} |L_{x(Q),i} \cap V_y|\\ 
         &\leq \sum_{y \in X_3} |L_{x(Q),i} \cap U_{f(y)}| && \text{by~\ref{item:making_a_td_natural:i}}\\
        &\leq \cLRS(t) \cdot |X_3| && \text{by~\ref{LRS:ltw}}\\
        &\leq \cLRS(t) \cdot \big(8\cLRS(t)+8\big)d && \text{by~\eqref{eq:size-X_3}.}
    \end{aligned}
    \label{eq:bound_on_layer_cap_Z}
    \end{equation}

\begin{figure}[tp]
    \centering
    \includegraphics{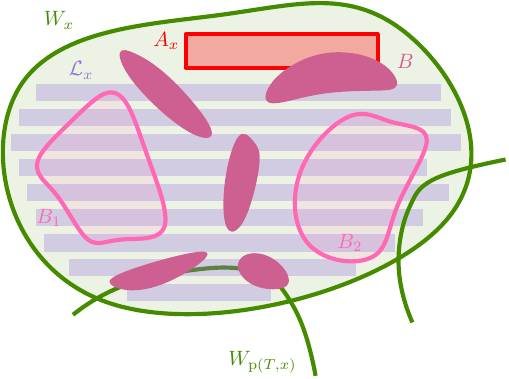}
    \caption{
        The set $Z$ is the union of $B$ and the members of $\calB$.
        We color the vertices in $B$ injectively.
        The role of the set $B$ is to separate members of $\calB$ so that if $B$ is intersected by a connected subgraph $H$ of $G_0$, we immediately get a center.
        In the figure, $B_1 = \bigcup_{y \in V(Q_1)} V_y \setminus B$ and $B_2 = \bigcup_{y \in V(Q_2)} V_y \setminus B$ for some $Q_1,Q_2 \in \calQ(X_3)$.
        If a connected subgraph $H$ of $G_0$, intersects both $B_1$ and $B_2$, it has to intersect $B$.
    }
    \label{fig:constructing-psi}
\end{figure}

    We define a coloring $\psi$ of $Z$ using at most $c_{\ref{lemma:Kt_free_graphs_have_centered_Helly_colorings}}(t)$ colors.
    See also~\Cref{fig:constructing-psi}.
    First, we color $B$ injectively, and then we color each member of $\calB$ also injectively avoiding colors in $\psi(B)$.
    In the first step we used at most $\cLRS(t)(4\cLRS(t)+2)d$ colors by~\eqref{eq:size-of-B} and in the second step, we used at most $\cLRS(t) \cdot \big(8\cLRS(t)+8\big)d$ colors by~\eqref{eq:bound_on_layer_cap_Z}, thus, $\psi$ is well-defined.

    We now show that $Z$ satisfy \ref{centered-Helly-hitting}--\ref{centered-Helly-restriction}.
    
    Recall that $X_0$ was chosen so that $X_0 \subseteq X_3$ and $\bigcup_{y \in X_0} V_y$ intersects every member of $\mathcal{F}$.
    Therefore,~\ref{centered-Helly-hitting} holds.
    Item~\ref{centered-Helly-components} holds
    by \Cref{lemma:increase_X_to_have_small_interfaces} since $Z = \bigcup_{y \in \LCA(S, X_2)} V_y$, and $\mathcal{V} = \big(S,(V_y \mid y \in V(S))\big)$ is a natural tree decomposition of $G_0$.

    It remains to prove~\ref{centered-Helly-restriction}.
    Let $P \subseteq Z$ and let $R = P \cup B$.
    First, $|R \setminus P| \leq |B| \leq \cLRS(t)(4\cLRS(t)+2)d.$
    Consider a connected subgraph $H$ of $G_0$ such that $V(H) \cap P \neq \emptyset$.
    If $V(H) \cap B \neq \emptyset$,
    then since every vertex of $B$ has a unique color, any vertex in $V(H) \cap B$ is a $\psi\vert_R$-center of $V(H) \cap P$, and so, a $(\phi \times \psi)$-center of $V(H) \cap Z$.
    Thus, we assume that $V(H) \cap B = \emptyset$.
    
    Since $V(H) \cap B = \emptyset$, it follows from~\Cref{claim:B_separates_the_Qs}, that $V(H) \cap Z \subseteq \bigcup_{z \in V(Q)} V_z$ for some $Q \in \mathcal{Q}(X_1)$.

    Additionally, by \Cref{claim:V(Q)-stays-in-one-W},
    \begin{align*}
        V(H) \cap Z &\subseteq W_x \setminus A_x && \textrm{if $x = r$,}\\
        V(H) \cap Z &\subseteq W_x \setminus (A_x \cup W_{\parent(T,x)}) && \textrm{if $x \neq r$.}
    \end{align*}
    Recall that $(L_{x,i} \mid i \in \mathbb{N})$ is a layering of $\torso_{G,\mathcal{W}}(W_x)-A_x$.
    Moreover, by definition of $\phi$,
    for every $u \in V(H) \cap Z$, $\phi(u) = i \bmod (q+1)$ where $i \in \NN$ is such that $u \in L_{x,i}$.

    Consider $u \in V(H) \cap R$.
    Let $i\in \NN$ be such that $u \in L_{x,i}$.
    If $u$ is a $(\phi \times \psi)$-center in $V(H) \cap R$, then we are done.
    Assume now that there exists $u' \in V(H) \cap R$ distinct from $u$
    such that $\phi(u)=\phi(u')$ and $\psi(u) = \psi(u')$.
    Let $j \in \NN$ be such that $u' \in L_{x,j}$.
    Without loss of generality, assume that $i \leq j$.
    By the definition of $\psi$, $j \neq i$,
    and by the definition of $\phi$, $|j-i| > q$.
    By \Cref{lemma:projections_on_a_torso_stay_connected},
    $V(H) \cap W_x$ induces 
    a connected subgraph of $\torso_{G,\mathcal{W}}(W_x)-A_x$.
    Since $(L_{x,k} \mid k \in \NN)$ is a layering of $\torso_{G,\mathcal{W}}(W_x)-A_x$,
    it follows that $V(H)$ intersects $L_{x,k}$ for every $k \in \{i,\dots, j-1\}$. 
    We deduce that $|\phi(V(H))| > q$.
    This proves that $\psi\vert_{V(H)}$ witnesses the fact that
    \[
        \cen_q(G,\phi,R) \leq \big(\cLRS(t)(4\cLRS(t)+2) + \cLRS(t) (8\cLRS(t)+8)\big) \cdot d
    \]
    Therefore, we obtain~\ref{centered-Helly-restriction}, which ends the proof.
\end{proof}

\subsection{The bounds}\label{sec:cen:the_bounds}
In the case of centered chromatic numbers, we handle the bounds involving treewidth and the general bounds in~\Cref{thm:centered} almost in the same way: the only distinction is applying either~\Cref{lemma:p1_good_colorings_when_bounded_tw} or~\Cref{lemma:Kt_free_graphs_have_centered_Helly_colorings}.
Note that this differs from how we proved~\Cref{thm:wcol} and~\Cref{thm:fragility}.
\Cref{thm:centered} will follow from~\Cref{thm:main_cen_St,thm:main_cen_Rt}.
Already developed material is almost enough to prove~\Cref{thm:main_cen_St}, where we exclude a graph from $\SRt_t$.
To complete the proof, we need to define an appropriate family of focused parameters related to centered colorings, fitting the framework from~\Cref{ssec:centered-nice}.
In order to prove~\Cref{thm:main_cen_Rt} (excluding a graph from $\Rt_t$), we also need to inspect the base cases more carefully.

Let $\ell$ be a positive integer.
Let \defin{$(G_{\ell,i} \mid i \in \posint)$} be a collection of all pairwise non-isomorphic (finite) $K_\ell$-minor-free graphs on pairwise disjoint vertex sets.
Let \defin{$(H_{\ell,i} \mid i \in \posint)$} be a collection of all pairwise non-isomorphic (finite) graphs of treewidth less than $\ell$ on pairwise disjoint vertex sets.

We define \defin{$\calK_\ell$} as the disjoint union of all the graphs $G_{\ell,i}$ for $i \in \posint$,
and we define \defin{$\mathcal{T}_\ell$} as the disjoint union of all the graphs $H_{\ell,i}$ for $i \in \posint$.
Note that the graphs $\calK_\ell$ and $\calT_\ell$ are infinite.

Let $q$ be a positive integer.
For each $i \in \posint$, we fix a $(q,c_{\ref{lemma:Kt_free_graphs_have_centered_Helly_colorings}}(\ell))$-good coloring \defin{$\varphi_{\ell,q,i}$} of $G_{\ell,i}$ using $q+1$ colors that exists by~\Cref{lemma:Kt_free_graphs_have_centered_Helly_colorings}, and we fix a $(q,2)$-good coloring \defin{$\psi_{\ell,q,i}$} of $H_{\ell,i} $ using $\ell$ colors that exists by~\Cref{lemma:p1_good_colorings_when_bounded_tw}.
We may assume that the colorings $\varphi_{\ell,q,i}$ for each $i \in \posint$ use the same set of $q+1$ colors, and that the colorings $\psi_{\ell,q,i}$ for each $i \in \posint$ use the same set of $\ell$ colors
Let \defin{$\Phi_{\ell,q}$} be the coloring of $\calK_\ell$ such that for every $i \in \posint$, we have $\Phi_{\ell,q}\vert_{V(G_{\ell,i})} = \varphi_{\ell,q,i}$.
Similarly, let \defin{$\Psi_{\ell,q}$} be the coloring of $\calT_\ell$ such that for every $i \in \posint$, we have $\Psi_{\ell,q}\vert_{V(H_{\ell,i} )} = \psi_{\ell,q,i}$.

Finally, we define \defin{$\mathcal{T}$} as the disjoint union of $\mathcal{T}_\ell$ for $\ell \in \posint$, and \defin{$\Psi_q$} as the coloring of $\mathcal{T}$ such that $\Psi_q\vert_{V(\mathcal{T}_\ell)} = \Psi_{q,\ell}$ for every $\ell \in \posint$.
\Cref{obs:qc-good-on-components} implies the following statement.

\begin{obs}\label{obs:combined-coloring-is-good}
    Let $\ell$ and $q$ be positive integers.
    For every finite subgraph $G$ of $\calK_\ell$, $\Phi_{\ell,q}\vert_{V(G)}$ is a $(q,c_{\ref{lemma:Kt_free_graphs_have_centered_Helly_colorings}}(\ell))$-good coloring of $G$ using $q+1$ colors.
    Similarly, for every finite subgraph $H$ of $\calT_\ell$, $\Psi_{\ell,q}\vert_{V(H)}$ is a $(q,2)$-good coloring of $H$ using $\ell$ colors.
\end{obs}

For every positive integer $q$,
for every graph $G$, and
for every $S \subseteq V(G)$, let
\[
    \defmath{\paramktmf_{\ell,q}(G,S)} =
    \begin{cases}
        \cen_q(G,\Phi_{\ell,q}\vert_{V(G)},S) & \textrm{if $V(G) \subseteq V(\mathcal{K}_\ell)$ and $E(G) \subseteq E(\mathcal{K}_\ell)$,} \\
        0 & \textrm{otherwise;}
    \end{cases}
\]
for every positive integer $\ell$,
and let
\[
    \defmath{\parambdtw_{q}(G,S)} =
    \begin{cases}
        \cen_q(G,\Psi_{q}\vert_{V(G)},S) & \textrm{if $V(G) \subseteq V(\mathcal{T})$ and $E(G) \subseteq E(\mathcal{T})$,} \\
        0 & \textrm{otherwise.}
    \end{cases}
\]
Note that $\parambdtw$ is well-defined since $\mathcal{T}$ is the disjoint union of $\mathcal{T}_\ell$ for $\ell \in \posint$.
Moreover, the families $(\paramktmf_{\ell,q} \mid q \in \posint)$ for every positive integer $\ell$, and $(\parambdtw_{q} \mid q \in \posint)$ are families of focused parameters.
We emphasize that these parameters are not invariant under graph-isomorphism.

\begin{theorem}\label{thm:main_cen_St}
    Let $t$ be an integer with $t \geq 2$.
    For every $X \in \SRt_t$, there exists an integer $c$ such that for every positive integer $q$ and for every $X$-minor-free graph $G$,
    \begin{align*}
        \cen_q(G) &\leq c \cdot (\tw(G)+1) \cdot q^{t-2} \\
    \intertext{and}
        \cen_q(G) &\leq c \cdot q^{t-1}.
    \end{align*}
\end{theorem}

\begin{proof}
    Let $\ell$ be a positive integer, let $\paramktmf_\ell = (\paramktmf_{\ell,q} \mid q \in \posint)$ and $\parambdtw = (\parambdtw_{\ell,q} \mid q \in \posint)$, and let $\param \in \{\paramktmf_\ell, \parambdtw\}$.
    By~\Cref{lemma:centered_colorings_are_nice}, the family of focused parameters $\param$ is nice.
    By~\Cref{obs:combined-coloring-is-good,lem:cen-good-const-edgeless-bounding} the function $q \mapsto 1$ is $(\param,\edgeless)$-bounding.
    Therefore, by \Cref{cor:SRt-edgeless-to-t}, 
    $q \mapsto q^{t-2}$ is $(\param, \SRt_{t})$-bounding.
    By~\Cref{lemma:par-bounding-to-bound}, this implies that for every $X \in \SRt_t$, there exists a constant $\beta_\param(X)$ such that for every positive integer $q$ and for every graph $G$, we have
    \begin{equation}
        \param_q(G,V(G)) \leq \beta_\param(X) \cdot q^{t-2}.\label{eq:param-bound}
    \end{equation}

    Let $X \in \SRt_t$, 
    let $q$ be a positive integer, 
    and let $G$ be an $X$-minor-free graph.
    In particular, $G$ is $K_{|V(X)|}$-minor-free, hence, there exists $i \in \posint$ such that $G$ is isomorphic to $G_{|V(X)|,i}$, which is one of the building blocks of $\calK_{|V(X)|}$.
    Recall that $\Phi_{|V(X)|,q}\vert_{V(G_{|V(X)|,i})} = \phi_{|V(X)|,i,q}$ uses at most $q+1$ colors.
    It follows that for $\ell = |V(X)|$, we have
    \begin{align*}
        \cen_q(G) = \cen_q(G_{\ell,i}) &\leq \cen_q(G_{\ell,i},\Phi_{\ell,q}\vert_{V(G_{\ell,i})},V(G_{\ell,i})) \cdot (q+1) && \text{by~\Cref{lem:ordered-to-normal}}\\
        &= \paramktmf_{\ell}(G_{\ell,i},V(G_{\ell,i})) \cdot (q+1)\\
        &\leq \beta_{\paramktmf_\ell}(X) \cdot q^{t-2} \cdot (q+1) \leq (2\beta_{\paramktmf_\ell}(X)) \cdot q^{t-1} && \text{by~\eqref{eq:param-bound}.}
    \end{align*}
    This gives the second bound in the statement of the theorem.
    For the first bound, note that there exists $i \in \posint$ such that $G$ is isomorphic to $H_{\tw(G)+1,i}$, a building block of $\calT_{\tw(G)+1}$.
    Recall that $\Psi_{q}\vert_{V(H_{\tw(G)+1,i})} = \Psi_{\tw(G)+1,i,q}$ uses at most $\tw(G)+1$ colors.
    It follows that with $\ell = \tw(G)+1$, we have
    \begin{align*}
        \cen_q(G) = \cen_q(H_{\ell,i}) &\leq \cen_q(H_{\ell,i},\Psi_{\ell,q}\vert_{V(H_{\ell,i})},V(H_{\ell,i})) \cdot \ell && \text{by~\Cref{lem:ordered-to-normal}}\\
        &= \parambdtw(H_{\ell,i},V(H_{\ell,i})) \cdot \ell\\
        &\leq \beta_{\parambdtw}(X) \cdot q^{t-2} \cdot \ell && \text{by~\eqref{eq:param-bound}.}\\
        &= \beta_{\parambdtw}(X) \cdot (\tw(G)+1) \cdot q^{t-2}.
    \end{align*}
    This concludes the proof of the theorem. 
\end{proof}

\begin{figure}[tp]
    \centering
    \includegraphics{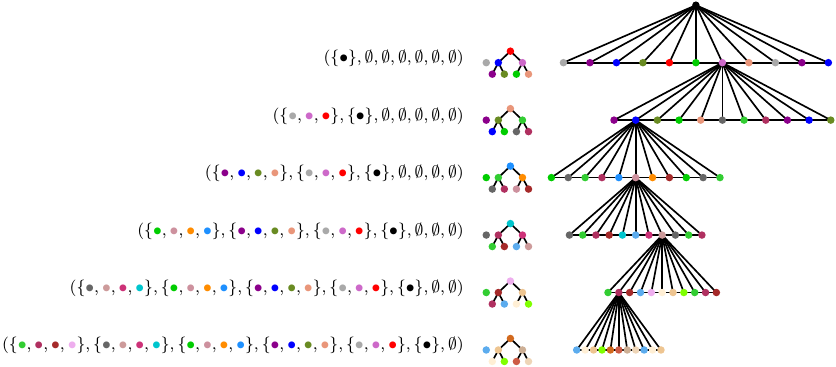}
    \caption{
        We present a sketch of the argument by Dębski et al.~\cite{Dbski2021} showing that outerplanar graphs have
        $q$-centered chromatic number in $\bigO(q \log q)$.
        Assume $q=2^k-1$ for some $k \in \NN$ ($k=3$ and $q=7$ in the picture).
        For simplicity, we consider a graph $G$ which is a \q{tree of fans}.
        Let $\Lambda$ be a set of $q(k+2)+1 = \bigO(q \log q)$ colors.
        We will construct from top to bottom a $q$-centered coloring $\psi$ of $G$ using only colors in $\Lambda$.
        We want $\psi$ to verify the following property: for every connected subgraph $H$ of $G$ using at most $q$ colors, the highest path intersecting $H$ in the tree of fans will contain a $\psi$-center of $V(H)$.
        We color the root of $G$ with an arbitrary color in $\Lambda$. 
        Now assume that a vertex $u$ of $G$ and all its ancestors are colored. 
        We proceed with a coloring of the path $P$ formed by the children of $u$ in the tree. 
        The coloring of layers above produces also some sets of forbidden colors to use on $P$. 
        This is encoded as a sequence $(\Lambda_1, \dots, \Lambda_q)$ of subsets of $\Lambda$, each of size at most $k+1$.
        The set $\Lambda_a$ for $a \in [q]$ represents the colors forbidden because of the $a$th previous layer.
        To color $P$, we use $q+1$ colors in $\Lambda \setminus \bigcup_{a \in [q]} \Lambda_a$.
        Among these colors, we choose one color to be special, and we arrange the $q=2^k-1$ other ones
        in a complete binary tree. Then, we color $P$ periodically with a pattern of length $q+1$ starting with
        the special color, followed by the $q$ colors in the auxiliary tree in the in-order.
        These auxiliary trees are depicted for some of the paths.
        Once $P$ is colored, it remains to color the subgraphs under $P$.
        Let $v \in V(P)$,
        and let $\tilde{\Lambda}_v$ be the set containing the special color and the ancestors of the color of $v$ in the auxiliary tree.
        Note that $|\tilde{\Lambda}_v| \leq k+1$.
        We apply the induction hypothesis
        to the subgraph under $v$
        with the sequence of sets of forbidden colors $(\tilde{\Lambda}_v, \Lambda_1, \dots, \Lambda_{q-1})$.
        It remains to show that the obtained coloring $\psi$ is $q$-centered.
        Let $H$ be a connected subgraph of $G$ using at most $q$ colors,
        and let $P$ be the highest path $P$ in this tree structure intersecting $V(H)$.
        Then $V(P) \cap V(H)$ induces a subpath of $P$.
        If $V(H) \cap V(P)$ contains a vertex with the special color,
        then this vertex is a $\psi$-center of $V(H)$, since $V(H)$ uses at most $q$ colors, and because the color of $v$ 
        is forbidden in every subtree below $P$ in the next $q$ layers.
        Otherwise, let $v \in V(H) \cap V(P)$ be such that the color $\psi(v)$ is the (unique) highest in the auxiliary tree.
        By construction, the color of $v$ is forbidden in the next $q$ layers of each subtree below $P$ intersecting $V(H)$,
        and so $v$ is a $\psi$-center of $V(H)$.
    }
    \label{fig:centered_coloring_tree_of_fans}
\end{figure}

In the reminder of this section,
we show the bounds of \Cref{thm:centered} for the case $X \in \Rt_t$.
The proof is similar to the one of \Cref{thm:main_cen_St},
except for the base case, which corresponds to excluding an apex-forest.
The base case is much more involved.
The core of this proof is inspired by a result of D\k{e}bski, Felsner, Micek, and Schröder~\cite{Dbski2021}
which asserts that outerplanar graphs have $q$-centered chromatic number in $\bigO(q \log q)$.
See \Cref{fig:centered_coloring_tree_of_fans} for a sketch of their argument.
Since outerplanar graphs exclude an apex-forest, namely $K_{2,3}$, this section can be seen as a generalization of their result.
As in \Cref{sec:common_base}, we also use ideas from Dujmovi\'c et al.'s paper~\cite{Dujmovi2023}.
We start by proving structural properties for graphs with no $\mathcal{F}$-rich model of a given star (\Cref{lemma:centered_coloring:base_case:excluding_a_star}).
Then, we extend it to graphs with no $\mathcal{F}$-rich model of a given forest (\Cref{lemma:centered_coloring:base_case:excluding_a_forest,lemma:centered_coloring:base_case:excluding_a_forest_small_interface}), and finally 
to graphs with no $\mathcal{F}$-rich model of a given apex-forest (\Cref{lemma:centered_coloring:base_case:excluding_an_apex_forest}).

\begin{lemma}\label{lemma:centered_coloring:base_case:excluding_a_star}
    Let $q$, $c$, and $d$ be positive integers,
    let $G$ be a connected graph,
    let $\phi$ be a $(q,c)$-good coloring of $G$,
    let $\mathcal{F}$ be a family of connected subgraphs of $G$ such that $G$ has no $\mathcal{F}$-rich model of $F_{2,d}$,
    and let $U \subseteq V(G)$ be such that $G[U]$ is connected.
    There exists $S \subseteq V(G)$, a path partition $(P_0, \dots, P_\ell)$ of $(G,S)$, and sets $R_1, \dots, R_\ell \subseteq V(G)$ such that
    \begin{enumerateOurAlph}
        \item $P_0 = U$; \label{item:centered_coloring:base_case:excluding_a_star:P0=U}
        \item $S \cap V(F) \neq \emptyset$ for every $F \in \mathcal{F}$; \label{item:centered_coloring:base_case:excluding_a_star:hitting}
        \item for every component $C$ of $G-S$, $N_G(V(C))$ intersects at most three components of $G-V(C)$; \label{item:centered_coloring:base_case:excluding_a_star:connectivity}
        \item $P_i \subseteq R_i$ for every $i \in [\ell]$; \label{item:centered_coloring:base_case:excluding_a_star:Pi_included_Ri}
        \item $\cen_q(G,\phi,R_i) \leq c \cdot d$ for every $i \in [\ell]$; and \label{item:centered_coloring:base_case:excluding_a_star:chi(Ri)_small}
        \item $|R_i \setminus P_i| \leq c \cdot d$ for every $i \in [\ell]$. \label{item:centered_coloring:base_case:excluding_a_star:Ri_minus_Pi_small}
    \end{enumerateOurAlph}
\end{lemma}

\begin{proof}
    We proceed by induction on $|V(G)\setminus U|$.
    If every member of $\mathcal{F}$ intersects $U$, then the result holds for $S=P_0=U$ and $\ell=0$. 
    Note that item~\ref{item:centered_coloring:base_case:excluding_a_star:connectivity} holds as $G[U]$ is connected. 
    Now suppose that $\mathcal{F}\vert_{G-U} \neq \emptyset$, and so, $V(G) \setminus U \neq \emptyset$.
    Let $\mathcal{F}'$ be the family of all the connected subgraphs $H$ of $G-U$ such that $N_G(V(H)) \cap U \neq \emptyset$ and
    there exists $F \in \mathcal{F}$ such that $F \subseteq H$.
    Since $G$ is connected, 
    $\mathcal{F}'$ is nonempty. 
    See an illustration of the proof in~\cref{fig:excluding-star-cen}.
    
    There are no $d+1$ pairwise disjoint members of $\mathcal{F}'$ as otherwise, if $H_1, \dots, H_{d+1}$ are $d+1$ pairwise disjoint members of $\mathcal{F}'$,
    then $(V(H_1), \dots, V(H_d), U \cup V(H_{d+1}))$ is an $\mathcal{F}$-rich model of $F_{2,d}$ in $G$, which is a contradiction.
    Hence, since $\phi$ is a $(q,c)$-good coloring of $G$,
    there exists $Z \subseteq V(G)$ such that
    \begin{enumerate}[label=(qc\arabic*')]
        \item $Z \cap V(F) \neq \emptyset$ for every $F \in \calF'$; \label{item:qc1'}
        \item for every component $C$ of $G-Z$, $N_{G}(V(C))$ intersects at most two components of $G-V(C)$; and \label{item:qc1''}
        \item for every $P \subseteq Z$, there exists $R \subseteq Z$ with \label{item:qcgood:extend}
            \begin{enumerate}
                \item $P \subseteq R$, \label{item:qcgood:extend:i}
                \item $\cen_q(G,\phi,R) \leq c \cdot d$, and\label{item:qcgood:extend:ii}
                \item $|R \setminus P| \leq c \cdot d$.\label{item:qcgood:extend:iii}
            \end{enumerate}
    \end{enumerate}

    \begin{figure}[tp]
    \centering
    \includegraphics{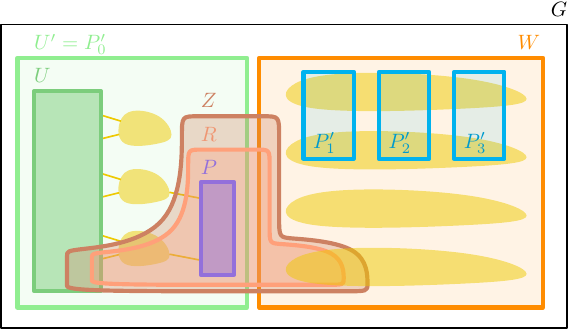}
    \caption{
        The objects defined and used in the proof of~\Cref{lemma:centered_coloring:base_case:excluding_a_star}.
    }
    \label{fig:excluding-star-cen}
\end{figure}

    Let $P$ be an inclusion-wise minimal subset of $Z \setminus U$ satisfying the following properties:
    \begin{enumerate}[label=(p\arabic*)]
        \item $P \cap V(F) \neq \emptyset$ for every $F \in \mathcal{F}'$ and \label{item:proof:centered_coloring:base_case:excluding_a_star:P_is_hitting}
        \item for every component $C$ of $G-U-P$ such that $\mathcal{F}\vert_C = \emptyset$, $N_{G}(V(C))$ intersects at most three components of $G-V(C)$. \label{item:proof:centered_coloring:base_case:excluding_a_star:P_has_a_nice_connectivity_property}
    \end{enumerate}

    To justify that $P$ is well-defined, we show that $Z - U$ satisfies~\ref{item:proof:centered_coloring:base_case:excluding_a_star:P_is_hitting} and~\ref{item:proof:centered_coloring:base_case:excluding_a_star:P_has_a_nice_connectivity_property}.
    Item~\ref{item:proof:centered_coloring:base_case:excluding_a_star:P_is_hitting} for $P=Z-U$ follows directly from~\ref{item:qc1'}.
    Let $C$ be a component of $G - U - Z$.
    There is a component $C'$ of $G - Z$ such that $C \subset C'$.
    By~\ref{item:qc1''}, $N_G(V(C'))$ intersects at most two components of $G - V(C')$. 
    These components are connected subgraphs of $G-V(C)$. 
    Moreover, $N_G(V(C)) \subseteq N_G(V(C'))\cup U$ and $G[U]$ is connected by the assumption of the lemma.
    Altogether, $N_G(V(C))$ is contained in at most three connected subgraphs of $G-V(C)$, which are contained in at most three components of $G-V(C)$.
    This shows~\ref{item:proof:centered_coloring:base_case:excluding_a_star:P_has_a_nice_connectivity_property} for $P = Z- U$ as claimed. 
    
    Let $W$ be the union of the vertex sets of all the components $C$ of $G-U-P$ such that $\mathcal{F}\vert_C \neq \emptyset$. 
    By~\ref{item:proof:centered_coloring:base_case:excluding_a_star:P_is_hitting}, 
    \begin{equation}
        N_G(W) \subseteq P.
        \label{eq:neighborhood-of-W-in-P}
    \end{equation}
    We claim that for every $u \in N_G(W)$, there is a \stpath{\{u\}}{U} in $G$ disjoint from $W$.
    By the minimality of $P$, either \ref{item:proof:centered_coloring:base_case:excluding_a_star:P_is_hitting} or \ref{item:proof:centered_coloring:base_case:excluding_a_star:P_has_a_nice_connectivity_property} is violated for $P - \{u\}$.
    First, we show that $P - \{u\}$ satisfies~\ref{item:proof:centered_coloring:base_case:excluding_a_star:P_has_a_nice_connectivity_property}.
    Let $C$ be a component of $G - U - (P - \{u\})$.
    If $u \notin V(C)$, then $C$ is a component of $G - U - P$, and so, by~\ref{item:proof:centered_coloring:base_case:excluding_a_star:P_has_a_nice_connectivity_property} for $P$, if $\calF\vert_C = \emptyset$, then $N_G(V(C))$ intersects at most three components of $G - V(C)$.
    Otherwise, we have $u \in V(C)$. 
    Since $u \in N_G(W)$, it follows that there is a component $C'$ of $G-U-P$, more precisely of $G[W]$, with $\mathcal{F}\vert_{C'} \neq \emptyset$ and $C' \subseteq C$.
    In particular, $\mathcal{F}\vert_C \neq \emptyset$.
    This yields~\ref{item:proof:centered_coloring:base_case:excluding_a_star:P_has_a_nice_connectivity_property} for $P - \{u\}$.
    We conclude that~\ref{item:proof:centered_coloring:base_case:excluding_a_star:P_is_hitting} is violated for $P - \{u\}$.
    That is, there exists $F \in \mathcal{F}'$ disjoint from $P \setminus \{u\}$.
    Then, $u \in V(F)$.
    Also, since $F \in \calF'$, $N_G(V(F)) \cap U \neq \emptyset$.
    Therefore, since $F$ is connected, there is a \stpath{\{u\}}{N_G(U)} $Q$ in $F$.
    By~\eqref{eq:neighborhood-of-W-in-P} and the fact that $V(F) \cap P = \{u\}$, the path $Q$ is disjoint from $W$.
    Finally, we can extend this path by one vertex in $U$ to be a \stpath{\{u\}}{U} in $G$.
    This is the desired path.

    Let $U' = V(G) \setminus W$.
    We claim that $G[U']$ is connected. 
    Consider a component $C$ of $G[U']$.
    Since $G$ is connected, $C$ contains a vertex $u$ of $N_G(W)$.
    By the previous claim, there exists a \stpath{\{u\}}{U} disjoint from $W$, i.e.\ contained in $G[U']$. 
    Thus, every component of $G[U']$ contains a vertex of $U$. 
    Since $G[U]$ is connected, we conclude that $G[U']$ is connected.

    Since $P$ hits all the members of $\mathcal{F'}$ (by~\ref{item:proof:centered_coloring:base_case:excluding_a_star:P_is_hitting}) and since $\mathcal{F}'$ is nonempty, $P$ is nonempty. 
    Therefore, $|V(G) \setminus U'| < |V(G) \setminus U|$ because $P\subseteq U' \setminus U$.
    Hence, by the induction hypothesis,
    there exists $S' \subseteq V(G)$, a path partition $(P'_0, \dots, P'_{\ell'})$ of $(G,S')$, and sets $R'_1, \dots, R'_{\ell'} \subseteq V(G)$ such that
    \begin{enumerateOurAlphPrim}
        \item $P'_0 = U'$; \label{item:centered_coloring:base_case:excluding_a_star:P0=U'}
        \item $S' \cap V(F) \neq \emptyset$ for every $F \in \mathcal{F}$; \label{item:centered_coloring:base_case:excluding_a_star:hitting'}
        \item for every component $C$ of $G-S'$, $N_G(V(C))$ intersects at most three components of $G-V(C)$;\label{item:centered_coloring:base_case:excluding_a_star:connectivity'}
        \item $P'_i \subseteq R'_i$ for every $i \in [\ell']$; \label{item:centered_coloring:base_case:excluding_a_star:Pi_included_Ri'}
        \item $\cen_q(G,\phi,R'_i) \leq c \cdot d$ for every $i \in [\ell']$; and \label{item:centered_coloring:base_case:excluding_a_star:chi(Ri)_small'}
        \item $|R'_i \setminus P'_i| \leq c \cdot d$ for every $i \in [\ell']$. \label{item:centered_coloring:base_case:excluding_a_star:Ri_minus_Pi_small'}
    \end{enumerateOurAlphPrim}
    Let $S = U \cup P \cup \textstyle\bigcup_{i \in [\ell']} P_i'$.
    Calling \ref{item:qcgood:extend}, 
    we fix $R \subseteq Z$ with $P \subseteq R$, $\cen_q(G,\phi,R) \leq c \cdot d$,
    and $|R \setminus P| \leq c \cdot d$.
    We set $\ell = \ell'+1$, and
    \[
        \begin{array}{l l}
            P_0 = U, & \\
            P_1  = P, \quad
            &R_1  = R, \\
             P_i = P'_{i-1}, \quad
            &R_i  = R'_{i-1} \quad \textrm{for every $i \in \{2, \dots, \ell\}$}.
        \end{array}
    \]
    To conclude the proof, we show that~\ref{item:centered_coloring:base_case:excluding_a_star:P0=U}--\ref{item:centered_coloring:base_case:excluding_a_star:Ri_minus_Pi_small} hold.
    By construction, 
    \ref{item:centered_coloring:base_case:excluding_a_star:P0=U} holds.
    
    Let $F \in \mathcal{F}$. By~\ref{item:centered_coloring:base_case:excluding_a_star:hitting'}, $S' \cap V(F) \neq \emptyset$, and so, there exists $i \in \{0,\dots,\ell'\}$ such that $P_i' \cap V(F)$.
    If $i \neq 0$, then $P_i' \subset S$, and thus, $S \cap V(F) \neq \emptyset$.
    If $i = 0$, then $U' \cap V(F) \neq \emptyset$.
    Recall that $U' = V(G) - W = U \cup P \cup \bigcup \{V(C) \mid C \in C\}$ where $\calC$ is the set of the components of $G - U - P$ such that $\calF\vert_C = \emptyset$.
    Therefore, $(U \cup P) \cap V(F) \neq \emptyset$, and so, $S \cap V(F) \neq \emptyset$.
    This concludes~\ref{item:centered_coloring:base_case:excluding_a_star:hitting}.

    Let $C$ be a component of $G-S$.
    If $V(C) \cap U' = \emptyset$, then $C$ is a component of $G-S'$ and so by~\ref{item:centered_coloring:base_case:excluding_a_star:connectivity'},
    $N_G(V(C))$ intersects at most three components of $G-V(C)$.
    Otherwise, $C$ is a component of $G-U-P$ with $\calF\vert_C = \emptyset$.
    Hence, by \ref{item:proof:centered_coloring:base_case:excluding_a_star:P_has_a_nice_connectivity_property},
    $N_{G}(V(C))$ intersects at most three components of $G-V(C)$.
    This proves \ref{item:centered_coloring:base_case:excluding_a_star:connectivity}.
    Items~\ref{item:centered_coloring:base_case:excluding_a_star:Pi_included_Ri}, \ref{item:centered_coloring:base_case:excluding_a_star:chi(Ri)_small}, and \ref{item:centered_coloring:base_case:excluding_a_star:Ri_minus_Pi_small} hold respectively by~\ref{item:centered_coloring:base_case:excluding_a_star:Pi_included_Ri'} and~\ref{item:qcgood:extend}.\ref{item:qcgood:extend:i}; \ref{item:centered_coloring:base_case:excluding_a_star:chi(Ri)_small'} and~\ref{item:qcgood:extend}.\ref{item:qcgood:extend:ii}; and \ref{item:centered_coloring:base_case:excluding_a_star:Ri_minus_Pi_small'} and \ref{item:qcgood:extend}.\ref{item:qcgood:extend:iii}.
    This ends the proof of the lemma.
\end{proof}

\begin{lemma}\label{lemma:centered_coloring:base_case:excluding_a_forest}
    Let $q,$ $c$, and $d$ be positive integers with $q \geq 2$ and let $h$ be a nonnegative integer,
    let $G$ be a graph,
    let $\phi$ be a $(q,c)$-good coloring of $G$,
    and let $\mathcal{F}$ be a family of connected subgraphs of $G$ such that $G$ has no $\mathcal{F}$-rich model of $F_{h+1,d}$.
    There exists
    a set $\Lambda$ of at most $12 c(hd+\binom{h}{2}) \cdot q\log q$ colors,
    a set $S \subseteq V(G)$, a coloring $\lambda \colon S \to \Lambda$,
    and a family $(\Lambda_C \mid \text{$C$ component of $G-S$})$ of subsets of $\Lambda$, each of size at most $20c (hd + \binom{h}{2}) \cdot \log q$ such that
    \begin{enumerateOurAlph}
        \item $S \cap V(F) \neq \emptyset$ for every $F \in \mathcal{F}$; \label{item:centered_coloring:base_case:excluding_a_forest:hitting}
        \item for every component $C$ of $G-S$, $N_G(V(C))$ intersects at most $8c (hd + \binom{h}{2}) \cdot \log q$ components of $G-V(C)$; and \label{item:centered_coloring:base_case:excluding_a_forest:connectivity}
        \item \label{item:centered_coloring:base_case:excluding_a_forest:center_among_log_colors}
            for every connected subgraph $H$ of $G$ intersecting $S$,
            \begin{enumerate}
                \item $|\phi(V(H))| > q$, or \label{item:c-i}
                \item $|\lambda(V(H) \cap S)| > q$, or\label{item:c-ii}
                \item there is a $(\phi \times \lambda)$-center $u$ of $V(H) \cap S$ such that for every component $C$ of $G - S$ intersecting $V(H)$, we have $\lambda(u) \in \Lambda_{C}$.\label{item:c-iii}
            \end{enumerate}
    \end{enumerateOurAlph}
\end{lemma}

\begin{proof}
    We proceed by induction on $h+|V(G)|$.
    When $h=0$, $F_{h+1,d}=K_1$ and so $\mathcal{F}$ is empty.
    Hence, setting $S = \emptyset$ and $\Lambda = \emptyset$ gives the assertion.
    Next suppose that $h>0$.

    First assume that $G$ is not connected.
    Let $\mathcal{C}$ be the family of all the components of $G$.
    By the induction hypothesis, and by possibly relabeling the obtained color sets,
    there exists
    a set $\Lambda$ of at most $20 c(hd+\binom{h}{2}) \cdot q\log q$ colors,
    such that for every $C \in \mathcal{C}$,
    there exists a set $S_C \subseteq V(C)$,
    a coloring $\lambda_C \colon S_C \to \Lambda$,
    and a family $(\Lambda_{C'} \mid \text{$C'$ component of $C-S_C$})$ of subsets of $\Lambda$, each of size at most $12c (hd + \binom{h}{2}) \cdot \log q$, 
    such that
    \begin{enumerateOurAlphPrim}
        \item $S_C \cap V(F) \neq \emptyset$ for every $F \in \mathcal{F}\vert_C$; \label{item:centered_coloring:base_case:excluding_a_forest:hitting:induction_call1}
        \item for every component $C'$ of $C-S_C$, $N_C(V(C'))$ intersects at most $8c (hd + \binom{h}{2}) \cdot \log q$ components of $C-V(C')$; and \label{item:centered_coloring:base_case:excluding_a_forest:connectivity:induction_call1}
        \item \label{item:centered_coloring:base_case:excluding_a_forest:center_among_log_colors:induction_call1}
            for every connected subgraph $H$ of $C$ intersecting $S_C$, 
            \begin{enumerate}
                \item $|\phi(V(H))| > q$, or
                \item $|\lambda_C(V(H) \cap S_C)| > q$, or
                \item there is a $(\phi \times \lambda_C)$-center $u$ of $V(H) \cap S_C$ such that for every component $C'$ of $C-S_C$ intersecting $V(H)$, we have $\lambda_C(u) \in \Lambda_{C'}$.
            \end{enumerate}
    \end{enumerateOurAlphPrim}

    Let
    \[
        S = \bigcup_{C \in \mathcal{C}} S_C,
    \]
    and let $\lambda \colon S \to \Lambda$ such that $\lambda \vert_{S_C} = \lambda_C$ for every $C \in \mathcal{C}$.
    Then, $S$, $\lambda$, and $(\Lambda_{C'} \mid \text{$C'$ component of $G-S$})$ satisfy the assertion of the lemma.
    Now suppose that $G$ is connected.
    
    Let $\mathcal{F}'$ be the family of all the connected subgraphs of $G$ containing an $\mathcal{F}$-rich model of $F_{h,d+1}$.

    We claim that there is no $\mathcal{F}'$-rich model of $F_{2,d}$ in $G$.
    Suppose for a contradiction that there is an $\mathcal{F}'$-rich model
    $\big(B_x \mid x \in V(F_{2,d})\big)$ of $F_{2,d}$ in $G$.
    We construct an $\mathcal{F}$-rich model of $F_{h+1,d}$ in $G$.
    Let $z$ the center of the star $F_{2,d}$
    and let $x_1, \dots, x_d$ be the leaves of this star.
    Let $i \in [d]$ and let $u_i \in B_{x_i}$ be a neighbor of some vertex of $B_z$ in $G$.
    Since $\big(B_x \mid x \in V(F_{2,d})\big)$ is an $\mathcal{F}'$-rich model of $F_{2,d}$,
    the graph $G[B_{x_i}]$ has an $\mathcal{F}$-rich model
    $\big(B_{i,x} \mid x \in V(F_{h,d+1})\big)$ of $F_{h,d+1}$.
    Since $G[B_{x_i}]$ is connected, we can assume that $u_i \in B_{i,y_i}$ for some $y_i \in V(F_{h,d+1})$.
    By \Cref{lemma:rooting_a_model_of_Fhd},
    there is an $\mathcal{F}$-rich model $\mathcal{M}_i$ of $F_{h,d}$ in $G[B_{x_i}]$ 
    whose branch set of the root contains a neighbor $u_i$ of $B_z$.
    It follows that the union of the models $\mathcal{M}_1, \dots, \mathcal{M}_d$,
    together with $B_z$ as the branch set of the root,
    yields an $\mathcal{F}$-rich model of $F_{h+1,d}$ in $G$, which is a contradiction.
    This proves that there is no $\mathcal{F}'$-rich model of $F_{2,d}$ in $G$.

    Hence, by \Cref{lemma:centered_coloring:base_case:excluding_a_star} applied to $U$, an arbitrary singleton of a vertex in $G$,
    there exists $S' \subseteq V(G)$, a path partition $(P_0, \dots, P_\ell)$ of $(G,S')$, and sets $R_1, \dots, R_\ell \subseteq V(G)$ such that
    \begin{enumerate}[label=\ref{lemma:centered_coloring:base_case:excluding_a_star}.(\alph*)]
        \item $P_0 = U$; \label{item:centered_coloring:base_case:excluding_a_star:P0=U:call}
        \item $S' \cap V(F) \neq \emptyset$ for every $F \in \mathcal{F}'$; 
            \label{item:centered_coloring:base_case:excluding_a_star:hitting:call}
        \item for every component $C$ of $G-S'$, $N_G(V(C))$ intersects at most three components of $G-V(C)$; 
            \label{item:centered_coloring:base_case:excluding_a_star:connectivity:call}
        \item $P_i \subseteq R_i$ for every $i \in [\ell]$; \label{item:centered_coloring:base_case:excluding_a_star:Pi_included_Ri:call}
        \item $\cen_q(G,\phi,R_i) \leq c \cdot d$ for every $i \in [\ell]$; and \label{item:centered_coloring:base_case:excluding_a_star:chi(Ri)_small:call}
        \item $|R_i \setminus P_i| \leq c \cdot d$. \label{item:centered_coloring:base_case:excluding_a_star:Ri_minus_Pi_small:call}
    \end{enumerate}

    We set $R_0=P_0=U$. 
    Note that \ref{item:centered_coloring:base_case:excluding_a_star:Pi_included_Ri:call}--\ref{item:centered_coloring:base_case:excluding_a_star:Ri_minus_Pi_small:call} also hold for $i=0$.

    Let $i \in \{0, \dots, \ell\}$.
    Let $\eta_i \colon R_i \to [c d]$ witness \ref{item:centered_coloring:base_case:excluding_a_star:chi(Ri)_small:call}, and 
    let $\psi_i\colon R_i \to [2c d]$ be such that $\psi_i\vert_{P_i} = \eta_i\vert_{P_i}$, $\psi_i\vert_{R_i \setminus P_i}$ is injective, 
    and $\psi_i(u) \neq \psi_i(v)$ for every $u \in P_i$ and $v \in R_i \setminus P_i$. 
    Such a coloring exists by \ref{item:centered_coloring:base_case:excluding_a_star:Ri_minus_Pi_small:call}.

    Let $\mathcal{C}$ be the family of all the components of $G-S'$.
    For every $C \in \mathcal{C}$, fix an index $i(C) \in [\ell]$ such that $N_G(V(C)) \subseteq P_{i(C)-1} \cup P_{i(C)}$.
    Next, for every $i \in [\ell]$, let
    \[
        W_i = P_{i-1} \cup P_i \cup \bigcup_{C \in \mathcal{C}, i(C)=i} V(C).
    \]
    Note that $(W_1, \dots, W_\ell)$ is a path decomposition of $G$.
    See the first part of~\Cref{fig:construction-R'-cen}.
    Also
    \begin{equation}
    W_{i-1} \cap W_{i+1} = \emptyset\qquad\textrm{for every $i \in [\ell-1]$}.
    \label{eq:Ws-disjoint}
    \end{equation}
    For every $u \in V(G)$, let 
    \[
        \mu(u) = 
        \begin{cases}
            i \mod (q+1) &\textrm{if $u \in S'$, where $i \in \{0, \dots, \ell\}$ is such that $u \in P_i$,} \\
            0 &\textrm{otherwise.}
        \end{cases}
    \]

    Let $s = \lceil \log (3q) \rceil$.
    For every $r \in \{0, \dots, s\}$, let 
    \[
        I_r = \set{i \in \{0, \dots, \ell\} \mid i \equiv 0 \bmod 2^r}.
    \]
    For each $r \in \{0,\dots,s\}$ and $i \in I_r$, we will build a set $R_i' \subset R_i$ such that $(R_i'\mid i\in I_r)$ are pairwise disjoint. 
    The construction is recursive starting from $r=s$ and going down to $r=0$. 
    When $R_i'$ are defined for all $i\in I_r$, we set 
    \[
    S_r = \bigcup \set{R'_i \mid i \in I_r}.
    \]
    We maintain the invariants
    \begin{align} 
      P_i\subseteq S_r&\text{ for every } i\in I_r, \label{eq:Pi_in_Sr} \\
      P_{i} \setminus S_{r+1} \subseteq R'_{i}&\text{ for every } i\in I_r, \label{eq:Pi_minus_Sr+1_in_R'i}
    \end{align}
    with the convention $S_{s+1} = \emptyset$.
    See an illustration of this construction in~\Cref{fig:construction-R'-cen}.

    We start with $r=s$.
    For every $i\in I_s$, let 
    \begin{align*}
        R'_i &= R_i \cap \bigcup\set{W_j\mid j \in \{0, \dots, \ell\} \cap \set{i-q+1,\ldots,i+q}}.
    \end{align*}
    When $a,b\in I_s$ and $a<b$, 
    we have $b-a\geq 2^s\geq 3q\geq2q+2$ since $q\geq 2$. 
    Thus, $(b - q)-(a+q+1)\geq 1$ and therefore, by~\eqref{eq:Ws-disjoint}, 
    the sets $R'_i$ for $i \in I_s$ are pairwise disjoint.
    Note also that $P_i\subseteq R_i \cap W_i\subseteq R'_i$ for every $i \in I_s$,
    and so \eqref{eq:Pi_in_Sr} holds.

    Next, let $r<s$ and suppose that $R'_i$ for $i \in I_{r+1}$, and $S_{r+1}$ has been constructed.
    For every $i \in I_r \setminus I_{r+1}$, let
    \begin{align*}
        R'_i &= \left(R_i\setminus S_{r+1}\right) \cap \bigcup\set{W_j\mid j\in\{0, \dots, \ell\}\cap\set{i-2^r+1,\ldots,i+2^r}}.
    \end{align*}
    When $a,b\in I_r \setminus I_{r+1}$ and $a<b$, we have 
    $b-a\geq 2^{r+1}$.
    Moreover, $R'_a \subseteq \bigcup\{W_j \mid j \in \{0, \dots, \ell\}, j \leq a+2^r\}$ and $R'_b \subseteq \bigcup\{W_j \mid j \in \{0, \dots, \ell\}, j \geq b-2^r+1\}$.
    Hence, by \eqref{eq:Ws-disjoint} for $r+1$,
    $R'_a \cap R'_b \subseteq W_{a+2^r} \cap W_{a+2^r+1} \subseteq P_{a+2^r} \subseteq S_{r+1}$, where the last inclusion follows from \eqref{eq:Pi_in_Sr} and the fact that $a+2^r \in I_{r+1}$.
    As $R'_a$ and $R'_b$ are both disjoint from $S_{r+1}$, we deduce that $R'_a \cap R'_b = \emptyset$.
    This proves that the sets $R'_i$ for $i \in I_r$ are pairwise disjoint.
    For every $i \in I_r \setminus I_{r+1}$,
    $P_i \setminus S_{r+1} \subseteq (R_i\setminus S_{r+1}) \cap W_i\subseteq R'_i$ so \eqref{eq:Pi_minus_Sr+1_in_R'i} holds.
    Moreover, $P_i \setminus S_{r+1} \subseteq S_r$.
    Since $S_{r+1} \subseteq S_r$, \eqref{eq:Pi_in_Sr} holds.

    \begin{figure}[tp]
        \centering
        \includegraphics{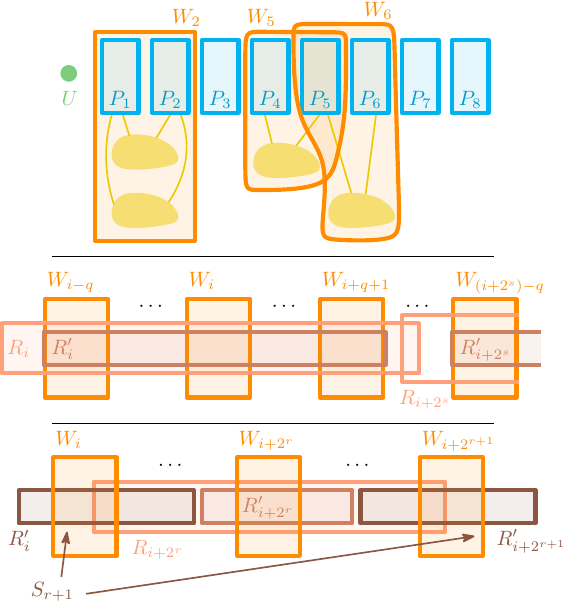}
        \caption{
            An illustration of objects in the proof of~\Cref{lemma:centered_coloring:base_case:excluding_a_forest}.
            In the first part, we show a few examples of what $W_i$ looks like.
            In the second part, we show the idea behind the construction of the sets $R_i'$ in the base case.
            Note that the sets $R_i$ are not necessarily disjoint from each other, while the sets $R_i'$ are.
            In the last part, we show the idea for constructing $R_{i+2^r}'$ given that $R_i'$ and $R_{i+2^{r+1}}'$ are constructed.
        }
        \label{fig:construction-R'-cen}
    \end{figure}

    Next, we define a coloring $\lambda_0$ of $S_0$ as follows.
    For convenience, let $I_{s+1} = \emptyset$.
    For every $r \in \{0, \dots, s\}$, for every $i \in I_r \setminus I_{r+1}$, and for every $u \in R_i'$, we set
    \[
        \lambda_0(u) = (r, \psi_i(u)).
    \]
    Now that $S_0$ and $\lambda_0$ are defined,
    we decompose $G-S_0$.
    Since $S'$ intersects every member of $\mathcal{F}'$ by~\ref{item:centered_coloring:base_case:excluding_a_star:hitting:call},
    and since $S' \subseteq S_0$,
    $G-S_0$ has no $\mathcal{F}$-rich model of $F_{h,d+1}$.
    Hence, by the induction hypothesis,
    there exists 
    a set $\tilde{\Lambda}$ of at most $20 c((h-1)(d+1)+\binom{h-1}{2}) \cdot q\log q$ colors,
    a set $\tilde{S} \subseteq V(G-S_0)$, 
    a coloring $\tilde{\lambda} \colon \tilde{S} \to \tilde{\Lambda}$,
    and a family $(\tilde{\Lambda}_C \mid \text{$C$ component of $G-S_0$})$ of subsets of $\Tilde{\Lambda}$, each of size at most $12c ((h-1)(d+1) + \binom{h-1}{2}) \cdot \log q$, 
    such that
    \begin{enumerate}[label=(\~{\alph*})]
        \item $\tilde{S} \cap V(F) \neq \emptyset$ for every $F \in \mathcal{F}\vert_{G-S_0}$; \label{item:centered_coloring:base_case:excluding_a_forest:hitting'}
        \item for every component $C$ of $G-S_0-\tilde{S}$, $N_{G-S_0}(V(C))$ intersects at most $8c ((h-1)(d+1) + \binom{h-1}{2}) \cdot \log q$ components of $G-S_0-V(C)$; and \label{item:centered_coloring:base_case:excluding_a_forest:connectivity'}
        \item \label{item:centered_coloring:base_case:excluding_a_forest:center_among_log_colors'}
            for every connected subgraph $H$ of $G - S_0$ intersecting $\tilde{S}$, 
            \begin{enumerate}
                \item $|\phi(V(H))| > q$, or \label{item:c-tilde:i}
                \item $|\tilde{\lambda}(V(H) \cap \tilde{S})| > q$, or\label{item:c-tilde:ii}
                \item there is a $(\phi \times \tilde{\lambda})$-center $u$ of $V(H) \cap \tilde{S}$ such that
                for every component $C$ of $G-S_0$ intersecting $V(H)$, we have
                $\tilde{\lambda}(u) \in \tilde{\Lambda}_{C}$.\label{item:c-tilde:iii}
            \end{enumerate}
    \end{enumerate}
    Up to relabeling the elements of $\tilde{\Lambda}$, we can assume that $\tilde{\Lambda}$ is disjoint from the image of $\mu \times \lambda_0$.

    \begin{figure}[tp]
        \centering
        \includegraphics{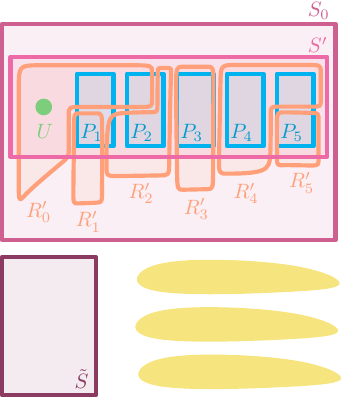}
        \caption{
            An illustration of objects in the proof of~\Cref{lemma:centered_coloring:base_case:excluding_a_forest}.
            Remember that $S_0 = \bigcup_{i \in \{0,\dots,\ell\}} R_i'$.
        }
        \label{fig:base-case-cen-Ss}
    \end{figure}

    Finally, let 
    \[
        S = S_0 \cup \tilde{S}.
    \]
    See~\Cref{fig:base-case-cen-Ss}.
    Also, for every $u \in S$, let
    \[
        \lambda(u) = 
        \begin{cases}
            (\mu(u),\lambda_0(u)) &\textrm{if $u \in S_0$,} \\
            \tilde{\lambda}(u)    &\textrm{if $u \in \tilde{S}$.}
        \end{cases}
    \]
    Let $\Lambda$ be the image of $\lambda$.
    First, observe that $s\leq 1 + \log(3q) \leq 1 + \log(4q) = 3 + \log q$ since $q\geq 2$. 
    Therefore,
    \begin{align*}
        |\Lambda| 
        &= |\lambda(S_0)| + \big|\lambda\big(\Tilde{S}\big)\big| \\
        &\leq |\mu(S_0)|\cdot |\lambda_0(S_0)| + \big|\Tilde{\Lambda}\big| \\
        &\leq (q+1) \cdot \big((s+1) \cdot 2cd\big) + 20c\left(\textstyle(h-1)(d+1)+\binom{h-1}{2}\right) \cdot q \log q \\
        &\leq (q+1) \cdot \big((4+\log q) \cdot 2c d\big) + 20c\left(\textstyle(h-1)(d+1)+\binom{h-1}{2}\right) \cdot q \log q \\
        &\leq 2q \cdot (5 \log q \cdot 2c  d) + 20c\left(\textstyle(h-1)(d+1)+\binom{h-1}{2}\right) \cdot q \log q \\
        &\leq 20cd \cdot q \log q + 20c\left(\textstyle(h-1)(d+1)+\binom{h-1}{2}\right) \cdot q \log q \\
        &\leq 20c\left(\textstyle hd + h-1+\binom{h-1}{2}\right) \cdot q \log q\\
        &\leq 20c\left(\textstyle hd+\binom{h}{2}\right) \cdot q \log q.
    \end{align*}

    For all $i\in\set{0,\ldots,\ell}$ and $r \in \{0, \dots, s\}$, let 
    \begin{align*}
        \ext{i}{r}{-} &= \max (\{j \in I_r \mid j < i\} \cup \{0\}) \text{ and}\\
        \ext{i}{r}{+} &= \min (\{j \in I_r \mid j \geq i\} \cup \{\ell\}).
    \end{align*}
    Moreover, for every $i\in\set{0,\ldots,\ell}$, let
    \[
        \Trace(i)=\set{\ext{i}{r}{\epsilon}\mid r\in\set{0,\ldots,s},\ \epsilon\in\set{-,+}}.
    \]

    \begin{claim}\label{claim:cen_excluding_a_forest:intermediate}
        Let $H$ be a connected subgraph of $G$ such that $S_0 \cap V(H) \neq \emptyset$, and let $r = \max \{t \in \{0, \dots, s\} \mid V(H) \cap S_t \neq \emptyset\}$.
        Then, either
        \begin{enumerate}[label=\normalfont (\arabic*)]
            \item $|\{i \in \{0, \dots, \ell\} \mid V(H) \cap P_i \neq \emptyset\}| > q$, or
            \item there is a unique $j \in I_r$ such that $V(H) \cap R'_j \neq \emptyset$ and for every $i \in [\ell]$, if $V(H) \cap W_i \neq \emptyset$, then $j \in \Trace(i)$.
        \end{enumerate}
    \end{claim}

    \begin{proofclaim}
        Suppose that $|\{i \in \{0, \dots, \ell\} \mid V(H) \cap P_i \neq \emptyset\}| \leq q$.
        Recall that $S_r = \bigcup_{j \in I_r} R'_j$.
        We show that there is unique $j \in I_r$ such that $V(H) \cap R_j' \neq \emptyset$.
        To this end, let $j_1,j_2 \in I_r$ be such that $V(H)$ intersects both $R'_{j_1}$ and $R'_{j_2}$, and $j_1 \leq j_2$.
        Suppose to the contrary that $j_1 < j_2$.

        First, consider the case of $r<s$.
        By the maximality of $r$, we have $j_1,j_2 \in I_r \setminus I_{r+1}$,
        which implies $j_2 - j_1 \geq 2^{r+1}$,
        and $j_1 + 2^r \in I_{r+1}$.
        However, since $(W_1, \dots, W_\ell)$ is a path decomposition of $G$,
        and since $V(H)$ intersects both $R'_{j_1} \subseteq \bigcup_{1 \leq i \leq j_1+2^r} W_i$ and $R'_{j_2} \subseteq \bigcup_{j_2-2^r+1 \leq i \leq \ell} W_i$,
        we deduce by \eqref{eq:Ws-disjoint} that $V(H)$ intersects
        \begin{align*}
            \left(\bigcup_{1 \leq i \leq j_1+2^r} W_i\right) \cap \left( \bigcup_{j_2-2^r+1 \leq i \leq \ell} W_i \right) &\subseteq \left(\bigcup_{1 \leq i \leq j_1+2^r} W_i\right) \cap \left( \bigcup_{j_1+2^r+1 \leq i \leq \ell} W_i \right)\\
            &= W_{j_1+2^r} \cap W_{j_1+2^r+1} = P_{j_1+2^r}.
        \end{align*}
        But since $j_1+2^r \in I_{r+1}$, we have $P_{j_1 + 2^r} \subseteq S_{r+1}$ by \eqref{eq:Pi_in_Sr},
        and so $V(H)$ intersects $S_{r+1}$,
        which contradicts the maximality of $r$.

        Now suppose that $r=s$.
        Since $j_1,j_2 \in I_s$, we have $j_2-j_1 \geq 2^s \geq 3q$.
        Again, since $(W_1, \dots, W_\ell)$ is a path decomposition of $G$,
        and since $V(H)$ intersects both $R'_{j_1} \subseteq \bigcup_{1 \leq i \leq j_1+q} W_i$ and $R'_{j_2} \subseteq \bigcup_{j_2-q+1 \leq i \leq \ell} W_i$,
        we deduce that $V(H)$ intersects
        $W_{i}$ for every integer $i$ with $j_1+q \leq i \leq j_2-q+1$.
        By the definition of $(W_1,\dots, W_\ell)$, and because
        $(P_0, \dots, P_\ell)$ is a path partition of $(G,S')$,
        it follows that $V(H)$ intersects $P_i$ for every $i$ with $j_1+q \leq i \leq j_2-q$,
        and so $|\{i \in \{0, \dots, \ell\} \mid V(H) \cap P_i \neq \emptyset\}| \geq (j_2-q)-(j_1+q)+1 = (j_2-j_1) - 2q+1 \geq 3q - 2q+1 = q+1 > q$,
        a contradiction.
        This concludes the proof of $j_1 = j_2$.

        In particular, we may fix the unique index $j \in I_r$ such that $V(H) \cap R'_j \neq  \emptyset$.
        Next, we show that 
        for every $i \in [\ell]$ such that $V(H) \cap W_i \neq \emptyset$, we have 
        $j \in \{\ext{i}{r}{-}, \ext{i}{r}{+}\}$.
        Since $\{\ext{i}{r}{-}, \ext{i}{r}{+}\} \subseteq \Trace(i)$ for every $i \in [\ell]$, this will prove the lemma.
        Let $i \in [\ell]$ such that $V(H) \cap W_i \neq \emptyset$,
        and suppose for contradiction that $j \not\in \{\ext{i}{r}{-}, \ext{i}{r}{+}\}$.
        Note that since $(P_0, \dots, P_\ell)$ is a path partition of $(G,S')$,
        and by the definitions of $(W_1, \dots, W_\ell)$, $\ext{i}{r}{-}$, and $\ext{i}{r}{+}$,
        we have $P_{\ext{i}{r}{-}} \cup P_{\ext{i}{r}{+}}$ separates $W_i$ and
        $\bigcup \{W_{i'} \mid i' \in [\ell], \text{ $i' \leq \ext{i}{r}{-}$ or $i' > \ext{i}{r}{+}$}\}$.
        However, since $j \not\in \{\ext{i}{r}{-}, \ext{i}{r}{+}\}$,
        we have $R'_j \subseteq \bigcup \{W_{i'} \mid i' \in [\ell], \text{ $i' \leq \ext{i}{r}{-}$ or $i' > \ext{i}{r}{+}$}\}$.
        Hence, since $V(H)$ intersects $R'_j$ and $H$ is connected,
        it follows that $V(H)$ intersects $P_{\ext{i}{r}{-}} \cup P_{\ext{i}{r}{+}}$.
        Then, since $V(H)$ is disjoint from $S_{r+1}$, and 
        because $P_{i'} \setminus S_{r+1} \subseteq R'_{i'}$ for every $i' \in I_r$ by \eqref{eq:Pi_minus_Sr+1_in_R'i},
        it follows that
        $V(H)$ intersects $R'_{\ext{i}{r}{-}} \cup R'_{\ext{i}{r}{+}}$. 
        By the uniqueness of $j$ we obtain that
        $j \in \{\ext{i}{r}{-}, \ext{i}{r}{+}\}$, a contradiction.
        Since $\{\ext{i}{r}{-}, \ext{i}{r}{+}\} \subseteq \Trace(i)$, this concludes the proof of the claim.
    \end{proofclaim}

    For every $i \in [\ell]$, let
    \[
        \Lambda_i = \bigcup \{\lambda(R'_j) \mid j \in \Trace(i)\}.
    \]

    \begin{claim}\label{claim:cen:apex_forest:size_Lambda_i}
        For every $i \in [\ell]$,
        \[
            |\Lambda_i| \leq cd(s+2) \leq 12 cd \cdot \log q.
        \]
    \end{claim}

    \begin{proofclaim}
        Let $j \in [\ell]$.
        Note that
        \begin{align*}
            |\lambda(R'_j)| 
            &\leq |\lambda(R'_j \setminus P_j)| + |\lambda(P_j)| \\
            &\leq |R'_j \setminus P_j| + |\psi_j(P_j)| \\
            &\leq |R_j \setminus P_j| + c d \\
            &\leq c d + c d \\
            &= 2 c d
        \end{align*}
        where the last inequality follows from \ref{item:centered_coloring:base_case:excluding_a_star:Ri_minus_Pi_small:call}.
        Hence $|\Lambda_i| \leq 2cd |\Trace(i)|$.

        Now we argue that $|\Trace(i)|\leq (s+2)$. 
        First, note that $|\Trace(i)\cap I_0|\leq 2$. 
        For every $r \in \{0, \dots, s-1\}$,
        $\{\ext{j}{r}{-},\ext{j}{r}{+}\}$ intersects $\{\ext{j}{r+1}{-},\ext{j}{r+1}{+}\}$. 
        Thus, for each $r\in \{0, \dots, s\}$, we have 
        $|\Trace(i)\cap \bigcup_{t\in\set{0,\ldots,r}} I_t| \leq 2 + r$. 
        In particular, $|\Trace(i)| \leq s+2$.
        Recall that $s\leq 3+\log q$.
        Therefore,
        \begin{align*}
            |\Lambda_i| 
            &\leq 2cd \cdot (s + 2) \\
            &\leq 2cd \cdot (5 + \log q)\\
            &\leq 12 cd \cdot \log q. \qedhere
        \end{align*}
    \end{proofclaim}

    \begin{claim}\label{claim:lemma:centered_coloring:base_case:excluding_a_forest}
        For every connected subgraph $H$ of $G$ intersecting $S_0$, 
        one of the following holds,
        \begin{enumerate}[label=(\arabic*)]
            \item $|\phi(V(H))| > q$, or \label{claim:item-i}
            \item $|\lambda(V(H) \cap S_0)| > q$, or \label{claim:item-ii}
            \item there is a $(\phi \times \lambda)$-center $u$ of $V(H) \cap S_0$ such that for every $i \in [\ell]$, if $W_i \cap V(H)\neq \emptyset$, then $\lambda(u) \in \Lambda_{i}$. \label{claim:item-iii}
        \end{enumerate}
    \end{claim}

    \begin{proofclaim}
        Let $H$ be a connected subgraph of $G$ intersecting $S_0$,
        and such that $|\phi(V(H))| \leq q$ and
        $|\lambda(V(H) \cap S_0)| \leq q$.
        Let $r = \max \{t \in \{0, \dots, s\} \mid V(H) \cap S_t \neq \emptyset\}$.
        By \Cref{claim:cen_excluding_a_forest:intermediate},
        either
        \begin{enumerate}[label=(\arabic*)]
            \item $|\{i \in \{0, \dots, \ell\} \mid V(H) \cap P_i \neq \emptyset\}| > q$, or
            \item there is a unique $j \in I_r$ such that $V(H) \cap R'_j \neq \emptyset$. Moreover, for every $i \in [\ell]$, if $V(H) \cap W_i \neq \emptyset$, then $j \in \Trace(i)$.
        \end{enumerate}
        Note that, since $(P_0, \dots, P_\ell)$ is a path partition of $(G,S')$,
        the set $\{i \in \{0, \dots, \ell\} \mid V(H) \cap P_i \neq \emptyset\}$ is an interval of integers.
        By the definition of $\mu$, this implies that $|\mu(V(H) \cap S')| \geq \min \{q+1, |\{i \in \{0, \dots, \ell\} \mid V(H) \cap P_i \neq \emptyset\}|\}$.
        Therefore,
        in the first case, 
        $|\lambda(V(H) \cap S_0)| \geq |\mu(V(H) \cap S')| > q$, a contradiction.
        Hence there is a unique $j \in I_r$ such that $V(H) \cap R'_j \neq \emptyset$,
        and for every $i \in [\ell]$, if $V(H) \cap W_i \neq \emptyset$,
        then $j \in \Trace(i)$.
        We fix this index $j \in I_r$.
        We will find a $(\phi \times \lambda)$-center $u$ of $V(H) \cap S_0$ which belongs to $R'_j$.
        First, by the definition of $\psi_j$,
        there is a $(\phi \times \eta_j)$-center $u_0$ of $V(H) \cap R_j$.
        It follows that $u_0$ is also a $(\psi_j\times\phi)$-center of $V(H) \cap R_j$.

        First, suppose that $u_0 \in R'_j$.
        We claim that $u_0$ is a $(\phi \times \lambda)$-center of $V(H) \cap S_0$.
        Indeed, for every $u' \in V(H) \cap S_0$ such that $(\phi \times \lambda)(u') = (\phi \times \lambda)(u_0)$,
        we have $u' \in R'_{j'}$ for some $j' \in I_r$.
        It follows from the uniqueness of $j$ that $j=j'$, and so $u' \in R'_j$.
        As $(\phi \times \psi_j)(u') = (\phi \times \psi_j)(u_0)$, and because $u_0$ is a $(\phi \times \psi_j)$-center of $V(H) \cap R'_j$,
        we deduce that $u'=u_0$.
        This shows that $u_0$ is a $(\phi \times \lambda)$-center of $V(H) \cap S_0$.

        Next, suppose that $u_0 \not\in R'_j$.
        We claim that $V(H)$ is disjoint from $P_j$.
        First, suppose $r<s$.
        Then, since $\{j-2^r,j+2^r\} \cap \{0, \dots, \ell\} \subseteq I_{r+1}$,
        we have by \eqref{eq:Pi_in_Sr} the inclusion $P_{j-2^r} \cup P_{j+2^r} \subseteq S_{r+1}$ (with the convention $P_i = \emptyset$ for every $i \not\in \{0, \dots, \ell\}$).
        By the maximality of $r$, $V(H)$ is disjoint from $S_{r+1}$,
        and so $V(H)$ is disjoint from $P_{j-2^r} \cup P_{j+2^r}$.
        Since $(P_0, \dots, P_\ell)$ is a path partition of $(G,S')$,
        and by the definition of $(W_1, \dots, W_\ell)$,
        we deduce that $V(H) \subseteq \bigcup \{W_i \mid j-2^r+1 \leq i \leq j+2^r\}$.
        But then, $V(H) \cap R'_j = V(H) \cap R_j$, which contradicts the fact that $u_0 \in R_j \setminus R'_j$.
        Hence $r=s$.
        Then, by the definition of $R'_j$,
        $u_0 \not\in \bigcup \{W_i \mid j-q+1 \leq i \leq j+q\}$.
        Hence, if $V(H)$ intersects $P_j$, 
        because $(P_0, \dots, P_\ell)$ is a path partition of $(G,S')$,
        and by the definition of $(W_1, \dots, W_\ell)$,
        $V(H)$ intersects either $P_i$ for each $i \in \{j-q, \dots, j\}$,
        or $P_i$ for each $i \in \{j, \dots, j+q\}$.
        In both cases,
        we have $|\{i' \in [\ell] \mid P_{i'} \cap V(H) \neq \emptyset\}| > q$,
        and so $|\lambda(V(H) \cap S_0)| \geq |\mu(V(H) \cap S')| > q$, a contradiction.
        This proves that $V(H)$ is disjoint from $P_j$.

        Then, let $u$ be an arbitrary vertex in $R'_j \cap V(H)$.
        We claim that $u$ is a $(\phi \times \lambda)$-center of $V(H) \cap S_0$.
        Let $u' \in V(H) \cap S_0$ be such that $(\phi \times \lambda)(u') = (\phi \times \lambda)(u)$.
        Then, $u' \in S_r \setminus S_{r+1}$, and so $u' \in R'_{j'}$ for some $j' \in I_r$.
        By the uniqueness of $j$, we deduce that $j'=j$ and so $u' \in R'_j$.
        Moreover, since $V(H)$ is disjoint from $P_j$, $u' \in R'_j \setminus P_j$.
        Since $\psi_j(u') = \psi_j(u)$,
        and because $\psi_j \vert_{R_j \setminus P_j}$ is injective,
        we deduce that $u=u'$.
        This proves that $u$ is a $(\phi \times \lambda)$-center of $V(H) \cap S_0$.

        Finally, 
        for every $i \in [\ell]$ such that $V(H) \cap W_i \neq \emptyset$,
        we have $j \in \Trace(i)$ and $u \in R'_j$, 
        so $\lambda(u) \in \lambda(R'_j)$.
        By the definition of $\Lambda_i$, we deduce that $\lambda(u) \in \Lambda_i$.
        This proves the claim.
    \end{proofclaim}

    Let $C$ be a component of $G-S$.
    Let $i$ be the unique index in $[\ell]$ such that $W_i$ contains $V(C)$.
    We set
    \[
        \Lambda_C = \Lambda_i \cup \tilde{\Lambda}_{C}.
    \]
    To conclude the proof, we show that \ref{item:centered_coloring:base_case:excluding_a_forest:hitting}--\ref{item:centered_coloring:base_case:excluding_a_forest:center_among_log_colors} hold.

    For every $F \in \mathcal{F}$, if $S_0 \cap V(F) = \emptyset$,
    then $\tilde{S} \cap V(F) \neq \emptyset$ by \ref{item:centered_coloring:base_case:excluding_a_forest:hitting'}.
    Therefore, $S \cap V(F) \neq \emptyset$, which proves \ref{item:centered_coloring:base_case:excluding_a_forest:hitting}.

    Next, we prove \ref{item:centered_coloring:base_case:excluding_a_forest:connectivity}.
    Let $C$ be a component of $G-S$.
    Let $C'$ be the component of $G-S'$ containing $C$.
    First, by \ref{item:centered_coloring:base_case:excluding_a_forest:connectivity'},
    $N_{G-S_0}(V(C)) = N_{G}(V(C)) \cap \tilde{S}$ intersects at most $8c ((h-1)(d+1) + \binom{h-1}{2}) \cdot \log q$ component of $G-S_0-V(C)$.
    Second, by \ref{item:centered_coloring:base_case:excluding_a_star:connectivity:call},
    $N_G(V(C'))$ intersects at most three components of $G-V(C')$.
    It is now enough to bound $|N_G(V(C)) \cap (S_0 \setminus S')|$.

    Let $r \in \{0, \dots, s\}$.
    Consider $H = G[V(C) \cup (N_G(V(C)) \cap (S_r \setminus S_{r+1}))]$ (with the convention $S_{s+1} = \emptyset$).
    By \Cref{claim:cen_excluding_a_forest:intermediate}, and because $H$ intersects at most two of the $P_{i'}$ for $i' \in \{0, \dots, \ell\}$ (namely $P_{i'-1}$ and $P_{i'}$),
    there is a unique $j \in I_r \setminus I_{r+1}$ (with the convention $I_{s+1} = \emptyset$) such that $H$ intersects $R'_j$. Therefore,
    \[
        |N_G(V(C)) \cap (S_r \setminus (S_{r+1} \cup S'))| \leq |R'_j \setminus S'| \leq |R'_j \setminus P_j| \leq cd.
    \]
    
    Altogether, this yields
    \[
        |N_G(V(C)) \cap (S_0 \setminus S')| \leq cd(s+1)
    \]
    and so the number of components of $G-V(C)$ intersecting $N_G(V(C))$ is at most
    \begin{align*}
        3 + cd (s+1) &+ 8c \left(\textstyle(h-1)(d+1) + \binom{h-1}{2}\right) \cdot \log q \\
        &\leq 3 + cd(4+\log q) + 8c \left(\textstyle(h-1)(d+1) + \binom{h-1}{2}\right) \cdot \log q \\
        &\leq 3 + cd\cdot 5\log q + 8c \left(\textstyle(h-1)(d+1) + \binom{h-1}{2}\right) \cdot \log q \\
        &\leq 8cd \cdot \log q + 8c \left(\textstyle(h-1)(d+1) + \binom{h-1}{2}\right) \cdot \log q \\
        &\leq 8c \left(\textstyle hd + \binom{h}{2}\right) \cdot \log q.
    \end{align*}
    This proves \ref{item:centered_coloring:base_case:excluding_a_forest:connectivity}.

    Let $C$ be a component of $G-S$. 
    By \Cref{claim:cen:apex_forest:size_Lambda_i},
    \begin{align*}
        |\Lambda_C|
        &\leq 12cd \cdot \log q + 12c \left(\textstyle(h-1)(d+1) + \binom{h-1}{2}\right) \cdot \log q \\
        &\leq 12c \left(\textstyle hd + \binom{h}{2}\right) \cdot \log q.
    \end{align*}
    We claim that~\ref{item:centered_coloring:base_case:excluding_a_forest:center_among_log_colors} is satisfied.
    Let $H$ be a connected subgraph of $G$ intersecting $S$.
    If $H$ is disjoint from $S_0$, then,
    by \ref{item:centered_coloring:base_case:excluding_a_forest:center_among_log_colors'}, one of \ref{item:centered_coloring:base_case:excluding_a_forest:center_among_log_colors'}.\ref{item:c-tilde:i}, \ref{item:centered_coloring:base_case:excluding_a_forest:center_among_log_colors'}.\ref{item:c-tilde:ii}, and~\ref{item:centered_coloring:base_case:excluding_a_forest:center_among_log_colors'}.\ref{item:c-tilde:iii} holds.
    If $|\phi(V(H))| > q$ (i.e.\ \ref{item:centered_coloring:base_case:excluding_a_forest:center_among_log_colors'}.\ref{item:c-tilde:i}), then~\ref{item:centered_coloring:base_case:excluding_a_forest:center_among_log_colors}.\ref{item:c-i} holds.
    If $|\tilde{\lambda}(V(H) \cap \tilde{S})| > q$ (i.e.\ \ref{item:centered_coloring:base_case:excluding_a_forest:center_among_log_colors'}.\ref{item:c-tilde:ii}), then 
    \[|\lambda(V(H) \cap S)| \geq |\tilde{\lambda}(V(H) \cap S)| = |\tilde{\lambda}(V(H) \cap \tilde{S})|  > q,\] 
    and so,~\ref{item:centered_coloring:base_case:excluding_a_forest:center_among_log_colors}.\ref{item:c-ii} holds.
    Finally, assume that there is a $(\phi \times \tilde{\lambda})$-center $u$ of $V(H) \cap \tilde{S}$ such that for every component $C$ of $G - S_0$ intersecting $V(H)$, $\tilde{\lambda}(u) \in \tilde{\Lambda}_{C}$ (i.e.\ \ref{item:centered_coloring:base_case:excluding_a_forest:center_among_log_colors'}.\ref{item:c-tilde:iii}).
    It follows that $u$ is a $(\phi \times \lambda)$-center $u$ of $V(H) \cap S$ such that for every component $C$ of $G$ intersecting $V(H)$, $\lambda(u) \in \Lambda_{C}$.
    This implies~\ref{item:centered_coloring:base_case:excluding_a_forest:center_among_log_colors}.\ref{item:c-iii}.

    Next, suppose that $V(H)$ intersects $S_0$.
    By \Cref{claim:lemma:centered_coloring:base_case:excluding_a_forest}, one of \ref{claim:lemma:centered_coloring:base_case:excluding_a_forest}.\ref{claim:item-i}, \ref{claim:lemma:centered_coloring:base_case:excluding_a_forest}.\ref{claim:item-ii}, and~\ref{claim:lemma:centered_coloring:base_case:excluding_a_forest}.\ref{claim:item-iii} holds.
    If $|\phi(V(H))| > q$ (i.e.\ \ref{claim:lemma:centered_coloring:base_case:excluding_a_forest}.\ref{claim:item-i}), then~\ref{item:centered_coloring:base_case:excluding_a_forest:center_among_log_colors}.\ref{item:c-i} holds.
    If $|\lambda(V(H) \cap S)| > q$ (i.e.\ \ref{claim:lemma:centered_coloring:base_case:excluding_a_forest}.\ref{claim:item-ii}), then $|\lambda(V(H) \cap S)| \geq |\lambda(V(H) \cap S_0)|  > q$, and so,~\ref{item:centered_coloring:base_case:excluding_a_forest:center_among_log_colors}.\ref{item:c-ii} holds.
    Finally, assume that there is a $(\phi \times \lambda)$-center $u$ of $V(H) \cap S_0$ such that for every $i \in [\ell]$ with $W_i \cap V(H) \neq \emptyset$, we have $\lambda(u) \in \Lambda_i$ (i.e.\ \ref{claim:lemma:centered_coloring:base_case:excluding_a_forest}.\ref{claim:item-iii}).
    Let $C$ be a component of $G - S$ intersecting $V(H)$.
    Let $i$ be the unique index in $[\ell]$ such that $C$ intersects $W_i$.
    Since $\lambda$ uses the colors of $\Lambda_i$ only on $S_0$, $u$ is $(\phi \times \lambda)$-center of $V(H) \cap S$ and $\lambda(u) \in \Lambda_i \subset \Lambda_C$.
    This implies~\ref{item:centered_coloring:base_case:excluding_a_forest:center_among_log_colors}.\ref{item:c-iii} and concludes the proof of the lemma.
\end{proof}

\begin{lemma}[{\cite[Proposition 7.2.1]{Diestel_book}}]\label{lemma:average_degree_Kk_minor_free}
    For every integer $k$ with $k \geq 3$, every $K_k$-minor-free graph $G$ has at most $2^{k-3} |V(G)|-1$ edges.
\end{lemma}

We now show a slightly modified version of
\Cref{lemma:centered_coloring:base_case:excluding_a_forest}
where in \ref{item:centered_coloring:base_case:excluding_a_forest:connectivity}, the number of components of $G-V(C)$ intersecting $N_G(V(C))$ is
now only a function of $k$, under the assumption that $G$ is $K_k$-minor-free for some fixed positive integer $k$.

\begin{lemma}\label{lemma:centered_coloring:base_case:excluding_a_forest_small_interface}
    Let $q$, $c$, $d$, and $k$ be positive integers with $q \geq 2$ and $k \geq 3$.
    Let $h$ be a nonnegative integer,
    let $G$ be a $K_k$-minor-free graph,
    let $\phi$ be a $(q,c)$-good coloring of $G$, and
    let $\mathcal{F}$ be a family of connected subgraphs of $G$ such that $G$ has no $\mathcal{F}$-rich model of $F_{h+1,d}$.
    There exists
    a set $\Lambda$ of at most $28c^2(hd+\binom{h}{2}) \cdot q\log q$ colors,
    a set $S \subseteq V(G)$, 
    a function $\lambda \colon S \to \Lambda$,
    and a family $(\Lambda_C \mid \text{$C$ component of $G-S$})$ of subsets of $\Lambda$, each of size at most $20c^2 (hd + \binom{h}{2}) \cdot \log q$ such that
    \begin{enumerateOurAlph}
        \item $S \cap V(F) \neq \emptyset$ for every $F \in \mathcal{F}$; \label{item:centered_coloring:base_case:excluding_a_forest_small_interface:hitting}
        \item for every component $C$ of $G-S$, $N_G(V(C))$ intersects at most $2^{k-1}$ components of $G-V(C)$; and \label{item:centered_coloring:base_case:excluding_a_forest_small_interface:connectivity}
        \item \label{item:centered_coloring:base_case:excluding_a_forest_small_interface:center_among_log_colors}
            for every connected subgraph $H$ of $G$ intersecting $S$, 
            \begin{enumerate}[label=\normalfont (\arabic*)]
                \item $|\phi(V(H))| > q$, or
                \item $|\lambda(V(H) \cap S)| > q$, or
                \item there is a $(\phi \times \lambda)$-center $u$ of $V(H) \cap S$ such that 
                for every component $C$ of $G-S$ intersecting $V(H)$,
                $\lambda(u) \in \Lambda_{C}$.
            \end{enumerate}
    \end{enumerateOurAlph}
\end{lemma}

\begin{proof}
    By \Cref{lemma:centered_coloring:base_case:excluding_a_forest},
    there exists
    a set $\Lambda_0$ of at most $20c(hd+\binom{h}{2}) \cdot q\log q$ colors,
    a set $S_0 \subseteq V(G)$, a function $\lambda_0 \colon S_0 \to \Lambda_0$,
    and a family $(\Lambda_{0,C} \mid \text{$C$ component of $G-S_0$})$ of subsets of $\Lambda_0$, each of size at most $12c (hd + \binom{h}{2}) \cdot \log q$,
    such that
    \begin{enumerate}[label=\ref{lemma:centered_coloring:base_case:excluding_a_forest}.(\alph*)]
        \item $S_0 \cap V(F) \neq \emptyset$ for every $F \in \mathcal{F}$; \label{item:centered_coloring:base_case:excluding_a_forest:hitting:call}
        \item for every component $C$ of $G-S_0$, $N_G(V(C))$ intersects at most $8c (hd + \binom{h}{2}) \cdot \log q$ components of $G-V(C)$; and \label{item:centered_coloring:base_case:excluding_a_forest:connectivity:call}
        \item \label{item:centered_coloring:base_case:excluding_a_forest:center_among_log_colors:call}
            for every connected subgraph $H$ of $G$ intersecting $S$, 
            \begin{enumerate}
                \item $|\phi(V(H))| > q$, or
                \item $|\lambda_0(V(H) \cap S)| > q$, or
                \item there is a $(\phi \times \lambda_0)$-center $u$ of $V(H) \cap S_0$ such that for every component $C$ of $G-S_0$ intersecting $V(H)$, $\lambda_0(u) \in \Lambda_{0,C}$.
            \end{enumerate}
    \end{enumerate}
    By relabeling the elements of $\Lambda_0$, we assume $\Lambda_0$ is disjoint from $\nonnegint$.

    Let $\mathcal{C}$ be the family of all the components of $G-S_0$.
    Let $C \in \mathcal{C}$.
    Let $\mathcal{F}'$ be the family of all the connected subgraphs $H$ of $C$ such that $N_G(V(H))$ intersects at least $2^{k-2}$ components of $G-V(C)$.
    Let $N = \big\lfloor 8c (hd + \binom{h}{2}) \cdot \log q \big\rfloor$.
    We claim that there are no $N+1$ pairwise disjoint members of $\mathcal{F}'$.
    Suppose for contradiction that $H_1, \dots, H_{N+1}$ are pairwise disjoint members of $\mathcal{F}'$.
    Let $A$ be the minor of $G$ obtained from $G$ by contracting every component of $G-V(C)$ having a neighbor in $V(C)$ in to a single vertex,
    and every $H_i$ into a single vertex.
    We denote by $u_1, \dots, u_m$ the vertices resulting the contractions of these components of $G-V(C)$,
    and by $v_1, \dots, v_{N+1}$ the vertices resulting from the contractions of $H_1, \dots, H_{N+1}$.
    By \ref{item:centered_coloring:base_case:excluding_a_forest:connectivity:call},
    $m \leq N+1$.
    Then, $A$ is a graph with $m+N+1 \leq 2(N+1)$ vertices, and every $v_i$ for $i \in [N+1]$ has at least $2^{k-2}$ neighbors in $\{u_1, \dots, u_m\}$.
    Hence $|E(A)| \geq 2^{k-2}(N+1)$. On the other hand, $A$ is a minor of $G$ and so is $K_k$-minor-free.
    By \Cref{lemma:average_degree_Kk_minor_free}, this implies that $|E(A)| \leq 2^{k-3}|V(A)|-1 < 2^{k-2}(N+1)$.
    This contradiction proves that there are no $N+1$ pairwise disjoint members of $\mathcal{F}'$.

    Hence, since $\phi$ is a $(q,c)$-good coloring,
    there exists a set $Z_C \subseteq V(C)$ and a coloring $\psi_C \colon Z_C \to [c \cdot N]$
    such that
    \begin{enumerate}
        \item $Z_C \cap V(F) \neq \emptyset$ for every $F \in \mathcal{F}'$; \label{item:proof:centered_coloring:base_case:excluding_a_forest_small_interface:ZC_hitting}
        \item for every component $C'$ of $C-Z$, $N_C(V(C'))$ intersects at most two components of $C-V(C')$; and \label{item:proof:centered_coloring:base_case:excluding_a_forest_small_interface:ZC_interface}
        \item for every connected subgraph $H$ of $C$, either $|\phi(V(H))| > q$, or there is a $(\phi \times \psi_C)$-center of $V(H) \cap Z$. \label{item:proof:centered_coloring:base_case:excluding_a_forest_small_interface:ZC_chip}
    \end{enumerate}

    Let
    \[
        S = S_0 \cup \bigcup_{C \in \mathcal{C}} Z_C,
    \]
    and for every $u \in S$, let
    \[
        \lambda(u) =
        \begin{cases}
            \lambda_0(u) &\textrm{if $u \in S_0$,} \\
            \psi_C(u) &\textrm{otherwise, where $C$ is the unique member of $\mathcal{C}$ containing $u$.}
        \end{cases}
    \]

    First note that $\lambda$ uses at most
    $20c (hd + \binom{h}{2}) \cdot q\log q + c \cdot N \leq 28c^2(hd+\binom{h}{2}) \cdot q\log q$ colors.

    We now prove \ref{item:centered_coloring:base_case:excluding_a_forest_small_interface:hitting}--\ref{item:centered_coloring:base_case:excluding_a_forest_small_interface:center_among_log_colors}.

    First, $S_0 \subseteq S$ and so $S \cap V(F) \neq \emptyset$ for every $F \in \mathcal{F}$ by \ref{item:centered_coloring:base_case:excluding_a_forest:hitting:call}.
    Hence \ref{item:centered_coloring:base_case:excluding_a_forest_small_interface:hitting} holds.

    Let $C'$ be a component of $G-S$.
    Let $C$ be the component of $G-S_0$ containing $C'$.
    Since $C'$ is disjoint from $Z_C$, $N_G(V(C'))$ intersects at most $2^{k-2}-1$ components of $G-V(C)$ by \ref{item:proof:centered_coloring:base_case:excluding_a_forest_small_interface:ZC_hitting}.
    Moreover, by \ref{item:proof:centered_coloring:base_case:excluding_a_forest_small_interface:ZC_interface}, $N_G(V(C'))$ intersects at most two components of $C-V(C')$ and thus of $G-V(C')$.
    Hence, $N_G(V(C'))$ intersects at most $2^{k-2}+1 \leq 2^{k-1}$ components of $G-V(C')$.
    This proves \ref{item:centered_coloring:base_case:excluding_a_forest_small_interface:connectivity}.

    Let $C$ be a component of $G-S$.
    Let $C'$ be the component of $G-S_0$ containing $C$.
    Let
    \[
        \Lambda_{C} = \Lambda_{0,C'} \cup [c \cdot N].
    \]
    Note that
    \begin{align*}
    |\Lambda_C| &\leq |\Lambda_{0,C'}|+c\cdot N \\
    &\leq 12c\left(dh + \binom{h}{2}\right) \log q + 8c^2 \left(dh + \binom{h}{2}\right) \log q\\
    &\leq 20c^2 \left(dh + \binom{h}{2}\right) \log q.
    \end{align*}
    Let $H$ be a connected subgraph of $G$ intersecting $S$ such that $|\phi(V(H))| \leq q$ and $|\lambda(V(H) \cap S)| \leq q$.
    Suppose that $C_1, \dots, C_m$ are the components of $G-S$ intersecting $V(H)$.
    For every $i \in [m]$, we denote by $C'_i$ the component of $G-S_0$ containing $V(C_i)$.
    If $H$ intersects $S_0$,
    then by \ref{item:centered_coloring:base_case:excluding_a_forest:center_among_log_colors:call}, 
    there is a $(\phi \times \lambda_0)$-center $u$ of $V(H) \cap S_0$ with $\lambda_0(u) \in \bigcup_{i \in [m]} \Lambda_{0,C'_i}$,
    and so $u$ is a $(\phi \times \lambda)$-center of $V(H) \cap S$ with $\lambda(u) \in \bigcup_{i \in [m]} \Lambda_{C_i}$.
    Now suppose that $H$ is disjoint from $S_0$.
    Then $C'_1 = \dots = C'_m$.
    By \ref{item:proof:centered_coloring:base_case:excluding_a_forest_small_interface:ZC_chip},
    there is a $(\phi \times \psi_{C'_1})$-center $u$ of $V(H) \cap Z_{C'_1}$,
    and so $u$ is a $(\phi \times \lambda)$-center of $V(H) \cap S$ with $\lambda(u) \in [c \cdot N] \subseteq \Lambda_{C_i}$ for every $i \in [m]$.
    This shows \ref{item:centered_coloring:base_case:excluding_a_forest_small_interface:center_among_log_colors},
    and concludes the proof of the lemma.
\end{proof}

Next, we apply iteratively \Cref{lemma:centered_coloring:base_case:excluding_a_forest_small_interface}
to obtain a coloring for graphs with no $\mathcal{F}$-rich model of a given apex-forest.
We use the same technique as in \Cref{claim:X-to-A(X)}, while building a suitable coloring.
Item \ref{item:centered_coloring:base_case:excluding_a_forest_small_interface:center_among_log_colors} in
\Cref{lemma:centered_coloring:base_case:excluding_a_forest_small_interface} allows us to reuse most colors in a component $C$ of $G-S$, except a few ones in the set $\Lambda_C$ of size $\bigO(\log q)$.
The idea is to maintain a small ``buffer'' of colors $(\Lambda_1, \dots, \Lambda_a)$ which are not yet reusable.
Here, $\Lambda_a$ for $a \in [q]$ is a set of the form $\Lambda_C$ obtained $a$ steps before the current induction call.
We will maintain the following property:
if a connected subgraph $H$, intersecting the component currently being processed, 
has as center with color in $\Lambda_a$, 
then $H$ already contains at least $a$ different colors.
Thus, when we reach $a=q+1$, we can reuse the colors in this set.

\begin{lemma}\label{lemma:centered_coloring:base_case:excluding_an_apex_forest}
    Let $c$ and $k$ be positive integers, and
    let $X$ be a forest.
    There exists an integer $c_{\ref{lemma:centered_coloring:base_case:excluding_an_apex_forest}}(c,k,X)$ such that the following holds.
    Let $q$ be an integer with $q \geq 2$,
    let $\Lambda$ be a set of at least $2c_{\ref{lemma:centered_coloring:base_case:excluding_an_apex_forest}}(c,k,X) \cdot q\log q$ colors.
    let $G$ be a connected $K_k$-minor-free graph,
    let $\phi$ be a $(q,c)$-good coloring of $G$, and
    let $\mathcal{F}$ be a family of connected subgraphs of $G$ such that $G$ has no $\mathcal{F}$-rich model of $K_1 \oplus X$.
    For every nonempty set $U \subseteq V(G)$, for every component $C$ of $G-U$, for every $\lambda_0 \colon U \to \Lambda$, and for all subsets $\Lambda_1, \dots, \Lambda_{q}$ of $\Lambda$ 
    with $|\Lambda_a| \leq c_{\ref{lemma:centered_coloring:base_case:excluding_an_apex_forest}}(c,k,X) \cdot \log q$ for every $a \in [q]$
    if
    \begin{enumerate}
        \item $N_G(V(C))$ intersects at most $2^{k-1}$ components of $G-V(C)$ and \label{lemma:centered_coloring:base_case:excluding_an_apex_forest:i}
        \item for every connected subgraph $H$ of $G$ intersecting both $U$ and $V(C)$ \label{lemma:centered_coloring:base_case:excluding_an_apex_forest:ii}
            \begin{enumerate}[label={\normalfont (\arabic*)}]
                \item $|\phi(V(H))|> q$, or \label{item:lemma:centered_coloring:base_case:excluding_an_apex_forest:ii:1}
                \item $|\lambda_0(V(H)\cap U)| >q$, or \label{item:lemma:centered_coloring:base_case:excluding_an_apex_forest:ii:2}
                \item there exists a $(\phi \times \lambda_0)$-center $u_H$ of $V(H) \cap U$ and $a \in [q]$ such that $\lambda_0(u_H) \in \Lambda_a$ and
                    $|\lambda_0(V(H) \cap U) \cap (\bigcup_{b \in [a]} \Lambda_b)| \geq a$ -- we fix such a vertex $u_H$;
                    \label{item:lemma:centered_coloring:base_case:excluding_an_apex_forest:ii:3}
            \end{enumerate}
    \end{enumerate}
    then there exist $S$ with $U\subseteq S \subseteq V(C) \cup U$ and a coloring $\lambda \colon S \to \Lambda$ extending $\lambda_0$ such that
    \begin{enumerateOurAlph}
        \item $S \cap V(F) \neq \emptyset$ for every $F \in \mathcal{F}\vert_C$; \label{item:centered_coloring:base_case:excluding_an_apex_forest:hitting}
        \item for every component $C'$ of $C-S$, $N_G(V(C'))$ intersects at most $2^k$ components of $G-V(C')$; and \label{item:centered_coloring:base_case:excluding_an_apex_forest:connectivity}
        \item for every connected subgraph $H$ of $G$ intersecting both $S$ and $V(C)$, \label{item:centered_coloring:base_case:excluding_an_apex_forest:chip}
            \begin{enumerate}[label={\normalfont (\arabic*)}]
                \item $|\phi(V(H))| > q$, or
                \item $|\lambda(V(H) \cap S)| > q$, or
                \item there is a $(\phi \times \lambda)$-center $u$ of $V(H) \cap S$ and moreover, if $V(H) \cap U \neq \emptyset$, then we may set $u=u_H$. \label{item-2-iii}
            \end{enumerate}
    \end{enumerateOurAlph}
\end{lemma}

\begin{figure}[tp]
    \centering
    \includegraphics{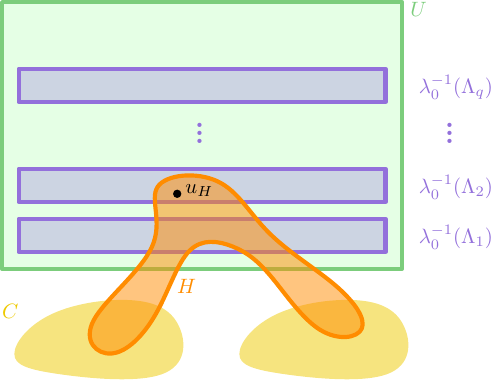}
    \caption{Illustration for the assumption \ref{lemma:centered_coloring:base_case:excluding_an_apex_forest:ii}.\ref{item-2-iii} in~\Cref{lemma:centered_coloring:base_case:excluding_an_apex_forest}.
    Intuitively, the set $\Lambda_a$ for $a \in [q]$ represents colors that can be taken by a center $u_H$ of a connected subgraph $H$ intersecting both $U$ and $V(C)$, but at a cost that $V(H)$ uses at least $a$ colors.}
    \label{fig:lemma:cen_base_case_apex_forest}
\end{figure}

\begin{proof}
    See \Cref{fig:lemma:cen_base_case_apex_forest}.
    By \Cref{lemma:Rt_has_coloring_elimination_property},
    the class of all forests has the coloring elimination property.
    Let $X'$ be a forest witnessing the coloring elimination property for $X \sqcup K_1$ and let $h$ and $d$ be positive integers with $X' \subseteq F_{h+1,d}$.
    Let
    \[
        \textstyle c_{\ref{lemma:centered_coloring:base_case:excluding_an_apex_forest}}(c,k,X) = 28 c^2\left(hd+\binom{h}{2}\right).
    \]

    A tuple $(U,C,\lambda_0, (\Lambda_a \mid a \in [q]), \mathcal{H}, (u_H \mid H \in \mathcal{H}))$ is an instance if $U \subseteq V(G)$ is a nonempty set, $C$ is a component of $G-U$, $\lambda_0 \colon U \to \Lambda$, and $\Lambda_1, \dots, \Lambda_{q}$ are subsets of $\Lambda$ 
    with $|\Lambda_a| \leq c_{\ref{lemma:centered_coloring:base_case:excluding_an_apex_forest}}(c,k,X) \cdot \log q$ for every $a \in [q]$, and these objects satisfy \ref{lemma:centered_coloring:base_case:excluding_an_apex_forest:i} and \ref{lemma:centered_coloring:base_case:excluding_an_apex_forest:ii}; additionally, $\calH$ is the family of all connected subgraphs $H$ of $G$ such that $V(H)\cap U \neq \emptyset$, $V(H)\cap V(C)\neq\emptyset$, $|\phi(V(H))| \leq q$, and $|\lambda_0(V(H) \cap U)| \leq q$; and finally $(u_H \mid H\in \calH)$ witnesses item~\ref{lemma:centered_coloring:base_case:excluding_an_apex_forest:ii}.\ref{item:lemma:centered_coloring:base_case:excluding_an_apex_forest:ii:3}.

    It suffices to prove that given an instance, there exist $S$ with $U\subseteq S \subseteq V(C) \cup U$ and a coloring $\lambda \colon S \to \Lambda$ extending $\lambda_0$ such that \ref{item:centered_coloring:base_case:excluding_an_apex_forest:hitting}--\ref{item:centered_coloring:base_case:excluding_an_apex_forest:chip} hold.
    We prove this statement by induction on $|V(C)|$.

    Let $(U,C,\lambda_0, (\Lambda_a \mid a \in [q]), \mathcal{H}, (u_H \mid H \in \mathcal{H}))$ be an instance.
    If $\mathcal{F}\vert_C = \emptyset$, then it suffices to set $S=U$ and $\lambda=\lambda_0$, thus, we suppose that $\mathcal{F}\vert_C \neq \emptyset$.
    Some objects considered in the proof are illustrated in~\Cref{fig:lemma:proof:cen_base_case_apex_forest}.

    Let $A_1, \dots, A_\ell$ be the components of $G-V(C)$.
    By assumption \ref{lemma:centered_coloring:base_case:excluding_an_apex_forest:i} for $\calI$, we have 
    $1 \leq \ell \leq 2^{k-1}$.
    Let $\mathcal{F}'$ be the family of all the connected subgraphs $H$ of $C$
    such that $H$ contains a member of $\mathcal{F}$ and intersects $N_G(U) = \bigcup_{i=1}^\ell N_G(V(A_i))$.
    We claim that there is no $\mathcal{F}'$-rich model of $X'$ in $C$.
    Suppose to the contrary that $C$ has an $\mathcal{F}'$-rich model $(D_x \mid x \in V(X'))$ of $X'$.
    For every $i \in [\ell]$, let $S_i$ be the set of all $x \in V(X')$ such that $D_x$ contains a vertex of $N_G(V(A_i))$.
    By the coloring elimination property,
    there exists $i \in [\ell]$ such that $X'$ contains an $S_i$-rooted model of $X \sqcup K_1$.
    It follows that $C$ contains an $\mathcal{F}$-rich model $(B_x \mid x \in V(X \sqcup K_1))$ of $X \sqcup K_1$ 
    such that every branch set has a neighbor in $A_i$.
    Let $x_0$ be the vertex of $K_1$ in $X \sqcup K_1$.
    Then the model $(B_x \mid x \in V(X))$ together with $B_{x_0} \cup A_i$ is an $\mathcal{F}$-rich model of $K_1 \oplus X$ in $G$, contradicts one of the assumptions.
    This shows that there is no $\mathcal{F}'$-rich model of $X'$ in $C$, as claimed.

    Let $\tilde{\Lambda} = \Lambda \setminus \bigcup_{a \in [q]} \Lambda_a$.
    Note that 
    \begin{align*}
    |\tilde{\Lambda}| \geq 2c_{\ref{lemma:centered_coloring:base_case:excluding_an_apex_forest}}(c,k,X) \cdot q\log q - q \cdot c_{\ref{lemma:centered_coloring:base_case:excluding_an_apex_forest}}(c,k,X) \cdot \log q 
    \geq c_{\ref{lemma:centered_coloring:base_case:excluding_an_apex_forest}}(c,k,X) \cdot q\log q.
    \end{align*}
    By \Cref{lemma:centered_coloring:base_case:excluding_a_forest_small_interface},
    there exists
    a set $\tilde{S} \subseteq V(C)$,
    a function $\tilde{\lambda} \colon \tilde{S} \to \tilde{\Lambda}$,
    and a family $(\tilde{\Lambda}_{C'} \mid \text{$C'$ component of $C-\tilde{S}$})$ of subsets of $\Tilde{\Lambda}$, 
    each of size at most $20c^2 (hd + \binom{h}{2}) \cdot \log q$,
    such that
    \begin{enumerate}[label=\ref{lemma:centered_coloring:base_case:excluding_a_forest_small_interface}.(\alph*)]
        \item $\tilde{S} \cap V(F) \neq \emptyset$ for every $F \in \mathcal{F}'$; 
            \label{item:centered_coloring:base_case:excluding_a_forest_small_interface:hitting:call}
        \item for every component $C'$ of $C-\tilde{S}$, $N_C(V(C'))$ intersects at most $2^{k-1}$ components of $C-V(C')$;  and 
            \label{item:centered_coloring:base_case:excluding_a_forest_small_interface:connectivity:call}
        \item \label{item:centered_coloring:base_case:excluding_a_forest_small_interface:center_among_log_colors:call}
            for every connected subgraph $H$ of $C$ intersecting $\tilde{S}$,
            \begin{enumerate}[label=(\arabic*)]
                \item $|\phi(V(H))| > q$, or
                \item $|\tilde{\lambda}(V(H) \cap \tilde{S})| > q$, or
                \item there is a $(\phi \times \tilde{\lambda})$-center $\tilde{u}_H$ of $V(H) \cap \tilde{S}$ 
                    such that for every component $C'$ of $C-\tilde{S}$ intersecting $V(H)$, $\tilde{\lambda}(\tilde{u}_H) \in \tilde{\Lambda}_{C'}$.
                    We fix such a vertex $\tilde{u}_H$.
            \end{enumerate}
    \end{enumerate}
    Since $\mathcal{F}\vert_C \neq \emptyset$ and $G$ is connected, $\tilde{S}$ is nonempty.
    Let 
    \[
        U' = U \cup \tilde{S},
    \]
    and let $\lambda'_0 \colon U' \to \Lambda$ be defined for every $u \in U'$ as
    \[
        \lambda'_0(u) = 
        \begin{cases}
            \lambda_0(u) &\textrm{if $u \in U$,} \\
            \tilde{\lambda}(u) &\textrm{if $u \in \tilde{S}$.}
        \end{cases}
    \]
    Let $H$ be a connected subgraph of $G$ intersecting $U'$ with $|\phi(V(H))| \leq q$ and $|\lambda'_0(V(H) \cap U')| \leq q$.
    In particular, $|\lambda_0(V(H) \cap U)| \leq q$ and $|\tilde{\lambda}(V(H) \cap \tilde{S})| \leq q$.
    If $V(H)\cap U \neq \emptyset$, then let $u'_H = u_H$,
    and if $V(H) \cap U = \emptyset$, then let $u'_H = \tilde{u}_H$.    
    \begin{figure}[tp]
        \centering
        \includegraphics{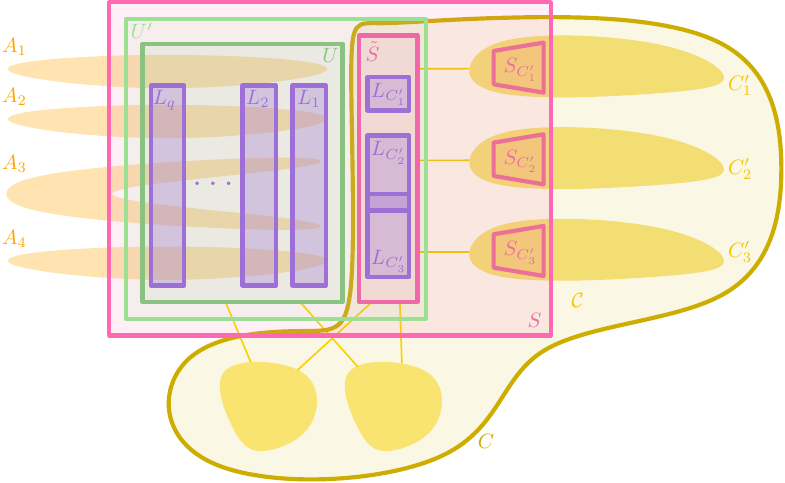}
        \caption{An illustration of some objects in the proof of~\Cref{lemma:centered_coloring:base_case:excluding_an_apex_forest}.
        Here, $L_a = \lambda_0^{-1}(\Lambda_a)$ for each $a \in [q]$ and $L_{C_i'} = \tilde{\lambda}^{-1}(\Lambda_{C_i}')$.
        Moreover, $\calC = \{C_1', C_2', C_3'\}$.
        }
        \label{fig:lemma:proof:cen_base_case_apex_forest}
    \end{figure}

    Let $\mathcal{C}$ be the family of all the components of $C-\tilde{S}$ that contain a member of $\mathcal{F}$.
    Let $C' \in \mathcal{C}$,
    and let $\mathcal{H}'$ be the family of all the connected subgraphs $H$ of $G$ intersecting both $V(C')$ and $U'$, and such that $|\phi(V(H))| \leq q$ and $|\lambda'_0(V(H) \cap U'| \leq q$.
    Let $(\Lambda'_1, \dots, \Lambda'_q) = (\tilde{\Lambda}_{C'}, \Lambda_1, \dots, \Lambda_{q-1})$.
   
    We now show that the tuple 
    \[
        \mathcal{I}' = \big(U',C', \lambda'_0, (\Lambda'_a \mid i \in [q]), \mathcal{H}', (u'_H \mid H \in \mathcal{H}')\big)
    \]
    is an instance. 
    
    Since $C' \not \in \mathcal{F}'$, we have $N_G(V(C')) \cap U = \emptyset$,
    and so $N_G(V(C')) \subset \tilde{S}$.
    Thus, $N_G(V(C'))$ intersects at most $2^{k-1}$ components of $G-V(C')$ 
    by \ref{item:centered_coloring:base_case:excluding_a_forest_small_interface:connectivity:call} for $\mathcal{I}$.
    This proves \ref{lemma:centered_coloring:base_case:excluding_an_apex_forest:i} for $\mathcal{I}'$.
    
    Let $H \in \mathcal{H}'$. 
    We assumed that \ref{lemma:centered_coloring:base_case:excluding_an_apex_forest:ii}.\ref{item:lemma:centered_coloring:base_case:excluding_an_apex_forest:ii:1} and \ref{lemma:centered_coloring:base_case:excluding_an_apex_forest:ii}.\ref{item:lemma:centered_coloring:base_case:excluding_an_apex_forest:ii:2} do not hold for $H$ in the context of the instance $\calI'$.
    Thus, the goal is to prove~\ref{lemma:centered_coloring:base_case:excluding_an_apex_forest:ii}.\ref{item:lemma:centered_coloring:base_case:excluding_an_apex_forest:ii:3}, and moreover, that it is witnessed by $u_H'$.
    Namely, we will show that $u_H'$ is a $(\phi\times\lambda_0')$-center of $V(H) \cap U'$ and there exists $a \in [q]$ such that $\lambda_0'(u_H') \in \Lambda_a'$ and $|\lambda_0'(V(H) \cap U') \cap (\bigcup_{b \in [a]} \Lambda_b')| \geq a$.
    
    First, suppose that $V(H)$ intersects $U$.
    Since $(u_H \mid H \in \mathcal{H}')$ witnesses \ref{lemma:centered_coloring:base_case:excluding_an_apex_forest:ii}.\ref{item:lemma:centered_coloring:base_case:excluding_an_apex_forest:ii:3} for $\mathcal{I}$, $u'_H = u_H$ is a $(\phi \times \lambda_0)$-center of $V(H) \cap U$ such that
    for some $a \in [q]$, we have $\lambda_0(u'_H) \in \Lambda_{a}$.
    Since the colors used by $\lambda_0$ and $\tilde{\lambda}$
    are pairwise distinct, $u'_H$ is a $(\phi \times \lambda'_0)$-center of $V(H) \cap U'$.
    Since $N_G(V(C')) \subseteq \tilde{S}$, the set $V(H)$ intersects $\tilde{S}$.
    Hence there is an edge $vw$ in $H$ with $v \in \Tilde{S}$ and $w \in V(C')$.
    Since this edge is a connected subgraph of $C$ intersecting both $\tilde{S}$ and $V(C')$, and as $q \geq 2$, by \ref{item:centered_coloring:base_case:excluding_a_forest_small_interface:center_among_log_colors:call},
    there is a $(\phi \times \tilde{\lambda})$-center of $\{v,w\} \cap \tilde{S} = \{v\}$ and $\tilde{\lambda}(v) \in \tilde{\Lambda}_{C'}$.
    Therefore, $\lambda'_0(V(H) \cap U') \cap \tilde{\Lambda}_{C'}$ is nonempty.
    By the assumption on $\calI$, $|\lambda_0(V(H) \cap U) \cap (\bigcup_{b \in [a]} \Lambda_b)| \geq a$.
    Since $\tilde{\Lambda}_{C'}$ is disjoint from $\bigcup_{b \in [a]} \Lambda_b$, we deduce that
    \begin{align*}
        |\lambda'_0(V(H) \cap U') \cap (\textstyle\bigcup_{b \in [a]} \Lambda_b \cup \tilde{\Lambda}_{C'})|
        \geq a+1.
    \end{align*}
    If $a=q$, then we have $|\lambda'_0(V(H) \cap U')| > q$, a contradiction.
    Otherwise, $|\lambda'_0(V(H) \cap U') \cap (\textstyle\bigcup_{b \in [a+1]} \Lambda'_b)|\geq a+1$, as desired.
    
    Next, suppose that $V(H)$ is disjoint from $U$.
    Then $H$ is a connected subgraph of $C$.
    Hence, by \ref{item:centered_coloring:base_case:excluding_a_forest_small_interface:center_among_log_colors:call},
    and because $|\phi(V(H))| \leq q$ and $|\tilde{\lambda}(V(H) \cap \tilde{S})| \leq q$,
    the vertex $u'_H$ is a $(\phi \times \tilde{\lambda})$-center of $V(H) \cap \tilde{S} = V(H) \cap U'$ such that $\lambda'_0(u'_{H}) \in \tilde{\Lambda}_{C'} = \Lambda'_1$.
    Clearly, $|\lambda'_0(V(H) \cap U) \cap (\textstyle \bigcup_{b \in [1]} \Lambda'_b)| \geq 1$.
    This proves \ref{lemma:centered_coloring:base_case:excluding_an_apex_forest:ii} for $\mathcal{I}'$, and so $\mathcal{I}'$ is an instance.

    Since $\mathcal{F}' \neq \emptyset$, the set $\tilde{S}$ is nonempty, and so $|V(C')| < |V(C)|$.
    Hence we can apply the induction hypothesis to the instance $\mathcal{I}'$.
    This gives a set $S_{C'}$ with $U' \subseteq S_{C'} \subseteq V(C') \cup U'$ and a coloring $\lambda_{C'} \colon S_{C'} \to \Lambda$ extending $\lambda'_0$ such that
    \begin{enumerateOurAlphPrim}
        \item $S_{C'} \cap V(F) \neq \emptyset$ for every $F \in \mathcal{F}\vert_{C'}$; \label{item:centered_coloring:base_case:excluding_an_apex_forest:hitting'}
        \item for every component $C''$ of $C'-S_{C'}$, $N_G(V(C''))$ intersects at most $2^k$ components of $G-V(C'')$; and \label{item:centered_coloring:base_case:excluding_an_apex_forest:connectivity'}
        \item for every connected subgraph $H$ of $G$ intersecting $V(C')$ and $S_{C'}$, \label{item:centered_coloring:base_case:excluding_an_apex_forest:chip'}
            \begin{enumerate}[label=(\arabic*)]
                \item $|\phi(V(H))| > q$, or
                \item $|\lambda_{C'}(V(H) \cap S_{C'})| >q$, or
                \item there is a $(\phi \times \lambda_{C'})$-center $u$ of $V(H) \cap S_{C'}$ and moreover, if $V(H) \cap U' \neq \emptyset$, then we may set $u=u_H'$.\label{item:centered_coloring:base_case:excluding_an_apex_forest:chip'-3}
            \end{enumerate}
    \end{enumerateOurAlphPrim}

    Let
    \[
        S = \bigcup_{C' \in \mathcal{C}} S_{C'}
    \]
    and let $\lambda\colon S \to \Lambda$ be defined as follows.
    For every $u \in S$, 
    we find $C' \in \mathcal{C}$ with $u \in S_{C'}$ and we set
    \[
        \lambda(u) = \lambda_{C'}(u).
    \]
    This is well defined since the colorings $\lambda_{C'}$ for $C' \in \mathcal{C}$ coincide on $U'$ as they extend $\lambda'_0$.
    We now show that \ref{item:centered_coloring:base_case:excluding_an_apex_forest:hitting}--\ref{item:centered_coloring:base_case:excluding_an_apex_forest:chip} hold.

    Let $F \in \mathcal{F} \vert_C$.
    If $U' \cap V(F) \neq \emptyset$ then $S \cap V(F) \neq \emptyset$.
    Otherwise, there exists $C' \in \mathcal{C}$ such that $F \in \calF\vert_{C'}$.
    Then, by \ref{item:centered_coloring:base_case:excluding_an_apex_forest:hitting'}, 
    $S_{C'} \cap V(F) \neq \emptyset$ and so $S \cap V(F) \neq \emptyset$.
    This proves \ref{item:centered_coloring:base_case:excluding_an_apex_forest:hitting}.

    Let $C''$ be a component of $C-S$.
    Since $C''$ is disjoint from $\tilde{S}$, either $N_G(V(C'')) \subseteq U \cup \tilde{S}$,
    or there exists $C' \in \mathcal{C}$ such that $C'' \subseteq C'$.
    In the first case,
    $N_G(V(C''))$ intersects at most $2^{k-1}$ components of $G-V(C)$ by~\ref{lemma:centered_coloring:base_case:excluding_an_apex_forest:i}, 
    and at most $2^{k-1}$ components of $C-\tilde{S}$ by \ref{item:centered_coloring:base_case:excluding_a_forest_small_interface:connectivity:call}.
    This implies that $N_G(V(C''))$ intersects at most $2^k$ components of $G-V(C'')$.
    In the second case, let $C' \in \mathcal{C}$ such that $C'' \subseteq C'$.
    Then, by \ref{item:centered_coloring:base_case:excluding_an_apex_forest:connectivity'},
    $N_G(V(C''))$ intersects at most $2^k$ components of $G-V(C'')$.
    This proves \ref{item:centered_coloring:base_case:excluding_an_apex_forest:connectivity}.

    Let $H$ be a connected subgraph of $G$ intersecting both $V(C)$ and $S$ such that
    $|\phi(V(H))| \leq q$ and $|\lambda(V(H) \cap S)| \leq q$.
    Let $C'_1, \dots, C'_m$ be the components of $C - \tilde{S}$ intersecting $V(H)$ (with possibly $m=0$).
    First, suppose that $V(H)$ intersects $U$.
    By~\ref{item:centered_coloring:base_case:excluding_an_apex_forest:chip'}.\ref{item:centered_coloring:base_case:excluding_an_apex_forest:chip'-3}, the vertex $u'_H = u_{H}$ is a  $(\phi \times \lambda_0)$-center of $V(H) \cap S_{C'_i}$ for every $i \in [m]$.
    Therefore, $u_H$ is a $(\phi \times \lambda)$-center of $V(H) \cap S$.
    Now suppose that $V(H)$ is disjoint from $U$.
    Assume that $V(H)$ intersects $\tilde{S}$.
    By \ref{item:centered_coloring:base_case:excluding_a_forest_small_interface:center_among_log_colors:call},
    the vertex $u'_{H}$ is a $(\phi \times \tilde{\lambda})$-center $u$ of $V(H) \cap S_{C'_i}$ for every $i \in [m]$,
    and so $u'_H$ is a $(\phi \times \lambda)$-center of $V(H) \cap S$.
    Finally, assume that $V(H)$ is disjoint from $\tilde{S} \cup U$.
    Then, there exists $C' \in \mathcal{C}$ such that $H \subseteq C'$.
    Therefore, by \ref{item:centered_coloring:base_case:excluding_an_apex_forest:chip'},
    there is a $(\phi \times \lambda_{C'})$-center $u$ of $V(H) \cap S_{C'}$.
    By the definition of $\lambda$, and because $V(H) \cap S_{C'} = V(H) \cap S$,
    we conclude that $u$ is a $(\phi \times \lambda)$-center of $V(H) \cap S$.
    This shows \ref{item:centered_coloring:base_case:excluding_an_apex_forest:chip} and concludes the proof of the lemma.
\end{proof}

\begin{lemma}\label{lemma:base_case_R2_cen_total}
    Let $c$ and $k$ be positive integers.
    There exists a constant $\alpha_{\ref{lemma:base_case_R2_cen_total}}(k)$ such that
    for every $X \in \Apex(\Rt_2)$,
    there exists a constant $c_{\ref{lemma:base_case_R2_cen_total}}(c,k,X)$ such that the following holds.
    For every connected $K_k$-minor-free graph $G$,
    for every integer $q$ with $q \geq 2$,
    for every $(q,c)$-good coloring $\phi$ of $G$,
    for every family $\mathcal{F}$ of connected subgraphs of $G$,
    if there is no $\mathcal{F}$-rich model of $X$ in $G$,
    then there exists $S \subseteq V(G)$
    such that
    \begin{enumerateOurAlph}
        \item $S \cap V(F) \neq \emptyset$ for every $F \in \mathcal{F}$; \label{item:lemma:base_case_R2_cen_total:hitting}
        \item for every component $C$ of $G-S$, $N_G(V(C))$ intersects at most $\alpha_{\ref{lemma:base_case_R2_cen_total}}(k)$ components of $G-V(C)$; and \label{item:lemma:base_case_R2_cen_total:components}
        \item $\cen_q(G,\phi,S) \leq c_{\ref{lemma:base_case_R2_cen_total}}(c,k,X) \cdot q \log q$. \label{item:lemma:base_case_R2_cen_total:chip}
    \end{enumerateOurAlph}
\end{lemma}

\begin{proof}
    We set $\alpha_{\ref{lemma:base_case_R2_cen_total}}(k) = 2^{k}$.
    Let $X \in \Apex(\Rt_2)$. 
    There exists $X' \in \Rt_2$ such that $X \subseteq K_1 \oplus X'$.
    Fix such a forest $X'$,
    and let $c_{\ref{lemma:base_case_R2_cen_total}}(c,k,X) = c_{\ref{lemma:centered_coloring:base_case:excluding_an_apex_forest}}(c,k,X')$.
    Let $q$ be an integer with $q \geq 2$,
    let $G$ be a connected $K_k$-minor-free graph,
    let $\mathcal{F}$ be a family of connected subgraphs of $G$
    such that there is no $\mathcal{F}$-rich model of $X$ in $G$.
    
    Let $\Lambda$ be a set of 
    $c_{\ref{lemma:centered_coloring:base_case:excluding_an_apex_forest}}(c,k,X') \cdot q \log q$ colors,
    let $y_0 \in \Lambda$,
    let $\Lambda_1 = \{y_0\}$ and $\Lambda_2 = \dots = \Lambda_q = \emptyset$.
    Let $u$ be a arbitrary vertex of $G$. 
    Let $U=\set{u}$, let $\lambda_0 \colon U \rightarrow \Lambda$ with $\lambda_0(u)=y_0$. 
    For every connected subgraph $H$ of $G$ intersecting $U$,
    let $u_H=u$.

    Let $\mathcal{C}$ be the family of all the components of $G-U$,
    and let $C \in \mathcal{C}$.
    By \Cref{lemma:centered_coloring:base_case:excluding_an_apex_forest}
    there exists $S_{C} \subseteq V(C) \cup U$ containing $U$ and a coloring $\lambda_{C} \colon S_{C} \to \Lambda$ extending $\lambda_0$ such that
    \begin{enumerate}[label=\ref{lemma:centered_coloring:base_case:excluding_an_apex_forest}.(\alph*)]
        \item $S_{C} \cap V(F) \neq \emptyset$ for every $F \in \mathcal{F}$; \label{item:centered_coloring:base_case:excluding_an_apex_forest:hitting:call}
        \item for every component $C'$ of $C-S_{C}$, $N_G(V(C'))$ intersects at most $2^k$ components of $G-V(C')$; and \label{item:centered_coloring:base_case:excluding_an_apex_forest:connectivity:call}
        \item for every connected subgraph $H$ of $G$ intersecting $V(C')$ and $S_{C}$, \label{item:centered_coloring:base_case:excluding_an_apex_forest:chip:call}
            either
            \begin{enumerate}[label=(\arabic*)]
                \item $|\phi(V(H))| > q$, or
                \item $|\lambda_{C}(V(H) \cap S)| > q$, or
                \item there is a $(\phi \times \lambda_C)$-center $v$ of $V(H) \cap S$ and moreover, if $V(H) \cap U \neq \emptyset$, then we may set $v=u$.
            \end{enumerate}
    \end{enumerate}
    Clearly, \ref{item:centered_coloring:base_case:excluding_an_apex_forest:hitting:call} and \ref{item:centered_coloring:base_case:excluding_an_apex_forest:connectivity:call} imply \ref{item:lemma:base_case_R2_cen_total:hitting}~and~\ref{item:lemma:base_case_R2_cen_total:components}, respectively.
    Moreover, the coloring $\lambda \colon \bigcup_{C \in \mathcal{C}} S_{C} \to \Lambda$
    defined by $\lambda(v) = \lambda_{C}(v)$ for every $C \in \mathcal{C}$ and $v \in S_{C}$,
    witnesses the fact that 
    $\cen_q(G,\phi,\bigcup_{C \in \mathcal{C}} S_{C}) 
    \leq c_{\ref{lemma:centered_coloring:base_case:excluding_an_apex_forest}}(c,k,X') \cdot q \log q = c_{\ref{lemma:base_case_R2_cen_total}}(c,k,X) \cdot q \log q$.
    This proves \ref{item:lemma:base_case_R2_cen_total:chip} and completes the proof.
\end{proof}

We are almost ready to wrap up the proof of \Cref{thm:centered}.
The general idea is exactly the same as for~\Cref{thm:main_cen_St}.
Recall the notation introduced at the beginning of this subsection.
We summarize the base case in the following statement.

\begin{corollary}\label{cor:base_case_cen}
     Let $\ell$ be a positive integer, let $\paramktmf_\ell = (\paramktmf_{\ell,q} \mid q \in \posint)$ and $\parambdtw = (\parambdtw_{\ell,q} \mid q \in \posint)$, and let $\param \in \{\paramktmf_\ell, \parambdtw\}$.
    The function $q\mapsto q\log (q+1)$ is $(\param,\Apex(\Rt_2))$-bounding for $\param$.
\end{corollary}

\begin{proof}
    Let $X \in \Apex(\Rt_2)$.
    Let $\calU$ be equal to $\calK_\ell$ when $\param = \paramktmf_\ell$ and to $\calT$ when $\param = \parambdtw$.
    Let $c = \max\{2,c_{\ref{lemma:Kt_free_graphs_have_centered_Helly_colorings}}$.
    For each graph $G$ for which it does not hold that $V(G) \subset V(\calU)$ and $E(G) \subset E(\calU)$, choosing the empty set as $S$ witnesses that $G$ is even $(q \mapsto 1, \param, X, 1, 0)$-good.
    Let $k$ be a positive integer and let $k'$ be such that all $K_k$-minor-free graphs $G$ have chromatic number at most $k'$ or in other words, $\cen_q(G) \leq k'$.
    Next, let 
        \[\beta(k,X) = \max\{k',c_{\ref{lemma:base_case_R2_cen_total}}(c,k,X)\}\]
    Let $G$ be a graph for which $V(G) \subset V(\calU)$ and $E(G) \subset E(\calU)$.
    We claim that $G$ is $(q\mapsto q\log (q+1),\param,X,\alpha_{\ref{lemma:base_case_R2_cen_total}}(k),\beta(k,X)$-good. 
    By~\Cref{obs:cen_components}, we may assume without loss of generality that $G$ is connected (see also the proof of~\Cref{cor:base_case_wcol}).
    Let $q$ be a positive integer and let $\calF$ be a family of connected subgraphs of $G$ such that $G$ has no $\calF$-rich model of $X$.
    Let $\Xi_q$ be equal to $\Phi_{\ell,q}$ when $\param = \paramktmf_\ell$ and to $\Psi_q$ when $\param = \parambdtw$.
    By assumption, $\Xi_q\vert_{V(G)}$ is $(q,c)$-good.
    Therefore, by~\Cref{lemma:base_case_R2_cen_total}, there exists $S \subset V(G)$ such that 
    \begin{enumerate}[label={\normalfont\ref{lemma:base_case_R2_cen_total}.(\makebox[\mywidth]{\alph*})}]
        \item $S \cap V(F) \neq \emptyset$ for every $F \in \mathcal{F}$; \label{item:lemma:base_case_R2_cen_total:hitting-cor-cor}
        \item for every component $C$ of $G-S$, $N_G(V(C))$ intersects at most $\alpha_{\ref{lemma:base_case_R2_cen_total}}(k)$ components of $G-V(C)$; and \label{item:lemma:base_case_R2_cen_total:components-cor-cor}
        \item $\cen_q(G,\Xi_q\vert_{V(G)},S) \leq c_{\ref{lemma:base_case_R2_cen_total}}(c,k,X) \cdot q \log (q+1)$. \label{item:lemma:base_case_R2_cen_total:chip-cor-cor}
    \end{enumerate}
    This directly implies
            \begin{enumerate}[label={\normalfont (g\arabic*$^\star$)}]
            \item $S \cap V(F) \neq \emptyset$ for every $F \in \mathcal{F}$; \label{item:lemma:fragility_rate_base_case:hit-cor-cor}
            \item for every component $C$ of $G-S$, $N_{G}(V(C))$ intersects at most $\alpha_{\ref{lemma:base_case_R2_cen_total}}(k)$ component of $G-V(C)$; and \label{item:lemma:fragility_rate_base_case:components-cor-cor}
            \item $\param_q(G,S) \leq \beta(k,X) \cdot q\log (q+1)$. \label{item:lemma:fragility_rate_base_case:param-cor-cor}
        \end{enumerate}
    Thus, $G$ is indeed $(q\mapsto q\log (q+1),\param,X,\alpha_{\ref{lemma:base_case_R2_cen_total}}(k),\beta(k,X)$-good, and so, $q\mapsto q\log (q+1)$ is $(\param,\Apex(\Rt_2))$-bounding for $\param$.
\end{proof}

\begin{theorem}\label{thm:main_cen_Rt}
    Let $t$ be an integer with $t \geq 3$, and let $X \in \Rt_t$.
    There is a constant $c$ such that
    for every $X$-minor-free graph $G$,
    for every integer $q$ with $q \geq 2$,
    \begin{align*}
        \cen_q(G) &\leq c \cdot (\tw(G)+1) \cdot q^{t-2} \log q, \\
    \intertext{and}
        \cen_q(G) &\leq c \cdot q^{t-1} \log q.
    \end{align*}
\end{theorem}
\begin{proof}
    Let $\ell$ be a positive integer, let $\paramktmf_\ell = (\paramktmf_{\ell,q} \mid q \in \posint)$ and $\parambdtw = (\parambdtw_{\ell,q} \mid q \in \posint)$, and let $\param \in \{\paramktmf_\ell, \parambdtw\}$.
    By~\Cref{lemma:centered_colorings_are_nice}, the family of focused parameters $\param$ is nice.
    We show by induction on $t$ (integers with $t \geq 3$) that $q \mapsto q^{t-2}\log (q+1)$ is $(\param, \Rt_t)$-bounding.
    Since $t\geq 3$, $\Rt_{t-1}$ is closed under disjoint union and leaf addition, and by \Cref{lemma:Rt_has_coloring_elimination_property}, $\Rt_{t-1}$ has the coloring elimination property.
    By \Cref{cor:base_case_cen},
    the function $q \mapsto q \log(q+1)$ is $(\param,\Apex(\Rt_2))$-bounding,
    and so is $(\param, \Rt_3)$-bounding by \Cref{thm:main_A(X)_to_T(X)}.
    This gives the statement for $t = 3$, so suppose that $t \geq 4$.
    By induction, $q \mapsto q^{t-3}\log (q+1)$ is $(\param, \Rt_{t-1})$-bounding.
    Therefore, by \Cref{thm:abstract_induction_main}, 
    $q \mapsto q^{t-2} \log (q+1)$ is $(\param, \Rt_{t})$-bounding.
    This completes the inductive proof that for every integer $t$ with $t\geq 3$, $q \mapsto q^{t-2}\log (q+1)$ is $(\param, \Rt_{t})$-bounding, as desired.
    By~\Cref{lemma:par-bounding-to-bound}, this implies that for every $X \in \Rt_t$, there exists a constant $\beta_\param(X)$ such that for every positive integer $q$ and for every graph $G$, we have
    \begin{equation}
        \param_q(G,V(G)) \leq \beta_\param(X) \cdot q^{t-2} \log (q+1).\label{eq:param-bound-Rt}
    \end{equation}

    Let $X \in \Rt_t$, 
    let $q$ be an integer with $q \geq 2$, 
    and let $G$ be an $X$-minor-free graph.
    In particular, $G$ is $K_{|V(X)|}$-minor-free, hence, there exists $i \in \posint$ such that $G$ is isomorphic to $G_{|V(X)|,i}$, which is one of the building blocks of $\calK_{|V(X)|}$.
    Recall that $\Phi_{|V(X)|,q}\vert_{V(G_{|V(X)|,i})} = \phi_{|V(X)|,i,q}$ uses at most $q+1$ colors.
    It follows that with $\ell = |V(X)|$, we have
    \begin{align*}
        \cen_q(G) = \cen_q(G_{\ell,i}) &\leq \cen_q(G_{\ell,i},\Phi_{\ell,q}\vert_{V(G_{\ell,i})},V(G_{\ell,i})) \cdot (q+1) && \text{by~\Cref{lem:ordered-to-normal}}\\
        &= \paramktmf_{\ell}(G_{\ell,i},V(G_{\ell,i})) \cdot (q+1)\\
        &\leq \beta_{\paramktmf_\ell}(X) \cdot q^{t-2} \log(q+1) \cdot (q+1)  && \text{by~\eqref{eq:param-bound-Rt},}\\
        &\leq (4\beta_{\paramktmf_\ell}(X)) \cdot q^{t-1} \log q.
    \end{align*}
    This gives the second bound in the statement of the theorem.
    For the first bound, note that there exists $i \in \posint$ such that $G$ is isomorphic to $H_{\tw(G)+1,i}$, a building block of $\calT_{\tw(G)+1}$.
    Recall that $\Psi_{q}\vert_{V(H_{\tw(G)+1,i})} = \psi_{\tw(G)+1,i,q}$ uses at most $\tw(G)+1$ colors.
    It follows that for $\ell = \tw(G)+1$, we have
    \begin{align*}
        \cen_q(G) = \cen_q(H_{\ell,i}) &\leq \cen_q(H_{\ell,i},\Psi_{\ell,q}\vert_{V(H_{\ell,i})},V(H_{\ell,i})) \cdot \ell && \text{by~\Cref{lem:ordered-to-normal}}\\
        &= \parambdtw(H_{\ell,i},V(H_{\ell,i})) \cdot \ell\\
        &\leq \beta_{\parambdtw}(X) \cdot q^{t-2} \log(q+1) \cdot \ell && \text{by~\eqref{eq:param-bound-Rt}.}\\
        &= 2\beta_{\parambdtw}(X) \cdot (\tw(G)+1) \cdot q^{t-2} \log q. && \qedhere 
    \end{align*}
\end{proof}

\bibliographystyle{plain}
\bibliography{biblio}

\clearpage

\appendix

\section{Lower bounds}
\label{sec:lower-bounds}

\subsection{Weak coloring numbers} \label{ssec:wcols}

In this section, we present a construction by Grohe, Kreutzer, and Siebertz~\cite{Grohe15}, which was initially presented as a lower bound for weak coloring numbers of graph of treewidth $t$.
This was later adapted to graph of simple treewidth at most $t$ by Joret and Micek~\cite{JM22}.
To show the link between these constructions and the classes $(\Rt_t)_{t \in \posint}, (\SRt_t)_{t \in \posint}$,
we elected to provide them in the current appendix.
The core of the argument is the following lemma.

\begin{lemma}\label{lemma:lower_bound:apply_T}
    Let $\mathcal{G}$ be a nonempty class of graphs, 
    and let
    $g(q) = \max_{G \in \mathcal{G}} \wcol_q(G)$ 
    for every nonnegative integer $q$.
    For every nonnegative integer $q$, there exists $G \in \Tree(\mathcal{G})$
    such that
    \[
        \wcol_q(G) \geq \sum_{i=0}^q g(i).
    \]
\end{lemma}

\begin{proof}
    We proceed by induction on $q$.
    First, suppose that $q=0$. 
    Then consider $G \in \mathcal{G}$ such that $\wcol_q(G) = g(0)$.
    Since $\mathcal{G} \subseteq \Tree(\mathcal{G})$,
    $G$ is as desired.
    
    Next, suppose that $q \geq 1$. 
    By the inductive hypothesis there exists $G_0 \in \Tree(\mathcal{G})$ with
    \[
        \wcol_{q-1}(G_0) \geq \sum_{i=0}^{q-1} g(i).
    \]
    Let $H \in \mathcal{G}$ be such that $\wcol_q(H) = g(q)$.
    Let $k = \sum_{i=0}^{q} g(i)$. 
    Let $G$ be a graph obtained from $G_0$ by adding for each $u \in V(G_0)$, 
    $k$ disjoint copies $H_{u,1}, \dots, H_{u,k}$ of $H$,
    and all possible edges between $u$ and $V(H_{u,i})$,
    for every $i \in [k]$.

    We claim that $G \in \Tree(\mathcal{G})$.
    Let $\big(T_0,(W_{0,x} \mid x \in V(T_0))\big)$
    be a rooted forest decomposition of $G_0$ witnessing the fact that
    $G_0 \in \Tree(\mathcal{G})$.
    Then, for every $u \in V(G_0)$, choose $x_u \in V(T_0)$ such that $u \in W_{0,x_u}$.
    Let $T$ be the rooted forest obtained from $T_0$ by adding,
    for every $u \in V(G_0)$ and for every $i \in [N]$,
    a leaf $x_{u,i}$ with parent $x_u$.
    Then, for every $x \in V(T)$, let
    \[
      W_x = 
      \begin{cases}
          W_{0,x} & \textrm{if $x \in V(T_0)$} \\
          \{u\} \cup V(H_{u,i}) &\textrm{if $x = x_{u,i}$ for $u \in V(G_0)$, $i \in [k]$.}
      \end{cases}
    \]
    It follows that $\big(T,(W_x \mid x \in V(T))\big)$ is a rooted forest
    decomposition of $G$ witnessing the fact that $G \in \Tree(\mathcal{G})$.

    It remains to show that $\wcol_q(G) \geq k$.
    Let $\sigma$ be an ordering of $V(G)$.
    Since $\wcol_{q-1}(G_0) \geq \sum_{i=0}^{q-1} g(i)$,
    there exists $v \in V(G_0)$ such that 
    \[
        \left|\WReach_{q-1}\left[G_0,\sigma\vert_{V(G_0)},v\right]\right| \geq \sum_{i=0}^{q-1} g(i).
    \]
    If, for every $i \in [k]$, there is a vertex $w$ in $V(H_{v,i})$ 
    such that $w <_\sigma v$, then
    $|\WReach_q[G,\sigma,v]| \geq k$.
    Otherwise, there exists $i \in [k]$ such that
    every vertex $w$ in $H_{v,i}$ is such that $w >_\sigma v$.
    Since $\wcol_q(H_{v,i}) = \wcol_q(H) = g(q)$, there exists
    $u \in V(H_{v,i})$ such that $\left|\WReach_q\left[H_{v,i},\sigma\vert_{V(H_{v,i})}, u\right]\right| \geq g(q)$.
    Then, 
    \[
        \WReach_{q-1}\left[G_0,\sigma\vert_{V(G_0)},v\right] \cup \WReach_q\left[H_{v,i},\sigma\vert_{V(H_{v,i})},u\right]
        \subseteq \WReach_{q}\left[G,\sigma\vert_{V(G)},u\right],
    \]
    and we conclude that $|\WReach_r[G,\sigma,u]| \geq \left(\sum_{i=0}^{q-1} g(i)\right) + g(q) = k$.
\end{proof}

\begin{corollary}\label{cor:lower-bound:wcol:rtd}
    For every integer $t$ with $t\geq1$,
    \begin{align*}
        \max_{G \in \Rt_t} \wcol_q(G) &= \Omega(q^{t-1}).
    \end{align*}
\end{corollary}

\begin{proof}
We show by induction on $t$ that
for every nonnegative integer $q$,
\[
    \max_{G \in \Rt_t} \wcol_q(G) \geq \binom{q+t-1}{t-1}.
\]
For $t=1$, recall that $K_1\in\calR_1$. Therefore, 
for all nonnegative integers $q$ we have $\wcol_q(K_1)=1 \geq \binom{q}{0}$, as desired.

For $t \geq 2$, by the induction hypothesis, 
for every nonnegative integer $q$,
\[
    \max_{G \in \Rt_{t-1}} \wcol_q(G) \geq \binom{q+t-2}{t-2}.
\]
Therefore, by \Cref{lemma:lower_bound:apply_T},
\begin{align*}
    \max_{G \in \Rt_t} \wcol_q(G)
    &\geq \sum_{i=0}^q \binom{i+t-2}{t-2} 
    = \binom{q+t-1}{t-1}. \qedhere
\end{align*}
\end{proof}

\begin{corollary}\label{cor:lower_bound_SRt_cen}\label{cor:lower-bound:wcol:srtd}
    For every positive integer $t$ with $t \geq 2$,
    \begin{align*}
        \max_{G \in \SRt_t} \wcol_q(G) &= \Omega(q^{t-2} \log q).
    \end{align*}
\end{corollary}

\begin{proof}
    We proceed by induction on $t$.
    For $t=2$,
    $\SRt_2$ contains every path, which have $q$th weak coloring number in $\Omega(\log q)$, see~\cite{JM22}.

    For $t \geq 3$, by the induction hypothesis,
    there exists $c > 0$ such that
    for every integer $q$ with $q \geq 2$,
    \[
    \max_{G \in \SRt_{t-1}} \wcol_q(G) \geq c \cdot q^{t-3} \log q.
    \]
    Therefore, by \Cref{lemma:lower_bound:apply_T}, for every integer $q$ with $q \geq 4$,
    \begin{align*}
        \max_{G \in \SRt_{t}} \wcol_q(G)
        &\geq c \cdot \sum_{i=1}^q i^{t-3} \log i \\
        &\geq c \cdot \sum_{q/2 \leq i \leq q} i^{t-3} \log i \\
        &\geq c \cdot \left\lfloor \frac{q-1}{2} \right\rfloor (q/2)^{t-3} \log(q/2) \\
        &\geq c \cdot \frac{q}{4} \cdot \frac{1}{2^{t-3}} \cdot q^{t-3} \cdot \frac{\log q}{2} \\
        &=    \frac{c}{2^t} \cdot q^{t-2} \log q. \qedhere
    \end{align*}
\end{proof}

\subsection{Centered chromatic numbers}

In this section, we prove lower bounds for centered chromatic numbers 
of graphs in $\Rt_t$ and $\SRt_t$.
This is a rewriting with our notations of the proof published in~\cite{Dbski2021}.

Let $q$ be a positive integer.
Let $G$ be a graph.
Let $K$ be a set of colors.
Let $\phi \colon V(G) \to K$ be a coloring
and $\Xi \colon V(G) \to 2^{2^K}$ be a function that maps each vertex of $G$ to a family of set of colors in $K$.
The pair $(\phi,\Xi)$ is a \defin{generalized $q$-centered coloring} if
for every connected subgraph $H$ of $G$, 
for every function $\xi \in \prod_{u \in V(H)} \Xi(u)$\footnote{Given a set $X$ and a family $(Y_x \mid x \in X)$ of sets, we denote by \defin{$\prod_{x \in X} Y_x$} the set of all the functions $f \colon X \to \bigcup_{x \in X} Y_x$ such that $f(x) \in Y_x$ for every $x \in X$.}
that maps each vertex $u\in V(H)$ to a set of colors in $\Xi(u)$,
either
\begin{enumerate}
    \item $|\phi(V(H)) \cup \bigcup \xi(V(H))| > q$ or
    \item there is a vertex $u \in V(H)$
    such that 
    $\phi(u) \not\in \phi(V(H)\setminus \{u\}) \cup \bigcup\xi(V(H))$.
    \label{item:lower_bound:cen:generalized_coloring:center}
\end{enumerate}
The number of colors used by $(\phi,\Xi)$ is the integer $|\phi(V(G)) \cup \bigcup_{u \in V(G)}\bigcup \Xi(u)|$.
Note that if $(\phi,u\mapsto \{\emptyset\})$ is a generalized $q$-centered coloring of $G$, then $\phi$ is a $q$-centered coloring of $G$.

Let $p$ be a nonnegative integer.
A mapping $\Xi \colon V(G) \to 2^{2^K}$ is said to be \defin{$p$-good} if
for every $u \in V(G)$,
$\emptyset \in \Xi(u)$, and
for every subset of colors $L \subseteq K$,
if there exists a set of colors $M \in \Xi(u)$ such that $L \subseteq M$,
then there exists a set of colors $M' \in \Xi(u)$ such that $L\subseteq M'$ and $|M' \setminus L| \leq p$.
Note that a $p$-good mapping is also a $(p+1)$-good mapping and every image of a $0$-good mapping is a family of sets that is closed under taking subsets.

Now, we define \defin{$\cen^{(k)}_{q,p}(G)$}
to be the minimum integer $\ell$ such that
there exists a generalized $q$-centered coloring $(\phi,\Xi)$
of $G$ using at most $k$ colors such that 
\begin{enumerate}
    \item $\Xi$ is $p$-good, and
    \item for every $u \in V(G)$, there exists $L \in \Xi(u)$ such that
    $|\phi(V(G)) \setminus L| \leq \ell$.
\end{enumerate}
By convention, we set $\cen^{(k)}_{q,p}(G) = \infty$
if there is no such integer $\ell$.
Note that for every nonnegative integer $p$, $\cen^{(k)}_{q,p+1}(G)\leq \cen^{(k)}_{q,p}(G)$.
Moreover,
if $\phi$ is a $q$-centered coloring of $G$,
then $(\phi, u \mapsto \{\emptyset\})$ is a generalized $q$-centered coloring of $G$, and $u \mapsto \{\emptyset\}$ is $0$-good.
Therefore,
for every positive integer $k$,
\[
\cen_q(G)\geq  \min \{k,\cen^{(k)}_{q,0}(G)\}.
\]

\begin{lemma}\label{lemma:lower_bound:cen:path}
    Let $k$ and $q$ be positive integers.
    For every positive integer $n$ with $n \leq 2q$,
    \[
        \cen_{4q,q}^{(k)}(P_n) \geq \lfloor \log n \rfloor.
    \]
\end{lemma}

\begin{proof}    
    It is enough to consider the case $n=2^\ell$ for some nonnegative integer $\ell$.
    Let $(\phi,\Xi)$ be a generalized $4q$-centered coloring of $P_n$ 
    using at most $k$ colors
    with $\Xi$ $q$-good,
    and let $N$ be the set of the colors of all the $\phi$-centers of $V(P_n)$, i.e.\ $N=\{\alpha\mid \alpha\in \phi(V(P_n)),|\phi^{-1}(\alpha)|=1\}$. 

    We will prove by induction on $\ell$ that
    there exists $u \in V(P_n)$ such that for every $L \in \Xi(u)$,
    \[
        |N \setminus L| \geq \ell.
    \]
    This is clear when $\ell=0$.
    Next, suppose that $\ell \geq 1$.

    Let $Q_1,Q_2$ be the subpaths of $P_n$
    induced by, respectively, the first $n/2 = 2^{\ell-1}$ vertices
    and last $n/2 = 2^{\ell-1}$ vertices of $P_n$.
    Let $N_1 = N \cap \phi(V(Q_1))$ and $N_2 = N \cap \phi(V(Q_2))$.
    By definition of $N$, $N_1$ and $N_2$ are disjoint.

    Suppose there exists $u_1 \in V(Q_1)$ and $L_1 \in \Xi(u_1)$ such that $N_2 \subseteq L_1$,
    and $u_2 \in V(Q_2)$ and $L_2 \in \Xi(u_2)$ such that $N_1 \subseteq L_2$.
    Since $\Xi$ is $q$-good, there exists
    $M_1 \in \Xi(u_1)$ and $M_2 \in \Xi(u_2)$
    such that $N_2 \subseteq M_1$, $N_1 \subseteq M_2$, $|M_1 \setminus N_2|\leq q$, and $|M_2 \setminus N_1| \leq q$.
    Let $\xi \in \prod_{u \in V(P_n)} \Xi(u)$
    be defined by
    \[
        \xi(u) =
        \begin{cases}
            M_1 &\textrm{if $u=u_1$,} \\
            M_2 &\textrm{if $u=u_2$,} \\
            \emptyset &\textrm{otherwise.} \\
        \end{cases}
    \]
    Since $(\phi,\Xi)$ is $4q$-centered, either
    \begin{enumerate}
        \item $|\phi(V(P_n)) \cup \bigcup \xi(V(P_n))| > 4q$,
            but then $|V(P_n)| + |M_1 \setminus N| + |M_2 \setminus N| > 4q$,
            which implies $n > 4q - |M_1\setminus N_2| - |M_2\setminus N_1| = 2q$, a contradiction; or
        \item there is a vertex $u$ of $V(P_n)$ such that $\phi(u)\notin\phi(V(P_n)-\{u\})\cup\bigcup\xi(V(P_n))$.
            In particular, $u$ must be a $\phi$-center of $P_n$, so 
            $\phi(u) \in N$. 
            However, $\phi(u)\in N = N_1 \cup N_2 \subseteq M_2 \cup M_1 \subseteq \bigcup \xi(V(P_n))$, a contradiction.
    \end{enumerate}
    This proves that there exists $i \in \{1,2\}$
    such that
    for every $u \in V(Q_i)$, for every $L \in \Xi(u)$,
    $N_{3-i} \setminus L \neq \emptyset$.
    We now call the induction hypothesis on $Q_i$,
    $\phi \vert_{V(Q_i)}, \Xi\vert_{V(Q_i)}$.
    Therefore, there exists $u \in V(Q_i)$
    such that for every $L \in \Xi(u)$,
    \[
        |N_i \setminus L| \geq \ell-1.
    \]
    Since $N_i$ and $N_{3-i}$ are disjoint,
    we deduce that
    for every $L \in \Xi(u)$,
    \[
        |N \setminus L| \geq |N_i \setminus L| + |N_{3-i} \setminus L| \geq \ell.
    \]
    This concludes the proof of the induction.
    
    Since $\Xi$ is $q$-good and there exists $u\in V(P_n)$ such that for every $L\in\Xi(u)$, $|\phi(V(P_n))\setminus L| \geq |N\setminus L|\geq \ell = \lfloor\log n\rfloor$, this also concludes the proof of the lemma. 
\end{proof}

\begin{lemma}\label{lemma:lower_bound:cen:induction_step}
    Let $k$ and $q$ be positive integers,
    let $\mathcal{G}$ be a nonempty class of graphs,
    and let $g(p)= \max_{G \in \mathcal{G}} \cen^{(k)}_{q,p}(G)$ for every nonnegative integer $p$.
    For every integer $p$ with $0 \leq p \leq q$,
    there exists $G \in \Tree(\mathcal{G})$ such that
    \[
        \cen^{(k)}_{q,p}(G) \geq \sum_{r=p}^{q} g(r).
    \]
\end{lemma}

\begin{proof}
    We proceed by induction on $q-p$.
    When $q=p$, consider $G \in \mathcal{G}$
    such that $\cen^{(k)}_{q,p}(G) = g(q)$.
    Since $\mathcal{G} \subseteq \Tree(\mathcal{G})$,
    $G$ is as desired.
    Next, suppose that $q - p \geq 1$.

    By the induction hypothesis,
    there exists $G_0 \in \Tree(\mathcal{G})$
    such that
    \[
        \cen^{(k)}_{q,p+1}(G_0) \geq \sum_{r=p+1}^{q} g(r).
    \]
    Moreover, by the definition of $g$,
    there exists $G' \in \mathcal{G}$
    such that
    \[
        \cen^{(k)}_{q,p}(G') = g(p).
    \]
    Let $G$ be the graph obtained from $G_0$ as follows.
    For every $u \in V(G_0)$,
    consider disjoint copies $G'_{u,i}$ of $G'$ for $i \in [2^k+1]$,
    and add every possible edge between $u$ and $\bigcup_{i = 1}^{2^k+1} V(G'_{u,i})$.

    First, we show that $G \in \Tree(\mathcal{G})$.
    Let $\big(T_0,(W_{0,x} \mid x \in V(T_0))\big)$
    be a rooted forest decomposition of $G_0$ witnessing the fact that
    $G_0 \in \Tree(\mathcal{G})$.
    Then, for every $u \in V(G_0)$, choose $x_u \in V(T_0)$ such that $u \in W_{0,x_u}$.
    Let $T$ be the rooted forest obtained from $T_0$ by adding,
    for every $u \in V(G_0)$ and for every $i \in [2^k+1]$,
    a leaf $x_{u,i}$ with parent a vertex $x_u$.
    Then, for every $x \in V(T)$, let
    \[
      W_x = 
      \begin{cases}
          W_{0,x} & \textrm{if $x \in V(T_0)$} \\
          \{u\} \cup V(G'_{u,i}) &\textrm{if $x = x_{u,i}$ for $u \in V(G_0)$, $i \in [2^k+1]$.}
      \end{cases}
    \]
    It follows that $\big(T,(W_a \mid a \in V(T))\big)$ is a rooted forest
    decomposition of $G$ witnessing the fact that $G \in \Tree(\mathcal{G})$.

    We now show that $\cen^{(k)}_{q,p}(G) \geq \sum_{r=p}^q g(r)$.
    
    Let $(\phi,\Xi)$ be a generalized $q$-centered coloring of $G$ 
    using at most $k$ colors and with $\Xi$ $p$-good.
    Let $u \in V(G_0)$.
    By the pigeonhole principle,
    there exists $i,j \in [2^k+1]$ such that $\phi(V(G'_{u,i})) = \phi(V(G'_{u,j}))$.
    Without loss of generality, assume $i=1$ and $j=2$.
    Since $\cen^{(k)}_{q,p}(G'_{u,1}) = g(p)$,
    and because $(\phi \vert_{V(G'_{u,1})},\Xi \vert_{V(G'_{u,1})})$
    is a generalized $q$-centered coloring of $G'_{u,1}$,
    we define 
    \[
    v(u) \in V(G'_{u,1}) \textrm{ such that for every } M\in\Xi(v(u)), |\phi(V(G'_{u,1})) \setminus M| \geq g(p).
    \]
    Let 
    \begin{align*}
        \Xi_0(u) &=
        \{\emptyset\} \cup \big\{L \cup M \mid \phi(v(u)) \in L \subseteq \phi(V(G'_{u,1})), M \in \Xi(v(u))\big\}\textrm{ and}\\
        \phi_0 &= \phi \vert_{V(G_0)}.
    \end{align*}

    \begin{claim}\label{claim:lower_bound_cen}
    $(\phi_0, \Xi_0)$ is a generalized $q$-centered coloring of $G_0$.
    \end{claim}

    \begin{proofclaim}
    Let $H_0$ be a connected subgraph of $G_0$ and
    let $\xi_0 \in \prod_{u \in V(H_0)} \Xi_0(u)$. 
    We need to prove that either
    \begin{enumerateOurAlph}
        \item $|\phi_0(V(H_0)) \cup \bigcup \xi_0(V(H_0))| > q$ or\label{item:106a}
    \item there is a vertex $u \in V(H_0)$
    such that 
    $\phi_0(u) \not\in \phi_0(V(H_0)\setminus \{u\}) \cup \bigcup\xi_0(V(H_0))$.
    \label{item:106b}
    \end{enumerateOurAlph}
    Consider a subgraph $H$ of $G$ obtained from $H_0$ as follows.
    Let $u \in V(H_0)$.
    If $\xi_0(u) = \emptyset$, then let $L(u) = M(u) = \emptyset$.
    If $\xi_0(u) \neq \emptyset$, then fix 
    $L(u) \subseteq \phi(V(G'_{u,1}))$ with $\phi(v(u)) \in L(u)$
    and $M(u) \in \Xi(v(u))$ such that 
    \[
        \xi_0(u) = L(u) \cup M(u).
    \]
    Since $\phi(V(G'_{u,1})) = \phi(V(G'_{u,2}))$,
    there exists $V_1 \subseteq V(G'_{u,1})$ and $V_2 \subseteq V(G'_{u,2})$
    such that $\phi(V_1) = \phi(V_2) = L(u)$, and $v(u) \in V_1$.
    Then, add the vertices in $V_1 \cup V_2$ to $H$ and for every $w\in V_1\cup V_2$, add the edge $uw$.

    We define $\xi \in \prod_{x \in V(H)} \Xi(x)$ as follows.
    Let
    \[
        \xi(x) = 
        \begin{cases}
            \emptyset & \textrm{if $x\in V(H_0)$,}\\
            M(u) &\textrm{if $x \notin V(H_0)$ and $x=v(u)$ for some $u\in V(H_0)$,} \\
            \emptyset &\textrm{otherwise.}
        \end{cases}
    \]
    This completes the definition of $H$ and $\xi$.
    Since $H_0$ is connected and all the added vertices are neighbors of vertices in $H_0$, we conclude that $H$ is connected as well. 

    Observe that
    \begin{align*}
        \phi(V(H)\setminus V(H_0)) & = \bigcup\{L(u)\mid u\in V(H_0)\},\\
        \bigcup\xi(V(H)) &= \bigcup\{M(u)\mid u\in V(H_0)\}, \textrm{ and}\\
        \bigcup\xi_0(V(H_0)) &= \bigcup\{L(u)\cup M(u)\mid u\in V(H_0)\}.
        \intertext{These three equations imply that}
        \phi_0(V(H_0)) \cup \bigcup \xi_0(V(H_0)) &= \phi(V(H)) \cup \bigcup \xi(V(H)).
    \end{align*}
    Since $(\phi,\Xi)$ is a generalized $q$-centered coloring of $G$,
    and because $H$ is a connected subgraph of $G$,
    we have either
    \begin{enumerate}
        \item $|\phi(V(H)) \cup \bigcup \xi(V(H))|  > q$; or
        \item there is a vertex $u\in V(H)$ such that $\phi(u)\notin\phi(V(H)\setminus\{u\})\cup\bigcup\xi(V(H))$. \label{item:claim-lower-bound-cen-ii}
    \end{enumerate}
    In the former case, we have
    $|\phi_0(V(H_0)) \cup \bigcup \xi_0(V(H_0))| = |\phi(V(H)) \cup \bigcup \xi(V(H))| > q$, so~\ref{item:106a} holds. 
    In the latter case, fix such a vertex $u \in V(H)$. 
    We now argue that $\phi_0(u) \not\in \phi_0(V(H_0) \setminus \{u\}) \cup \bigcup \xi_0(V(H_0))$.
    Since the color of every vertex in $V(H) \setminus V(H_0)$ is repeated at least twice under $\phi$, we conclude that $u \in V(H_0)$. 
    Since $u\in V(H_0)$,
    \begin{align*}
    \phi(V(H)\setminus\{u\}) &= \phi_0(V(H_0)\setminus\{u\})\cup \phi(V(H)\setminus V(H_0))\\
         &= \phi_0(V(H_0)\setminus\{u\})\cup \bigcup\{L(v)\mid v\in V(H_0)\}.\\
    \intertext{Thus, we have }
    \phi(V(H)\setminus\{u\})\cup&\bigcup \xi(V(H)) \\
        &= \phi_0(V(H_0)\setminus\{u\})\cup \bigcup\{L(v)\mid v\in V(H_0)\}\cup \bigcup \{M(v)\mid v\in V(H_0)\}\\
        &= \phi_0(V(H_0)\setminus\{u\})\cup\xi_0(V(H_0)).
    \end{align*}
    Finally, using \ref{item:claim-lower-bound-cen-ii},
    we conclude that 
    \[
        \phi_0(u)=\phi(u)\notin \phi(V(H)\setminus\{u\})\cup\bigcup\xi(V(H)) = \phi_0(V(H_0)\setminus\{u\})\cup\xi_0(V(H_0)),
    \]
    and so \ref{item:106b} holds.
    This proves that $(\phi_0, \Xi_0)$ is a generalized $q$-centered coloring of $G_0$.\qedhere
    \end{proofclaim}

    \begin{claim}\label{claim:p-good}
        $\Xi_0$ is $(p+1)$-good.
    \end{claim}

    \begin{proofclaim}
    For every set of colors $N$,
    if $N \subseteq L \cup M$ for some $L \subseteq \phi(V(G'_{u,1}))$ with $\phi(v(u)) \in L$ and $M \in \Xi(v(u))$,
    then, because $\Xi$ is $p$-good, there exists $M' \in \Xi(v(u))$
    such that $N \cap M \subseteq M'$ and $|M' \setminus (N \cap M)| \leq p$,
    and it follows that $M_0 = (L \cap N) \cup \{\phi(v(u))\} \cup M'$
    is a member of $\Xi_0(v(u))$
    such that $N \subseteq M_0$ and $|M_0 \setminus N| \leq p+1$.
    This proves the claim.
    \end{proofclaim}

    Since $\cen^{(k)}_{q,p+1}(G_0) \geq \sum_{r=p+1}^q g(r)$ and by Claims \ref{claim:lower_bound_cen} and \ref{claim:p-good},
    we deduce that there exists $u \in V(G_0)$ such that for every $M_0 \in \Xi_0(u)$,
    \[
        |\phi(V(G_0)) \setminus M_0| \geq \sum_{r=p+1}^q g(r).
    \]
    Therefore,
    for every $M \in \Xi(v(u))$,
    apply this inequality for $M_0 = \phi(V(G'_{u,1})) \cup M$, we deduce that
    \begin{align*}
        |\phi(V(G_0)) \setminus (\phi(V(G'_{u,1})) \cup M)| 
        &\geq \sum_{r=p+1}^q g(r). \\
    \intertext{Moreover, by the definition of $v(u)$,}
        |\phi(V(G'_{u,1})) \setminus M| 
        &\geq g(p)
    \end{align*}
    Since the set of colors 
    $\phi(V(G_0)) \setminus (\phi(V(G'_{u,1})) \cup M)$
    and $\phi(V(G'_{u,1})) \setminus M$ are disjoint,
    we deduce
    \[
        |\phi(V(G)) \setminus M| \geq g(p) + \sum_{r=p+1}^q g(r) 
        = \sum_{r=p}^q g(r).
    \]
    This proves that 
    $\cen^{(k)}_{q,p}(G) \geq \sum_{r=p}^q g(r)$ and concludes the proof of the lemma.
\end{proof}

\begin{corollary}\label{cor:lower-bound:cen:rtd}
    For every positive integer $t$,
    \begin{align*}
        \max_{G \in \Rt_t} \cen_q(G) &= \Omega(q^{t-1}).
    \end{align*}
\end{corollary}

\begin{proof}
    Let $q$ and $t$ be positive integers.
    Let $k = \binom{q-1+t-1}{t-1}-1$.
    
    \begin{claim}
    For every positive integer $s$,
    for every nonnegative integer $p$ with $p < q$,
    \[
        \max_{G \in \Rt_s} \cen^{(k)}_{q,p}(G) \geq \binom{q-1 +s-1-p}{s-1}.
    \]
    \end{claim}

    \begin{proofclaim}
    We proceed by induction on $s$.
    
    First suppose $s=1$.
    Note that $\binom{q-1 + s-1-p}{s-1} = 1$.
    Consider $K_1 \in \Rt_1$,
    and let $(\phi, \Xi)$ be a generalized $q$-centered coloring of $K_1$ such that $\Xi$ is $p$-good.
    We denote by $u$ the unique vertex of $K_1$.
    Then, consider $L \in \Xi(u)$. 
    If $\phi(u) \in L$,
    then, since $\Xi$ is $p$-good, there exists $M \in \Xi(u)$ such that $\phi(u) \in M$ and $|M| \leq p+1 \leq q$.
    But then,
    the subgraph $K_1$ equipped with $\xi \colon u \mapsto \{M\}$ uses at most $q$ color and $\phi(u)\in \xi(\{u\})$.
    This contradicts the fact that $(\phi,\Xi)$ is a generalized $q$-centered coloring of $K_1$.
    Therefore, $\phi(u) \not\in L$, and so $|\phi(V(K_1)) \setminus L| \geq 1$.
    This proves that $\cen^{(k)}_{q,p}(K_1) \geq 1$ and concludes the case $s=1$.
    
    Now, suppose that $s>1$.
    By the fact that $\Rt_{s} = \Tree(\Rt_{s-1})$,
    \Cref{lemma:lower_bound:cen:induction_step}, 
    and the induction hypothesis,
    we deduce that
    \begin{align*}
        \max_{G \in \Rt_s} \cen^{(k)}_{q,p}(G)
        &\geq \sum_{r=p}^q \max_{G \in \Rt_{s-1}} \cen^{(k)}_{q,r}(G) \\
        &\geq \sum_{r=p}^{q-1} \max_{G \in \Rt_{s-1}} \cen^{(k)}_{q,r}(G) \\
        &\geq \sum_{r=p}^{q-1} \binom{q-1+s-2-r}{s-2} 
        = \binom{q-1+s-1-p}{s-1}.
    \end{align*}
    This proves the induction.
    \end{proofclaim}

    By the previous claim, we consider a graph $G \in \Rt_t$
    such that $\cen^{(k)}_{q,0}(G) \geq \binom{q-1+t-1}{t-1}$.
    Suppose by contradiction that $G$ admits a $q$-centered coloring $\phi$
    using at most $k$ colors.
    Then $(\phi, v \mapsto \{\emptyset\})$ is a generalized $q$-centered coloring of $G$
    and $v \mapsto \{\emptyset\}$ is $0$-good.
    Therefore, there exists $u \in V(G)$, and $L \in \{\emptyset\}$ such that $|\phi(V(G)) \setminus L| = |\phi(V(G))| \geq \binom{q-1+t-1}{t-1}$,
    contradicting the fact that $\phi$ uses at most $k$ colors.
    Thus, every $q$-centered coloring of $G$ uses at least $\binom{q-1+t-1}{t-1}$
    colors.
    Since $\binom{q-1+t-1}{t-1} = \Omega(q^{t-1})$, this proves the corollary.
\end{proof}

\begin{corollary}\label{cor:lower-bound:cen:srtd}
    For every integer $t$ with $t \geq 2$,
    \begin{align*}
        \max_{G \in \SRt_t} \cen_q(G) &= \Omega(q^{t-2} \log q).
    \end{align*}
\end{corollary}

\begin{proof}
    Let $q$ and $t$ be positive integers.
    Let $k = \binom{q+t-2}{t-2}\lfloor\log q\rfloor-1$.

    \begin{claim}
    For every integer $s$ with $s\geq 2$,
    for every nonnegative integer $p$ with $p\leq q$,
    \[
        \max_{G \in \SRt_s} \cen^{(k)}_{4q,p}(G) \geq \binom{q + s-2-p}{s-2} \log q.
    \]
    \end{claim}

    \begin{proofclaim}
    We proceed by induction on $s$.
    If $s=2$, then $\binom{q + s-2-p}{s-2} = 1$,
    and this inequality follows from
    \Cref{lemma:lower_bound:cen:path} and the fact that $\cen^{(k)}_{4q,p}(G)\geq \cen^{(k)}_{4q,q}(G)$ for every nonnegative integer $p\leq q$.
    Next, suppose that $s>2$.
    
    By \Cref{lemma:lower_bound:cen:induction_step},
    because $\SRt_{s} = \Tree(\SRt_{s-1})$, 
    and applying the induction hypothesis,
    we deduce
    \begin{align*}
        \max_{G \in \SRt_s} \cen^{(k)}_{4q,p}(G)
        &\geq \sum_{r=p}^{4q} \max_{G \in \SRt_{s-1}} \cen^{(k)}_{4q,r}(G) \\
        &\geq \sum_{r=p}^{q} \max_{G \in \SRt_{s-1}} \cen^{(k)}_{4q,r}(G) \\
        &\geq \sum_{r=p}^q \binom{q+s-3-r}{s-3}\log q = \binom{q+s-2-p}{s-2}\log q.
    \end{align*}
    This concludes the induction.
    \end{proofclaim}

    By the previous claim, we consider a graph $G \in \SRt_t$
    such that $\cen^{(k)}_{4q,0}(G) \geq \binom{q+t-2}{t-2}\log q$.
    Suppose by contradiction that $G$ admits a $4q$-centered coloring $\phi$
    using at most $k$ colors.
    Then $(\phi,v \mapsto \{\emptyset\})$ is a generalized $4q$-centered coloring of $G$
    and $v \mapsto \{\emptyset\}$ is $0$-good.
    Therefore, there  exist $u \in V(G)$ and $L \in \{\emptyset\}$ such that 
    $|\phi(V(G)) \setminus L| = |\phi(V(G))| \geq \binom{q+t-2}{t-2} \log q > k$,
    contradicting the fact that $\phi$ uses at most $k$ colors.
    Therefore, every $4q$-centered coloring of $G$ uses at least $\binom{q+t-1}{t-1}\lfloor\log q\rfloor$
    colors.
    Since $\binom{\lfloor q/4 \rfloor+t-2}{t-2}\lfloor \log(\lfloor q/4 \rfloor)\rfloor = \Omega(q^{t-2} \log q)$, this proves the corollary.
\end{proof}

\subsection{Fractional treedepth fragility rates}
 
In this section, we prove lower bounds for fractional treedepth fragility rates
of graphs in $\mathcal{R}_t$ and $\mathcal{S}_t$ 
(see \Cref{cor:lower-bound:frates:rtd,cor:lower-bound:frates:srtd}).
This was essentially proved by Dvořák and Sereni~\cite{DS_2020},
but we elected to provide full proofs (following the same approach)
since their statements are only in terms of treewidth and simple treewidth.

For every graph $G$ and every function $\ttt \colon V(G) \to \nonnegint$, we define \defin{$\td(G,\ttt)$} recursively as
\begin{enumerate}
    \item $\td(G,\ttt) = 0$ if $G$ is the null graph;
    \item $\td(G,\ttt) = \max\{\td(C,\ttt\vert_{V(C)})\mid C \textrm{ component of } G\}$ if $G$ is nonnull and not connected; and
    \item $\td(G,\ttt) = 1+\min_{u \in V(G)} \max\{\ttt(u), \td(G-u,\ttt\vert_{V(G-u)})\}$ otherwise.
\end{enumerate}

A straightforward induction shows that
for every graph $G$, for every subgraph $G'$ of $G$,
for every $\ttt \colon V(G) \to \nonnegint$,
\[
    \td(G,\ttt) \geq \td(G',\ttt\vert_{V(G')}).
\]
First, we prove a few other properties of $\td(\cdot,\cdot)$.

\begin{lemma}\label{lemma:wieghted_treedepth:edge_contraction}
    Let $G$ be a graph,
    let $\ttt \colon V(G) \to \nonnegint$,
    let $U \subseteq V(G)$ with $G[U]$ connected,
    and let $u \in U$.
    Let $G'$ be the graph obtained from $G$ by identifying the vertices in $U$ into 
    a single vertex $w$.
    Let $\ttt'\colon V(G') \to \nonnegint$
    be such that $\ttt'(w) \leq \ttt(u)$ and $\ttt'(x) \leq \ttt(x)$ for every $x \in V(G) \setminus U$.
    We have
    \[
        \td(G,\ttt) \geq \td(G',\ttt').
    \]
\end{lemma}

\begin{proof}
    We proceed by induction on $|V(G)|$.
    The result is clear if $V(G) = U$. Next, suppose that $V(G) \setminus U \neq \emptyset$.
    
    Suppose that $G$ is not connected.
    Consider the families $\mathcal{C}$ and $\mathcal{C}'$ of all components of $G$ and $G'$ respectively.
    Let $C_0 \in \mathcal{C}$ be the component of $G$ containing $U$.
    This component corresponds to a component $C'_0 \in \mathcal{C}'$ of $G'$.
    By the induction hypothesis, $\td(C_0, \ttt\vert_{V(C_0)}) \geq \td(C'_0, \ttt'\vert_{V(C'_0)})$.
    Moreover, $\mathcal{C} \setminus \{C_0\} = \mathcal{C}'\setminus \{C'_0\}$.
    Therefore, 
    \[
        \td(G,\ttt) = \max_{C \in \mathcal{C}} \td(C,\ttt\vert_{V(C)}) \geq \max_{C' \in \mathcal{C}'} \td(C',\ttt'\vert_{V(C')}) = \td(G',\ttt').
    \]

    Now, suppose that $G$ is connected.
    Consider $z \in V(G)$ such that
    $\td(G,\ttt) = 1+\max\{\ttt(z),\td(G-z,\ttt\vert_{V(G-z)})\}$.
    
    First, suppose that $z \not\in U$.
    Then, by the induction hypothesis,
    $\td(G-z,\ttt\vert_{V(G-z)}) \geq \td(G'-z, \ttt'\vert_{V(G'-z)})$,
    and since $\td(G',\ttt') \leq 1+\max\{\ttt'(z),\td(G'-z,\ttt'\vert_{V(G'-z)})\}$ and $\ttt'(z)\leq\ttt(z)$,
    we deduce that $\td(G,\ttt) \geq \td(G',\ttt')$.
    
    Next, suppose that $z \in U$.
    Since $G'-w \subseteq G-z$ and $\ttt\vert_{V(G-U)}(v)\geq \ttt'\vert_{V(G'-w)}(v)$ for every $v\in V(G'-w)$,
    $\td(G'-w,\ttt'\vert_{V(G'-w)}) \leq \td(G-z, \ttt\vert_{V(G-z)})$.
    Moreover, if $z=u$, then $\ttt(z) \geq \ttt'(w)$,
    and if $z \neq u$, then $\td(G-z,\ttt\vert_{V(G-z)}) \geq 1+\ttt(u)\geq 1+\ttt'(w)$.
    In both cases, $1+\ttt'(w) \leq \max\{1+\ttt(z), \td(G-z,\ttt\vert_{G-z})\}$,
    and it follows that $\td(G,\ttt) = 1+\max\{\ttt(z),\td(G-z,\ttt\vert_{V(G-z)})\} \geq 1+\max\{\ttt'(w),\td(G'-w,\ttt'\vert_{V(G'-w)})\} \geq \td(G',\ttt')$.
\end{proof}

\begin{lemma}\label{lemma:weighted_treedepth:model}
Let $G$ and $H$ be graphs and
let $(B_x\mid x\in V(H))$ be a model of $H$ in $G$.
Let $\ttt\colon V(G) \to \nonnegint$ and $\ttt'\colon V(H) \to \nonnegint$ be such that for every $x\in V(H)$, $\max\{\ttt(u)\mid u\in B_x\} \geq \ttt'(x)$.
We have
\[
\td(G,\ttt)\geq \td(H,\ttt').
\]    
\end{lemma}

\begin{proof}
Let $V(H)=\{x_1,\dots,x_{|V(H)|}\}$.
Let $G_0=G$.
For every $i\in [|V(H)|]$, let $G_i$ be obtained from $G_{i-1}$ by identifying the vertices in $B_{x_i}$ into a single vertex $x_i$,
let $u_i\in B_{x_i}$ be such that $\ttt'(x_i)\leq \ttt(u_i)$,
and let $\ttt'_i\colon V(G_i)\to \nonnegint$ be such that $\ttt'_i(x_i) = \ttt'(x_i)$ and $\ttt'_i=\ttt(v)$ for every $v\in V(G_i)\setminus B_{x_i}$.
By applying iteratively \Cref{lemma:wieghted_treedepth:edge_contraction} for every $i\in [|V(H)|]$, we have
\[
\td(G,\ttt) \geq \td(G_1,\ttt'_1)\geq \dots\geq \td(G_{|V(H)|},\ttt'_{|V(H)|})=\td(H,\ttt').\qedhere
\]
\end{proof}

\begin{lemma}\label{lemma:compose_treedepth}
    Let $G_0$ be a graph.
    Let $\ttt_0 \colon V(G_0) \to \nonnegint$.
    Let $G_{u,i}$ be a family of pairwise disjoint
    graphs, all of them disjoint from $G_0$, for $u \in V(G_0)$ and $i \in \{1,2\}$.
    For every $u \in V(G_0)$ and $i \in \{1,2\}$,
    let $\ttt_{u,i} \colon V(G) \to \nonnegint$.
    Let $G$ be a graph with vertex set $V(G_0) \cup \bigcup_{u \in V(G_0), i \in \{1,2\}} V(G_{u,i})$
    and such that $E(G_0) \subseteq E(G)$, $E(G_{u,i}) \subseteq E(G)$, and $G[\{u\} \cup V(G_{u,i})]$ is connected
    for all $u \in V(G_0)$ and $i \in \{1,2\}$.
    Let $\ttt \colon V(G_0) \to \NN$ be such that
    $\ttt(v) \geq \ttt_{u,i}(v)$ for all $i \in \{1,2\}$ and $v \in V(G_{u,i})$.
    If $\td(G_{u,i}, \ttt_{u,i}) \geq \ttt_0(u)$ for all $u \in V(G_0)$ and $i \in \{i,2\}$,
    then
    \[
        \td(G,\ttt) \geq \td(G_0,\ttt_0).
    \]
\end{lemma}

\begin{proof}
    For every $u \in V(G)$, for every $i \in \{1,2\}$,
    for every $v \in V(G_{u,i})$,
    by replacing $\ttt_{u,i}(v)$ with $\min\{\ttt_{u,i}(v), \ttt_0(u)\}$,
    we maintain that $\td(G_{u,i},\ttt_{u,i})\geq \ttt_0(u)$
    so we suppose that $\ttt(v) = \ttt_{u,i}(v) \leq \ttt_0(u)$.
    
    We now proceed by induction on $|V(G)|$.
    If $G$ is the null graph, then
    $G_0$ is null as well, and so $\td(G,\ttt) = 0 = \td(G_0,\ttt_0)$.
    Now, suppose that $G$ is nonnull.

    Suppose that $G$ is not connected.
    By definition of $G$, $G_0$ is not connected.
    Let $\mathcal{C}_0$ be the family of all components of $G_0$
    and let $C_0 \in \mathcal{C}_0$.
    Let $G' = G[\bigcup_{u \in V(C_0)} (\{u\} \cup V(G_{u,1}) \cup V(G_{u,2}))]$.
    Since $G'$ is a subgraph of $G$, $\td(G,\ttt) \geq \td(G',\ttt\vert_{V(G')})$.
    Hence, by the induction hypothesis, 
    \[
        \td(G,\ttt) \geq \td(G',\ttt\vert_{V(G')}) \geq \td(C_0, \ttt_0\vert_{V(C_0)}).
    \]
    Therefore,
    \[
        \td(G,\ttt) \geq \max_{C_0 \in \mathcal{C}_0} \td(C_0, \ttt\vert_{V(C_0)}) = \td(G_0,\ttt_0).
    \]
    This concludes the case where $G$ is not connected.
    
    Suppose that $G$ is connected.
    By definition of $G$, $G_0$ is connected.
    Consider $v \in V(G)$ such that
    \[
        \td(G,\ttt) = 1 + \max\{\ttt(v), \td(G-v, \ttt\vert_{V(G-v)})\}.
    \]
    Let $u \in V(G_0)$
    and $i \in \{1,2\}$ such that $v \in V(G_{u,i}) \cup \{u\}$.
    Observe that $G_{u,3-i} \subseteq G - v$, and
    so $\td(G-v,\ttt\vert_{V(G-v)}) \geq \td(G_{u,3-i},\ttt_{u,3-i}) \geq \ttt_0(u) \geq \ttt(v)$.
    Therefore, $\td(G,\ttt) = 1 + \td(G-v,\ttt\vert_{V(G-v)})$.
    Moreover, by the induction hypothesis
    applied to $G' = G - (\{u\} \cup V(G_{u,1}) \cup V(G_{u,2}))$,
    we have $\td(G', \ttt\vert_{V(G')}) \geq \td(G_0-u,\ttt_0\vert_{V(G_0-u)})$.
    Since $G-v$ contains both $G_{u,3-i}$ and $G'$ as subgraphs,
    we conclude that
    \begin{align*}
        \td(G-v, \ttt\vert_{V(G-v)}) 
        &\geq \max\{\td(G_{u,3-i}, \ttt_{u,3-i}), \td(G', \ttt\vert_{V(G')})\}\\
        &\geq \max\{\ttt_0(u), \td(G_0 - u, \ttt\vert_{V(G_0-u)})\}.
    \end{align*}
    It follows that
    \begin{align*}
        \td(G,\ttt) 
        &= 1+\td(G-v,\ttt\vert_{V(G-v)}) \\
        &\geq 1+\max\{\ttt_0(u), \td(G_0 - u, \ttt\vert_{V(G_0-u)})\}\\
        &\geq \td(G_0,\ttt_0) && \textrm{ since $G_0$ is connected}.
    \end{align*}
    This concludes the proof of the lemma.
\end{proof}

We denote by \defin{$\nonnegrat$} the set of nonnegative rationals.
Let $G$ be a graph.
Let $w\colon V(G) \to \nonnegrat$.
For every $U \subseteq V(G)$,
we write $\defmath{w(U)} = \sum_{u \in U} w(u)$.

Let $q$ be a positive integer.
We denote by \defin{$\frate'_q(G)$} the maximum integer $k$ such that
there exists $w \colon V(G) \to \nonnegrat$ such that
for every $\tau \in \nonnegint$,
for every $\ttt \colon V(G) \to \nonnegint$,
for every $Y \subseteq V(G)$,
if
\begin{enumerate}
    \item $w(Y) \leq \frac{1}{q} w(V(G))$, and
    \item $w(\{u \in V(G) \mid \ttt(u) < \tau\}) \leq \frac{1}{2} w(V(G))$,
\end{enumerate}
then there exists disjoint subsets $S_1,S_2 \subseteq V(G) \setminus Y$
such that 
\[
    \td(G[S_i],\ttt\vert_{S_i}) \geq \tau + k
\]
for every $i \in \{1,2\}$.

Now, we show that $\frate'_q$ is a lower bound on $\frate_q$.
\begin{lemma}\label{lemma:lower_bound_fragility:frate'_and_frate}
    For every graph $G$ and for every positive integer $q$,
    \[
        \frate_q(G) \geq \frate'_q(G).
    \]
\end{lemma}

\begin{proof}
    Let $G$ be a graph,
    let $q$ be a positive integer,
    and let $k = \frate'_q(G)$.
    Let $\rvar{Y}$ be a random variable over subsets of $V(G)$
    equipped with a $q$-thin distribution.
    We want to show that there exists a value $Y$ of $\rvar{Y}$ such that $\td(G\setminus Y)\geq \frate'_q(G)$.
    
    Let $w\colon V(G) \to \nonnegrat$ witness $\frate'_q(G) \geq k$.
    The expected value of $w(\rvar{Y})$
    is at most $\frac{1}{q}w(V(G))$ since $\rvar{Y}$ is equipped with a $q$-thin distribution.
    So, we fix $Y \subseteq V(G)$ of $\rvar{Y}$
    with $w(Y) \leq \frac{1}{q}w(V(G))$.
    Then, considering $\ttt \colon V(G) \to \nonnegint$ constant to $0$ and $\tau=0$,
    we have $S \subseteq V(G)\setminus Y$ with 
    $\td(G[S],\ttt\vert_{S}) \geq \frate'_q(G)$.
    In particular, $\td(G-Y) = \td(G-Y,\ttt\vert_{V(G-Y)}) \geq \td(G[S],\ttt\vert_S) \geq k$.
    This proves that $\frate_q(G) \geq k$.
\end{proof}

\begin{lemma}\label{lemma:lower_bound_fragility:paths}
    For every integer $q$ with $q\geq 5$,
    \[
        \frate'_q(P_{q-1}) \geq \log(q) - 2.
    \]
\end{lemma}

\begin{proof}
    Let $q$ be a positive integer.
    Let $k =\lceil \log (1+\lfloor(q-1)/4\rfloor) \rceil$.
    Let $w \colon V(P_{q-1}) \to \nonnegrat$ be constant to $1$.
    Let $\ttt \colon V(P_{q-1}) \to \nonnegint$,
    and let $\tau$ be a nonnegative integer such that
    $w(\{u \in V(P_{q-1}) \mid \ttt(u) < \tau\}) \leq \frac{q-1}{2}$.
    Every $Y\subseteq V(P_{q-1})$ such that $w(Y)\leq w(V(P_q))/q = \frac{q-1}{q}$, is empty.
    So, we only have to show that 
    there exists disjoint subsets $S_1,S_2 \subseteq V(P_{q-1})$ such that
    $\td(G[S_i],\ttt\vert_{S_i}) \geq \tau + k$ for every $i \in \{1,2\}$.

    Since $w(\{u \in V(P_q) \mid \ttt(u) \geq \tau\}) \geq \frac{q-1}{2}$
    and by the definition of $w$,
    we have that $|\{u \in V(P_{q-1}) \mid \ttt(u) \geq \tau\}| \geq \frac{q-1}{2}$.
    Consider disjoint subpaths $Q_1,Q_2$ of $P_{q-1}$ such that 
    $|\{v \in V(Q_i) \mid \ttt(u) \geq \tau\}| \geq \left\lfloor \frac{q-1}{4}\right\rfloor$ for each $i \in \{1,2\}$.
    Let $\ttt'\colon V(P_{\lfloor (q-1)/4 \rfloor}) \to \nonnegint$ be constant to $\tau$.
    Thus, for every $i\in\{1,2\}$, $Q_i$ contains a model $(B_x\mid x\in V(P_{\lfloor (q-1)/4 \rfloor}))$ such that $\ttt'\leq \max\{\ttt(u)\mid u\in B_x\}$ for every $x\in V(P_{\lfloor (q-1)/4 \rfloor})$.
    Therefore, by \Cref{lemma:weighted_treedepth:model},
    for each $i \in \{1,2\}$,
    \[
        \td(Q_i,\ttt\vert_{V(Q_i)}) \geq \td(P_{\lfloor (q-1)/4 \rfloor}, \ttt').
    \]
    Then, a simple induction shows that for every $q\geq 5$,
    \begin{align*}
        \td(P_{\lfloor (q-1)/4 \rfloor}, \ttt') 
        &\geq \tau + \td(P_{\lfloor(q-1)/4 \rfloor}) \\
        &\geq \tau + \lceil \log (1+\lfloor(q-1)/4\rfloor) \rceil =\tau +k.
    \end{align*}

    It follows that for $S_1 = V(Q_1)$ and $S_2 = V(Q_2)$,
    we have $\td(G[S_i], \ttt\vert_{S_i}) \geq \tau + k$ for every $i \in \{1,2\}$.
    Since $k = \lceil \log (1+\lfloor(q-1)/4\rfloor) \rceil \geq \log(q/4) = \log q - 2$,
    this proves the lemma.
\end{proof}

\begin{lemma}\label{lemma:lower_bound:fragility:applying_T}
    Let $\mathcal{G}$ be a nonempty class of graphs,
    and let $g(q) = \max_{G \in \mathcal{G}} \frate'_q(G)$ for every positive integer $q$.
    For every positive integer $q$,
    there exists $G \in \Tree(\mathcal{G})$ such that
    \[
        \frate'_{8q}(G) \geq qg(2q). 
    \]
\end{lemma}

\begin{proof}
    Let $q$ be a positive integer,
    let $G'$ be a graph in $\mathcal{G}$ such that 
    $\frate'_{2q}(G') \geq g(2q)$,
    and let $w' \colon V(G') \to \nonnegrat$ witness this fact.
    Observe that, by the definition of $\frate'$, we have $\frate'_{2q}(G'-v)=\frate'_{2q}(G')$ for every $v \in V(G')$ with $w'(v)=0$.
    Therefore, removing these vertices from $v$,
    we assume that $w'$ takes only positive values.
    By possibly replacing $w'$ by $\frac{1}{w'(V(G'))} \cdot w'$, 
    we also assume $w'(V(G')) = 1$.

    Let $h$ be a positive integer
    and let $G_h$ be the graph defined by
    \begin{align*}
        V(G_h) &= \bigcup_{i=0}^{h} V(G')^i \\
    \intertext{and}
        E(G_h) &= \bigcup_{i=1}^{h} 
                \Big(
                    \left\{(u_1, \dots, u_{i-1}) (u_1, \dots, u_{i}) \mid (u_1, \dots, u_{i}) \in V(G')^i\right\} \\
             &\hspace{8.5mm}
                \cup
                    \left\{(u_1, \dots, u_{i-1},v)(u_1, \dots, u_{i-1},w) \mid (u_1, \dots, u_{i-1}) \in V(G')^{i-1}, vw \in E(G')\right\}
                \Big).
    \end{align*}
    
    \begin{figure}[tp]
        \centering
        \includegraphics{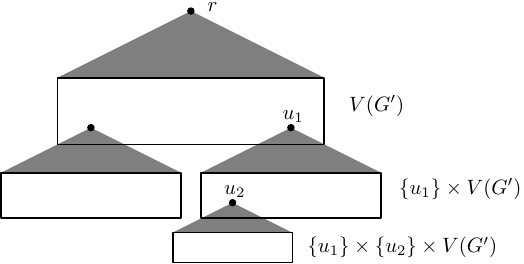}
        \caption{Illustration for the definition of $G_h$.}
        \label{fig:frate_construction}
    \end{figure}
    
    We denote by $r$ the unique vertex of $G_h$ in $V(G')^0$.
    See \Cref{fig:frate_construction}.

    Let $T_h$ be the tree on $V(G_h)$ defined by
    \[
    E(T_h) = \bigcup_{i=1}^{h} \left\{(u_1, \dots, u_{i-1}) (u_1, \dots, u_{i}) \mid (u_1, \dots, u_{i}) \in V(G')^i\right\}.
    \]
    We root $T_h$ in $r$.

    Recall that $G' \in \mathcal{G}$.
    For every $x \in V(T_h)$,
    let
    \[
        W_x =
        \begin{cases}
        \{x\} \cup \{(x_1, \dots, x_i, v) \mid v \in V(G')\} & \textrm{for $1\leq i\leq h-1$ and $x=(x_1,\dots,x_i)$,}\\
        \{x\} & \textrm{ otherwise.}
        \end{cases}
    \]
    The pair $\big(T_h,(W_x \mid x \in V(T_h))\big)$ is a tree decomposition of $G$ witnessing the fact that $G_h \in \Tree(\mathcal{G})$.

    For every $u$, let $G_{h,u} = G_h[V(\subtree{T_h}{u})]$.
    Note that $G_{h,u}$ is isomorphic to $G_{h-i}$ where $i$ is the height of $u$ in $T_h$.

    Let $w\colon V(G_h) \to \nonnegrat$ be defined as
    \[
    w(u) =
    \begin{cases}
        1 & \textrm{if $u=r$,}\\
        \prod_{j=1}^i w'(u_j) & \textrm{if $u = (u_1, \dots, u_i)$.}
    \end{cases}
    \]
    Observe that
    \begin{enumerate}
        \item for every $u\in V(G_h)$, $w(u)=\max\{w(v)\mid v\in V(G_{h,u})\}$; \label{obs:w_max}
        \item for every $i\in [h]$, $w(V(G')^i)=\sum_{(u_1, \dots, u_i)\in V(G')^i}\prod_{j=1}^i w'(u_j) =(\sum_{v\in V(G')}w'(v))^i=1$ so $w(V(G_h))=h+1$; and \label{obs:w_1}
        \item for every $u \in V(G_h)$,
            $w(u) = w(N_{G_{h,u}}(u))$. \label{obs:w_neighbors}
    \end{enumerate}
    
    \begin{claim}\label{claim:lower_bound_fragility:finding_S1_S2}
        Let $p$ be a nonnegative integer,
        let $h$ be a positive integer,
        let $s \in V(G_h)$,
        let $S \subseteq V(G_{h,s})$ such that $s \in S$ and $G_h[S]$ connected,
        let $\ttt \colon S \to \nonnegint$,
        and let $\tau \in \nonnegint$.
        If $w(\{u \in S \mid \ttt(u) \geq \tau\}) > 2(q+p) \cdot w(s)$, then
        there exist disjoint subsets $S_0, S_1, S_2 \subseteq S$
        such that
        \begin{enumerate}[label=(\alph*)]
            \item $s \in S_0$; \label{item:claim:lower_bound_fragility:finding_S1_S2:i}
            \item $G_h[S_a]$ is connected for every $a \in \{0,1,2\}$; \label{item:claim:lower_bound_fragility:finding_S1_S2:ii}
            \item there is an edge between $S_0$ and $S_a$ in $G_h$, for every $a \in \{1,2\}$; and \label{item:claim:lower_bound_fragility:finding_S1_S2:iii}
            \item $\td(G_h[S_a], \ttt\vert_{S_a}) \geq \tau + p g(2q)$ for every $a \in \{1,2\}$. \label{item:claim:lower_bound_fragility:finding_S1_S2:iv}
        \end{enumerate}
    \end{claim}

    \begin{proofclaim}
        Let $d$ be the height of $s$ in $G$.
        We proceed by induction on $p + (h - d)$.

        If $d=h$, then $S = \{s\}$,
        and so $w(\{u \in S \mid \ttt(u) \geq \tau\}) \leq w(s)$.
        It follows that the claim is vacuously true.
        Next, suppose that $d<h$ and $w(\{u \in S \mid \ttt(u) \geq \tau\}) > 2(q+p) \cdot w(s)$.
        
        Let $N = N_{G_{h,s}}(s)$, i.e.\ the set of all the children of $s$ in $T_h$,
        and let $N_S = N \cap S$.
        For every $v \in N_S$, let $S_v = S \cap V(G_{h,v})$.
        Since $G_h[S]$ is connected, for every $v \in N_S$, $G_h[S_v]$ is connected.
        
        If there exists $v \in N_S$ such that $w(\{u \in S_v \mid \ttt(u) \geq \tau\}) > 2(q+p) \cdot w(v)$,
        then, since the height of $v$ is larger than the height of $s$,
        we can call the induction hypothesis on $S_v$.
        This gives three disjoint sets $S^0_0, S^0_1, S^0_2 \subseteq S_v$
        satisfying \ref{item:claim:lower_bound_fragility:finding_S1_S2:i}--\ref{item:claim:lower_bound_fragility:finding_S1_S2:iv}.
        Then $(S_0,S_1,S_2) = (S^0_0 \cup \{s\}, S^0_1, S^0_2)$ are as wanted.

        Now assume for every $v\in N_S$,
        \begin{equation}\label{eq:lower_bound_fragility:claim:small_subtrees}
            w(\{u \in S_v \mid \ttt(u) \geq \tau\}) \leq 2(q+p) \cdot w(v).
        \end{equation}

        We deduce that
        \begin{align*}
            2(q+p) \cdot w(s)
            &< w(\{u \in S \mid \ttt(u) \geq \tau\}) \\
            &= w(s) + \sum_{v \in N_S} w(\{v \in S_v \mid \ttt(u) \geq \tau\}) \\
            &\leq w(s) + \sum_{v \in N_S} 2(q+p) \cdot w(v) && \textrm{by \eqref{eq:lower_bound_fragility:claim:small_subtrees}} \\
            &= w(s) + 2(q+p) \cdot w(N_S).
        \end{align*}
        Since $w(s) = w(N)$ by \ref{obs:w_neighbors},
        \begin{equation}\label{eq:lower_bound_fragility:claim:B_is_big}
            w(N_S) > \left(1-\frac{1}{2(q+p)}\right) \cdot w(N).
        \end{equation}

        Let 
        \begin{align*}   
        L_S &= \big\{v \in N_S \mid w(\{u \in S_v \mid \ttt(u) \geq \tau\}) \leq 2(q+p-1) \cdot w(v)\big\} \textrm{ and}\\ 
        M_S &= \big\{v \in N_S \mid w(\{u \in S_v \mid \ttt(u) \geq \tau\}) > 2(q+p-1) \cdot w(v)\big\}.
        \end{align*}
        Observe that $L_S$ and $M_S$ form a partition of $N_S$.
        
        Using \eqref{eq:lower_bound_fragility:claim:small_subtrees}, we deduce
        \begin{align*}
            2(q+p) \cdot w(s) 
            &< w(\{u \in S \mid \ttt(u) \geq \tau\}) \\
            &= w(s) + \sum_{v \in L_S} w(\{u \in S_v \mid \ttt(u) \geq \tau\}) + \sum_{v \in M_S} w(\{u \in S_v \mid \ttt(u) \geq \tau\}) \\
            &\leq w(s) + 2(q+p-1) \cdot w(L_S) + 2(q+p) \cdot w(M_S) \\
            &= w(s) + + 2(q+p-1) \cdot w(L_S) + 2(q+p) \cdot (w(N_S) -w(L_S))\\
            &= w(s) + 2(q+p) \cdot w(N_S) - 2 w(L_S) \\
            &\leq (2q+2p+1) \cdot w(s) - 2w(L_S) \qquad \textrm{ since $w(N_S)\leq w(N)=w(s)$ by \ref{obs:w_neighbors}}.
        \end{align*}
        Therefore, since $w(s) = w(N)$ by \ref{obs:w_neighbors},
        \begin{equation}\label{eq:lower_bound_fragility:claim:C_is_small}
            w(L_S) < \frac{1}{2} w(N).
        \end{equation}
        
        Recall that $G[N]$ is isomorphic to $G'$ and observe that for every $u\in N$, $w(u) = w(s)w'(u_{d+1})$ where $u=(s_1,\dots,s_d,u_{d+1})$.
        Since $w'$ witnesses the fact that $\frate'_{2q}(G') \geq g(2q)$, $w\vert_{N}$ witnesses the fact that $\frate'_{2q}(G[N]) \geq g(2q)$.
        By \eqref{eq:lower_bound_fragility:claim:B_is_big}, the set $Y = N \setminus N_S$ satisfies $w(Y) < \frac{1}{2(q+p)} w(N) \leq \frac{1}{2q} w(N)$.
        
        Let
        \begin{align*}
        \tau' &= \tau + (p-1)g(2q)\\
        \intertext{and $\ttt' \colon N \to \nonnegint$ defined by}
        \ttt'(v) &=
        \begin{cases}
            \tau' &\textrm{if $v \in M_S$,} \\
            0 &\textrm{otherwise.}
        \end{cases}
        \end{align*}
        
        By \eqref{eq:lower_bound_fragility:claim:C_is_small}, $w(\{u \in N_S \mid \ttt'(u) < \tau'\})  = w(L_S) < \frac{1}{2} w(N)$.
        By definition of $\frate'$,
        we take disjoint subsets $S'_1, S'_2 \subseteq N\setminus Y = N_S$ such that for every $a\in\{1,2\}$,
        \[
            \td(G_h[S'_a], \ttt'\vert_{S'_a}) \geq \tau' + g(2q) = \tau + p g(2q).
        \]
        By definition of $\td(\cdot,\cdot)$, we can assume that $G_h[S'_a]$ is connected for every $a\in\{1,2\}$.
        
        If $p=0$,
        then, $(S_0,S_1,S_2) = (\{s\}, S'_1,S'_2)$ satisfy \ref{item:claim:lower_bound_fragility:finding_S1_S2:i}--\ref{item:claim:lower_bound_fragility:finding_S1_S2:iv} and we are done.
        Now suppose $p \geq 1$.

        Let $v \in M_S$.
        By definition of $M_S$, $w(\{u \in S_v \mid \ttt(u) \geq \tau\}) > 2(q+p-1) w(v)$.
        By the induction hypothesis,
        there exist disjoint subsets $S_{v,0},S_{v,1},S_{v,2} \subseteq S_v$ such that
        \begin{enumerate}[label=(\alph*')]
            \item $v \in S_{v,0}$; \label{item:claim:lower_bound_fragility:finding_S1_S2:i:call}
            \item $G_h[S_{v,a}]$ connected for every $a \in \{0,1,2\}$; \label{item:claim:lower_bound_fragility:finding_S1_S2:ii:call}
            \item there is an edge between $S_{v,0}$ and $S_{v,a}$ in $G_h$, for every $a \in \{1,2\}$; and \label{item:claim:lower_bound_fragility:finding_S1_S2:iii:call}
            \item $\td(G_h[S_{v,a}], \ttt\vert_{S_{v,a}}) \geq \tau + (p-1)g(2q) = \tau'$ for every $a \in \{1,2\}$.\label{item:claim:lower_bound_fragility:finding_S1_S2:iv:call}
        \end{enumerate}

        Let $S_a = S'_a \cup \bigcup_{v\in S'_a\cap M_S}(S_{v,0}\cup S_{v,1}\cup S_{v,2})$ for every $a\in\{1,2\}$.
        Let $S_0 = \{s\}$ so \ref{item:claim:lower_bound_fragility:finding_S1_S2:i} holds.
        By \ref{item:claim:lower_bound_fragility:finding_S1_S2:i:call}, \ref{item:claim:lower_bound_fragility:finding_S1_S2:ii:call}, and the fact that $S'_a$ is connected for every $a\in\{1,2\}$, \ref{item:claim:lower_bound_fragility:finding_S1_S2:ii} holds.
        Since $S'_a\subseteq N$ for every $a\in\{1,2\}$, \ref{item:claim:lower_bound_fragility:finding_S1_S2:iii} holds.

        Let $a \in \{1,2\}$.
        Let $H$ be the graph obtained from $G_h[S_a]$ by contracting, for each $v \in S'_a$, $S_{v,0}$
        into a single vertex that we still call $v$.
        Thus, $\ttt\vert_{V(H)}$ is well defined.
        By \Cref{lemma:weighted_treedepth:model},
        $\td(G_h[S_a], \ttt\vert_{S_a}) \geq \td(H,\ttt\vert_{V(H)})$.

        Let $H_0 = H[S'_a]$ and for all $v\in V(H_0) = S'_a$ and $b \in \{1,2\}$, let 
        \[
            H_{v,b} =
            \begin{cases}
                H[S_{v,b}] & \textrm{if $v \in M_S$,} \\
                \textrm{the null graph} & \textrm{otherwise.}
            \end{cases}
        \]
        By \ref{item:claim:lower_bound_fragility:finding_S1_S2:iv:call},
        we have $\td(H_{v,b},\ttt\vert_{V(H_{u,b})}) \geq \ttt'(v)$ for all $v \in S'_a$ and $b \in \{1,2\}$.
        Therefore, by \Cref{lemma:compose_treedepth},
        we have $\td(H,\ttt\vert_{V(H)}) \geq \td(G_h[S'_a], \ttt'\vert_{S'_a})$.
        Therefore, $\td(G_h[S_a], \ttt\vert_{S_a}) \geq \tau + pg(2q)$.
        This shows \ref{item:claim:lower_bound_fragility:finding_S1_S2:iv} and
        concludes the proof of the lemma.
    \end{proofclaim}

    Now, we exhibit a graph in $\Tree(\mathcal{G})$
    with $\frate'_{8q}(G) \geq qg(2q)$.
    We take $G=G_h$ for $h = 8q-1$.

    Let $\tau \in \nonnegint$,
    let $\ttt \colon V(G) \to \nonnegint$,
    and let $Y \subseteq V(G)$ such that
    \begin{enumerate}
        \item $w(Y) \leq \frac{1}{8q} w(V(G))$, and
        \item $w(\{u \in V(G) \mid \ttt(u) < \tau\}) \leq \frac{1}{2} w(V(G))$.
    \end{enumerate}

    For every $y \in Y \cup \{r\}$, let $S_y$ be the vertex set of the component of $y$ in $G_{h,y}[V(G_{h,y} - (Y\cup\{r\})) \cup \{y\}]$.
    Note that $\bigcup_{y \in Y \cup \{r\}} S_y = V(G)$,
    and so
    \[
        \frac{w(V(G))}{w(Y\cup\{r\})} \leq \sum_{y \in Y\cup\{r\}} \frac{w(y)}{w(Y\cup \{r\})}\cdot\frac{w(S_y)}{w(y)}.
    \]
    Therefore,
    there exists $y \in Y \cup \{r\}$
    such that $\frac{w(S_y)}{w(y)} \geq \frac{w(V(G))}{w(Y\cup \{r\})}$.
    Fix such $y$ in $Y\cup \{r\}$.
    By \ref{obs:w_1} and the fact that $w(Y)\leq \frac{1}{8q}w(V(G))$,
    \[
        w(S_y) 
        \geq \frac{w(V(G))}{w(r) + w(Y)} \cdot w(y) 
        \geq \frac{h+1}{1 + \frac{h+1}{8q}} \cdot w(y) 
        = \frac{8q}{1+\frac{8q}{8q}}\cdot w(y)
        = 4q \cdot w(y).
    \]
    Therefore, by \Cref{claim:lower_bound_fragility:finding_S1_S2} applied for $p=q$,
    there exist disjoint $S_{1},S_{2} \subseteq S_y \setminus \{y\}$
    such that for every $a\in\{1,2\}$
    \[
        \td(G[S_{a}], \ttt\vert_{S_{a}}) \geq \tau + q g(2q).
    \]
    Then, the sets $S_1,S_2$
    are as wanted in the definition of $\frate'$.
    Since $Y,\ttt,\tau$ were arbitrary,
    this shows that $\frate'_{8q}(G) \geq q g(2q)$.
\end{proof}

\begin{corollary}\label{cor:lower-bound:frates:rtd}
    For every positive integer $t$,
    \begin{align*}
        \max_{G \in \Rt_t} \frate_q(G) &= \Omega(q^{t-1}).
    \end{align*}
\end{corollary}

\begin{proof}
    Let $t$ be a positive integer.
    We show by induction on $t$ that
    there exists $c>0$ such that for every integer $q$ with $q\geq 16$,
    \[
        \max_{G \in \Rt_t} \frate'_{q}(G) \geq c \cdot q^{t-1}.
    \]
    
    For $t=1$, this is clear for $G = K_1$ and $c=1$.
    Next, suppose that $t \geq 2$.

    By the induction hypothesis,
    there exists $c' > 0$ such that
    $\max_{G \in \Rt_{t-1}} \frate'_{q}(G) \geq c' \cdot q^{t-2}$
    for every integer $q$ with $q\geq 16$.
    Let $q$ be an integer with $q\geq 16$.
    Since $\Rt_t = \Tree(\Rt_{t-1})$,
    by \Cref{lemma:lower_bound:fragility:applying_T},
    \begin{align*}
        \max_{G \in \Rt_t} \frate'_{q}(G) 
        &\geq \lfloor q/8 \rfloor \cdot c' \cdot (2\lfloor q/8 \rfloor)^{t-2} \\
        &\geq c' \cdot q/16 \cdot (q/8)^{t-2} \\
        &= \frac{c'}{2^{3t-2}} \cdot q^{t-1}\\
        &= c\cdot q^{t-1} \textrm{ for }c = \frac{c'}{2^{3t-2}}.
    \end{align*}
    This concludes the induction,
    which, by \Cref{lemma:lower_bound_fragility:frate'_and_frate}, implies the corollary.
\end{proof}

\begin{corollary}\label{cor:lower-bound:frates:srtd}
    For every integer $t$ with $t \geq 2$,
    \begin{align*}
        \max_{G \in \SRt_t} \frate_q(G) &= \Omega(q^{t-2} \log q).
    \end{align*}
\end{corollary}

\begin{proof}
    Let $t$ be an integer with $t \geq 2$.
    We show by induction on $t$ that
    there exists $c > 0$ such that
    \[
        \max_{G \in \SRt_t} \frate'_{q}(G) \geq c \cdot q^{t-2} \log q
    \]
    for every integer $q$ with $q \geq 16$.
    When $t=2$, this is a consequence of \Cref{lemma:lower_bound_fragility:paths}.
    Next, suppose that $t \geq 3$.

    By the induction hypothesis,
    there exists $c' > 0$ such that
    $\max_{G \in \Rt_t} \frate'_{q}(G) \geq c' \cdot q^{t-3} \log q$
    for every integer $q$ with $q \geq 16$.
    Let $q$ be an integer with $q \geq 16$.
    Since $\SRt_t = \Tree(\SRt_{t-1})$,
    by \Cref{lemma:lower_bound:fragility:applying_T},
    \begin{align*}
        \max_{G \in \SRt_t} \frate'_{q}(G) 
        &\geq \lfloor q/8 \rfloor \cdot c' \cdot (2\lfloor q/8\rfloor)^{t-3} \log(2\lfloor q/8 \rfloor)\\
        &\geq c' \cdot q/16 \cdot (q/8)^{t-3} \log(q/8)\\
        &\geq c' \cdot q/16 \cdot (q/8)^{t-3} \log(q)/4 \\
        &= \frac{c'}{2^{3t-3}} \cdot q^{t-2} \log q\\
        &= c\cdot q^{t-2}\log q \textrm{ for }c = \frac{c'}{2^{3t-3}}.
    \end{align*}
    This concludes the induction,
    which, by \Cref{lemma:lower_bound_fragility:frate'_and_frate}, implies the corollary.
\end{proof}

\section{\texorpdfstring{$2$}{2}-treedepth and rooted \texorpdfstring{$2$}{2}-treedepth}
\label{sec:td2_and_rtd2}

In this section, we give a few facts about rooted $2$-treedepth.
In particular, we show that it is tied with $2$-treedepth introduced by Huynh, Joret, Micek, Seweryn, and Wollan~\cite{HJMSW22}.
In fact, these parameters are the same up to a multiplicative constant of $2$, see~\Cref{lemma:rtd_2-and-td_2}.
Moreover, the bounds that we give are optimal.
The \defin{$2$-treedepth} of a graph $G$, denoted by \defin{$\td_2(\cdot)$}, is defined recursively as
\begin{enumerate}[label={\normalfont (t\arabic*)}]
    \item $\td_2(G) = 0$ if $G$ is the null graph;
    \item $\td_2(G) = \max_{B} \td_2(B)$ over all the blocks\footnote{A \defin{block} of a graph $G$ is a maximal subgraph of $G$ which is either $2$-connected, a single edge, or an isolated vertex.} $B$ of $G$, if $G$ has more than one block; and
    \item $\td_2(G) = 1 + \min_{u \in V(G)} \td_2(G-u)$ if $G$ has exactly one block.
\end{enumerate}

Let $G$ be a graph and let $k$ be a nonnegative integer.
A \defin{separation of order $k$} of $G$ is a pair $(A, B)$ of subgraphs of $G$ such that $A\cup B = G$,
$E(A \cap B) = \emptyset$, and $|V (A \cap B)| = k$.

Our new parameter, rooted $2$-treedepth can be equivalently defined in a similar way to $2$-treedepth.
For this, we need the following two lemmas.

\begin{lemma}\label{lemma:rtd_and_separations}
    Let $G$ be a graph.
    For every separation $(A,B)$ of $G$ of order at most one,
    \[
        \rtd_2(G) \leq \max \left\{\rtd_2(A), |V(A) \cap V(B)| + \rtd_2(B \setminus V(A))\right\}.
    \]
\end{lemma}

\begin{proof}
    Consider a separation $(A,B)$ of $G$ of order at most one.
    Let 
    \[t = \max\left\{\rtd_2(A),\rtd_2\big(B \setminus V(A)\big) + |V(A) \cap V(B)|\right\}.\]
    Suppose that $G$ has at least one edge, and so $t \geq 2$.
    
    First consider the case where $(A,B)$ is of order $0$.
    By the definition of $t$, we have $A,B \in \Rt_t$.
    Let $\big(F_1, (W_{1,x} \mid x \in V(F_1))\big)$ and $\big(F_2, (W_{2,x} \mid x \in V(F_2))\big)$,
    be rooted forest decompositions of $A$ and $B$ witnessing the fact that $A,B \in \Rt_t = \Tree(\Rt_{t-1})$ respectively.
    Then $\big(F_1 \cup F_2, (W_x \mid x \in V(F_1 \cup F_2))\big)$,
    where $W_x = W_{i,x}$ for every $x \in V(F_i)$, for every $i \in \{1,2\}$,
    is a rooted forest decomposition of $G$ witnessing the fact that $G \in \Tree(\Rt_{t-1}) = \Rt_t$.
    
    Next, suppose that $(A,B)$ is of order $1$.
    Let $u$ be the unique vertex in $V(A) \cap V(B)$.
    By the definition of $t$, we have $A \in \Rt_t$ and $B-u \in \Rt_{t-1}$.
    Let $\big(F, (W_x \mid x \in V(F))\big)$ be a rooted forest decomposition of $A$ witnessing the fact that $A \in \Rt_t = \Tree(\Rt_{t-1})$.
    There exists $y \in V(F)$ such that $u \in W_y$.
    Let $F'$ be obtained from $F$ by adding a fresh leaf $z$ with parent $y$.
    Then, let 
    \[
        W'_x = 
        \begin{cases}
            W_x &\textrm{if $x \in V(F)$,} \\
            V(B) &\textrm{if $x=z$.}
        \end{cases}
    \]
    It follows that $\big(F', (W'_x \mid x \in V(F'))\big)$ is a rooted forest decomposition
    witnessing the fact that $G \in \Tree(\Rt_{t-1}) = \Rt_t$.
    This proves the observation.
\end{proof}

\begin{lemma}\label{lemma:rtd_and_separations2}
    Let $G$ be a graph with at least two vertices.
    There is a separation $(A,B)$ of $G$ of order at most one with $V(A) \neq \emptyset$ and $V(B) \setminus V(A) \neq \emptyset$ such that
    \[
        \rtd_2(G) \geq \max \left\{\rtd_2(A), |V(A) \cap V(B)| + \rtd_2(B - V(A))\right\}.
    \]
\end{lemma}

\begin{proof}
    Let $t = \rtd_2(G)$.
    Let $\big(T,(W_x \mid x \in V(T))\big)$ be a rooted forest decomposition of $G$
    witnessing the fact that $G \in \Tree(\Rt_{t-1}) = \Rt_t$.
    Consider such a forest decomposition with $|V(T)|$ minimal.
    In particular, $W_x \neq \emptyset$ for every $x \in V(T)$
    and $W_x \not\subseteq W_y$ for every $xy \in E(T)$.

    If $T$ is not connected, then consider a component $S$ of $T$.
    Let $A = G[\bigcup_{x \in V(S)} W_x]$ and $B = G - V(A)$.
    Then $(A,B)$ is a separation of $G$ of order $0$ such that $V(A)\neq \emptyset$ and $V(B)\setminus V(A)\neq \emptyset$.
    Moreover,
    $\big(S, (W_x \mid x \in V(S))\big)$ witnesses the fact that $S \in \Tree(\Rt_{t-1}) = \Rt_t$
    and $\big(T - V(S), (W_x \mid x \in V(T - V(S)))\big)$ witnesses the fact that
    $B \in \Tree(\Rt_{t-1}) = \Rt_t$.
    
    Next, suppose that $T$ is connected.
    If $T$ has only one vertex, then $V(G)$ has at most one vertex, a contradiction.
    Hence $T$ has at least two vertices.
    Therefore, $T$ has a leaf $y$.
    Let $A = G[\bigcup_{x \in V(T-y)} W_x]$ and $B = G[W_y]$.
    Since $W_y$ is not included in $W_{\parent(T,y)}$,
    $V(B) \setminus V(A) \neq \emptyset$,
    and since $W_{\parent(T,y)} \neq \emptyset$, $V(A) \neq \emptyset$.
    Then $B - V(A) \in \Rt_{t-1}$ by the choice of $\big(T,(W_x \mid x \in V(T))\big)$.
    Moreover, $(T - y, (W_x \mid x \in V(T-y))\big)$ witnesses the fact that $A \in \Tree(\Rt_{t-1}) = \Rt_t$.
    Since $|V(T)|$ is minimal and $T$ is connected, $\big(T,(W_x \mid x \in V(T))\big)$ has adhesion exactly one by definition.
    Therefore, $|V(A)\cap V(B)|=\left|G[W_y\cap W_{\parent(T,y)}]\right|=1$.
    This concludes the proof.
\end{proof}

\Cref{lemma:rtd_and_separations,lemma:rtd_and_separations2} yield the following recursive definition of rooted $2$-treedepth.
For every graph $G$,
\begin{enumerate}[label={\normalfont (r\arabic*)}]
    \item $\rtd_2(G) = 0$ if $G$ is the null graph; \label{def:rtd2:item:null-graph}
    \item $\rtd_2(G) = 1$ if $G$ is a one vertex graph; and otherwise \label{def:rtd2:item:one-vertex}
    \item $\rtd_2(G)$ is the minimum of $\max\left\{\rtd_2(A),\rtd_2\big(B \setminus V(A)\big) + |V(A) \cap V(B)|\right\}$ over all separations $(A,B)$ of $G$ of order at most one with $V(A) \neq \emptyset$ and $V(B) \setminus V(A) \neq \emptyset$. \label{def:rtd2:item:separations}
\end{enumerate}

The following properties are direct consequences of the previous definition.
For every graph $G$,
\begin{enumerate}[resume*]
    \item $\rtd_2(G)$ is the maximum of $\rtd_2(C)$ over all components $C$ of $G$ when $G$ is not connected;\label{rtd2:item:components}
    \item $\rtd_2(G)$ is the minimum of $\rtd_2(G-v)+1$ over all vertices $v$ of $G$, when $G$ consists of one block;\label{rtd2:item:one-block}
    \item $\rtd_2(G)$ is the minimum of $\max\{\rtd_2(A),\rtd_2(B \setminus V(A)) + 1\}$ over all separations $(A,B)$ of $G$ of order one with $V(A) \cap V(B)$ consisting of a cut-vertex, when $G$ is connected and consists of more than one block; \label{rtd2:item:more-blocks}
    \item $\rtd_2(G) \leq 1+\rtd_2(G-u)$ for every vertex $u$;\label{rtd2:item:removing-vertex}
\end{enumerate}

Item~\ref{def:rtd2:item:separations} in the definition can be in fact strengthened in the following way.

\begin{lemma} \label{lemma:rtd2:separations-stronger}
For every graph $G$ with at least two vertices, $\rtd_2(G)$ is the minimum of $\max\left\{\rtd_2(A),\rtd_2\big(B \setminus V(A)\big) + |V(A) \cap V(B)|\right\}$ over all separations $(A,B)$ of $G$ of order at most one with $V(A) \neq \emptyset$ and $V(B) \setminus V(A) \neq \emptyset$ such that $B$ is a block.
\end{lemma}

\begin{proof}
Indeed, first, assume that $G$ is connected. 
Consider $(A,B)$ a separation of $G$ such that $\rtd_2(G) = \max\left\{\rtd_2(A),\rtd_2\big(B \setminus V(A)\big) + |V(A) \cap V(B)|\right\}$, $V(A) \neq \emptyset$, and $V(B) \setminus V(A) \neq \emptyset$.
Since $G$ is connected, $(A,B)$ has order 1.
Let $(A',B')$ be a separation of $G$ of order at most one with $V(A') \neq \emptyset$ and $V(B') \setminus V(A') \neq \emptyset$ such that $B'$ is a block and $A \subset A'$ and $B' \subset B$.
By \ref{def:rtd2:item:separations}, $\rtd_2(G) \leq \max\left\{\rtd_2(A'),|V(A')\cap V(B')|+\rtd_2\big(B' \setminus V(A'))\right\}$.
To prove the reverse inequality, first observe that $\rtd_2(A') \leq \rtd_2(G)$ because $A'\subseteq G$ and $\rtd_2(B' - V(A')) \leq \rtd_2(B-V(A))$ because $B'-V(A')\subseteq B-V(A)$.
Since $|V(A')\cap V(B')|\leq 1=|V(A)\cap V(B)|$, $\rtd_2(B' - V(A')) + |V(A')\cap V(B')| \leq \rtd_2(B-V(A)) + |V(A)\cap V(B)|\leq \rtd_2(G)$.
Therefore, connected graphs with at least two vertices satisfy the lemma.

Next, suppose that $G$ is not connected. 
By \ref{rtd2:item:components}, $\rtd_2(G)=\rtd_2(C)$ for some component $C$ of $G$.
If $C$ is a block, then $(G-V(C),C)$ is a separation witnessing that $G$ satisfies the statement of the lemma.
If $C$ consists of more than one block, then $C$ satisfies the lemma.
Therefore, we consider a separation $(A,B)$ of $C$ such that $\rtd_2(C) = \max\left\{\rtd_2(A),\rtd_2\big(B \setminus V(A)\big) + |V(A) \cap V(B)|\right\}$, $V(A) \neq \emptyset$, $V(B) \setminus V(A) \neq \emptyset$, and $B$ is a block.
Either $\rtd_2(B-V(A))+|V(A)\cap V(B)| = \rtd_2(C) = \rtd_2(G)$,
or $\rtd_2(A)=\rtd_2(C)$ and $\rtd_2((G-V(C))\cup A) = \max\{\rtd_2(G-V(C)),\rtd_2(A)\} = \rtd_2(C) = \rtd_2(G)$.
Thus, $((G-V(C))\cup A,B)$ witnesses that $G$ satisfies the lemma.
This concludes the proof.
\end{proof}

\begin{lemma}\label{lemma:rtd_2-and-td_2}
    For every graph $G$ with at least one edge,
    \[
    \td_2(G) \leq \rtd_2(G) \leq 2 \td_2(G) - 2.
    \]
\end{lemma}

\begin{proof}
    First, we prove that $\td_2(G) \leq \rtd_2(G)$ for every graph $G$.
    We proceed by induction on $|V(G)|$.
    When $G$ is a null graph, we have $\td_2(G) = \rtd_2(G) = 0$ and when $G$ is a one-vertex graph, we have $\td_2(G) = \rtd_2(G) = 1$.
    Thus, we assume that $|V(G)| \geq 2$.
    If $G$ consists of one block, then by~\ref{rtd2:item:one-block} and induction hypothesis,
        \[\td_2(G) = \min_{v\in V(G)} \td_2(G - v) + 1 \leq \min_{v\in V(G)} \rtd_2(G - v) + 1 = \rtd_2(G).\]
    If $G$ consists of blocks $B_1,\dots,B_k$ for $k > 1$, then by the induction hypothesis,
        \[\td_2(G) = \max_{i \in [k]} \td_2(B_i) \leq \max_{i \in [k]} \rtd_2(B_i) \leq \rtd_2(G).\]

    Now, we prove the other inequality for every graph $G$ with at least one edge. 
    We again proceed by induction on $|V(G)|$.
    If $\td_2(G)=2$, then $G$ is a forest with at least one edge, and so as mentioned earlier $\rtd_2(G) = \td_2(G) = 2$.
    Now assume that $\td_2(G) \geq 3$, and so in particular $|V(G)| \geq 3$, and that the result holds for smaller graphs.
    In particular, for every graph $H$ with $|V(H)|<|V(G)|$, either $H$ has no edge and so $\rtd_2(H)=\td_2(H)=1$,
    or $\rtd_2(H) \leq 2\td_2(H)-2$. In both cases, $\rtd_2(H) \leq \max\{1,2\td_2(H)-2\}$.
    By \Cref{lemma:rtd2:separations-stronger},
    there is a separation $(A,B)$ of $G$ of order at most one such that $\rtd_2(G) = \max\{\rtd_2(A),\rtd_2(B \setminus V(A)) + |V(A) \cap V(B)|\}$, $V(B) \setminus V(A) \neq \emptyset$, $V(A) \neq \emptyset$, and $B$ is a block of $G$.
    If $|V(A) \cap V(B)| = 0$, then $B-V(A) = B$ and so
    \begin{align*}
        \rtd_2(G) 
        &= \max\{\rtd_2(A),\rtd_2(B)\} \\
        &\leq \max\{\max\{1,2\td_2(A)-2\},\max\{1,2\td_2(B)-2\}\} \\
        &\leq \max\{1,2\max\{\td_2(A),\td_2(B)\} - 2\} \\
        &\leq 2\td_2(G)-2.
    \end{align*}
    Therefore, we assume that $|V(A) \cap V(B)|=1$ and $V(A) \cap V(B) = \{u\}$.
    There exists $v \in V(B)$ such that $\td_2(B-v) = \td_2(B)-1$.
    Then, by~\ref{rtd2:item:removing-vertex},
    \begin{align*}
        \rtd_2(B-u) &\leq \rtd_2(B-v) + 1\\
       3 &\leq \max\{1,2\td_2(B-v)-2\} + 1 \\
        &\leq \max\{2,2\td_2(B-v)-1\} \\
        &\leq \max\{2, 2\td_2(B)-3\}.
    \end{align*}
    Finally, since $\td_2(G) \geq 3$,
    \begin{align*}
        \rtd_2(G) &= \max\{\rtd_2(A),\rtd_2(B \setminus u) + 1\}\\
        &\leq \max\{\max\{1,2\td_2(A)-2\}, \max\{2, 2\td_2(B)-3\}+1\}\\
        &\leq \max\{3,2\max\{\td_2(A),\td_2(B)\} - 2\} \\
        &\leq \max\{3,2\td_2(G)-2\}\\
        &\leq 2\td_2(G) - 2.\qedhere
    \end{align*}    
\end{proof}

The bounds in \Cref{lemma:rtd_2-and-td_2} are tight. 
Indeed, for every positive integer $n$, we have $\td_2(K_n) = \rtd_2(K_n) = n$, which witnesses that the first inequality is tight.
For the second one, see \Cref{lemma:construction_rtd2_td2}, which we precede with a simple observation.

\begin{lemma}\label{lemma:add_apex_increases_rtd2}
    For every graph $G$, 
    \[\rtd_2(K_1 \oplus G) = 1 + \rtd_2(G).\]
\end{lemma}

\begin{proof}
    Let $G$ be a graph and let $s$ the vertex of $K_1$ in $K_1 \oplus G$.
    By \ref{rtd2:item:removing-vertex}, $\rtd_2(K_1 \oplus G) \leq 1 + \rtd_2(G)$.
    For the other inequality, we proceed by induction on $\rtd_2(G)$.
    The assertion is clear when $G$ is the null graph, thus, assume that $G$ is not the null graph.
    If $G$ is not connected, then $\rtd_2(G) = \rtd_2(C)$ for some component $C$ of $G$ by \ref{rtd2:item:components}, and since $\rtd_2(K_1 \oplus G) \geq \rtd_2(K_1 \oplus C)$, it suffices to show $\rtd_2(K_1 \oplus C) \geq 1 + \rtd_2(C)$.
    Therefore, we assume that $G$ is connected.
    Since $K_1 \oplus G$ is also connected, there is a separation $(A,B)$ of $K_1 \oplus G$ of order one such that $\rtd_2(G) = \max\{\rtd_2(A),\rtd_2(B \setminus V(A)) + 1\}$, $V(B) \setminus V(A) \neq \emptyset$, and $V(A) \neq \emptyset$ by \ref{def:rtd2:item:separations}.
    Since $s$ is adjacent to all other vertices in $K_1 \oplus G$, the only possibility is that $V(A) = \{s\}$ and $V(B) = V(K_1 \oplus G)$.
    It follows that $B-V(A)=G$ thus, $\rtd_2(K_1 \oplus G) \geq 1 + \rtd_2(B \setminus V(A)) = 1 + \rtd_2(G)$.
\end{proof}

\begin{lemma}\label{lemma:construction_rtd2_td2}
    For every integer $k$ with $k \geq 2$, there is a graph $G$ with
    $\td_2(G)\leq k$ and $\rtd_2(G)\geq 2k-2$.
\end{lemma}

\begin{figure}[tp]
    \centering 
    \includegraphics[scale=1]{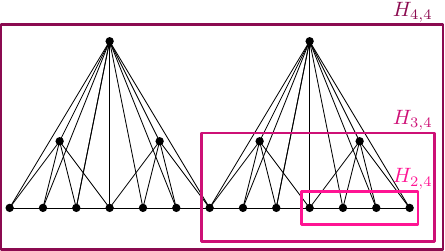} 
    \caption{The proof of \Cref{lemma:construction_rtd2_td2} implies that $\td_2(H_{4,4}) \leq 4$ and $\rtd_2(H_{4,4}) \geq 6$.}
    \label{fig:lower_bound_rtd2_td2}
\end{figure}

\begin{proof}
    We define inductively graphs $H_{k,\ell}$ with two distinguished vertices $u_{k,\ell}$ and $v_{k,\ell}$ for every integers $k,\ell$ with $k,\ell \geq 2$.
    For $k=2$, $H_{k,\ell}$ is a path on $\ell$ vertices and $u_{k,\ell},v_{k,\ell}$ are its endpoints.
    For $k \geq 3$, $H_{k,\ell}$ is obtained from two disjoint copies $H_1,H_2$ of $K_1\oplus H_{k-1,\ell}$ by identifying the copy of $v_{k-1,\ell}$ in $H_1$ with the copy of $u_{k-1,\ell}$ in $H_2$. The vertices $u_{k,\ell}, v_{k,\ell}$ are then respectively the copy of $u_{k-1,\ell}$ in $H_1$ and the copy of $v_{k-1,\ell}$ in $H_2$.
    See Figure~\ref{fig:lower_bound_rtd2_td2}.
    
    By induction on $k$, we show that $\td_2(H_{k,\ell}) \leq k$ and $\rtd_2(H_{k,\ell}) \geq 2k-2$ for all integers $k,\ell$ with $\ell \geq k \geq 2$.
    When $k=2$, $H_{2,\ell}$ is a path on at least two vertices and so $\td_2(H_{2,\ell}) = \rtd_2(H_{2,\ell}) = 2$.
    Now suppose that $k\geq 3$.
    First, observe that $H_{k,\ell}$ has exactly two blocks $H_1,H_2$, both isomorphic to $K_1 \oplus H_{k-1,\ell}$.
    Hence, $\td_2(H_{k,\ell}) \leq \td_2(K_1 \oplus H_{k-1,\ell}) \leq 1 + \td_2(H_{k-1,\ell}) \leq k$ by the induction hypothesis.
    Let $v$ be the unique cut-vertex of $H_{k,\ell}$.
    Since $H_{k,\ell}$ is connected and consists of more than one block, by~\ref{rtd2:item:more-blocks}, there is a separation $(A,B)$ of $H_{k,\ell}$ such that $V(A) \cap V(B) = \{v\}$ and $\rtd_2(H_{k,\ell}) = \max\{\rtd_2(A), \rtd_2(B-v)+1\}$.
    It follows that the graph $B-v$ contains $K_1 \oplus H_{k-1,\ell-1}$ as a subgraph, and so, applying \Cref{lemma:add_apex_increases_rtd2},
    \[\rtd_2(H_{k,\ell-1}) \geq \rtd_2(B-v) + 1 \geq \rtd_2(K_1 \oplus H_{k-1,\ell-1}) + 1 \geq \rtd_2(H_{k-1,\ell-1}) + 2 \geq 2k-2.\]
    This concludes the proof of the lemma.
\end{proof}

\section{Bounding weak coloring numbers in terms of pathwidth} \label{ssec:pathwidth}

For completeness, we provide a proof of a logarithmic bound on weak coloring numbers of graphs of bounded pathwidth.

Let $G$ be a graph.
A \defin{path decomposition} of $G$ is a sequence $(W_1, \dots, W_m)$ of subsets of vertices such that
\begin{enumerate}
    \item for every $u \in V(G)$, $\{i \in [m] \mid u\in W_i\}$ is a nonempty interval of integers and
    \item for every edge $uv \in E(G)$, there exists $i \in [m]$ such that $u,v \in W_i$.
\end{enumerate}
The \defin{width} of this path decomposition is $\max_{i \in [m]} |W_i|-1$.
The \defin{pathwidth} of $G$,
denoted by \defin{$\pw(G)$}, is the minimum width of a path decomposition of $G$.

\begin{lemma}\label{lemma:wcol_path_log}
    Let $G$ be a graph.
    For every positive integer $q$,
    \[
        \wcol_q(G) \leq \binom{\pw(G)+2}{2} (\lceil\log q\rceil + 2).
    \]
\end{lemma}

\begin{proof}
    We proceed by induction on $\pw(G)$.
    If $\pw(G) = -1$, then $G$ is the null graph and $\wcol_q(G)=-\infty$ for all positive integers $q$.
    Suppose that $\pw(G) \geq 0$.
    Let $(W_1, \dots, W_m)$ be a path decomposition of $G$ of width $\pw(G)$,
    and let $q$ be a positive integer.
    By convention, let $W_{m+1} = \emptyset$.
    Let $I \subseteq [m]$ be inclusion-wise maximal set such that the bags $W_i$ for $i \in I$
    are pairwise disjoint.
    Let $S = \bigcup_{i \in I} W_i$.
    First, observe that for every $i \in [m]$, by maximality of $I$, $W_i \cap S \neq \emptyset$, 
    and so $\pw(G-S) < \pw(G)$ as witnessed by the path decomposition $(W_1 \setminus S, \dots, W_m \setminus S)$.
    Second, for every component $C$ of $G-S$, there exists $i(C),j(C) \in I$ with $i(C) \leq j(C)$, $\{i(C)+1, \dots, j(C)-1\} \cap I = \emptyset$ and such that $N_G(V(C)) \subseteq W_{i(C)} \cup W_{j(C)}$. 
    We claim that $\wcol_q(G,S) \leq (\pw(G)+1)(\lceil\log q\rceil + 2)$.

    To see that, consider the path $Q$ with vertex set $I$, in the natural order.
    By \Cref{lemma:wcol_paths}, there is an ordering $\sigma_Q$ of $I$ such that $\wcol_q(Q,\sigma_Q) \leq \lceil\log q\rceil + 2$.
    Let $\sigma$ be an ordering of $S$ such that for every $i,j \in I$ with $i <_{\sigma_Q} j$,
    for every $u \in W_i$ and for every $v \in W_j$, 
    we have $u <_\sigma v$.
    Let $u \in V(G)$.
    If $u \in S$, then let $i(u) \in I$ such that $u \in W_{i(u)}$.
    If $u \not\in S$, consider the component $C$ of $u$ in $G-S$ and let $i(u) = i(C)$.
    Let $v \in \WReach_q[G,S,\sigma,u]$ and let $(u_1, \dots, u_\ell,v)$ be a path witnessing this fact.
    Then, $(i(u_1), \dots, i(u_\ell),i(v))$ is a walk in $Q$ whose minimum according to $\sigma_Q$ is $i(v)$.
    Therefore, $i(v) \in \WReach_q[Q,\sigma_Q,i(u)]$.
    This proves that
    \(
        \WReach_q[G,S,\sigma,u] \subseteq \bigcup \{W_i \mid i \in \WReach_q[Q,\sigma_Q,i(u)]\}.
    \)
    Since $|W_i| \leq \pw(G)+1$ for every $i \in [m]$,
    this proves that $\wcol_q(G,S) \leq (\pw(G)+1)(\lceil\log q)\rceil + 2)$.

    Finally, by \Cref{obs:wcol_union} and the induction hypothesis applied to $G-S$, we deduce
    \begin{align*}
        \wcol_q(G) 
        &\leq \wcol_q(G,S) + \wcol_q(G-S) \\
        &\leq (\pw(G)+1)(\lceil\log q\rceil + 2) + \binom{\pw(G)+1}{2} (\lceil\log q\rceil + 2) \\
        &\leq \binom{\pw(G)+2}{2} (\lceil\log q\rceil + 2). \qedhere
    \end{align*}
\end{proof}
\end{document}